\documentclass[reqno,12pt]{amsbook}
\usepackage{amsmath,amsthm,amscd,amssymb}
\numberwithin{equation}{section}
\makeatletter

\def\th@plain{%
\let\thm@indent\noindent
\thm@headfont{\bfseries}
\let\thmhead\thmhead@plain
\let\swappedhead\swappedhead@plain
\thm@preskip.5\baselineskip\@plus.2\baselineskip\@minus.2\baselineskip
\thm@postskip\thm@preskip
\slshape%\itshape %
}

\def\th@remark{%
\let\thm@indent\noindent
\thm@headfont{\bfseries}
\let\thmhead\thmhead@plain
\let\swappedhead\swappedhead@plain
\thm@preskip.5\baselineskip\@plus.2\baselineskip\@minus.2\baselineskip
\thm@postskip\thm@preskip
\upshape
}

\makeatother

\theoremstyle{plain}
\newtheorem{Theorem}{Theorem}[section]
\newtheorem{Corollary}[Theorem]{Corollary}
\newtheorem{Lemma}[Theorem]{Lemma}
\newtheorem{Lemma-Definition}[Theorem]{Lemma-Definition}
\newtheorem{Proposition}[Theorem]{Proposition}

\theoremstyle{plain} 
\newtheorem{ThmA}{Theorem A}

\theoremstyle{plain} 
\newtheorem{ThmB}{Theorem B}

\theoremstyle{plain} 
\newtheorem{ThmC}{Theorem C}

\theoremstyle{plain} 
\newtheorem{ThmD}{Theorem D}

\theoremstyle{plain} 
\newtheorem{ThmE}{Theorem E}

\theoremstyle{plain} 
\newtheorem{ThmF}{Theorem F}

\theoremstyle{plain} 
\newtheorem{ThmG}{Theorem G}

\theoremstyle{plain} 
\newtheorem{ThmH}{Theorem H}

\theoremstyle{remark}
\newtheorem{Remark}[Theorem]{\bfit Remark}

\newtheorem{Definition}[Theorem]{\bfit Definition}
\newtheorem{Definition-Notation}[Theorem]{\bfit Definition-Notation}
\newtheorem{Example}[Theorem]{\bfit Example}
\newtheorem{Examples}[Theorem]{\bfit Examples}
\newtheorem{Example-Definition}[Theorem]{\bfit Example-Definition}
\newtheorem{Notation}[Theorem]{\bfit Notation}
\newtheorem{Problem}[Theorem]{\bfit Problem}

\newtheorem{Observation}{Observation}

\newcommand{\Def}{\stackrel{\mathrm{def}}{=\!\!=}}

\catcode`\@=11

 \def\AMSTeXfeatures{\Plainheads 
   \let\current@vert=\AMS@vert}

 \def\Plainheads{\sh@ftdiam=0.05em
   \getlabeldims
   \let\vshaftfill=\plnvsolidfill
   \let\hshaftfill=\plnhsolidfill
   \let\th@rhead=\plnrhead
   \let\th@lhead=\plnlhead
   \let\th@dnhead=\plndnhead
   \let\th@uphead=\plnuphead}
 
 \def\glet{\global\let}

 \def\LaTeXfeatures{\catcode`\@=11
   \ifx\@clnwd\undefined \nol@g
      \input ltxcode.tex \dol@g \fi
   \ltxheads \let\current@vert=\new@vert
   \providelto \catcode`\@=\active}

 \def\nol@g{\def\wlog{\edef\garbage}}
 \def\dol@g{\let\wlog=\wl@g} \let\wl@g=\wlog
 \nol@g % Silences allocations with \newbox, etc.

 \newbox\ltobox
 \def\providelto{{\setbox\z@=
   \hbox{$\to$}\minharrlen=\wd\z@
   \global\setbox\ltobox=\hbox{$\activeat>>>$}}
   \def\lto{\mathrel{\copy\ltobox}}}

 \def\ltxheads{\sh@ftdiam=\@wholewidth
   \getlabeldims
   \let\vshaftfill= \ltxvsolidfill
   \let\hshaftfill=\ltxhsolidfill
   \let\th@rhead=\ltxrhead
   \let\th@lhead=\ltxlhead
   \let\th@dnhead=\ltxdnhead
   \let\th@uphead=\ltxuphead}
 {\catcode`\@=\active
   \gdef@#1{\csname #1\string@at\endcsname}
   \glet\activeat=@}
 \def\def@#1{\expandafter\def\csname #1@at\endcsname}

 \def@>#1>#2>{\@rrow R{#1}{#2}}
 \def@<#1<#2<{\@rrow L{#1}{#2}}
 \def@ V#1V#2V{\@rrow V{#1}{#2}}
 \def@ A#1A#2A{\@rrow A{#1}{#2}}
 \def@/#1/#2/#3/{\@rrow{#1}{#2}{#3}}
 \def@.{\ifodd\row\ifmmode\noharrow
     \else\leavevmode.\spacefactor3000 \fi
   \else\novarrow\fi}
 \def@={\ifodd\row\harrow\hequalfill{}{}%
   \else\varrow\vequalfill{}{}\fi}
 \def@:#1{\ifx=#1\harrow\deffill{}{}%
   \else\leavevmode\null:#1\fi}
 \def@|{\current@vert}
  \def\AMS@vert{\varrow\vequalfill{}{}}
  \def\new@vert#1|#2|{\ifodd\row
   \let\nextarrow\vertexvarrow
   \else\let\nextarrow\varrow\fi
   \nextarrow\vshaftfill{#1}{#2}}
 \def@-{\ifmmode\let\next\hl@ne
   \else\let\next\AMSatdash \fi \next}
  \def\hl@ne#1-#2-{\harrow\hshaftfill{#1}{#2}}
  \def\AMSatdash{\let\next\relax\leavevmode
    \def\next@{\ifx\next-%
      \def\next-{\futurelet\next\nextii@}%
     \else\def\next{\hbox{-}}\fi\next}%
    \def\nextii@{\ifx\next-\def\next-{\hbox{---}}%
      \else\def\next{\hbox{--}}\fi\next}%
    \futurelet\next\next@}
 \def@(#1){\tweenarrows{#1}}
 \def@[#1]{\setsp@n#1\relax\activeat}
 \def\fiberbox{\hbox{$\vcenter{\hr@le\hbox{\vr@le
   \kern1ex\vbox{\kern1.2ex}\vr@le}\hr@le}$}}
  \def\hr@le{\hrule height \sh@ftdiam}
  \def\vr@le{\vrule width \sh@ftdiam}
 \def@+#1+#2+#3+{\ifodd\row \harrow{#1}{#2}{#3}%
   \else \varrow{#1}{#2}{#3}\fi}

 \def\Dnarrfill{\vequalfill\Dnhe@d}
 \def\Uparrfill{\Uphe@d\vequalfill}
 \def\ontofill{\rtarrfill\kern-0.3em %2\he@dwd
   \th@rhead\kern 0.3em} %new def

 \def\rtarrfill{\hshaftfill\th@rhead}
 \def\ltarrfill{\th@lhead\hshaftfill}
 \def\dnarrfill{\vshaftfill\th@dnhead}
 \def\uparrfill{\th@uphead\vshaftfill}
 \def\hequalfill{\plnhfill=}
 \def\deffill{:\plnhfill=}
 \def\plnvextfill#1{\setbox\z@
   \hbox{\the\textfont3 #1}%
   \dimen@=\dp\z@\advance\dimen@\ht\z@
   \copy\z@ \kern-\dimen@ %-\dp\z@
   \cleaders\copy\z@ \vfill
   \kern-\dimen@ %-\dp\z@
   \box\z@}
 \def\plnhfill#1{$\m@th\mkern-1.5mu\mathord#1\mkern-6mu
    \cleaders\hbox{$\mkern-2mu\mathord#1\mkern-2mu$}\hfill
    \mkern-6mu\mathord#1\mkern-1.5mu$}
 \def\vequalfill{\plnvextfill{\char'167}}
 \def\plnvsolidfill{\plnvextfill{\char'077}}
 \def\plnhsolidfill{\plnhfill-}
 \def\ltxhsolidfill{\leaders\hrule height\topofshaft depth\botofshaft
   \hfill}
 \def\ltxvsolidfill{\leaders\vrule width\sh@ftdiam\vfill}
 \def\hdashfill{\hd@sh\wd@sh
   \xleaders \hbox{\wd@sh\hd@sh\wd@sh}\hfill
   \wd@sh\hd@sh}
 \def\vdashfill{\vd@sh\wd@sh
   \xleaders \vbox{\wd@sh\vd@sh\wd@sh}\vfill
   \wd@sh\vd@sh}
 \def\dashed{\ifinmeasureCD\else
    \ifodd\row\option{\let\hshaftfill=\hdashfill}%
   \else\option{\let\vshaftfill=\vdashfill}\fi\fi}
 \newdimen\CDstrutht  \newdimen\CDstrutdp
 \newdimen\CDstrutlen \CDstrutlen=\CDstrutht
   \advance\CDstrutlen by \CDstrutdp
 \def\CDstrut{\vrule
   height \ifnum\row=1 \z@\else\CDstrutht \fi
   depth \ifnum\row=\numrows \z@ \else\CDstrutdp \fi
   width\z@}
 \newdimen\CDarrsurr \CDarrsurr=0.375em
 \newdimen\CDdashlen
 \newdimen\CDvarrlen \CDvarrlen=1.5\baselineskip
 \newdimen\minharrlen %Used outside CD's
  \setbox\z@\hbox{$\longrightarrow$} \minharrlen=\wd\z@
 \newdimen\minCDharrlen \minCDharrlen=2.5em %825 2pc %2.5pc
\newdimen \minc@lwd
\def\findminc@lwd{\minc@lwd=2\CDarrsurr
  \advance\minc@lwd\minCDharrlen}
 \newdimen\sh@ftdiam

 \newdimen\labelsurr \labelsurr=1.25 em

\newcount\sp@ncnt \sp@ncnt=\@ne
\newcount\sp@ncnt@ \sp@ncnt@=\@ne
\newdimen\@rrwd \newdimen\@rrdp

 \def\adjustbot#1{\option{\advance\@rrdp#1\relax}}
\def\pushvertex#1{\global\p@shlen#1\relax
   \global\let\maybepush=\dopush}

 \newdimen\p@shlen \p@shlen=\z@

 \let\maybepush=\relax
 \def\dopush{\ifinmeasureCD %omitted by accident
   \advance\locdimen by -\p@shlen %AL
   \else\advance \@rrwd by -\p@shlen \fi %AL
   \global\let\maybepush=\relax \global\p@shlen=\z@\relax}
 \def\span@ne{\global\sp@ncnt=\@ne\relax}
 \def\setsp@n#1#2{\global\sp@ncnt=#1\relax
   \ifx\relax#2\relax\else\global\sp@ncnt@=#2\relax\fi}

 \def\plnrhead{\llap{$\rightarrow\mkern-1.5mu$}}
 \def\plnlhead{\rlap{$\mkern-1.5mu\leftarrow$}}

 \def\clap#1{\hbox to \z@{\hss #1\hss}}

 \def\plndnhead{\hbox{\the\textfont3 \char'171}}
 \def\plnuphead{\hbox{\the\textfont3 \char'170}}
 \def\Dnhe@d{\hbox{\the\textfont3 \char'177}}
 \def\Uphe@d{\hbox{\the\textfont3 \char'176}}

 \def\ltxrhead{\raise\@xisheight
   \llap{\smash{\@linefnt\@getrarrow(1,0)}}}
 \def\ltxlhead{\raise\@xisheight
   \rlap{\@linefnt\@getlarrow(-1,0)}}
 \def\ltxuphead{\setbox\z@=\rlap{%
   \kern\@halfwidth\@linefnt\char'66}%
   \copy\z@\kern-\ht\z@}
 \def\ltxdnhead{\setbox\z@=\rlap{%
   \kern\@halfwidth\@linefnt\char'77}%
   \ht\z@=\z@\box\z@}

 \def\wd@sh{\kern0.5\CDdashlen}
 \def\hd@sh{\vrule height\topofshaft depth\botofshaft
    width\CDdashlen}
 \def\vd@sh{\hrule height\CDdashlen
   depth\z@ width\sh@ftdiam}
\def\xylist{14{3434}13{2414}12{1723}%
  23{1413}34{1153}11{0867}43{0707}%
  32{0580}21{0414}31{0291}41{0}}
\newcount\tgtcnt@
\def\find@xyargs{\dimen@=\@rrdp
  \advance\dimen@ by \CDstrutlen
  \tgtcnt@=\dimen@ \dimen@=\@rrwd %\relax
  \divide\dimen@ by \@m %\relax
  \divide \tgtcnt@ by \dimen@ %\relax
  \expandafter\testxy\xylist\relax
  \unitlength=\@xarg\@rrdp
  \divide\unitlength by\@yarg\relax}
\def\testxy#1#2#3{\ifnum\tgtcnt@>#3
    \@xarg=#1\relax \@yarg=#2\relax
    \let\next=\ignorerest
  \else\let\next\testxy\fi\next}
\def\ignorerest#1\relax{\relax}

\let\scalefactor=\@ne
\def\SWarrow{\find@xyargs\vector
  (-\@xarg,-\@yarg)\scalefactor\hskip-\wd\@linechar}
\def\NWarrow{\find@xyargs\vector
  (-\@xarg,\@yarg)\scalefactor\hskip-\wd\@linechar}
\def\NEarrow{\find@xyargs\vector
  (\@xarg,\@yarg)\scalefactor}
\def\SEarrow{\find@xyargs\vector
  (\@xarg,-\@yarg)\scalefactor}
\def\rightupline{\find@xyargs\@linelen=\scalefactor
     \unitlength\@sline}
\def\rightdownline{\find@xyargs\@yarg=-\@yarg\relax
     \@linelen=\scalefactor\unitlength\@sline}

\def\Sim{\ifodd\row\setbox\z@=\hbox{$\sim$}\dimen@=\ht\z@
 \advance\dimen@ by -\@xisheight
  \vbox{\box\z@\kern-\@xisheight\kern\dimen@}%
  \else\hbox{$\wr$}\fi}

\def\harrow#1#2#3{\inmeasureCDtrue\findminarrwd
  {#2}{#3}{\sp@ncnt\minharrlen}\inmeasureCDfalse\span@ne
  \mathrel{\hbox{\options\hplace{#1}\ulabel{#2}\dlabel{#3}}}}

\def\noharrow{\harrow\hfill{}{}}
\def\vertexvarrow#1#2#3{\findarrdp \@rrwd=\z@ \setsp@n\@ne\@ne
  \vbox to \z@{\kern-1.2\CDstrutht
  \rlap{\options\vplace{#1}\llabel{#2}\rlabel{#3}}\vss}}

\newif\ifinmeasureCD
\def\measurelabel#1{\setbox\z@
  \hbox{$\scriptstyle#1\kern\labelsurr$}%
  \ifdim\wd\z@>\@rrwd \@rrwd=\wd\z@\fi}
\def\findminarrwd#1#2#3{\@rrwd=#3\relax
   \measurelabel{#1}\measurelabel{#2}}
\def\findCDarrwd#1#2{\@rrwd=\minCDharrlen
   \measurelabel{#1}\measurelabel{#2}%
  }

\newcount\row \row=\@ne \newcount\col \col=\@ne %96 initialized
 \newcount\numrows 
\numrows=\@ne
 \newcount\numcols
\newcount\arrspan \newdimen\vrtxhalfwd  \newbox\tempbox

\def\DANABUG{\advance\col by \@ne
 \@rrwd=\minCDharrlen
  \advance\@rrwd by \vrtxhalfwd
  \advance\@rrwd by \CDarrsurr
  \ifnum\col>\numcols \numcols=\col
     \newlocdimen{col\the\col}\locdimen=\@rrwd %AL
  \else \ifdim\@rrwd>\c@l \c@l=\@rrwd\fi\fi}

\def\drop#1\\{%\noharrow %caused by DANABUG
  \findvrtxhalfsum\DANABUG\advance\row by 2 \measureinit}

\def\measureinit{\col=\@ne \vrtxhalfwd=-\CDarrsurr\arrspan=\@ne\@rrwd=\z@
   \setbox\tempbox=\hbox\bgroup$}
\def\measure{%CR \bgroup
  \let\harrow\measureCDarrow
  \let\CDCR=\measureCR %CR
   \findminc@lwd 
  \inmeasureCDtrue
  \row=\@ne \numcols=\z@ \measureinit}

\def\endmeasure{\findvrtxhalfsum\DANABUG
  \numrows=\row %CR \egroup
  \inmeasureCDfalse}

\def\newlocdimen#1{\advance\dimenc@unt by \@ne
  \ifnum\dimenc@unt<\insc@unt
     \else\errmessage{No room for the CD}\fi
  \dimendef\locdimen=\dimenc@unt
  \expandafter\dimendef\csname#1\endcsname=\dimenc@unt}

 \def\r@wc@l{\csname row\the\row col\the\col\endcsname}
 \def\c@l{\csname col\the\col\endcsname}

 \def\findvrtxhalfsum{$\egroup
  \newlocdimen{row\the\row col\the\col}%%AL
  \locdimen=\vrtxhalfwd %AL
  \vrtxhalfwd=0.5\wd\tempbox %\maybes@ve %8231
  \advance\vrtxhalfwd by \CDarrsurr
  \advance\locdimen by \vrtxhalfwd %AL
  \advance\@rrwd by \locdimen %AL
  \maybepush
  \divide\@rrwd by \arrspan\relax
  \ifdim\@rrwd<\minc@lwd
    \ifnum\col>\@ne \@rrwd=\minc@lwd\fi \fi
  \loop %94 \edef\c@l{\csname col\the\col\endcsname}
    \ifnum\col>\numcols \numcols=\col
       \newlocdimen{col\the\col}% %AL
       \locdimen=\@rrwd %AL
    \else \ifdim\@rrwd>\c@l \c@l=\@rrwd\fi \fi
   \ifnum\arrspan>\@ne
      \advance\arrspan by -1 \advance\col by \@ne
  \repeat }

 \def\measureCDarrow#1#2#3{\findvrtxhalfsum
   \arrspan=\sp@ncnt\relax\global\sp@ncnt=1\relax
   \advance\col by \@ne
   \findCDarrwd{#2}{#3}%
   \setbox\tempbox=\hbox\bgroup$}

 \newcount\dr@tn \dr@tn=\z@
 \def\locate#1:#2{\ifinmeasureCD\else
   \count@=-#1
   \multiply\count@ by 2
   \advance\count@ by #2
   \dimen@=\count@\@rrwd
   \ifnum\dr@tn=\@ne\relax \else\dimen@=-\dimen@ \fi
   \dimen@i=\@rrdp
   \ifnum\dr@tn>\z@\advance\dimen@i by \CDstrutlen \fi
   \dimen@i=\count@\dimen@i
   \count@=#2 \multiply\count@ by 2
   \divide\dimen@ by \count@
   \divide\dimen@i by \count@
   \lift\dimen@i\nudge\dimen@\fi}
\def\betweenCDrows{\advance\row by \@ne \col=\@ne
\options}

\def\hbegin{\hbox\bgroup\kern\c@l \kern-\r@wc@l$}
\def\hend{$\glet\maybepush\relax \CDstrut\egroup}
\def\vbegin{\setbox\tempbox=\hbox\bgroup$}
\def\vend{$\egroup\ht\tempbox=\z@\dp\tempbox\CDvarrlen
  \box\tempbox}
\def\setCD{\let\harrow=\setCDarrow
  \let\CDCR=\setCR %CR
  \row=\@ne \col=\@ne \hbegin}
\let\endsetCD=\hend %AT (Assume CD ends with hmaterial)

\def\findarrwd{\@rrwd=\z@ \count@=\col \advance\count@ by\sp@ncnt
  \loop\ifnum\count@>\col \advance\count@ by -1
      \advance\@rrwd by\csname col\the\count@\endcsname\repeat}
\def\setCDarrow#1#2#3{\kern\CDarrsurr\advance\col by \@ne
  \findarrwd \advance\@rrwd by -\r@wc@l  
  \@rrdp=\z@ %&0211 (It might be used by \locate).
  \maybepush
  \advance\col by -\@ne \advance\col by \sp@ncnt \span@ne
  \hbox to \@rrwd{\options
   \@rrwd=\scalefactor\@rrwd\hss
   \hplace{#1}\ulabel{#2}\dlabel{#3}\hss}%
   \kern\CDarrsurr}

\newdimen\labspacei %96 use subscript min of TeXbook 13a, p.444
\newdimen\labspaceii %96 Note many letters stick down below their boxes.

\newdimen\@xisheight
  \@xisheight=\the\fontdimen22\textfont2
\newdimen\labelskip
  \labelskip=\the\fontdimen10\textfont3 %2pt
\newdimen\topofshaft
\newdimen\botofshaft
\newdimen\botofulabel
\newdimen\topofdlabel
\def\getlabeldims{
  \topofshaft=0.5\sh@ftdiam
  \botofshaft=\topofshaft
  \advance\topofshaft by \@xisheight  
  \advance\botofshaft by -\@xisheight  
  \botofulabel=\topofshaft
  \advance\botofulabel by \labelskip
  \topofdlabel=\botofshaft
  \advance\topofdlabel by \labelskip}

\def\ulabel{\ifnum\row=\@ne\let\next\ulabeli
   \else\let\next\ulabellap\fi\next}
\def\ulabeli#1{\vbox{
  \clap{\kern-\@rrwd$\scriptstyle#1$}%
  \kern\botofulabel}\maybeoffset}
\def\ulabellap#1{\vbox to \z@{\vss
  \clap{\kern-\@rrwd$\scriptstyle#1$}%
  \kern\botofulabel}\maybeoffset}
\def\dlabel{\ifnum\row=\numrows\let\next\dlabeli
   \else\let\next\dlabellap\fi\next}
\def\dlabeli#1{\vtop{\kern\topofdlabel
  \clap{\kern-\@rrwd$\scriptstyle#1$}%
  }\maybeoffset}
\def\dlabellap#1{\vbox to \z@{\kern\topofdlabel
  \clap{\kern-\@rrwd$\scriptstyle#1$}%
  \vss}\maybeoffset}
\def\rlabel#1{\vbox to \z@{\vss
  \rlap{\kern\labelskip$\scriptstyle#1$}%
  \vss\kern-\@rrdp}\maybeoffset}
\def\llabel#1{\vbox to \z@{\vss
  \llap{$\scriptstyle#1$\kern\labelskip}%
  \vss\kern-\@rrdp}\maybeoffset}
\def\swlabel#1{\vtop{\kern0.5\@rrdp
  \llap{$\scriptstyle#1$\kern\labelskip\kern-0.5\@rrwd}
  }\maybeoffset}
\def\nwlabel#1{\vbox{
  \llap{$\scriptstyle#1$\kern\labelskip\kern-0.5\@rrwd}%
  \kern-0.5\@rrdp}\maybeoffset}
\def\selabel#1{\vtop{\kern0.5\@rrdp
  \rlap{\kern0.5\@rrwd\kern\labelskip$\scriptstyle#1$}%
  }\maybeoffset}
\def\nelabel#1{\vbox{
  \rlap{\kern0.5\@rrwd\kern\labelskip$\scriptstyle#1$}%
  \kern-0.5\@rrdp}\maybeoffset}
\def\cplace#1{\vbox to \z@{\vss
  \clap{$#1$\kern-\@rrwd}%
  \kern-\@rrdp\vss}\maybeoffset}
\def\hplace#1{\hbox to \@rrwd{#1}\maybeoffset}
\def\vplace#1{\clap{\vbox to \z@{#1\kern-\@rrdp}}\maybeoffset}
\newdimen\nudgeamount \nudgeamount=\z@
\newdimen\liftamount \liftamount=\z@
\let\maybeoffset\relax
\newbox\offsetbox \newdimen\lastheight
\def\dooffset{%assumes that \lastbox is a <box> set in horiz. mode
  \setbox\offsetbox=\lastbox \lastheight=\ht\offsetbox 
  \setbox\offsetbox=\vbox{\kern-\liftamount\box\offsetbox}%
  \ht\offsetbox=\lastheight
  \kern\nudgeamount\box\offsetbox\kern-\nudgeamount
  \global\nudgeamount=\z@ \global\liftamount=\z@
  \glet\maybeoffset=\relax}
\def\nudge#1{\ifinmeasureCD\else
  \global\advance\nudgeamount#1\relax
  \global\let\maybeoffset\dooffset\fi}
\def\lift#1{\ifinmeasureCD\else
  \global\advance\liftamount#1\relax
  \global\let\maybeoffset\dooffset\fi}

\def\findarrdp{\@rrdp=\CDvarrlen
  \ifnum\sp@ncnt@>1
    \advance\@rrdp by \CDstrutlen
    \multiply\@rrdp by \sp@ncnt@
    \advance\@rrdp by -\CDstrutlen \fi
 }

\def\varrow#1#2#3{\ifnum\sp@ncnt>\@ne 
     \sp@ncnt@=\sp@ncnt\relax\fi
  \findarrdp \@rrwd=\z@ %&0211 It might be used by \locate
  \kern\c@l
   \hbox to \z@{\options
   \@rrdp=\scalefactor\@rrdp
    \hss\vplace{#1}\llabel{#2}\rlabel{#3}\hss}%
  \global\advance\col by \@ne \setsp@n\@ne\@ne
  }

\def\novarrow{\varrow\vfill{}{}}

\def\tweenarrows#1{\findarrwd \findarrdp \setsp@n\@ne\@ne
  \rlap{\options\cplace{#1}}}

\def\usarrow #1#2#3{\dr@tn=\@ne
  \findarrwd \findarrdp \setsp@n\@ne\@ne 
  \rlap{\options\cplace{#1}\nwlabel{#2}\selabel{#3}}%
  \dr@tn=\z@}
\def\dsarrow #1#2#3{\dr@tn=\tw@
  \findarrwd \findarrdp \setsp@n\@ne\@ne 
  \rlap{\options\cplace{#1}\swlabel{#2}\nelabel{#3}}%
  \dr@tn=\z@}
 \def\@rrow#1{\csname #1@rrow\endcsname}
 \def\R@rrow{\harrow \rtarrfill}
 \def\L@rrow{\harrow \ltarrfill}
 \def\V@rrow{\varrow \dnarrfill}
 \def\A@rrow{\varrow \uparrfill}
 \def\SE@rrow{\dsarrow \SEarrow}
 \def\NW@rrow{\dsarrow \NWarrow}
 \def\SW@rrow{\usarrow \SWarrow}
 \def\NE@rrow{\usarrow \NEarrow}
 \def\DS@rrow{\dsarrow \dnslope}
 \def\US@rrow{\usarrow \upslope}
 \def\upslope{\find@xyargs
       \@linelen=\unitlength\@sline}
 \def\dnslope{\find@xyargs\@yarg=-\@yarg\relax
       \@linelen=\unitlength\@sline}

\newtoks\optionlist 
\optionlist={}
\let\options\relax
\def\dooptions{\the\optionlist\global\optionlist={}%
  \glet\options=\relax}
\def\option#1{\ifinmeasureCD\else
  \glet\options=\dooptions
  \global\optionlist=\expandafter{\the\optionlist\relax#1}\fi}
\def\wider#1{\ifinmeasureCD\else
   \option{\advance\@rrwd by #1}\fi}
\def\deeper#1{\ifinmeasureCD\else
   \option{\advance\@rrdp by #1}\fi}
{\def\\{\global\let\sptoken= }\\ }%now \sptoken is a spacetoken

\def\CR{\futurelet\nexttok\testCR}
\def\testCR{\ifx\nexttok\sptoken
   \let\next\eatspaceCR\else\let\next\CDCR\fi\next}
\def\eatspaceCR#1 {\CR}
\def\measureCR{\ifx\nexttok\endmeasure\let\nextCR\relax
    \else\let\nextCR\drop\fi\nextCR}
\def\setCR{\ifodd\row
  \ifx\nexttok\endsetCD\else\hend\betweenCDrows\vbegin\fi
  \else\vend\betweenCDrows\hbegin\fi}

\countdef\dimenc@unt=11
\def\CD#1\endCD{%CRAL
   \begingroup\let\\=\CR
  \m@th\offinterlineskip
   \measure#1\endmeasure\null\,\vcenter{\setCD#1\endsetCD}\,
   \endgroup
    }

\ifx\@clnwd\undefined \nol@g\else\catcode`\ =14\relax\fi
 \font\@linefnt=line10 
 \newcount\@tempcnta
 \newcount\@tempcntb
 \newdimen\@tempdima
 \newdimen\@tempdimb
 \newdimen\@wholewidth
 \newdimen\@halfwidth
   \@wholewidth\fontdimen8\@linefnt \@halfwidth .5\@wholewidth
 \newdimen\unitlength
 \newcount\@xarg
 \newcount\@yarg
 \newcount\@yyarg
 \newbox\@linechar
 \newdimen\@linelen
 \newdimen\@clnwd
 \newdimen\@clnht
 \newif\if@negarg
 
 \def\@whilenoop#1{}

 \def\@whiledim#1\do #2{\ifdim #1\relax#2\@iwhiledim{#1\relax#2}\fi}

 \def\@iwhiledim#1{\ifdim #1\let\@nextwhile=\@iwhiledim 
         \else\let\@nextwhile=\@whilenoop\fi\@nextwhile{#1}}

 \def\@sline{\ifnum\@xarg< 0 \@negargtrue \@xarg -\@xarg \@yyarg -\@yarg
   \else \@negargfalse \@yyarg \@yarg \fi
 \ifnum \@yyarg >0 \@tempcnta\@yyarg \else \@tempcnta -\@yyarg \fi
 \ifnum\@tempcnta>6 \@badlinearg\@tempcnta0 \fi
 \ifnum\@xarg>6 \@badlinearg\@xarg 1 \fi
 \setbox\@linechar\hbox{\@linefnt\@getlinechar(\@xarg,\@yyarg)}%
 \ifnum \@yarg >0 \let\@upordown\raise \@clnht\z@
    \else\let\@upordown\lower \@clnht \ht\@linechar\fi
 \@clnwd=\wd\@linechar
 \if@negarg \hskip -\wd\@linechar \def\@tempa{\hskip -2\wd\@linechar}\else
      \let\@tempa\relax \fi
 \@whiledim \@clnwd <\@linelen \do
   {\@upordown\@clnht\copy\@linechar
    \@tempa
    \advance\@clnht \ht\@linechar
    \advance\@clnwd \wd\@linechar}%
 \advance\@clnht -\ht\@linechar
 \advance\@clnwd -\wd\@linechar
 \@tempdima\@linelen\advance\@tempdima -\@clnwd
 \@tempdimb\@tempdima\advance\@tempdimb -\wd\@linechar
 \if@negarg \hskip -\@tempdimb \else \hskip \@tempdimb \fi
 \multiply\@tempdima \@m
 \@tempcnta \@tempdima \@tempdima \wd\@linechar \divide\@tempcnta \@tempdima
 \@tempdima \ht\@linechar \multiply\@tempdima \@tempcnta
 \divide\@tempdima \@m
 \advance\@clnht \@tempdima
 \ifdim \@linelen <\wd\@linechar
    \hskip \wd\@linechar
   \else\@upordown\@clnht\copy\@linechar\fi}
 
 \def\@getlinechar(#1,#2){\@tempcnta#1\relax\multiply\@tempcnta 8
 \advance\@tempcnta -9 \ifnum #2>0 \advance\@tempcnta #2\relax\else
 \advance\@tempcnta -#2\relax\advance\@tempcnta 64 \fi
 \char\@tempcnta}
 
 \def\vector(#1,#2)#3{\@xarg #1\relax \@yarg #2\relax
 \@tempcnta \ifnum\@xarg<0 -\@xarg\else\@xarg\fi
 \ifnum\@tempcnta<5\relax
 \@linelen=#3\unitlength
 \ifnum\@xarg =0 \@vvector 
   \else \ifnum\@yarg =0 \@hvector \else \@svector\fi
 \fi
 \else\@badlinearg\fi}
 
 \def\@svector{\@sline
 \@tempcnta\@yarg \ifnum\@tempcnta <0 \@tempcnta=-\@tempcnta\fi
 \ifnum\@tempcnta <5
   \hskip -\wd\@linechar
   \@upordown\@clnht \hbox{\@linefnt  \if@negarg 
   \@getlarrow(\@xarg,\@yyarg) \else \@getrarrow(\@xarg,\@yyarg) \fi}%
 \else\@badlinearg\fi}
 
 \def\@getlarrow(#1,#2){\ifnum #2 =\z@ \@tempcnta='33\else
 \@tempcnta=#1\relax\multiply\@tempcnta \sixt@@n \advance\@tempcnta
 -9 \@tempcntb=#2\relax\multiply\@tempcntb \tw@
 \ifnum \@tempcntb >0 \advance\@tempcnta \@tempcntb\relax
 \else\advance\@tempcnta -\@tempcntb\advance\@tempcnta 64
 \fi\fi\char\@tempcnta}
 
 \def\@getrarrow(#1,#2){\@tempcntb=#2\relax
 \ifnum\@tempcntb < 0 \@tempcntb=-\@tempcntb\relax\fi
 \ifcase \@tempcntb\relax \@tempcnta='55 \or 
 \ifnum #1<3 \@tempcnta=#1\relax\multiply\@tempcnta
 24 \advance\@tempcnta -6 \else \ifnum #1=3 \@tempcnta=49
 \else\@tempcnta=58 \fi\fi\or 
 \ifnum #1<3 \@tempcnta=#1\relax\multiply\@tempcnta
 24 \advance\@tempcnta -3 \else \@tempcnta=51\fi\or 
 \@tempcnta=#1\relax\multiply\@tempcnta
 \sixt@@n \advance\@tempcnta -\tw@ \else
 \@tempcnta=#1\relax\multiply\@tempcnta
 \sixt@@n \advance\@tempcnta 7 \fi\ifnum #2<0 \advance\@tempcnta 64 \fi
 \char\@tempcnta}
\catcode`\ =10

\dol@g %925
\catcode`\@=\active
\LaTeXfeatures

\LaTeXfeatures
\DeclareMathOperator{\Fix}{Fix}

\DeclareMathOperator{\GCD}{G.C.D.}

\DeclareMathOperator{\Hom}{Hom}

\DeclareMathOperator{\id}{id}
\DeclareMathOperator{\Img}{Im}
\DeclareMathOperator{\Inv}{Inv}

\DeclareMathOperator{\Ker}{Ker}
\DeclareMathOperator{\LCM}{L.C.M.}

\DeclareMathOperator{\modl}{mod}
\DeclareMathOperator{\ord}{ord}

\DeclareMathOperator{\supp}{supp}

\font\scaps=cmcsc9
\font\caps=cmcsc10
\font\bcaps=cmcsc10 scaled \magstep1

\font\bfit=cmbxti10 scaled \magstep1

\font\twelvegtc=eufm10 scaled 1200

\font\ninegtc=eufm9
\font\sevengtc=eufm7

\newfam\gtcfam

\textfont\gtcfam=\twelvegtc
\scriptfont\gtcfam=\ninegtc
\scriptscriptfont\gtcfam=\sevengtc

\def\bold{\mathbf}

\def\modo#1{\left| #1 \right|}

\def\eqnum#1{\eqno (#1)}
\begin{document}

\title{Braids and Permutations}

\author{Vladimir Lin} 

\begin{abstract}
{E. Artin described all irreducible representations
of the braid group $B_k$ to the symmetric group ${\mathbf S}(k)$.
We strengthen some of his results and, moreover, exhibit a complete picture
of homomorphisms $\psi\colon B_k\to{\mathbf S}(n)$
for $n\le 2k$. We show that the image $\Img\psi$ of $\psi$
is cyclic whenever either \ $(*)$ $n<k\ne 4$ or \ $(**)$ $\psi$
is irreducible and $6<k<n<2k$. For $k>6$ there exist, 
up to conjugation, exactly $3$ irreducible representations
$B_k\to{\mathbf S}(2k)$ with non-cyclic images but
they all are imprimitive. We use these results to prove 
that for $n<k\ne 4$ the image of any braid homomorphism
$\varphi\colon B_k\to B_n$ is cyclic, whereas any endomorphism
$\varphi$ of $B_k$ with non-cyclic image preserves the pure braid group
$PB_k$. We prove also that for $k>4$ the intersection $PB_k\cap B'_k$
of $PB_k$ with the commutator subgroup $B'_k=[B_k,\, B_k]$
is a completely characteristic subgroup of $B'_k$.}
\end{abstract}

\maketitle

\vskip3cm

\thanks{
\noindent Some results of this paper were obtained
and the draft version \cite{Lin96b} was written
during my stay at the Max-Planck-Institut f\"ur Mathematik in Bonn
in 1996. I am deeply grateful to MPI for hospitality.}

\frontmatter

%\centerline{\bcaps preface}%\twotnbf
%\bigskip

%%%%%%%%%%%%%%%%%%%%
\tableofcontents

\mainmatter

\markboth{%\textsc 
\bcaps {Vladimir Lin}}
{\textsc {Braids and permutations}}

%%%%%%%%%%%%%%%%%%%%%%%%%%%%%%%%%%%%%%%%%%%%%%%%%%%%%%%%%%%%%%%%%%%%%%%%%%%

\newpage

%Sec. 1
\section{Introduction}\label{sec: Introduction} 

In the middle $1970$'s I announced some results on braid
homomorphisms and representations of braids by permutations 
\cite{Lin72a,Lin74,Lin79}.
The proofs were never published for the following two reasons.
First, certain proofs were based on a straightforward
modification of Artin's methods\footnote{Perhaps with the exception of
with exception of a systematical exploitation of the fact that for $k>4$
the commutator subgroup of the braid group $B_k$ is a perfect group
\cite{GorLin69}, which seemingly was not noticed by Artin and hence
may be considered as a rather new element.} \cite{Art47b}.
On the other hand, some other proofs contained very long 
combinatorial computations and looked just ugly;
I always felt that they might be improved.
Step by step I found that this in fact can be done by
a new approach, which involves certain cohomology of braid groups.
This approach leads also to stronger results. The new proofs are still lengthy
but involve less combinatorics and seem more suitable for publication.

My interest to braid homomorphisms and representains of braids by permutations
was motivated by the fact that these things are closely related
to algebraic equations with function coefficients, algebraic functions
(in particular, to the 13th Hilbert problem for such functions),
configuration spaces and their holomorphic maps etc. Some information
about such connections may be found in the papers
\cite{Arn70a,Arn70b,Arn70c}, \cite{GorLin69,GorLin74},
\cite{Kal75,Kal76,Kal93}, \cite{Kurth97},
\cite{Lin71,Lin72b,Lin79,Lin96a,Lin96b,Lin03,Lin04a}.
One can also think (and I would certainly be on his side!)
that no motivation is needed at all,
since both braids and permutations are in any ``short list" of
important mathematical notions.

\subsection {Notation and some definitions} 
\label{Ss: Notation and definitions}

\noindent For the reader convenience we start with some notation
and definitions, which are used throughout the paper.\footnote{The main
results of the paper are stated in Sec. \ref{Ss: Main results} 
below; to understand their formulations it suffices to take a look at
Definition \ref{Def: sets, groups, homomorphisms}{\bf (f)}
({\em abelian, cyclic} and {\em integral} \ homomorphisms),
Definition \ref{Def: symmetric groups}{\bf (c)} 
({\em transitive} and {\em primitive}
homomorphisms), Sec. \ref{Ss: Pure braid group} 
(the {\em pure braid group}),
and Sec. \ref{Ss: Commutator subgroup B'(k)}{\bf (c)}
(the {\em pure commutator subgroup}).} 

\begin{Definition}%[{\caps sets, groups, homomorphisms}] 
\label{Def: sets, groups, homomorphisms}
\noindent {\bf a)} {\bfit General notation.} 
$\#\,\Gamma$ is the cardinality of a set $\Gamma$ and
$\ord g$ denotes the order of an element $g\ne 1$ of a group $G$.
We write $g\sim h$ whenever the elements $g,h\in G$ are conjugate.

${\mathbb F}_m$ denotes a {\sl free group of rank} $m$
($m\in{\mathbb Z}_+$ or $m=\infty = \# \mathbb N$); in particular,
${\mathbb F}_1\cong{\mathbb Z}$; we write simply ${\mathbb F}$
if a specific value of $m$ does not matter.
\index{${\mathbb F}_m$;\, $\mathbb F$\hfill}
\vskip0.2cm

{\bf b)} {\bfit Commutator subgroup; perfect groups;
residually finite groups.} 
\index{Group!perfect\hfill}
A group $G$ is called {\em perfect} if it
coincides with its commutator subgroup $G'=[G,G]$. 
A quotient group of a perfect group is perfect;
a perfect group does not possess non-trivial homomorphisms
to any abelian group.

A group $G$ is called {\em residually finite}
\index{Group!residually finite\hfill}
if homomorphisms to finite groups separate elements of $G$,
or, which is the same, the intersection of all finite index
subgroups $H\subseteq G$ is $\{1\}$. Any free group is residually finite.
The following theorem is due to A. I. Malcev \cite{Mal40}:
\vskip0.3cm

\noindent{\bcaps Malcev Theorem.}
\index{Malcev Theorem\hfill}
\index{Theorem!of Malcev\hfill}
{\sl Any semidirect product of finitely generated residually finite groups is a residually 
finite group.}
\vskip0.3cm

{\bf c)} {\bfit Hopfian groups.}
A group $G$ is called {\em Hopfian}
\index{Group!Hopfian\hfill}
if any epimorphism $G\to G$ is an automorphism.
{\sl Any finitely generated residually finite group is
Hopfian} (see, for instance, \cite{HNeum}, Sec. 41.44, p. 151).
\vskip0.2cm

{\bf d)} {\bfit Braid-like couples.}
\index{Braid-like couples\hfill}
A pair of elements $g,h$ of a group $G$ is called {\em braid-like}
if $gh\ne hg$ and $ghg = hgh$; in this case we write $g\infty h$.
Clearly, $g\infty h$ implies $g\sim h$.
\vskip0.2cm

{\bf e)} {\bfit Conjugate homomorphisms.}
\index{Conjugate homomorphisms\hfill}
Group homomorphisms $\phi,\psi\colon G\to H$ are {\em conjugate}
if there is $h\in H$ such that
$\psi(g)=h\phi(g)h^{-1}$ for all $g\in G$; in this case
we write $\phi\sim\psi$; ``$\sim$" is an equivalence
relation in $\Hom(G,H)$.
\vskip0.2cm

{\bf f)} {\bfit Abelian, cyclic and integral homomorphisms.}
\index{Homomorphism!abelian\hfill}
\index{Homomorphism!cyclic\hfill}
\index{Homomorphism!integral\hfill}
A group homomorphism $\phi\colon G\to H$ is said to be
{\em abelian}, {\em cyclic} or {\em integral} respectively
if its image $\Img\phi = \phi (G)$ is an abelian,
cyclic or torsion free cyclic subgroup of $H$.
We include the trivial homomorphism
in each of these three classes. 
\hfill $\bigcirc$
\end{Definition}

\begin{Remark}\label{Rmk: If G' is finitely generated then
any homomorphism G to F is integral} 
{\sl If the commutator subgroup $G'$ of a group $G$
is finitely generated then any homomorphism $\phi\colon G\to{\mathbb F}$
to a free group ${\mathbb F}$ is integral.}
Indeed, $H=\phi(G)\subseteq{\mathbb F}$
is a free group of some rank $r$; the image
$\phi(G')=H'\cong({\mathbb F}_{r})'$ must be finitely generated,
which cannot happen unless $r\le 1$.
\hfill $\square$
\end{Remark}

\begin{Definition}[{\caps symmetric groups}]
\label{Def: symmetric groups}
{\bf a)} {\bfit Invariant sets; fixed points.}
All permutations of a set $\Gamma$ form the
{\em symmetric group} ${\mathbf S}(\Gamma)$, which is
regarded as acting on $\Gamma$ from the left hand side.

For $H\subseteq{\mathbf S}(\Gamma)$ a subset
$\Sigma\subseteq\Gamma$ is $H$-{\em invariant} if
$S(\Sigma)=\Sigma$ for all $S\in H$; we denote by $\Inv H$ the
family of all {\sl non-trivial} (i. e., $\ne\varnothing$ and
$\ne\Gamma$) $H$-invariant subsets of $\Gamma$;
\index{$\Inv H$; \ $\Inv_r H$\hfill}
$\Inv_r H$ consists of all $\Sigma\in\Inv H$ such that
$\#\Sigma=r$, $1\le r<\#\Gamma$ (if $H$ consists of a
single $S\in{\mathbf S}(\Gamma)$ we write $\Inv S$ and $\Inv_r S$
instead of $\Inv \{S\}$ and $\Inv_r \{S\}$ respectively); the restriction
$S|\Sigma$ of $S$ to a set $\Sigma\in\Inv S$ is
regarded as an element of ${\mathbf S}(\Sigma)$.
$\Fix S$ denotes the set $\{\gamma\in\Gamma|\, S(\gamma)=\gamma\}$
of all fixed points of $S$; the complement $\supp S=\Gamma\setminus\Fix S$
is called the {\em support} of $S$. For $\Sigma\subseteq\Gamma$ we
identify ${\mathbf S}(\Sigma)$ with the subgroup 
$\{S\in{\mathbf S}(\Gamma)|\,\supp S\subseteq\Sigma\}
\subseteq{\mathbf S}(\Gamma)$. $S,S'\in{\mathbf S}(\Gamma)$ are
{\em disjoint} if $\supp S\cap\supp S'=\varnothing$.

${\mathbf S}(n)$ denotes the symmetric
group of degree $n$, that is, the permutation group
${\mathbf S}({\boldsymbol\Delta}_n)$ of the $n$ point set
${\boldsymbol\Delta}_n=\{1,2,\ldots ,n\}$.
The {\em alternating subgroup} ${\mathbf A}(n)\subset{\mathbf S}(n)$
consists of all even permutations $S\in{\mathbf S}(n)$ and coincides
with the commutator subgroup ${\mathbf S}'(n)$; for $n>4$ the group
${\mathbf A}(n)$ is perfect.
\vskip0.2cm

{\bf b)} {\bfit Cyclic types, r-components.}
\index{Cyclic type\hfill}
For $A,B\in{\mathbf S}(n)$ we write $A\preccurlyeq B$
\index{$\preccurlyeq$\hfill}
whenever the cyclic decomposition of $B$ contains all cycles
that occur in the cyclic decomposition of $A$.
Let $A = C_1\cdots C_q$ be the cyclic
decomposition of $A\in{\mathbf S}(n)$ and
$r_i\ge 2$ be the length of the cycle
$C_i$ \ ($1\le i\le q$); the {\sl unordered} $q$-tuple of
the natural numbers $[r_1,\ldots ,r_q]$ is called
the {\em cyclic type} of $A$ and is denoted by $[A]$
(each $r_i$ occurs in $[A]$ as many times
as many $r_i$-cycles occur in the cyclic decomposition of $A$).
Clearly, $\ord A = \LCM(r_1,\ldots r_q)$
(the least common multiple of $r_1,...,r_q$).

For $A\in{\mathbf S}(n)$ and a natural $r\ge 2$ we denote by
\index{r-component of a permutation\hfill}
\index{${\mathfrak C}_r(A)$\hfill}
${\mathfrak C}_r(A)$ the set of all $r$-cycles that occur in the cyclic
decomposition of $A$; we call this set the $r$-{\em component}
of $A$. The set $\Fix A$ is called the {\em degenerate component}
of $A$.
\vskip0.2cm

{\bf c)} {\bfit Transitive and primitive homomorphisms.}
\index{Homomorphism!transitive\hfill}
\index{Homomorphism!primitive\hfill}
A group homomorphism $\psi\colon G\to{\mathbf S}(n)$ is said to
be {\em transitive} (respectively {\em intransitive,
primitive, imprimitive}), if its image $\psi (G)$ is
a transitive (respectively intransitive, primitive,
imprimitive) subgroup of the symmetric group ${\mathbf S}(n)$.
\vskip0.2cm

{\bf d)} {\bfit Disjoint products. Reductions of homomorphisms
$G\to{\mathbf S}(n)$.}
\index{Disjoint product\hfill}
Given a disjoint decomposition 
${\boldsymbol\Delta}_n=D_1\cup\cdots\cup D_q$, \ $\# D_j=n_j$,
we have the corresponding embedding
${\mathbf S}(D_1)\times\cdots\times{\mathbf S}(D_q)
\hookrightarrow{\mathbf S}(n)$.
For a family of group homomorphisms
$\psi_j\colon G\to{\mathbf S}(D_j)\cong{\mathbf S}(n_j)$, \
$j=1,...,q$, we define the {\em disjoint product}
$$
\psi=\psi_1\times\cdots\times\psi_q\colon G\to
{\mathbf S}(D_1)\times\cdots\times{\mathbf S}(D_q)
\hookrightarrow{\mathbf S}(n)
$$
by $\psi(g)=\psi_1(g)\cdots\psi_q(g)\in{\mathbf S}(n)$, \ $g\in G$.

Let $\psi\colon G\to{\mathbf S}(n)$ be a group homomorphism, $H = \Img\psi$
and $\Sigma\in\Inv H$ (for instance, $\Sigma$ may be an $H$-orbit).
Define the homomorphisms
$$
\phi_\Sigma\colon H\ni S\mapsto S|\Sigma\in{\mathbf S}(\Sigma)
\quad\text{and}\quad
\psi_{\Sigma}=\phi_\Sigma\circ\psi\colon G\stackrel{\psi}{\longrightarrow} H 
\stackrel{\phi_\Sigma}{\longrightarrow} {\mathbf S}(\Sigma)\,;
$$
the composition $\psi_{\Sigma}=\phi_\Sigma\circ\psi$
\index{Reduction of homomorphism\hfill}
is called the {\em reduction of $\psi$ to $\Sigma$}; it is
transitive if and only if $\Sigma$ is
an $(\Img\psi)$-orbit. Any homomorphism $\psi$ is the disjoint
product of its reductions to all $(\Img\psi)$-orbits
(this is just the decomposition of the representation $\psi$
into the direct sum of irreducible representations).
The following simple observation is used throughout the paper:

\begin{Observation} {\sl A group homomorphism
$\psi\colon G\to{\mathbf S}(n)$ is abelian if and only if
all its reductions $\psi_\Sigma$ to $(\Img \psi)$-orbits
$\Sigma\subseteq{\boldsymbol\Delta}_n=\{1,...,n\}$ are abelian.}
\hfill $\square$
\end{Observation}

\end{Definition}

\subsection{Canonical presentation of the braid group $B_k$}
\label{Canonical presentation of B(k)}
Recall that the {\em canonical presentation}
\index{Canonical presentation of $B_k$\hfill}
of the Artin braid group $B_k$ on $k$ strings
involves $k-1$ {\em ''canonical"} generators $\sigma_{1},...,\sigma_{k-1}$
and the defining system of relations
%&
\begin{eqnarray}
&&\sigma_i\sigma_j
=\sigma_j\sigma_i\qquad\qquad\qquad \  (\modo{i-j}\ge 2)\,,
                            \label{eq: commut}\\
&&\sigma_i\sigma_{i+1}\sigma_i
=\sigma_{i+1}\sigma_i\sigma_{i+1}\qquad (1\le i\le k-2)\,.
                                             \label{eq: braid like}
\end{eqnarray}
%&

\subsection{Torsion}\label{Ss: Torsion}
\index{Torsion\hfill}
The next theorem was proved first by E. Fadell and
L. Neuwirth \cite{FadNeu62} by a topological argument based on
the existence of a finite dimensional Eilenberg-MacLane $K(B_k,1)$ space
(for instance, the {\em configuration space}
${\mathcal C}^k({\mathbb C})$ consisting of all $k$ point subsets
$Q\subset{\mathbb C}$ is such a space);
an algebraic proof was suggested by Joan L. Dyer \cite{Dye80}.
\vskip0.3cm

\noindent{\bcaps Fadell-Neuwirth Theorem.}
\index{Fadell-Neuwirth Theorem\hfill}
\index{Theorem!of Fadell-Neuwirth\hfill}
{\sl The Artin braid group $B_k$ is torsion free.}
\hfill $\square$
\vskip0.3cm

\noindent It follows from (\ref{eq: commut}), (\ref{eq: braid like}) that
$B_k/ B'_k\cong{\mathbb Z}$.
This fact and Fadell-Neuwirth Theorem imply: 
\vskip0.2cm

\begin{Corollary}\label{Crl: abelian and cyclic braid homomorphisms}
\noindent Any abelian homomorphism $\phi$ of $B_k$ is cyclic
and satisfies
$\phi(\sigma _{1})=\phi(\sigma _{2})=\cdots =\phi(\sigma _{k-1})$.
Any abelian homomorphism $B_k\to B_n$ is integral.
\hfill $\square$
\end{Corollary}

\subsection{Special presentation of $B_k$}
\label{Ss: Special presentation of B(k)}
For $1\le i<j\le k$, we put
%&
\begin{equation}\label{alpha beta}
\alpha_{ij} = \sigma_i\sigma_{i+1}\cdots \sigma_{j-1}\,, \ \
\beta_{ij} = \alpha_{ij}\sigma_i,\qquad
\alpha = \alpha_{1k} = \sigma_1\sigma_2\cdots \sigma_{k-1}\,, \ \
\beta = \alpha \sigma_1\,.
\end{equation}
%&
It is easily seen that
%&
\begin{eqnarray}
\sigma_{i+1} &=& \alpha\sigma_i\alpha^{-1}
\qquad\qquad\qquad \ \  
(1 \le  i \le  k-2)\,,\label{conjugation by alpha}\\
\sigma_i \ \ \ &=& \alpha^{i-1}\sigma_1\alpha^{-(i-1)}
\qquad\qquad(1\le i\le k-1)\,, \label{conjugation by alpha power}
\end{eqnarray}
%&
%&
\begin{equation}\label{conjugation by alphaij}
\hskip20pt\left\{
\aligned
\alpha_{ij}\sigma_m &
                   =\sigma_m\alpha_{ij}{\phantom{\sigma_{m+1}}}
\hskip35pt \text {for} \ \ m < i-1 \ \ \ \text {or} \ \ m > j\,,\\
\alpha_{ij}\sigma_m &=\sigma_{m+1}\alpha_{ij}{\phantom{\sigma_m}}
\hskip35pt \text {for} \ \ i\le m\le j-2\,,\\ 
\alpha_{ij}^q\sigma_m &= \sigma_{m+q}\alpha_{ij}^q{\phantom{\sigma_m}}
\hskip35pt \text {for} \ \ i\le m\le m+q\le j-1\,.
\endaligned \right.
\end{equation}
%&
\vskip0.1cm

\noindent 
Relations (\ref{conjugation by alpha}), (\ref{conjugation by alphaij}) 
imply that for $1\le q\le j-i$
$$
\beta_{ij}^q=(\alpha_{ij}\sigma_i)\cdot
\underbrace{(\alpha_{ij}\sigma_i)
     \cdots (\alpha_{ij}\sigma_i)}_{q-1 \ \text{times}}
=\alpha_{ij}\sigma_i\cdot
\sigma_{i+1}\cdots \sigma_{i+q-1}\alpha_{ij}^{q-1}\,;
$$
for $q=j-i$ this shows that
$$
\beta_{ij}^{j-i}=\alpha_{ij}\sigma_i
\sigma_{i+1}\cdots \sigma_{j-1}\alpha_{ij}^{j-i-1}=
\alpha_{ij}^{j-i+1}\,.
$$ 
Moreover, for $m=i$ relations (\ref{conjugation by alphaij}) 
may be written as
$$
\sigma_{i+q} = \alpha_{ij}^q \sigma_i \alpha_{ij}^{-q}
= \alpha_{ij}^{q-1} \beta_{ij} \alpha_{ij}^{-q}\quad
(i\le i+q\le j-1)\,.
$$ 
Therefore 
%&
\begin{equation}
\label{alpha(ij)pow(j-i+1)=beta(ij)pow(j-i)}
\aligned
\alpha_{ij}^{j-i+1} &= \beta_{ij}^{j-i}
\qquad\qquad\qquad &&\text {for} \ \
1\le i<j\le k\,,\\
\sigma_{i+q} \ \ \ &= \alpha_{ij}^{q-1} \beta_{ij} \alpha_{ij}^{-q}
&&\text {for} \ \ 0\le q\le j-i-1\,.
\endaligned
\end{equation}
%&
In particular, these relations show that
{\sl the element $\alpha_{ij}^{j-i+1} = \beta_{ij}^{j-i}$
commutes with all elements} $\sigma_i,\sigma_{i+1},...,\sigma_{j-1}$.

For $i=1$ and $j=k$ relations (\ref{alpha(ij)pow(j-i+1)=beta(ij)pow(j-i)})
take the form
%&
\begin{equation}\label{alpha pow k=beta pow(k-1)}
\aligned
\alpha^k \ \ \ &= \beta^{k-1}\,,\\
\sigma_{1+q} &= \alpha^{q-1} \beta \alpha^{-q}
\qquad \text {for} \ \ 0\le q\le k-2\,,
\endaligned
\end{equation}
%&
which shows that {\sl the elements $\alpha$, $\beta$
generate the whole group $B_k$, and the element
$\alpha^k = \beta^{k-1}$ is central in $B_k$.}
The defining system of relations for the generators
$\alpha$, $\beta $ is as follows:
%&
\begin{eqnarray}
\beta\alpha^{i-1}\beta &=&
\alpha^{i}\beta\alpha^{-(i+1)}\beta\alpha^{i}
\qquad (2\le i\le [k/2])\,,\label{special presentation1}\\ 
\alpha^{k} \qquad &=& \beta^{k-1} \label{special presentation2}
\end{eqnarray}
%&
(see, for instance, \cite{Kle}).
The presentation of $B_k$ given by 
(\ref{special presentation1}), (\ref{special presentation2})
is called {\em special.}
\index{Special presentation of $B_k$\hfill}
Notice that the elements $\sigma_1,\alpha$ also generate the whole
group $B_k$ since $\beta  = \alpha\sigma_1$.

\begin{Definition}\label{Def: special system of generators}
\index{Special system of generators in $B_k$\hfill}
A pair of elements $a,b\in B_k$ is said to be a
{\em special system of generators} if there exists an automorphism
$\psi$ of $B_k$ such that $\psi(\alpha)=a$ and $\psi(\beta)=b$.
If $\{a,b\}$ is such a pair then the elements
%&
\begin{equation*}
s_{i} = \psi(\sigma_{i}) = a^{i-2}ba^{-(i-1)}\qquad (1\le i\le k-1)
%\eqnum{1.11}
\end{equation*}
%&
also form a system of generators of $B_k$ that satisfy
relations (\ref{eq: commut}), (\ref{eq: braid like}); 
\index{Standard system of generators in $B_k$\hfill}
we call such a system of generators a {\em standard} one.
\hfill $\bigcirc$
\end{Definition}

\subsection{Pure braid group}\label{Ss: Pure braid group}
The {\em canonical epimorphism}
\index{Canonical epimorphism $\mu$\hfill}
$\mu=\mu_k\colon B_k\to{\mathbf S}(k)$
is defined by $\mu(\sigma_i)=(i,i+1)\in{\mathbf S}(k)$
\ ($1\le i\le k-1)$.
Its kernel $\Ker\mu = {PB_k}\subset B_k$
is the normal subgroup generated (as a normal subgroup) by the elements
$\sigma_1^2,\ldots ,\sigma_{k-1}^2$; it 
is called the {\em pure braid group}.
\index{Pure braid group\hfill}
\index{$PB_k$\hfill}

A presentation of the pure braid group 
${PB_k}$ was found first by
W. Burau \cite{Bur32} (see also \cite{Mar45},\cite{Bir74}).
Below we formulate some well-known properties of ${PB_k}$
in the form presented in \cite{Mar45}.

The group ${PB_k}$ is generated by the elements
$s_{i,j}\in B_k$ \ ($1\le i<j\le k$)
defined by the recurrent relations
$$
s_{i,i+1} = \sigma_i^2\qquad {\text {and}}\qquad
s_{i,j+1} = \sigma _{j}s_{i,j}\sigma ^{-1}_{j} \ \
{\text {for}} \ \ 1\le i<j<k.
$$
The elements $s_{i,j}$ are called the {\em canonical
generators} of ${PB_k}$. For $1<r<k$ denote by
$s'_{i,j}$, \ $1\le i<j\le r$, the canonical
generators of $PB_r$. The mapping $\xi_{k,r}$ of the generators
%&
\begin{equation}\label{homomorphism xi}
\aligned
\xi_{k,r}(s_{i,j}) &= 1 \qquad &&\text {if} \ \ 1\le i<j \ \ {\text {and}}
\ \ r<j\le k\,,\\
\xi_{k,r}(s_{i,j}) &= s'_{i,j} &&\text {if} \ \ 1\le i<j\le r\,,
\endaligned
\end{equation}
%&
defines an epimorphism $\xi_{k,r}\colon {PB_k}\to PB_r$.
Its kernel coincides with the subgroup $PB_k^{(r)}\subset {PB_k}$
generated by the elements $s_{i,j}$ with $j>r$.
\vskip0.3cm

\noindent{\bcaps Markov Theorem.} 
\index{Markov Theorem\hfill}
\index{Theorem!of Markov\hfill}
{\sl The normal subgroups
$PB_k^{(r)}\vartriangleleft {PB_k}$ fit into the normal series
$$
\{1\} = PB_k^{(k)}\subset PB_k^{(k-1)}\subset\cdots\subset PB_k^{(2)}
\subset PB_k^{(1)} = PB_k
$$
\index{Normal series in ${PB_k}$\hfill}
\index{$PB_k^{(r)}$\hfill}
such that $PB_k^{(r)}/PB_k^{(r+1)}\cong{\mathbb F}_r$
for $1\le r\le k-1$.}
\hfill $\square$
\vskip0.3cm

\noindent Each group $PB_k^{(r)}$ is finitely generated;
Markov Theorem implies the following two corollaries.

\begin{Corollary}\label{Crl: Perfect group has no homomorphisms to PB(k)}  
\noindent A perfect group does not possess
non-trivial homomorphisms to the pure braid group ${PB_k}$.
\end{Corollary}
%\vskip0.3cm

\begin{proof} 
It suffices to show that every non-trivial subgroup
$H\subseteq{PB_k}$ admits a non-trivial homomorphism to ${\mathbb Z}$.
Clearly $H\subseteq PB_k^{(r)}$ and $H\nsubseteq PB_k^{(r+1)}$
for a certain $r$, \ $1\le r\le k-1$. Projecting $H$ to the quotient group
$PB_k^{(r)}/PB_k^{(r+1)}\cong {\mathbb F}_r$ we
obtain a non-trivial free subgroup 
$\widetilde H\subseteq PB_k^{(r)}/PB_k^{(r+1)}$, which admits non-trivial
homomorphisms to ${\mathbb Z}$; hence $H$ itself
has such homomorphisms.
\end{proof}

\begin{Corollary}\label{Crl: B(k) is residually finite and Hopfian} 
\noindent The group $B_k$ is residually
finite; any finitely generated subgroup of $B_k$ is Hopfian.
\end{Corollary}

\begin{proof} By Markov Theorem, each $PB_k^{(r)}$ is a semidirect
product of the finitely generated groups $PB_k^{(r+1)}$
and ${\mathbb F}_r$. By induction, Malcev Theorem implies
that the group $PB_k$ is residually finite. Any finite index subgroup
$H\subseteq PB_k$ is also a finite index subgroup in $B_k$;
hence $B_k$ is residually finite and every finitely generated
subgroup of $B_k$ is Hopfian.
\end{proof}

\subsection{Center}\label{Ss: Center}
\index{Center\hfill}
Denote by ${CB_k}$ $(k\ge 2)$ the infinite cyclic subgroup
in $B_k$ generated by the element
$A_k = \alpha^k=(\sigma_1\sigma_2\cdots \sigma_{k-1})^k$.
Since $\mu(\alpha)=(1,2)(2,3)\cdots (k-1,k)
=(1,2,...,k)$,
we have $\mu(A_k)=1$; hence $CB_k\subseteq PB_k$.
Clearly, $CB_2 = PB_2$. W.-L. Chow \cite{Chow48} proved
that for $k\ge 3$ the subgroup $CB_k$ coincides with
the center of the braid group $B_k$ (see also \cite{Bohn47}).

\subsection{Transitive homomorphisms $B_k\to{\mathbf S}(k)$}
\label{Ss: Transitive homomorphisms B(k) to S(k)}
Any transitive abelian homomorphism
$\psi\colon B_k\to{\mathbf S}(n)$ is cyclic and conjugate
to the homomorphism $\psi_0$ defined by
$$
\psi_0(\sigma_1) =
\psi_0(\sigma_2) =\ldots
= \psi_0(\sigma_{k-1}) = (1,2,\ldots ,n);
$$
in particular $[\psi(\sigma_i)]=[n]$ for all $i=1,...,k-1$.
The following classical theorem of E. Artin \cite{Art47b} describes all
non-cyclic transitive homomorphisms of the group $B_k$ to
the symmetric group ${\mathbf S}(k)$. (See also Remark 3.2.)
\vskip0.3cm

\noindent{\bcaps Artin Theorem.}
\index{Artin Theorem\hfill|phantom\hfill}
\index{Theorem!of Artin\hfill}
\index{Theorem!on homomorphisms!$B_k\to{\mathbf S}(k)$\hfill}
\index{Homomorphisms!$B_k\to{\mathbf S}(k)$\hfill}
{\sl Let $\psi\colon B_k\to{\mathbf S}(k)$ be a non-cyclic
transitive homomorphism.
\vskip0.2cm

$a)$ If $k\ne 4$ and $k\ne 6$ then $\psi$ is conjugate to the
canonical epimorphism $\mu=\mu_k$.
\vskip0.2cm

$b)$ If $k = 6$ and $\psi\not\sim\mu_6$
then $\psi$ is conjugate to the homomorphism $\nu_6$ 
\index{Homomorphism $\nu_6$\hfill}
\index{Artin homomorphism $\nu_6$\hfill}
\index{$\nu_6$\hfill}
defined by
%&
\begin{equation}\label{Artin homomorphism vu6}
\nu_6(\sigma_1) = (1,2)(3,4)(5,6),
\qquad \nu_6(\alpha) = (1,2,3)(4,5).
\end{equation}
%&

$c)$ If $k = 4$ and $\psi\not\sim\mu_4$ 
then $\psi$ is conjugate to one of the 
three homomorphisms $\nu_{4,1}$, $\nu_{4,2}$ and $\nu_{4,3}$
\index{Homomorphisms!$\nu_{4,1}$, \ $\nu_{4,2}$, \ $\nu_{4,3}$\hfill}
\index{Artin homomorphisms $\nu_{4,1}$, $\nu_{4,2}$, $\nu_{4,3}$\hfill}
\index{$\nu_{4,1}$, $\nu_{4,2}$, $\nu_{4,3}$\hfill}
defined by
%&
\begin{equation}\label{Artin homomorphisms vu41-43}
\aligned
\nu_{4,1}(\sigma_1) &= (1,2,3,4),\qquad \\
\nu_{4,2}(\sigma_1) &= (1,3,2,4),\qquad \\
\nu_{4,3}(\sigma_1) &= (1,2,3), \qquad  \\
\endaligned
\aligned
\nu_{4,1}(\alpha) &= (1,2);\qquad \\
\nu_{4,2}(\alpha) &= (1,2,3,4);\qquad \\
\nu_{4,3}(\alpha) &= (1,2)(3,4);\qquad \\
\endaligned
\aligned
\lbrack\nu_{4,1}(\sigma_3) &= \nu_{4,1}(\sigma_1)\rbrack\\
\lbrack\nu_{4,2}(\sigma_3) &= \nu_{4,2}(\sigma_1^{-1})\rbrack\\
\lbrack\nu_{4,3}(\sigma_3) &= \nu_{4,3}(\sigma_1)\rbrack.\\
\endaligned
\end{equation}
%&

$d)$ $\psi$ is surjective except of the case
when $k=4$, $\psi\sim\nu_{4,3}$ and $\Img\psi={\mathbf A}(4)$.}
\hfill $\bigcirc$

\subsection{Commutator subgroup $B'_k$}
\label{Ss: Commutator subgroup B'(k)}{\bf a)} {\bfit Canonical integral 
projection.} We have already noticed that
$B_k/B'_k\cong{\mathbb Z}$; the homomorphism $\chi\colon B_k\to{\mathbb Z}$
defined by
$$
\chi(\sigma_1) =\ldots = \chi(\sigma_{k-1})=1\in{\mathbb Z}
$$
is called the {\em canonical integral projection} of $B_k$.
Clearly $\Ker\chi=B'_k$.

\begin{Remark}\label{Rmk: kernel of abelian braid homomorphism} 
{\sl If $G$ is a torsion free group and $\phi\colon B_k\to G$ is a non-trivial
abelian homomorphism then $\Ker\phi = B'_k$.}
(Clearly $B'_k\subseteq \Ker\phi$; this
inclusion cannot be strict, for otherwise
$\Img\phi\cong B_k/\Ker\phi$ would be a non-trivial proper quotient group
of the group $B_k/B'_k\cong{\mathbb Z}$, which
is impossible since $G$ is torsion free.)
\hfill $\square$
\end{Remark}

{\bf b)} {\bfit Presentation of $B'_k$.} 
First, $B'_2 = \{1\}$. The following theorem
(see \cite{GorLin69} for the proof) supplies us with a presentation
of the commutator subgroup $B'_k\subset B_k$ 
for $k\ge 3$. This presentation, in turn, implies
certain crucial properties of $B'_k$ (see statements $(c)$,$(d)$
of the theorem and Corollary \ref{Crl: subgroups of B3 and B4}).
\vskip0.3cm

\noindent{\bcaps Gorin-Lin Theorem}
\index{Gorin-Lin Theorem\hfill}
\index{Theorem!of Gorin-Lin\hfill}
$a)$ {\sl $B'_3$ is a free group of rank $2$ with the 
free base\footnote{We write the generators
of $B'_k$ in terms of the canonical generators 
$\sigma_1,...,\sigma_{k-1}$ of $B_k$.}
%&
\begin{equation}\label{generators of B'_3}
u = \sigma_2\sigma_1^{-1}, \ \ \
v = \sigma_1\sigma_2\sigma_1^{-2}.
\end{equation}
%&

$b)$ \index{Presentation of $B_k'$\hfill}
For $k>3$ the group $B'_k$ has a finite
presentation with the generators
%&
\begin{equation}\label{generators of commutator} %@
\aligned
\hskip13pt u &= \sigma_2\sigma_1^{-1}, \\
\hskip13pt w &= \sigma_2\sigma_3\sigma_1^{-1}\sigma_2^{-1},\qquad \\
\endaligned
\aligned
v &= \sigma_1\sigma_2\sigma_1^{-2}, \\
c_{i} &= \sigma_{i+2}\sigma_1^{-1}\qquad (1\le i\le k-3)
\endaligned
%\eqnum{1.13}
\end{equation}
%&
and with the following defining system of relations:
\footnote{For $k=4$ the generators $c_{i}$ with $i \ge 2$
and relations (\ref{1.18}) - (\ref{1.21}) do not appear.}
\begin{eqnarray}
&&uc_1 u^{-1} \hskip7pt = \hskip2pt w\,,\label{1.14}\\
%\eqnum{1.14}
&&uwu^{-1} \hskip7pt = \hskip2pt w^2 c_1^{-1}w\,,\label{1.15}\\
%\eqnum{1.15}
&&vc_1 v^{-1} \hskip7pt = \hskip2pt c_1^{-1}w\,,\label{1.16}\\
%\eqnum{1.16}
&&vwv^{-1} \hskip8pt = \hskip2pt (c_1^{-1} w)^3 c_1^{-2} w\,,\label{1.17}\\
%\eqnum{1.17}
&&uc_i \hskip25pt = \hskip2pt c_i v \hskip80pt (2\le i\le k-3)\,,\label{1.18}\\
%\eqnum{1.18}
&&vc_i \hskip25pt = \hskip2pt c_i u^{-1}v 
\hskip63pt (2\le i\le k-3)\,,\label{1.19}\\
%\eqnum{1.19}
&&c_i c_j \hskip22pt = \hskip2pt c_j c_i
\hskip78pt (1\le i<j-1\le k-4)\,,\label{1.20}\\
%\eqnum{1.20}
&&c_i c_{i+1}c_i \hskip4pt = \hskip2pt c_{i+1}c_i c_{i+1}
\hskip48pt (1 \le  i \le  k-4)\,.\label{1.21}
%\eqnum{1.21}
\end{eqnarray}
%\vskip0.2cm

$c)$ The subgroup $\mathbf T\subset B'_4$
generated by $w=\sigma_2\sigma_3\sigma_1^{-1}\sigma_2^{-1}$ and 
$c_1=\sigma_3\sigma_1^{-1}$ is a free group of rank $2$.
It coincides with the intersection of the lower central series
of the group $B'_4$ and $B'_4/{\mathbf T}\cong{\mathbb F}_2$.
\vskip0.2cm

$d)$ For $k>4$ the group $B'_k$ is perfect.}
\hfill $\square$
\vskip0.3cm

\noindent Since $B_k/B'_k\cong{\mathbb Z}$,
statements $(a)$ and $(c)$ of this theorem imply:

\begin{Corollary}\label{Crl: subgroups of B3 and B4} 
\noindent The groups $B_3$ and $B_4$
admit finite normal series with free quotient groups. Thus,
these braid groups contain no perfect subgroups.
For $k\le 4$ any non-trivial subgroup $G\subseteq B_k$
possesses non-trivial homomorphisms $G\to{\mathbb Z}$.
%\hfill $\square$
\end{Corollary}

\begin{Remark}\label{Gorin's relation} 
In 1967 I found certain identities in $B_k$, $k>6$,
which led to the first proof of the fact that for such $k$
the group $B'_k$ is perfect; during a long time,
the only known proof for smaller $k$ was based
on the above presentation
(\ref{generators of commutator})-(\ref{1.21})
of $B'_k$, which was derived in \cite{GorLin69}.

In the middle 1980's E. A. Gorin discovered the following
beautiful relation, which holds for any $k\ge 4$:
%&
\begin{equation}\label{eq: Gorin relation}
\index{Gorin's relation\hfill}
\sigma_3\sigma_1^{-1}
= (\sigma_1\sigma_2)^{-1}\cdot
\left[\sigma_3\sigma_1^{-1},\sigma_1\sigma_2^{-1}\right]
\cdot (\sigma_1\sigma_2)\,,
\end{equation}
%&
where
$\left[\sigma_3\sigma_1^{-1},\sigma_1\sigma_2^{-1}\right]
=\left(\sigma_3\sigma_1^{-1}\right)^{-1}
\cdot\left(\sigma_1\sigma_2^{-1}\right)^{-1}
\cdot\left(\sigma_3\sigma_1^{-1}\right)
\cdot\left(\sigma_1\sigma_2^{-1}\right)$ is the
commutator of $g_1=\sigma_3\sigma_1^{-1}$ and
$g_2=\sigma_1\sigma_2^{-1}$. Clearly $g_1,g_2\in B'_k$
and (\ref{eq: Gorin relation}) shows that 
$\sigma_3\sigma_1^{-1}\in[B'_k,B'_k]$. Hence
the normal subgroup $N\vartriangleleft B_k$ generated
(as a normal subgroup) by the element $\sigma_3\sigma_1^{-1}$
is contained in $[B'_k,B'_k]$. For $k>4$ it follows readily from
(\ref{eq: commut}) and (\ref{eq: braid like})
that $N$ contains the whole $B'_k$. Indeed, the presentation of $B_k/N$
involves the generators $\sigma_1,...,\sigma_{k-1}$ 
and the defining system of relations consisting of 
(\ref{eq: commut}), (\ref{eq: braid like}) and the additional
relation $\sigma_3\sigma_1^{-1}=1$; since $k>4$ 
relations (\ref{eq: commut}) and $\sigma_3\sigma_1^{-1}=1$ imply 
$\sigma_3\sigma_4=\sigma_4\sigma_3$; in view of 
(\ref{eq: braid like}), this shows that $\sigma_3=\sigma_4$
and finally $\sigma_1=\sigma_2=...=\sigma_{k-1}$.
Thus, $B_k/N\cong{\mathbb Z}$ and $N\supseteq B'_k$.
Thereby $B'_k\subseteq N\subseteq [B'_k,B'_k]$ and $B'_k$ is perfect. 
\hfill $\square$
\end{Remark}

\begin{Remark}\label{Rmk: embeddings lambda(k,m) of B(k-2) to B(k)}
\index{Embeddings $\lambda_{k,m}$ and $\lambda_k'$\hfill}
\index{Embedding $\lambda_k'$\hfill}
\index{$\lambda_{k,m}$\hfill}
\index{$\lambda_k'$\hfill}
Assume that $k\ge 4$ and denote
the canonical generators in $B_{k-2}$ and
$B_k$ by $s_i$ and $\sigma_j$ respectively. Since
$\sigma_1$ commutes with $\sigma_3,...,\sigma_{k-1}$,
for any integer $m$ the mapping
$$
s_i\mapsto\lambda_{k,m}(s_i)\Def\sigma_{i+2}\sigma_1^{-m} \ \ \ 
\text{for} \ \ 1\le i\le k-3
$$
defines a homomorphism 
$$
\lambda_{k,m}\colon B_{k-2}\to B_k\,.
$$
It is well known that $\lambda_{k,m}$ {\sl is an embedding}
(geometrically this is evident).

Furthermore, relations (\ref{1.20}), (\ref{1.21}) 
for the generators
$c_i= \sigma_{i+2}\sigma_1^{-1}$ in $B'_k$ 
show that the mapping
$$
s_i\mapsto\lambda_k'(s_i)\Def c_i \ \ \ 
\text{for} \ \ 1\le i\le k-3
$$
defines a homomorphism
$$
\lambda_k'\colon B_{k-2}\to  B'_k\,.
$$
The composition of $\lambda_k'$ with the natural
embedding $B'_k\hookrightarrow B_k$
coincides with $\lambda_{k,1}$; hence $\lambda_k'$
{\sl is an embedding}. 
We shall use this embedding in
Sec. \ref{Ss: 7.2. Certain homomorphisms of the commutator subgroup B'(k)}.
\hfill $\bigcirc$
\end{Remark}

{\bf c)} {\bfit Canonical homomorphism
$B'_k\to{\mathbf S}(k)$; pure commutator subgroup.}
The presentation of $B'_k$ 
given by Gorin-Lin Theorem is called {\em canonical.}
The restriction $\mu'$ of the canonical epimorphism $\mu$
to the commutator subgroup $B'_k$ of the group $B_k$,
$$
\mu' = \mu|_{B'_k}\colon B'_k\to{\mathbf S}(k),
$$
\index{Canonical homomorphism $B'_k\to{\mathbf S}(k)$\hfill}
is called {\em the canonical homomorphism} of
$B'_k$ to ${\mathbf S}(k)$. Its image
coincides with the alternating subgroup
${\mathbf A}(k)\subset {\mathbf S}(k)$, and its kernel
coincides with the normal subgroup $J_k=PB_k\cap B'_k$
of the group $B_k$. The latter normal subgroup
$$
J_k = PB_k\cap B'_k
$$
\index{Pure commutator subgroup $J_k$\hfill}
\index{$J_k$\hfill}
is called the {\it pure commutator subgroup} of the braid group
$B_k$. It is easily checked that
%&
\begin{equation}\label{1.22}
\aligned
\mu' (u) &= (1,3,2)\,, \quad
\mu' (v) = (1,2,3)\,, \quad
\mu' (w) = (1,3)(2,4)\,, \\
\mu' (c_i) &= (1,2)(i+2,i+3)\qquad (1 \le  i \le  k-1)\,.
\endaligned
%\eqnum{1.22}
\end{equation}
%&

\subsection{Main results}\label{Ss: Main results}
Here we formulate main results of the paper.

\begin{ThmA}\label{Thm A}  
If $k>4$ and $n<k$, then

$a)$ each homomorphism $B_k\to{\mathbf S}(n)$ 
is cyclic;
\index{Theorem!on homomorphisms!$B_k\to{\mathbf S}(n)$ for $n<k$\hfill}

$b)$ each homomorphism $B_k\to B_n$
is integral;
\index{Theorem!on homomorphisms!$B_k\to B_n$ for $n<k$\hfill}

$c)$ the commutator subgroup $B'_k$
of the group $B_k$ does not possess non-trivial homomorphisms
to the groups ${\mathbf S}(n)$ and $B_n$.
\index{Theorem!on homomorphisms!$B_k'\to{\mathbf S}(n)$ for $n<k$\hfill}
\index{Theorem!on homomorphisms!$B_k'\to B_n$ for $n<k$\hfill}
\end{ThmA}

\noindent E. Artin \cite{Art47b} proved that {\sl the pure braid group ${PB_k}$
is a characteristic subgroup of the braid group $B_k$}, that is,
{\sl $\phi(PB_k)=PB_k$ for any automorphism $\phi$ of $B_k$}.
Our next result shows that for $k\ne 4$ the subgroup $PB_k$
possesses in fact much stronger invariance properties.

\begin{ThmB}\label{Thm B}  
\index{Theorem!on endomorphisms of $B_k$ for $k\ne 4$\hfill}
If $k\ne 4$, then $\phi ({PB_k})\subseteq {PB_k}$, \
$\phi^{-1}({PB_k}) = {PB_k}$ \ and $\Ker\phi\subseteq {PB_k}$
for any non-integral endomorphism
$\phi\colon B_k\to B_k$.
\end{ThmB}

\noindent Let $\nu_6'$
\index{Homomorphism $\nu_6'$\hfill}
\index{$\nu_6'$\hfill}
denote the restriction of Artin's homomorphism
$\nu_6\colon B_6\to{\mathbf S}(6)$
to the commutator subgroup $B'_6$ of $B_6$.

\begin{ThmC}\label{Thm C}  
\index{Theorem!on homomorphisms!$B_k'\to{\mathbf S}(k)$\hfill}
Let $k>4$ and let 
$\psi\colon B'_k\to{\mathbf S}(k)$
be a non-trivial homomorphism. Then
either $\psi\sim \mu_k'$ 
{\rm (which may happen for any $k$)} or $k=6$ and
$\psi\sim \nu_6'$. In particular,
$\Img\psi=\mathbf A(k)$ and 
$\Ker\psi=J_k={PB_k}\cap B'_k$.
\end{ThmC}

\begin{ThmD}\label{Thm D}  
\index{Theorem!on pure commutator subgroup $J_k$\hfill}
\index{Pure commutator subgroup $J_k$\hfill}
\index{$J_k$\hfill}
Let $k>4$. Then the pure commutator subgroup $J_k = {PB_k}\cap B'_k$ 
is a completely characteristic subgroup of
the group $B'_k$, that is, $\phi(J_k)\subseteq J_k$
for each endomorphism \ $\phi\colon B'_k\to B'_k$.
Moreover, $\phi^{-1}(J_k) = J_k$ for every non-trivial
endomorphism $\phi$ of the group $B'_k$.
\end{ThmD}

\begin{ThmE}\label{Thm E}  
$a)$ For $k>5$ every transitive homomorphism
$B_k\to{\mathbf S}(k+1)$ is cyclic.

$b)$ For $k>4$ every transitive homomorphism
$B_k\to{\mathbf S}(k+2)$ is cyclic.
\end{ThmE}

\begin{Definition}[{\caps model and standard homomorphisms
$B_k\!\to{\mathbf S}(2k)$}] 
\label{Def: model and standard homomorphisms B(k) to S(2k)}
\index{Model homomorphisms $B_k\to{\mathbf S}(2k)$\hfill}
\index{Standard homomorphisms $B_k\to{\mathbf S}(2k)$\hfill}
\index{Homomorphisms!$B_k\to{\mathbf S}(2k)$!model\hfill}
\index{Homomorphisms!$B_k\to{\mathbf S}(2k)$!standard\hfill}
The following three homomorphisms
$\varphi_1,\varphi_2,\varphi_3\colon B_k\to{\mathbf S}(2k)$
are called the {\em model} ones:
$$
\aligned
&\varphi_1(\sigma_i) = 
\underbrace{(2i-1,2i+2,2i,2i+1)}_{\text {$4$-cycle}};\\
%&\\
&\varphi_2(\sigma_i)=(1,2)\cdot\cdot\cdot (2i-3,2i-2)
\underbrace{(2i-1,2i+1)(2i,2i+2)}_{\text {two transpositions}}\times\\
&\hskip2.9in
\times
(2i+3,2i+4)\cdot\cdot\cdot(2k-1,2k);\\
%\\
&\varphi_3(\sigma_i)= (1,2)\cdot\cdot\cdot (2i-3,2i-2)
\underbrace{(2i-1,2i+2,2i,2i+1)}_{\text {$4$-cycle}}\times\\
&\hskip2.9in
\times(2i+3,2i+4)\cdot\cdot\cdot (2k-1,2k);\\
\endaligned
$$
in each of the above formulae $i=1,...,k-1$.
A homomorphism \ $\psi\colon B_k\to{\mathbf S}(2k)$ \
is said to be {\em standard} if it is conjugate to one of
the three model homomorphisms \ $\varphi_1$, \ 
$\varphi_2$, \ $\varphi_3$. \hfill $\bigcirc$
\end{Definition}

\begin{ThmF}\label{Thm F}
\index{Theorem!on homomorphisms!
$B_k\to{\mathbf S}(n)$ for $6<k<n<2k$\hfill}
$a)$ For $6<k<n<2k$ every transitive
homomorphism $B_k\to{\mathbf S}(n)$ is cyclic.

$b)$ 
\index{Standard homomorphisms $B_k\to{\mathbf S}(2k)$\hfill}
\index{Homomorphisms!$B_k\to{\mathbf S}(2k)$!standard\hfill}
\index{Theorem!on homomorphisms!$B_k\to{\mathbf S}(2k)$ for $k>6$\hfill}
For $k>6$ every non-cyclic transitive homomorphism
$\psi\colon B_k\to{\mathbf S}(2k)$
is standard.\footnote{For $k>8$ this
was proved in \cite{Lin96b}; the cases $k=7$ and $k=8$
were treated by S. Orevkov \cite{Ore98}.}
\end{ThmF}

\begin{Remark}
\label{Rmk: finite overings of Cn(C)}
As we mentioned in Sec. \ref{Ss: Torsion}, the Artin braid group
$B_k$ is the fundamental group of the non-ordered $k$ point configuration
space ${\mathcal C}^k({\mathbb C})$ of the complex affine line ${\mathbb C}$
(in fact ${\mathcal C}^k({\mathbb C})=K(B_k,1)$).
Theorems A and F provide us with a complete description
of homomorphisms $B_k\to{\mathbf S}(n)$ for $n\le 2k$ (at least when $k>6$)
and thereby with a complete description of all unramified
$n$ coverings of the configuration space ${\mathcal C}^k({\mathbb C})$
of degree $n\le 2k$. See \cite{Lin03,Lin04a} for more details.
\hfill $\bigcirc$
\end{Remark}

\begin{ThmG}\label{Thm G}  
Let $k>4$ and $n<2k$. Then:

$a)$ any transitive imprimitive homomorphism
$\psi\colon B_k\to{\mathbf S}(n)$ is cyclic;

$b)$ any transitive homomorphism
$\psi'\colon B'_k\to{\mathbf S}(n)$ is primitive.
\end{ThmG}

\begin{Definition}[{\caps special homomorphisms}] 
\label{Def: special homomorphism}
\index{Special homomorphisms\hfill}
\index{Homomorphism!special\hfill}
Given a special system of generators $\{a,b\}$ in $Bm$,
we denote by ${\mathcal H}_{m}(a,b)$ the set
of all elements of the form $g^{-1}a^p g$ or $g^{-1}b^p g$, 
where $g$ runs over $B_m$ and $p$ runs over ${\mathbb Z}$.
A homomorphism $\phi\colon B_k\to B_n$
is said to be {\em special} if $\phi({\mathcal H}_k(a,b))\subseteq
{\mathcal H}_n(a',b')$ for some choice of special systems of
generators $a,b\in B_k$ and $a' ,b'\in B_n$.
\hfill $\bigcirc$
\end{Definition}

\noindent{\bcaps Murasugi Theorem} (\cite{Mur82}).
\index{Murasugi Theorem\hfill}
\index{Theorem!of Murasugi\hfill}
{\sl A braid $h$ belongs to ${\mathcal H}_m(a,b)$
if and only if $h$ is an element of finite order modulo the
center $CB_m$ of the group $B_m$.}
\vskip0.3cm

\noindent This theorem implies that {\sl the set 
${\mathcal H}_m={\mathcal H}_m(a,b)$ does not depend
on a choice of a special system of generators $a,b\in B_m$, and
a homomorphism $\phi\colon B_k\to B_n$
is special if and only if for any element $g\in B_k$
of finite order modulo $CB_k$ its image
$\phi(g)\in B_n$ is an element of finite order modulo
${CB_n}$} (in fact, we do not use this result in this paper;
it is convenient, however, to keep it in mind).
\vskip0.2cm

\noindent It was stated in \cite{Lin71} (see also \cite{Lin79}
and \cite{Lin03,Lin04a}, Sec. 9, for the proof) that
{\sl for every holomorphic
mapping $f\colon{\mathbf G}_k\to{\mathbf G}_n$ the induced
homomorphism of the fundamental groups
$$
f_*\colon B_k\cong\pi_1({\mathbf G}_k)
\to\pi_1({\mathbf G}_n)\cong B_n
$$ 
is special}. This is a reason to study such homomorphisms.

\begin{Definition-Notation}\label{Def: four progressions} 
\index{Arithmetic progressions $P^i(k)$\hfill}
\index{$P^i(k)$\hfill}
Let $P(k)$ be the union
of the four increasing infinite arithmetic progressions
$$
P^i(k)=\left\{\left. p^i_j(k)=p^i_1(k)+(j-1)d(k)\right|\
j\in\mathbb N\right\}\,,
\quad 1\le i\le 4\,,
$$
with the same difference $d(k) = k(k-1)$ with the initial terms
$$
\quad\qquad 
p^1_1(k)= k, \ \ p^2_1(k)= k(k-1), \ \ p^3_1(k)= k(k-1)+1, \ \
{\text {and}} \ \ p^4_1(k)= (k-1)^{2}
$$
respectively.
\end{Definition-Notation}

\begin{ThmH}\label{ThmH}  
\index{Special homomorphisms\hfill}
\index{Homomorphism!special\hfill}
$a)$ Let $k\ne 4$ and $n\notin P(k)$. 
Then every special homomorphism $B_k\to B_n$ is integral.

$b)$ 
\index{Non-abelian special homomorphisms\hfill}
\index{Special homomorphisms!non-abelian\hfill}
For any $k\ge 3$ and any $n\in P^1(k)\cup P^2(k)$
there exists a non-abelian special homomorphism
$B_k\to B_n$.
\end{ThmH}

\noindent Theorem H$(a)$ has the following immediate 

\begin{Corollary}\label{Crl: non-cyclic special homomorphisms}
If $k\ne 4$ and there is a non-cyclic special homomorphism
$B_k\to B_n$ then $k(k-1)$ divides $n(n-1)$.
\hfill $\bigcirc$
\end{Corollary}

\subsection{Strategy and location of proofs} 
\label{Ss: How and where}
Here I try to describe, step by step,
a strategy of the proofs of main results
and give an information about the proofs' geography.
If the reader is not interested in these rather long preliminary notes,
he may pass to the next section or even to Sec. \ref{sec: B(k) to S(n), n<=k}
\vskip0.2cm

\noindent For any couple $(k,n)$ the set $\Hom(B_k,{\mathbf S}(n))$
is finite (certainly $\#\Hom(B_k,{\mathbf S}(n))\le h(n)=(n!-1)n!$)
and all homomorphisms $B_k\to{\mathbf S}(n)$ may be found by a
straightforward combinatorial computation.
However for large $n$ this way is unrealistic
(say $h(10)\approx 1.3\cdot 10^{15}$) 
\vskip0.3cm

\noindent Any homomorphism $\psi\colon B_k\to{\mathbf S}(n)$
is a disjoint product of transitive
homomorphisms $\psi_j\colon B_k\to{\mathbf S}(n_j)$,
\ $\sum n_j=n$. If $\psi$ is non-cyclic, at least one of the
homomorphisms $\psi_j$ must be non-cyclic. Taking into account this fact,
we try to describe all possible cyclic types of the permutation
$\widehat\sigma_1=\psi(\sigma_1)\in{\mathbf S}(n)$
for a transitive (or non-cyclic and transitive) homomorphism
$\psi\colon B_k\to{\mathbf S}(n)$.
(Notice that all $\widehat\sigma_i=\psi(\sigma_i)$
are conjugate to each other and hence they are of the same cyclic type.)
\vskip0.3cm
%\pagebreak

\noindent{\bfit Fixed points and primes.
Transitive homomorphisms $B_k\to{\mathbf S}(n)$ 
for $n=k$ and $n=k+1$.} 
The following very remarkable result is due to Artin:
\vskip0.3cm

\noindent{\bcaps Artin Fixed Point Lemma.} 
\index{Artin Fixed Point Lemma\hfill|phantom\hfill}
\index{Artin Lemma!on fixed points of $\widehat\sigma_1$\hfill}
\index{Lemma!of Artin on fixed points \hfill}
\index{Fixed points and primes\hfill}
{\sl If $k>4$ and there is a prime $p>2$ such that $n/2 < p\le k-2$,
then for any non-cyclic transitive homomorphism
$\psi\colon B_k\to{\mathbf S}(n)$ the permutation
$\widehat\sigma_1$ has at least $k-2$ fixed points.}
\vskip0.3cm

\noindent Artin treated only the case $n=k$;
however his proof does not in fact depend
on the latter assumption (see Lemma \ref{Lm: Artin Fixed Point Lemma}
below). 

A famous theorem of P. L. Chebyshev
\index{Chebyshev Theorem\hfill}
\index{Theorem!of Chebyshev\hfill}
ensures the existence of a required prime $p$ whenever
eather $k>4$ and $n<k$ or $6\ne k>4$ and $n\le k$.
If $6\ne k=n>4$, Artin Fixed Point Lemma shows
that {\sl all permutations $\widehat\sigma_i$ must be transpositions},
and the whole rest of the proof of Artin Theorem $(a)$ takes a few words.
Moreover, for any $k>6$ there is a prime $p$ such that $(k+1)/2<p\le k-2$.
Hence the inequality $\#\Fix\widehat\sigma_1\ge k-2$ holds true
for any $k>6$ and any non-cyclic transitive homomorphism
$\psi\colon B_k\to{\mathbf S}(k+1)$, which yields
Theorem E$(a)$ whenever $k>6$
\index{Theorem!on homomorphisms $B_k\to{\mathbf S}(k+1)$\hfill} 
\index{Homomorphisms!$B_k\to{\mathbf S}(k+1)$\hfill} 
(in Theorem \ref{Thm: homomorphisms B(k) to S(k+1)}$(a)$
the case $k=6$ is treated as well).
\vskip0.3cm

\noindent{\bfit Homomorphisms $\psi\colon B_k\to{\mathbf S}(n)$, \
$n<k$. An improvement of Artin Theorem.}
\index{Homomorphisms!$B_k\to{\mathbf S}(n)$ for $n<k$\hfill} 
We represent $\psi$ as a disjoint product of transitive
homomorphisms $\psi_j$. If some $\psi_j$ is non-cyclic,
then, by Artin Fixed Point Lemma, all $\psi_j(\sigma_i)$
are transpositions; hence $\psi_j$ cannot be transitive and a contradiction
ensues. This proves Theorem A$(a)$
(see Theorem \ref{Thm: homomorphisms B(k) to S(n), n<k}$(a)$). 
Combining the latter theorem with Artin Theorem we show that
{\sl for $k\ne 4$ any non-surjective homomorphism
$B_k\to{\mathbf S}(k)$ is cyclic} 
(see Lemma \ref{Lm: non-surjective homomorphism B(k) to S(k) is cyclic}).
This implies (see Theorem \ref{Thm: Improvement of Artin Theorem})
that {\sl for $k\ne 4,\,6$ any non-cyclic homomorphism
$B_k\to{\mathbf S}(k)$ is conjugate
to the canonical epimorphism $\mu_k$}, which is a useful improvement
\index{Artin Theorem!improvement of\hfill}
\index{Theorem!of Artin!improvement of\hfill}
\index{Improvement of Artin Theorem\hfill}
of Artin Theorem.
\vskip0.3cm

\noindent{\bfit Homomorphisms $\psi\colon B'_k\to{\mathbf S}(n)$.}
\index{Homomorphisms!$B'_k\to{\mathbf S}(n)$\hfill}
Take the restriction $\phi$ of a homomorphism $\psi$ to the subgroup
$B\cong  B_{k-2}$ of $B'_k$ generated by the elements
$c_i=\sigma_{i+2}\sigma_1^{-1}$, \ $1\le i\le k-3$
\index{Embedding $\lambda_k'$\hfill}
(compare to Remark \ref{Rmk: embeddings lambda(k,m) of B(k-2) to B(k)}).
Using Theorem A$(a)$ and
Lemma \ref{Lm: non-surjective homomorphism B(k) to S(k) is cyclic}
mentioned above, we show that for $k>4$ and $n<k$ the homomorphism $\phi$
must be cyclic. On the other hand,
in Lemma \ref{Lm: composition B(k-2) to Bk' to H} we deduce from relations
(\ref{1.14})-(\ref{1.21}) that the original homomorphism $\psi$ is trivial
whenever $\phi$ is cyclic. This proves the statement of Theorem A$(c)$
concerning homomorphisms $ B'_k\to{\mathbf S}(n)$, which is an essential
strengthening of Theorem A$(a)$; moreover, this leads to a proof of Theorem G.
Part $(b)$ of the latter theorem asserting that transitive homomorphisms
$B'_k\to{\mathbf S}(n)$ are primitive occurs very helpful, since
it allows us to apply Jordan Theorem on primitive permutation groups.
(See Theorem
\ref{Thm: Bk' has no homomorphisms to S(n) and Bn when k>4 and n<k},
Lemma \ref{Lm: when transitive group homomorphism B to S(n) is primitive or not} and 
Proposition
\ref{Prp: transitive imprimitive homomorphisms Bk to S(n) is cyclic if n<2k}.)
\vskip0.3cm

\noindent{\bfit Small supports.} Let $\psi\colon B_k\to{\mathbf S}(n)$ 
be a homomorphism. Artin proved that a cycle of length $>n/2$ cannot occur
in the cyclic decomposition of $\widehat\sigma_1=\psi(\sigma_1)$ whenever
$k>4$ and $\psi$ is non-cyclic and transitive (see Lemma
\ref{Lm: Artin Lemma on cyclic decomposition of hatsigma(1)}$(a)$ for the proof).
\index{Artin Lemma!on cyclic decomposition of $\widehat\sigma_1$\hfill}
On the other hand, in a quite elementary
Lemma \ref{Lm: condition no cycles of the same length restricts support} 
we show that if no different cycles of the same length $\ge 2$ occur
in the cyclic decomposition of $\widehat\sigma_1$, then
the permutations $\widehat\sigma_i=\psi(\sigma_i)$
and $\widehat\sigma_j=\psi(\sigma_j)$ with $|i-j|\ge 2$ 
must be disjoint, which implies that $\#\supp\widehat\sigma_1\le n/E(k/2)$
($E(x)$ denotes the integral part of $x$).
\index{$E(x)$\hfill}
That is, in such a case the support
of $\widehat\sigma_1$ is relatively small;
say for $k\ge 6$ and $n\le 2k$ we have $\#\supp\,\widehat\sigma_1\le 4$.
The study of such homomorphisms is by no means complicated. In fact
for $6<k<n\le 2k$ 
Lemma \ref{Lm: homomorphisms Bk to S(n) for n<=2k and supp(hatsigma(1))<6}
exhibits an explicit description
of all non-cyclic transitive homomorphisms $\psi\colon B_k\to{\mathbf S}(n)$
with $\#\supp\,\widehat\sigma_1\le 5$. In particular, we prove that
for $k,n$ as above such a homomorphism does exist only if $n=2k$
and $\widehat\sigma_1$ is a $4$-cycle.
\vskip0.3cm

\noindent{\bfit Question:} {\em What to do if $n$ is near $2k$?} 
The above methods do not go too far from the original Artin's ideas
(the main innovation is that Gorin-Lin Theorem applies systematically).
With some exceptions, Artin Fixed Point Lemma still works when 
$n$ slightly exceeds $k+1$. For instance, if $n=k+2$ and $8,\,12\ne k\ge 7$
or $n=k+3$ and $11,\,12\ne k\ge 9$, then it follows from
the Finsler inequality $\pi (2m) - \pi (m) > m/(3\log (2m))$
%\footnote{$\pi (x)$ is the number of all primes $p\le x$.}
(see \cite{Fin45,Trost}) that there is
a prime $p$ on the interval $(n/2,k-2]$.
However, if $n$ is near $2k$, these methods hardly work,
since there is no hope to find a prime on
a rather short interval $(n/2,k-2]$. For such $n$,
even if we were lucky to get the inequality 
$\#\Fix\,\widehat\sigma_1\ge k-2$,
the support of $\widehat\sigma_1$
may still contain about $k+2$ points, and we cannot
come to any immediate conclusion. 
{\bfit Ultimate answer:} {\em Look for new ideas}.
\vskip0.3cm

\noindent{\bfit Homomorphism $\Omega$; cohomology.}
\index{Homomorphism!$\Omega$\hfill}
\index{Cohomology\hfill|phantom\hfill}
For the reasons explained above, we should mainly handle
the case when $\psi\colon B_k\to{\mathbf S}(n)$ is
a non-cyclic transitive homomorphism and 
the cyclic decomposition of $\widehat\sigma_1=\psi(\sigma_1)$ contains
several different cycles of a certain length $r\ge 2$.
A simple idea described below occurs crucial (see Sections
\ref{sec: Retraction of homomorphisms; homomorphisms and cohomology}, 
\ref{sec: Omega-homomorphisms and cohomology: some computations}).
\vskip0.2cm

\noindent Suppose that for some $r\ge 2$ the $r$-component
$\mathfrak C_r$ of $\widehat\sigma_1$ consists of $t\ge 2$ \ $r$-cycles
$C_1,...,C_t$. Since $\widehat\sigma_3,...,\widehat\sigma_{k-1}$ 
commute with $\widehat\sigma_1$,
the conjugation of $\widehat\sigma_1$ by any of the elements 
$\widehat\sigma_i$ ($3\le i\le k-1$) induces
a certain permutation of the $r$-cycles $C_1,...,C_t$. This gives rise
to a homomorphism
$$
\Omega\colon B_{k-2}\to{\mathbf S}(\mathfrak C_r)\cong {\mathbf S}(t),
$$
where $B_{k-2}$ is the subgroup of $B_k$ generated
by $\sigma_3,...,\sigma_{k-1}$.
\vskip0.2cm

\noindent Furthermore, the support 
$$
\Sigma=\supp{\mathfrak C}_r \Def \bigcup_{i=1}^t \supp C_i 
\subseteq{\boldsymbol\Delta}_n\qquad (\#\Sigma=rt)
$$
of the $r$-component $\mathfrak C_r$ is invariant with respect to each 
$\widehat\sigma_i$, $3\le i\le k-1$. Let 
$$
\phi=\psi|_{ B_{k-2}}\colon B_{k-2}\to{\mathbf S}(n)
$$ 
be the restriction of the original homomorphism $\psi$
to the subgroup $B_{k-2}\subset B_k$ generated by $\sigma_3,...,\sigma_{k-1}$.
the subset $\Sigma\subset{\boldsymbol\Delta}_n$ is ($\Img\phi$)-invariant.
Consider the reduction\footnote{See Definition 
\ref{Def: symmetric groups}{\bf {(d)}}.} %of $\phi$
$$
\varphi_{_\Sigma}\colon B_{k-2}\to 
{\mathbf S}(\Sigma)\cong {\mathbf S}(rt)\subseteq{\mathbf S}(n)
$$
of $\phi$ to $\Sigma$. It is easily seen that the homomorphism $\Omega$
and the homomorphism $\varphi_{_\Sigma}$ 
(which still keeps a certain information about the original 
homomorphism $\psi$) fit in a commutative diagram of the form
%&
\begin{equation}\label{eq: CD with Omega}
\CD
@. @.  B_{k-2}@. @. \\
@. @.  @V{\varphi_{_{\Sigma}}}VV @/SE//{\,\Omega}/ @. \\ 
1@>>>H@>>>G@>>{\pi}>{\mathbf S}({\mathfrak C}_r)@>>>1\,.
\endCD
\end{equation}
%&
Here $H\cong ({\mathbb Z}/r{\mathbb Z})^t$ is the 
subgroup of ${\mathbf S}(\Sigma)$ generated by the $r$-cycles
$C_1,...,C_t$, \ $G$ is the centralizer of the element
$\mathcal C=C_1\cdots C_t$ in ${\mathbf S}(\Sigma)$,
and the horizontal line in (\ref{eq: CD with Omega})
is an exact sequence with a certain splitting 
$\rho\colon {\mathbf S}({\mathfrak C}_r)\cong {\mathbf S}(t)\to G$.
The latter exact sequence is, in a sense, {\bfit universal}, 
meaning that it is ``the same" sequence for all homomorphisms
$\psi\colon B_k\to{\mathbf S}(n)$ 
that possess an $r$-component ${\mathfrak C}$ of length
$\#{\mathfrak C}=t$.
\vskip0.2cm

\noindent The complementary set 
$\Sigma'={\boldsymbol\Delta}_n\setminus\Sigma$ is 
($\Img \phi$)-invariant as well, and we can take the reduction
$\varphi_{_{\Sigma'}}$ of $\phi$ to $\Sigma'$.
The fact that the commutator subgroup $B'_{k-2}$ for $k>6$
is perfect implies the following properties of the
homomorphism $\Omega$ (see Lemma 
\ref{Lm: if psi: B_k to S(n) is non-cyclic and varphi(Sigma') is abelian then phi, varphi(Sigma) and Omega are non-abelian} 
and Theorem
\ref{Thm: if Omega is cyclic then every Omega-homomorphism is cyclic}$(a)$):
\vskip0.2cm 

\noindent{\sl Let $k>6$. Suppose that either $(i)$ $\psi$ is non-cyclic but 
$\varphi_{_{\Sigma'}}$ is abelian or $(ii)$ $\varphi_{_{\Sigma}}$ is non-cyclic.
Then $\Omega$ must be non-cyclic.}
\vskip0.2cm

\noindent Since $k-2<k$ and $t\le n/r\le n/2\ll n$, 
we may hope to handle homomorphisms $B_{k-2}\to{\mathbf S}(t)$
(say by an appropriate induction hypothesis) and thereby recognize
the properties of $\Omega$. Assuming that $\Omega$ is already known
(at least to some extent),  we can try to recover the homomorphism
$\varphi_{_{\Sigma}}$ as far as possible from the above commutative diagram.
This is certainly a homological problem. Indeed, the homomorphism
$\Omega$ and the splitting $\rho$ give rise to an action $T$ of
$B_{k-2}$ on the abelian normal subgroup $H\vartriangleleft G$. 
We show (Proposition \ref{Prp: Omega-homomorphisms and cocycles})
that {\sl there is a natural bijection 
between the cohomology group $H_T^1(B_{k-2},H)$ and
the set of the classes of $H$-conjugate homomorphisms
$\varphi\colon B_{k-2}\to G$ that satisfy
the commutativity relation $\pi\circ\varphi=\Omega$}.
The action $T$ and the corresponding cohomology
$H_T^1(B_{k-2},H)$ can be computed explicitly
in many cases that we are interested in (see 
Sections \ref{sec: Retraction of homomorphisms; homomorphisms and cohomology}
and \ref{sec: Omega-homomorphisms and cohomology: some computations}).
As a result, we obtain a description of all possible homomorphisms
$\varphi_{_{\Sigma}}$ (up to conjugation). This description
puts certain strict restrains in the original homomorphism $\psi$.
\vskip0.2cm 

\noindent Furthermore, $\widehat\sigma_{k-1}$ is conjugate to
$\widehat\sigma_1$ and hence the $r$-component ${\mathfrak C}_r^*$
of $\widehat\sigma_{k-1}$ is also of length $\#{\mathfrak C}_r^*=t$. 
The conjugation by the element $\widehat\alpha^{k-2}
=\psi(\sigma_1\cdots\sigma_{k-1})^{k-2}$ induces a bijection
$\mathfrak C_r\stackrel{\cong}{\longrightarrow}\mathfrak C_r^*$
and an isomorphism ${\mathbf S}(\mathfrak C_r)\stackrel{\cong}{\longrightarrow}
{\mathbf S}(\mathfrak C_r^*)$.
All permutations $\widehat\sigma_1,...,\widehat\sigma_{k-3}$
commute with $\widehat\sigma_{k-1}$; a construction completely 
similar to the previous one gives rise to a homomorphism
$$
\Omega^*\colon B^*_{k-2}\to{\mathbf S}({\mathfrak C}_r^*)
\cong{\mathbf S}(t)\,,
$$
where $B^*_{k-2}\cong  B_{k-2}$ is the subgroup of
$B_k$ generated by $\sigma_1,...,\sigma_{k-3}$.
The following simple observation occurs very helpful:
{\sl under the identification
$B_{k-2}\cong B^*_{k-2}$ given by $\sigma_{i+2}\mapsto\sigma_i$,
\ $1\le i\le k-3$, the homomorphism $\Omega^*$ coincides with
the above homomorphism $\Omega$} (Lemma \ref{Lm: Omega=Omega*}).
Thus, if a certain assumtion about a homomorphism
$\psi\colon B_k\to{\mathbf S}(n)$ leads to the conclusion
that $\Omega(\sigma)\ne\Omega^*(\sigma)$ for some $\sigma\in B_{k-2}$,
then this assumption is wrong. 
\vskip0.3cm

\noindent{\bfit Down with long components! Fixed points whithout primes.}
Using the cohomology approach (and also Lemma \ref{Lm: Omega=Omega*}
mentioned above), we prove that {\sl for a non-cyclic homomorphism
$\psi\colon B_k\to{\mathbf S}(n)$
the permutation $\widehat\sigma_1$ cannot have an $r$-component
of length $t>k-3$ whenever $6<k<n<2k$}
(see Lemma \ref{Lm: no components of length t>k-3},
which is one of the crucial technical results).
On the other hand, with help of Theorem A$(a)$, we prove in
Lemma \ref{Lm1: homomorphism with all components of hatsigma(1) of length<=k-3}
the following statement: {\sl if $k>6$ and all the components
of $\widehat\sigma_1$ $($including the degenerate component
$\Fix \widehat\sigma_1)$ are of length at most $k-3$, then the homomorphism
$\psi$ must be cyclic}. Combining these two results, we
prove the following analog of Artin Fixed Point Lemma
(see Corollary \ref{Crl: k-2 fixed points of psi(sigma(1)) for 6<k<n<2k}):
\vskip0.2cm

\noindent {\sl If $6<k<n<2k$, then $\#\Fix \widehat\sigma_1\ge k-2$
for any non-cyclic homomorphism $\psi$.}
\vskip0.2cm

\noindent This leads to Theorem E$(b)$ (at least for $k>6$;
the cases $k=5,\,6$ may be handeled as well; see 
Propositions \ref{Prp: B(5) to S(7)},\, 
\ref{Prp: B(6) to S(8)} and
Theorem \ref{Thm: homomorphisms B(k) to S(k+2)}$(a)$). 
\vskip0.3cm

\noindent{\bfit Homomorphisms $B_k\to{\mathbf S}(n)$, \ $k<n\le 2k$.}
\index{Homomorphisms!$B_k\to{\mathbf S}(n)$, $k<n\le 2k$\hfill}
Theorem F$(a)$ is proven by induction on $k$.
To get a base of induction, we study first
the cases $k=7,\,8$ (see Lemmas \ref{Lm: B(7) to S(n) for 7<n<14}
and \ref{Lm: B(8) to S(n), 8<n<16}). Then, assuming that $6<k<n<2k$ and
$\psi\colon B_k\to{\mathbf S}(n)$ is
a non-cyclic transitive homomorphism, we take the restriction
$\phi\colon B_{k-2}\to{\mathbf S}(n)$ of
$\psi$ to the subgroup $B_{k-2}\subset B_k$
generated by $\sigma_3,...,\sigma_{k-1}$ and the reductions
$\varphi\colon B_{k-2}\to{\mathbf S}(\Sigma)$
and $\varphi'\colon B_{k-2}\to{\mathbf S}(\Sigma')$
of $\phi$ to the mutually complementary
$\psi(B_{k-2})$-invariant sets $\Sigma=\supp\,\widehat\sigma_1$
and $\Sigma'={\boldsymbol\Delta}_n\setminus\Sigma$ respectively. 
So $\phi$ is the disjoint product of $\varphi$ and $\varphi'$.
In
Lemma \ref{Lm2: homomorphism with all components of hatsigma(1) of length<=k-3},
which is the main technical tool of the induction,
we show that {\sl $\varphi$ is trivial, $\varphi'$ is
non-cyclic and $\phi=\varphi'$}. The proof of this lemma
involves the homomorphisms $\Omega_{\mathfrak C},\,\Omega^*_{\mathfrak C}$ 
(corresponding to every non-degenerate component 
$\mathfrak C$ of $\widehat\sigma_1$) and
Theorem
\ref{Thm: if Omega is cyclic then every Omega-homomorphism is cyclic}$(a)$
mentioned above. 
Finally, assuming, by the inductive hypothesis, the existence of a natural
$m$ such that any transitive homomorphism
$\psi\colon B_k\to{\mathbf S}(n)$ is cyclic
whenever $k,n$ satisfy $6<k\le m$ and $k<n<2k$,
we prove that the same conclusion
must be true whenever $6<k\le m+2$ and $k<n<2k$.
The justification of this inductive step involves
Lemma \ref{Lm2: homomorphism with all components of hatsigma(1) of length<=k-3},
Corollary \ref{Crl: k-2 fixed points of psi(sigma(1)) for 6<k<n<2k},
Artin Theorem, Theorem G, and Jordan Theorem on primitive permutation groups.
(See Theorem \ref{Thm: homomorphisms B(k) to S(n), 6<k<n<2k}.)
\vskip0.2cm

\noindent The treatment of the case $n=2k$ is more sophisticated,
since the permutation $\widehat\sigma_1=\psi(\sigma_1)$
corresponding to a non-cyclic transitive homomorphism
$\psi\colon B_k\to{\mathbf S}(2k)$
may possess a $2$-component of length $t>k-3$. However, with help of 
an appropriate cohomology and Theorem F$(a)$,
we show that in the latter case either $t=k$ and $\psi$ is conjugate
to the model homomorphism $\varphi_2$
or $t=k-2$ and $\psi$ is conjugate to $\varphi_3$ (see Lemma
\ref{Lm: psi B(k) to S(2k) with psi(sigma(1)) having component of length>k-3}).
Combining this property with Theorem A$(a)$,
Artin Theorem, and Theorem F$(a)$, we prove Theorem F$(b)$
(see Theorem
\ref{Thm: non-cyclic transitive B(k) to S(2k) is standard for k>6}). 
\vskip0.3cm

\noindent{\bfit Homomorphisms $\phi\colon B_k\to B_n$.} 
\index{Homomorphisms!$B_k\to B_n$, $k>n$\hfill}
Assume that $k>4$ and consider the composition
$\psi = \mu\circ\phi\colon B_k\to{\mathbf S}(n)$ of $\phi$ with the canonical
epimorphism $\mu\colon B_n\to{\mathbf S}(n)$.
If $\psi$ is cyclic, then $\psi(B'_k)\subseteq\Ker\mu = {PB_k}$.
Since $ B'_k$ is perfect,
\index{Gorin-Lin Theorem\hfill}
\index{Theorem!of Gorin-Lin\hfill}
\index{Markov Theorem\hfill}
\index{Theorem!of Markov\hfill}
Markov Theorem implies that $\psi(B'_k)=\{1\}$ and $\phi$ is integral.
Combining this simple observation with Theorem A$(a)$,
we obtain Theorem A$(b)$. In a similar way we complete
the proof of Theorem A$(c)$ (we already have commented on
the absence of non-trivial homomorphisms
$ B'_k\to{\mathbf S}(n)$ for $k>4$ and $n<k$).
\vskip0.2cm

\noindent To prove Theorem B, we use
Lemma \ref{Lm: non-surjective homomorphism B(k) to S(k) is cyclic}
cited above and show that if
$k\ne 4$ then for any non-integral endomorphism $\phi$ of $B_k$
the composition $\psi = \mu\circ\phi$ is surjective and,
therefore, non-cyclic and transitive. By Artin Theorem,
\index{Artin Theorem\hfill}
\index{Theorem!of Artin\hfill}
$\Ker\psi = {PB_k}$, which implies
$\phi^{-1}({PB_k})={PB_k}$; the rest assertions
of the theorem follow readily from this fact.
\vskip0.3cm

\noindent{\bfit Homomorphisms $\psi\colon B'_k\to{\mathbf S}(k)$.
Endomorphisms of $ B'_k$.}
\index{Homomorphisms!$B'_k\to{\mathbf S}(k)$\hfill}
\index{Endomorphisms of $B'_k$\hfill}
Let $6\ne k>4$. Restricting a non-trivial homomorphism
$\psi$ to the subgroup $B_{k-2}\subset B'_k$
generated by the elements $c_i=\sigma_{i+2}\sigma_1^{-1}$,
\ $1\le i\le k-3$, we obtain a homomorphism
$\phi\colon B_{k-2}\to{\mathbf S}(k)$. Then we apply
Lemma \ref{Lm: composition B(k-2) to Bk' to H}
mentioned above together with Theorem A$(a)$ and Theorem E to show that 
$\psi$ is {\em tame}, meaning that the permutation group
$\phi(B_{k-2})\subset{\mathbf S}(k)$ has
an orbit $Q\subset{\boldsymbol\Delta}_k$ of length $k-2$
(see Definition \ref{Def: tame homomorphism psi: B'(k) to S(k)}
and Lemma \ref{Lm: tame and non-tame homomorphisms B'(k) to S(k)}). 
By Artin Theorem, the reduction $\phi_Q$ of $\phi$ to $Q$
is conjugate to the canonical epimorphism
$\mu_{k-2}\colon B_{k-2}\to{\mathbf S}(Q)\cong{\mathbf S}(k-2)$. 
Hence, without loss of generality, we may assume that
$\widehat c_i=\psi(c_i)=(1,2)(i+2,i+3)$ for all $i=1,...,k-3$.
Using this property and the defining relations (\ref{1.14})-(\ref{1.21}),
we show (by a straightforward computation) that $\psi\sim\mu_k'$;
this proves Theorem C (see Theorem \ref{Thm: homomorphisms B'(k) to S(k)},
where the case $k=6$ is treated as well).
In view of Markov and Gorin-Lin theorems, Theorem D follows immediately
from Theorem C (see Theorem
\ref{Thm: J(k) is a completely characteristic subgroup of B'(k)}).
\vskip0.3cm

\noindent{\bfit Special homomorphisms $\phi\colon B_k 
\to B_n$.} The proof of Theorem H is mainly based upon Lemma
\ref{Lm: k|m if psi(alpha m)<->psi(beta),(k-1)|m if psi(beta m)<->psi(alpha)},
which provides us with some arithmetical properties
of non-cyclic homomorphisms of braid groups
(see Theorem
\ref{Thm: arithmetical properties of special homomorphisms B(k) to B(n)}).

\newpage

%Sec. 2
\section{Auxiliary results}\label{sec: Auxiliary results}

To facilitate the exposition of the main proofs,
I collected here a number of various elementary statements.
Some of them apply just occasionally but
some other are involved throughout the paper.
I suggest to take a look at the corresponding parts of this section
when (and only if) the reader feels such a need.

\subsection{Three algebraic lemmas}
\label{Ss: 2.1. Three algebraic lemmas}
Recall that a word in the variables
$a_1,...,a_s$ and $a_1^{-1},...,a_s^{-1}$ is said
to be {\it reduced} if it does not contain a part of the form
$a_{i}a^{-1}_{i}$ or $a^{-1}_{i}a_{i}$ ($1\le i\le s$).

\begin{Lemma}
\label{Lm: if reduced f(u,uinv,v,vinv)=1 in a free group then uv=vu} 
Let $f(x,x^{-1},y,y^{-1})$ be a nonempty
reduced word in variables $x,y$ and $x^{-1},y^{-1}$.
If the elements $u$ and $v$ of a group $G$ satisfy
$f(u,u^{-1},v,v^{-1}) = 1$, then for any homomorphism
$\phi\colon G\to {\mathbb F}$ to a free group the elements
$\widehat u = \phi (u)$ and $\widehat v = \phi (v)$ commute.
\end{Lemma}

\begin{proof} $\widehat{u}$ and $\widehat{v}$
generate a free subgroup $H\subseteq {\mathbb F}$
of rank $r\le 2$. Were $r=2$,  
$\{\widehat{u},\widehat{v}\}$ would be a free base of $H$,
which is impossible since
$f(\widehat u,\,\widehat u^{-1},\,\widehat v,\,\widehat v^{-1})
= 1$. Hence $r\le 1$ and $H$ is commutative.
\end{proof}

\begin{Lemma}
\label{Lm: non-zero det(M2-I)->integral psi is trivial}
Let $\Psi\colon H\to G$ be a group
homomorphism and $K\subseteq H$ be a subgroup with two generators $u,v$.
Let $\sigma\in H$, $x,y\in K'=[K,K]$ and a matrix $M=\begin{pmatrix}
p & q\\
s & t
\end{pmatrix}$ $(p,q,s,t\in{\mathbb Z})$ be such that 
$\det M^2=1$ and $\det (M^2-I)\ne 0$, whereas
%&
\begin{equation}\label{2.1}
\sigma u\sigma^{-1} = u^{p}v^{q}x\quad
\text{and}\quad \sigma v\sigma^{-1} = u^{s}v^{t}y\,.
\end{equation}
%&
If the restriction $\psi=\Psi|K\colon K\to G$ is integral
then it is trivial.
\end{Lemma}

\begin{proof}
Put $\widehat{u}=\psi(u)$, $\widehat{v}=\psi (v)$,
$\widehat{\sigma} = \Psi(\sigma)$; it suffices to show that
$\widehat{u} = \widehat{v} = 1$. 
Since $\psi$ is integral, the subgroup
$\widehat{K}=\psi(K)\subseteq G$ generated by
$\widehat{u}$ and $\widehat{v}$ 
either is trivial and there is nothing to prove
or $\widehat{K}\cong{\mathbb Z}$.
In the latter case we still have $\psi(x) = \psi(y) = 1$,
for $x,y\in K'$, and relations (\ref{2.1}) show that
%&
\begin{equation}\label{2.2}
\widehat{\sigma}\widehat{u}\widehat{\sigma}^{-1}
= \widehat{u}^{p}\widehat{v}^{q}\quad\text{and}\quad
\widehat{\sigma}\widehat{v}\widehat{\sigma}^{-1}
= \widehat{u}^{s}\widehat{v}^{t}\,.
\end{equation}
Clearly $\det M = \pm 1$; hence (\ref{2.2}) imply that
the conjugation by the element $\widehat{\sigma}\in G$
determines an automorphism $S$ of the subgroup $\widehat{K}\cong{\mathbb Z}$.
Any automorphism of ${\mathbb Z}$ is involutive; that is,
$S^2 = \id$. Combining this fact with relations (\ref{2.2})
and passing to the additive notation (which is natural for
$\widehat u,\widehat v\in \widehat K\cong{\mathbb Z}$), 
we see that the vector $\begin{pmatrix} \widehat u\\ \widehat v\end{pmatrix}
\in{\mathbb Z}^2$ satisfies the system of linear equations
$$
\begin{pmatrix} \widehat u\\ \widehat v\end{pmatrix}
= S^2\begin{pmatrix} \widehat u\\ \widehat v\end{pmatrix} =
M^2\begin{pmatrix} \widehat u\\ \widehat v\end{pmatrix}\,.
$$
Since $\det(M^2-I)\ne 0$, this implies $\widehat u=\widehat v=0$
and concludes the proof.
\end{proof}

\begin{Example}\label{Exmp 2.1}
Let $H=G=\left<a,\sigma|\, \sigma a\sigma^{-1} = a^2\right>$.
The elements $u=a^2$ and $v=a^3$ generate the subgroup $K\cong{\mathbb Z}$.
The natural embedding $\psi\colon K\hookrightarrow H=G$
is a non-trivial integral homomorphism that extends to $H$. Clearly
$\sigma u\sigma^{-1}=u^2$, $\sigma v\sigma^{-1}=v^2$;
for the corresponding matrix $M=\begin{pmatrix}
2 & 0\\
0 & 2
\end{pmatrix}$ we have $\det M^2 = 16$ and $\det(M^2-I)=9$,
which shows that the condition $\det M^2 = 1$
in Lemma \ref{Lm: non-zero det(M2-I)->integral psi is trivial}
is essential and cannot be replaced with $\det M\ne 0$.
\hfill $\bigcirc$
\end{Example}

\begin{Notation}\label{Not: GCD and residue}
\index{$\GCD(m,n)$\hfill}
\index{$\modo{N}_r$\hfill}
The greatest common divisor of $m,n\in{\mathbb Z}$
is denoted by $\GCD(m,n)$ or simply $(m,n)$.
For natural $N,r$ we denote by $|N|_r$ the residue of $N$
modulo $r$ $(0\le |N|_r\le r-1)$.
\vskip0.2cm

\noindent We use also the notation $a\infty b$ for a braid-like couple $a,b$
(see Definition \ref{Def: sets, groups, homomorphisms}{\bf (d)}).
\hfill $\bigcirc$
\end{Notation}

\begin{Lemma}
\label{Lm: braid-like couples (a,b) and (aq,b)} 
Let $q\ge 2$ be a natural number and let
$\nu =\GCD(q+1,4)$. Assume that a braid-like couple
$a,b$ in a group $G$ satisfy also the condition $b\infty a^q$. Then
$$
a^{\nu (q-1)} = b^{\nu (q-1)} = 1.
$$
\end{Lemma}

\begin{proof} 
The conditions $b\infty a$ and $b\infty a^{q}$
imply that
%&
\begin{equation}\label{2.3}
a^{q-1}b = b^{q}a^{q-1},
\end{equation}
%&
%&
\begin{equation}\label{2.4}
ba^{q-1} = a^{q-1}b^{q},
\end{equation}
%&
and also
$$
\aligned
ab^{2(q-1)}a^{-1}&= b^{-1}a^{2(q-1)}b = b^{-1}a^{q-1}a^{q-1}b\\
&=b^{-1}a^{q-1}b^qa^{q-1}=
b^{-1}ba^{q-1}a^{q-1}=a^{2(q-1)}\,.
\endaligned
$$
Hence $b^{2(q-1)}=a^{2(q-1)}$ and $a$ commutes with $b^{2(q-1)}$. 
It follows from (\ref{2.3}) and (\ref{2.4}) that
$b^{q+1}a^{q-1} = ba^{q-1}b=a^{q-1}b^{q+1}$, i. e., $a^{q-1}$
commutes with $b^{q+1}$; therefore, the element $c=a^{q-1}$
commutes with $b^{q+1}$ and $b^{2(q-1)}$. This implies that
$c$ commutes also with $b^d$, where $d=\GCD(q+1,2(q-1))$
(since $d=(q+1)m+2(q-1)n$ for suitable $m,n$). But
$\GCD(q+1,2(q-1))=\GCD(q+1,4)=\nu$, and thus
%&
\begin{equation}\label{2.5}
a^{q-1} \ \text{commuts with} \ b^\nu\,.
\end{equation}
%&
Since $\nu$ divides $q+1$, we have $q=r\nu + (\nu-1)$
for some non-negative $r\in{\mathbb Z}$.
Taking into account (\ref{2.5}) and using one time (\ref{2.4})
and $\nu-1$ times (\ref{2.3}), we obtain
$$
ba^{q-1} = a^{q-1}b^q = a^{q-1}b^{r\nu}b^{\nu-1} =
b^{r\nu}a^{q-1}b^{\nu-1}=
b^{r\nu}b^{(\nu-1)q}a^{q-1}=b^{r\nu + (\nu-1)q}a^{q-1}\,.
$$
Hence $b^{\nu(q-1)} = b^{r\nu + (\nu-1)q - 1} = 1$.
The latter relation and $a\infty b$ imply $a^{\nu (q-1)} = 1$. 
\end{proof}

\begin{Example}\label{Exmp 2.2} 
Let $G = {\mathbf S}(8)$, \
$a = (1,2,3,4,5,6,7,8)$, $b = (1,7,6,8,5,3,2,4)$ and $q = 3$.
Then
$b\infty a$, \ $b\infty a^{3}$, \ $\nu=\GCD(q+1,4)=4$ and
$$
\nu (q-1) = 8 =\ord a = \ord b\,.
$$
This shows that the result of
Lemma \ref{Lm: braid-like couples (a,b) and (aq,b)} is sharp.
\hfill $\bigcirc$
\end{Example}

\subsection{Some properties of permutations}
\label{subsec 2.1. Some properties of permutations}
Here we prove some simple lemmas about permutations.

\begin{Lemma}[cf. \cite{Art47b}]
\label{Lm: invariant sets and components of commuting permutations}
Suppose that $A,B\in{\mathbf S}(n)$ and $AB = BA$. Then:

$a)$ the set $\supp A$ is $B$-invariant;

$b)$ if for some $r$, \ $2\le r\le n$, the $r$-component
of $A$ consists of a single $r$-cycle $C$,
then the set $\supp C$ is $B$-invariant and
$B|\supp C = C^q$ for some integer $q$, \ $0\le q<r$.
\end{Lemma}

\begin{proof}
$a)$ Let $i\in\supp A$; then $A(i)\ne i$ and
$AB(i)=BA(i)\ne B(i)$; hence $B(i)\in\supp A$.
\vskip0.2cm

$b)$ Let $A=CD_1\cdots D_t$ be the cyclic decomposition of $A$.
Then
$$
CD_1\cdots D_t = A = B^{-1}AB=B^{-1}CD_1\cdots D_tB
= B^{-1}CB\cdot B^{-1}D_1B\cdots B^{-1}D_tB\,.
$$
Since $C$ is the only $r$-cycle in the
cyclic decomposition of $A$, we have $C = B^{-1}CB$, and
$(a)$ implies that the set $\supp C$
is $B$-invariant. Let $C=(i_0,i_1,...,i_{r-1})$. Then $B(i_0)=i_q$
for some $q$, \ $0\le q\le r-1$. Let us prove that
$B|\supp C = C^q$. Since $C^q(i_s)=i_{\modo{s+q}_r}$,
we should show that
$B(i_s)=i_{|s+q|_r}$ for all $s=0,1,...,r-1$.
The proof is by induction on $s$ with the case $s=0$ clear
($0\le q\le r-1$ and $\modo{q}_r=q$).
Assume that for some $k$, $1\le k\le r-1$, we have
$B(i_s)=i_{\modo{s+q}_r}$ for all $s=0,...,k-1$.
Then
$$
\quad B(i_k)= BA(i_{k-1})=AB(i_{k-1})
=A(i_{\modo{k-1+q}_r})=
i_{\modo{\modo{k-1+q}_r+1}_r}=i_{\modo{k+q}_r}\,,
$$
which completes the proof.
\end{proof}

\begin{Lemma}
\label{Lm: a single invariant set of length r} 
Assume that $A\in{\mathbf S}(n)$ and for some
$r\in{\mathbb N}$ \ $(1\le r < n)$ the family $\Inv_r A$ consists of a single
set $\Sigma$. The following statements hold:

$a)$ If $C\in{\mathbf S}(n)$ and $D = CAC^{-1}$, then $\Inv_r D=\{C(\Sigma)\}$.

$b)$ If $B\in{\mathbf S}(n)$ and $AB=BA$, then $B(\Sigma)=\Sigma$.

$c)$ If $B\in{\mathbf S}(n)$, $AB=BA$ and $B\sim A$, then $\Inv_r B=\{\Sigma\}$.
\end{Lemma}

\begin{proof}
$a)$ Since $DC(\Sigma)=CA(\Sigma)=C(\Sigma)$,
we have $C(\Sigma)\in\Inv_r D$. If $\Sigma'\in\Inv_r D$,
then $AC^{-1}(\Sigma')=C^{-1}D(\Sigma')=C^{-1}(\Sigma')$, so that
$C^{-1}(\Sigma')\in \Inv_r A$; hence $C^{-1}(\Sigma') = \Sigma$.
\vskip0.2cm

$b)$ $A = BAB^{-1}$ and $(a)$ imply
$B(\Sigma)\in\Inv_r A=\{\Sigma\}$ so that $B(\Sigma)=\Sigma$.
\vskip0.2cm

$c)$ $B = CAC^{-1}$ for a certain $C\in{\mathbf S}(n)$.
By $(a)$ and $(b)$, we have $\Inv_r B=\{C(\Sigma)\}$,
$\Sigma\in\Inv_r B=\{C(\Sigma)\}$ and hence $C(\Sigma)=\Sigma$. 
\end{proof}

\begin{Lemma}
\label{Lm: braid-like couple of 3 cycles} 
$\#(\supp A\cap\supp B)=2$ for any braid-like
couple of $3$-cycles $A,B\in{\mathbf S}(n)$.
\end{Lemma}

\begin{proof}
Since $AB\ne BA$, we have
$\supp A\cap\supp B\ne\varnothing$. Moreover,
$\supp A\ne\supp B$ (every $3$-cycles with the same support
commute to each other). Finally, if $\#(\supp A\cap\supp B)=1$,
a simple computation shows that $ABA\ne BAB$.
\end{proof}

\noindent The next lemma is evident.

\begin{Lemma}
\label{Lm: commuting 3 cycles} 
If $AC=CA$, \ $[A]=[C]=[3]$, and
$\supp A\cap\supp C\ne\varnothing$,
then $\supp A = \supp C$ and $C=A^{q}$, where either $q=1$ or $q=2$.
\end{Lemma}

\begin{Lemma}\label{Lm: braid-like triple of 3 cycles} 
Assume that $[A]=[B]=[C]=[3]$, \
$ABA = BAB$, $BCB = CBC$ and $AC=CA$. Then $A = C$.
\end{Lemma}

\begin{proof}
If $BA=AB$ or $BC=CB$, the assumed relations imply $A=B=C$.
Thus, we may suppose that $B\infty A$ and $B\infty C$. 
Then Lemma \ref{Lm: braid-like couple of 3 cycles} implies
$\#(\supp\,B\cap\supp A)=\#(\supp\,B\cap\supp\,C)=2$.
Hence the commuting $3$-cycles
$A$ and $C$ are not disjoint and, by Lemma \ref{Lm: commuting 3 cycles},
we have $C=A^{q}$ with $q=1$ or $q=2$. If $q=2$, then 
$\GCD(q+1,4)=1$, \ $q-1=1$, \ $B\infty A$, \ $B\infty A^q$ and
Lemma \ref{Lm: braid-like couples (a,b) and (aq,b)} implies $A=B=1$,
which contradicts our assumptions. Thus, $q=1$ and $C=A$.
\end{proof}

\begin{Lemma}
\label{Lm: p cycles in S(2p)} 
Let $p$ be a prime number and let $A\in{\mathbf S}(2p)$ be a product
of two disjoint $p$-cycles $B = (b_{0},b_{1},\ldots ,b_{p-1})$ and
$C = (c_{0},c_{1},\ldots ,c_{p-1})$. If $D\in{\mathbf S}(2p)$
commutes with $A$, then the following three cases may only occur:
%Let $p$ be a prime number.
%Assume that $A,D\in{\mathbf S}(2p)$, $AD=DA$ and $A=BC$, where
%$B = (b_{0},b_{1},\ldots ,b_{p-1})$ and
%$C = (c_{0},c_{1},\ldots ,c_{p-1})$
%are disjoint $p$-cycles. Then the following three cases may only occur:

\begin{itemize}

\item[$(i)$] $D = B^m C^n$ for certain $m,n$ with $0\le m,n < p$;

\item[$(ii)$] $D$ is a product of $p$ disjoint transpositions
$D_{i} = (b_i,c_{|i+r|_p})$, \ $0\le i\le p-1$, 
where $r$ does not depend on $i$ and $0\le r<p$;

{\item[$(iii)$] $D$ is a $2p$-cycle of the form}
$$
{\phantom{(iiiii)}}(b_0, c_r, b_{|r+s|_p}, c_{|r+(r+s)|_p},
b_{|2(r+s)|_p}, c_{|r+2(r+s)|_p},
\ldots ,b_{|(p-1)(r+s)|_p},c_{|r+(p-1)(r+s)|_p})
$$
\vskip0.2cm

\item[] and $A=D^{2q}$, where $0\le r,s < p$, \ $|r+s|_p\ne 0$,
and $q$ is defined by the conditions $1\le q<p$, \ $|q(r+s)|_p=1$.

\end{itemize}
\end{Lemma}

\begin{proof}
Since $D$ commutes with $A=BC$, we have
$DBD^{-1}\cdot DCD^{-1} = BC$. Clearly $DBD^{-1}$ and
$DCD^{-1}$ are disjoint $p$-cycles. So either $DBD^{-1} = B$
and $DCD^{-1} = C$ or $DBD^{-1} = C$ and $DCD^{-1} = B$.
\vskip0.2cm

\noindent In the first case 
Lemma \ref{Lm: invariant sets and components of commuting permutations}$(b)$
implies that $D=B^m C^n$, where $0\le m,n<p$. 
\vskip0.2cm

\noindent In the second case there are
uniquely determined integers $r,s$ such that
$D(b_0)=c_r$, $D(c_0)=b_s$ and $0\le r,s<p$. Consequently,
$D(b_i)=c_{\modo{i+r}_p}$ and $D(c_j)=b_{\modo{j+s}_p}$
for all $i,j$, \ $0\le i,j<p$.

If $\modo{r+s}_p=0$, then 
$D^2(b_i) = D(c_{\modo{i+r}_p}) = b_{\modo{i+r+s}_p} = b_i$
and $D$ must coincide with the product or $p$ disjoint transpositions
$D_i=(b_i,c_{\modo{i+r}_p})$.

Finally, if $\modo{r+s}_p\ne 0$, then $\modo{t(r+s)}_p\ne 0$
for any $t$ with $1\le t<p$ (since $p$ is prime), and therefore
$D$ must coincide with the $2p$-cycle exhibited in the
formulation of the lemma. The rest assertions related to this case are evident.
\end{proof}

\noindent The next lemma will be used in Sec.
\ref{sec: Homomorphisms B'(k) to S(k) and endomorphisms of B'(k)}

\begin{Lemma}
\label{Lm: 3- and 5-cycles}
$a)$ Assume that $A,B\in{\mathbf S}(k)$ are
$3$-cycles. Then at least one of the permutations $AB$, $A^{-1}B$
is not a $3$-cycle.

$b)$ Assume that $A,B\in{\mathbf S}(5)$ are
$5$-cycles. Then at least one of the permutations $AB$, $A^{-1}B$, $A^{-2}B$,
$B^2AB$ is not a $5$-cycle.
\end{Lemma}

\begin{proof} 
$(a)$ is trivial. To check $(b)$, suppose 
that $A=(a,b,c,d,e)$ and $B,AB\in{\mathbf S}(5)$ are $5$-cycles. 
Then $B$ must be one of the following eight $5$-cycles:
$$
A,\ A^2,\ A^3,\ (a,b,d,e,c), \ (a,b,e,c,d),\ (a,c,d,b,e),\ 
(a,d,e,b,c),\ (a,d,b,c,e).
$$ 
The condition $[B^2AB]=[5]$ eliminates
all the cycles from this list but $A$ and $A^2$. Finally, for $B=A$
we have $A^{-1}B=1$; and if $B=A^2$, then $A^{-2}B=1$.
\end{proof}

\noindent In the following three lemmas $A,B\in{\mathbf S}(n)$;
we omit proofs which are straightforward.

\begin{Lemma}
\label{Lm: commuting permutations of type [2,2]} 
Let $A\ne B$, $AB=BA$, $[A]=[B]= [2,2]$ and $\supp A\cap\sup B\ne\varnothing$.
Then either $(i)$ $\supp A = \supp B$
but the cyclic decompositions of $A$ and $B$ contain no common
transpositions or $(ii)$ $\supp A$ and $\supp B$ have exactly
two common points whose transposition 
occurs in the cyclic decompositions of both $A$ and $B$.
\hfill $\square$
\end{Lemma}

\begin{Lemma}
\label{Lm: braid-like couples of type [2,2]}  
Let $[A]=[B]=[2,2]$ and $A\infty B$.
Then either $(i)$ $\supp A$ and $\supp B$ have exactly
three common points and the transposition of certain two of them
occurs in the cyclic decompositions of both $A$ and $B$ 
or $(ii)$ $\supp A$ and $\supp B$ have exactly
two common points and their transposition occurs in neither
of the cyclic decompositions of $A$ and $B$.
\hfill $\square$
\end{Lemma}

\begin{Lemma}
\label{Lm: braid-like couples of 4-cycles} 
Let $[A]=[B]=[4]$ and $A\infty B$.
Then either $(i)$ $\supp A = \supp B$ and $B$ may be obtained from $A$ by a
transposition of two neighboring symbols in $A$ or $(ii)$ 
$\supp A$ and $\supp B$ have exactly two common symbols which are
neighboring neither in $A$ nor in $B$.
\hfill $\square$
\end{Lemma}

\subsection{Elementary properties of braid homomorphisms}
\label{Ss: 2.2. Some elementary properties of braid homomorphisms}
Lemmas
\ref{Lm: if psi(sigma(i)) and psi(sigma(i+1)) commute then psi is cyclic}-\ref{Lm: if psi(alpha(ij)r) and psi(sigma(i)) commute then psi(sigma(1))=psi(sigma(3))}
concerning a group homomorphism $\psi\colon B_k\to H$ \ ($k\ge 3$)
are contained in \cite{Art47b} (the latter one in a slightly weaker form);
for the completeness of the exposition we give the proofs.

\begin{Lemma} 
\label{Lm: if psi(sigma(i)) and psi(sigma(i+1)) commute then psi is cyclic} 
Assume that for some $i$, \ $1\le i\le k-2$,
the elements $\psi(\sigma_i)$ and $\psi(\sigma_{i+1})$ commute.
Then the homomorphism $\psi$ is cyclic.
\end{Lemma}

\begin{proof}
Relation (\ref{eq: braid like}) implies
$\psi(\sigma_i)=\psi(\sigma_{i+1})$. In view of (\ref{eq: commut}), it
follows that $\psi(\sigma_{i-1})$ commutes with $\psi(\sigma_i)$,
and $\psi(\sigma_{i+1})$ commutes with $\psi(\sigma_{i+2})$. Hence,
we have $\psi(\sigma_{i-1})=\psi(\sigma_i)=
\psi(\sigma_{i+1})=\psi(\sigma_{i+2})$. Proceeding with this process,
we obtain that all elements $\psi(\sigma_i)$, \ $1\le i\le k-1$,
coincide.
\end{proof}

\begin{Lemma}
\label{Lm: if k ne 4 and psi(sigma(i))=psi(sigma(j)) for i ne j then psi is cyclic} 
If $\psi(\sigma_i)=\psi(\sigma_j)$ for some $i<j$ and 
either $k\ne 4$ or $k=4$ but $j-i\ne 2$, then $\psi$ is cyclic.
\end{Lemma}

\begin{proof}
If $k=3$ or $j=i+1$, then $\psi$ is cyclic
by Lemma
\ref{Lm: if psi(sigma(i)) and psi(sigma(i+1)) commute then psi is cyclic}.
So we may assume that $k>3$
and $j\ge i+2$. If $j>i+2$, then (\ref{eq: commut}) and the assumption
$\psi(\sigma_i)=\psi(\sigma_j)$
show that $\psi(\sigma_i)$ commutes
with $\psi(\sigma_{i+1})$, and $\psi$ is
cyclic (Lemma \ref{Lm: if psi(sigma(i)) and psi(sigma(i+1)) commute then psi is cyclic}). Finally, if $k>4$ and $\psi(\sigma_i)=
\psi(\sigma_{i+2})$, then either $i>1$ or $i=1$ and $i+2=3<k-1$;
in the first case $\psi(\sigma_{i-1})$ commutes
with $\psi(\sigma_i)$; in the second case
$\psi(\sigma_3)$ commutes
with $\psi(\sigma_4)$; by
Lemma \ref{Lm: if psi(sigma(i)) and psi(sigma(i+1)) commute then psi is cyclic},
$\psi$ is cyclic.
\end{proof}

\begin{Lemma}
\label{Lm: if psi(alpha(ij)r) and psi(sigma(i)) commute then psi(sigma(1))=psi(sigma(3))}
Assume that for some
$i,j \ \ (1\le i < j-1\le k-1)$ there exists a natural
$r\not\equiv 0 \, (\modl\,j-i+1)$ such that $\psi(\alpha_{ij}^r)$
commutes with $\psi(\sigma_i)$.
Then $\psi(\sigma_1)=\psi(\sigma_3)$; if in addition $k\ne 4$,
then $\psi$ is cyclic.
\end{Lemma}

\begin{proof}
Let $q=|r|_{j-i+1}$ (see Notation \ref{Not: GCD and residue}).
It was already noticed that relations
(\ref{alpha(ij)pow(j-i+1)=beta(ij)pow(j-i)})
imply that $\alpha_{ij}^{j-i+1}$ commutes with
$\sigma_i,\sigma_{i+1},...,\sigma_{j-1}$;
therefore, it follows from our assumption that 
%&
\begin{equation}\label{2.6}
\psi(\alpha^{q}_{ij})\psi(\sigma_{i})=\psi(\sigma_{i})\psi(\alpha^{q}_{ij})\,.
\end{equation}
%&

If $q=j-i$, then $q\equiv -1\,(\modl\,{j-i+1})$ and (\ref{2.6}) shows that
$\psi(\alpha^{-1}_{ij})$ commutes with $\psi(\sigma_i)$. But
then $\psi(\alpha_{ij})$ also commutes with $\psi(\sigma_i)$
and (\ref{conjugation by alphaij}) implies $\psi(\sigma_i)=\psi(\sigma_{i+1})$;
by Lemma
\ref{Lm: if k ne 4 and psi(sigma(i))=psi(sigma(j)) for i ne j then psi is cyclic},
$\psi$ is cyclic.

Assume now that $q\le j-i-1$, i. e., $i\le i+q\le j-1$; then the third of
relations (\ref{conjugation by alphaij}) (with $m=i$)
implies that $\alpha^q_{ij}\sigma_i=\sigma_{i+q}\alpha^q_{ij}$.
Combining this with \ref{2.6}, we obtain
$\psi(\sigma_i)=\psi(\sigma_{i+q})$.
If $q\ne 2$ or $k\ne 4$, the homomorphism $\psi$ is cyclic
(Lemma
\ref{Lm: if k ne 4 and psi(sigma(i))=psi(sigma(j)) for i ne j then psi is cyclic});
if $q=2$ and $k=4$, then $i=1$ and $\psi(\sigma_1)=\psi(\sigma_3)$.
\end{proof}

\noindent Lemma
\ref{Lm: if psi(sigma(i)) and psi(sigma(i+1)) commute then psi is cyclic},
Lemma \ref
{Lm: if psi(alpha(ij)r) and psi(sigma(i)) commute then psi(sigma(1))=psi(sigma(3))}
and relations (\ref{eq: braid like}) 
imply the following corollary:

\begin{Corollary}
\label{Crl: if psi is non-cyclic then psi(sigma(i)) infty psi(sigma(i+1))} 
Let $\psi\colon B_k\to H$ be a non-cyclic homomorphism. Then
$$
\psi(\sigma_i)\infty\psi(\sigma_{i+1}) \ \
{\text {for}} \ \ 1\le i\le k-2\qquad {\text{and}} \qquad
\psi(\alpha_{ij})\ne 1 \ \
{\text {for}} \ \ 1\le i<j\le k\,.
$$
%If \ $1\le i<j\le k$, then $\psi(\alpha_{ij})\ne 1$.
Moreover, if the group $H$ is finite then \
$\ord \psi(\alpha_{ij})\equiv 0\,(\modl\, {j-i+1})$ \
whenever either $1\le i<j-1\le k-1\ne 3$ \
or \ $1\le i\le 2$, \ $j=i+2$ and $k=4$.
\end{Corollary}

\begin{Lemma}
\label{Lm: k|m if psi(alpha m)<->psi(beta),(k-1)|m if psi(beta m)<->psi(alpha)}
Let $\psi\colon B_k\to H$ be a non-cyclic group homomorphism and $m$ be
a natural number.
\vskip0.2cm

$a)$ If $\psi(\alpha^m)$ and $\psi(\beta)$ commute,
then either $k$ divides $m$ or $k=4$, $\GCD(m,4)=2$
and $\psi(\sigma_1)=\psi(\sigma_3)$.
\vskip0.2cm

$b)$ If $\psi(\alpha)$ and $\psi(\beta^m)$ commute, then $k-1$ divides $m$.
\end{Lemma}

\begin{proof} 
Set $\widehat\alpha=\psi(\alpha)$, \ $\widehat\beta=\psi(\beta)$.
Since $\alpha^{k} = \beta^{k-1}$, the element
$\widehat\alpha^k = \widehat\beta^{k-1}$
commutes with both $\widehat\alpha$ and $\widehat\beta$.
Notice that $k>2$ and $\widehat\alpha$ does not commute with
$\widehat\beta$ (for $\psi$ is non-cyclic).

$a)$ Assume that $k$ does not divide $m$; then $\nu\Def\GCD(m,k)<k$.
Since both $\widehat\alpha^m$ and $\widehat\alpha^k$ commute
with $\widehat\beta$, the element $\widehat\alpha^\nu$ also commutes with
$\widehat\beta$ and hence $\nu\ge 2$.
Moreover, $\nu\le k-2$, since $k>2$.
Relations (\ref{alpha beta}), (\ref{alpha pow k=beta pow(k-1)}) show now that
$$
\psi(\sigma_{\nu+1}) = \psi(\alpha^{\nu-1}\beta\alpha^{-\nu})
=\widehat\alpha^{\nu-1}\widehat\beta\widehat\alpha^{-\nu}
=\widehat\alpha^{-1}\widehat\beta = \psi(\alpha^{-1}\beta)
= \psi(\sigma_1),
$$
and Lemma
\ref{Lm: if k ne 4 and psi(sigma(i))=psi(sigma(j)) for i ne j then psi is cyclic}
implies that $k=4$, \ $\nu=2$ and $\psi(\sigma_3) = \psi(\sigma_1)$.
\vskip0.2cm

$b)$ Assume that $k-1$ does not divide $m$ and set
$\mu = \GCD(m,k-1)$. Since $\psi$ is non-cyclic
and $\widehat\beta^\mu$ commutes with $\widehat\alpha$, we have
%&
\begin{equation}\label{2.7} 
2\le\mu\le k-2\,.
\end{equation}
%& 
It follows from relations (\ref{conjugation by alphaij}) that
%&
\begin{equation}\label{eq: temporary 1}
\sigma_1\sigma_2\cdots \sigma_\mu\cdot \sigma_1 
= \sigma_2\cdot \sigma_1\sigma_2\cdots \sigma_\mu\,;
\end{equation}
%&
expressing $\sigma_1,...,\sigma_\mu$ in terms of 
$\alpha$ and $\beta$ according to (\ref{alpha pow k=beta pow(k-1)}), 
we can rewrite (\ref{eq: temporary 1}) in the form
$$
\aligned
(\alpha^{-1}\beta\cdot \beta\alpha^{-1}\cdot \alpha^1\beta&\alpha^{-2}
\cdot \alpha^2\beta\alpha^{-3}\cdots 
\alpha^{\mu-2}\beta\alpha^{-\mu+1})\cdot \alpha^{-1}\beta \\
&=\beta\alpha^{-1}\cdot 
(\alpha^{-1}\beta\cdot \beta\alpha^{-1}\cdot \alpha^1\beta\alpha^{-2}
\cdot \alpha^2\beta\alpha^{-3}\cdots 
\alpha^{\mu-2}\beta\alpha^{-\mu+1})\,,
\endaligned
$$
which leads to the relation
$\alpha^{-1}\beta^\mu\alpha^{-\mu}\beta = 
\beta\alpha^{-2}\beta^\mu\alpha^{-\mu+1}$.
Hence
$$
\widehat\alpha^{-1}\widehat\beta^\mu\widehat\alpha^{-\mu}\widehat\beta
= \widehat\beta\widehat\alpha^{-2}\widehat\beta^\mu\widehat\alpha^{-\mu+1}.
$$
Since $\widehat\beta^\mu$ commutes with
$\widehat\alpha$, it follows from the
latter relation that
$\widehat\alpha^{-(\mu+1)}\widehat\beta
= \widehat\beta\widehat\alpha^{-(\mu+1)}$,
that is, $\psi(\alpha^{\mu+1})$ commutes with $\psi(\beta)$. 
Because of (\ref{2.7}), $k$ cannot divide $\mu+1$; 
it follows from statement $(a)$ that $k=4$
and $\GCD(\mu+1,4)=2$ so that $\mu$ must be odd. 
But for $k=4$ (\ref{2.7}) means that $\mu=2$, and a contradiction ensues. 
\end{proof}

\begin{Remark}\label{Rmk 2.1} For any homomorphism
$\psi\colon B_k\to H$,
it follows from relations (\ref{conjugation by alpha power}) that
$\psi (\sigma_i)\sim\psi(\sigma_j) \ \ (1 \le  i,j < k)$. 
In particular, if $H={\mathbf S}(n)$, then all permutations
$\psi(\sigma_i)$ have the same cyclic type, that is,
$[\psi(\sigma_1)]=\ldots =[\psi(\sigma_{k-1})]$.\hfill$\bigcirc$
\end{Remark}
\vskip0.3cm

\noindent In the following lemmas
\ref{Lm: if |Inv(r,psi(sigma(1)))|=1 then psi is transitive}-\ref{Lm: if [psi(sigma(1))]=[3] or |Fix(psi(sigma(1))|>n-4 then psi is intransitive}
we assume that $k\ge 3$ and 
consider a homomorphism $\psi\colon B_k\to{\mathbf S}(n)$.

\begin{Lemma}
\label{Lm: if |Inv(r,psi(sigma(1)))|=1 then psi is transitive}
If $k>4$ and for some $r$, \ $1\le r<n$, the family 
$\Inv_r (\widehat\sigma_1)$ consists of a single set $\Sigma$, then
the homomorphism $\psi$ is intransitive.
\end{Lemma}

\begin{proof}
For $i\ne 2$, the mutually conjugate elements $\widehat\sigma_i$
and $\widehat\sigma_1$ commute;
by Lemma \ref{Lm: a single invariant set of length r}$(c)$,
$\Sigma$ is the only $\widehat\sigma_i$-invariant
set of cardinality $r$. It is certainly so for $i=4$;
since $\widehat\sigma_2$ commutes with $\widehat\sigma_4$, we see that
$\Sigma\in\Inv(\widehat\sigma_2)$. Thus, $\Sigma\in\Inv(\widehat\sigma_i)$
for all $i=1,\ldots ,k-1$ and $\psi$ is intransitive. 
\end{proof}

\noindent Parts $(a)$ and $(b)$ of the following lemma were proved
in \cite{Art47b} (for $n=k$).

\begin{Lemma}[\bcaps Artin Lemma on cyclic decomposition of $\widehat\sigma_1$]
\label{Lm: Artin Lemma on cyclic decomposition of hatsigma(1)}
\index{Artin Lemma!on cyclic decomposition of $\widehat\sigma_1$\hfill}
Let $k>4$ and let 
$\psi\colon B_k\to{\mathbf S}(n)$ be transitive.
Assume that the cyclic decomposition of $\widehat\sigma_1=\psi(\sigma_1)$
contains an $r$-cycle $C$.
\vskip0.2cm

$a)$ If $r>n/2$, then $r=n$ and $\psi$ is cyclic.
\vskip0.2cm

$b)$ If $n$ is even and $r=n/2$, then the cyclic decomposition
of $\widehat\sigma_1$ is of the form $\widehat\sigma_1= B_1\cdots B_s C$,
where all cycles $B_i$ are of the same length $t$, \ $2\le t\le r$, 
and $r=st$.
\vskip0.2cm

$c)$ If $n=2p$ with a prime $p$, then either 
$r<p$ or $r=2p=n$; in the latter case $\psi$ is cyclic.
\end{Lemma}

\begin{proof}
$a)$ Let $\Sigma=\supp C$ so that $\#\Sigma=r$.
It follows from the assumption $r>n/2$ that $\Sigma$ is the only 
$\widehat\sigma_1$-invariant set of cardinality $r$;
since $\psi$ is transitive, Lemma
\ref{Lm: k|m if psi(alpha m)<->psi(beta),(k-1)|m if psi(beta m)<->psi(alpha)} 
implies that $\Sigma={\boldsymbol\Delta}_n$, 
\ $r=n$, and $\widehat\sigma_1=C$. The permutations
$\widehat\sigma_3$ and $\widehat\sigma_4$ commute with
$\widehat\sigma_1$; hence each of them is a power
of the cycle $C$ and $\psi$ is cyclic (Lemma
\ref{Lm: if psi(sigma(i)) and psi(sigma(i+1)) commute then psi is cyclic}).
\vskip0.2cm

$b)$ If $\widehat\sigma_1=BC$, where $B$ is an $r$-cycle, then
assertion $(b)$ holds (with $s=1$, \ $t=r$).
So we may assume that $C$ is the only $r$-cycle in the cyclic
decomposition of $\widehat\sigma_1$. In this case for each $i\ne 2$
the set $\Sigma=\supp C$ is $\widehat\sigma_i$-invariant
and the restriction $\widehat\sigma_i|\Sigma=C^{q_i}$ for some integer
$q_{i}$, \ $0\le q_i<r$
(Lemma \ref{Lm: invariant sets and components of commuting permutations}$(b)$). 

Set $s=\GCD(q_4,r)$; let us show that $1<s<r$. If $s=1$, then $C^{q_4}$ 
is the only $r$-cycle in the cyclic decomposition of
$\widehat\sigma_4$; since $\widehat\sigma_2$ commutes
with $\widehat\sigma_4$, the set $\Sigma$ is also $\widehat\sigma_2$-invariant,
which contradicts the transitivity of $\psi$. If $s=r$, then 
$q_4=0$ and the restriction $\widehat\sigma_4|\Sigma=C^{q_4}=\id_\Sigma$
so that $\Sigma\subseteq\Fix(\widehat\sigma_4)$.
However, the cyclic decomposition of $\widehat\sigma_4$
must contain some $r$-cycle (for $\widehat\sigma_4\sim \widehat\sigma_1$).
Consequently, $\Sigma=\Fix(\widehat\sigma_4)$ and
hence $\Sigma\in\Inv(\widehat\sigma_2)$,
which contradicts the transitivity of $\psi$. 
Thus, $1<s<r$, \ $r = st$, \ $2\le t<r$ and $C^{q_{4}}$ is a
product of $s$ disjoint $t$-cycles. 

Since $\widehat\sigma_1\sim\widehat\sigma_4$ and
$C\preccurlyeq\widehat\sigma_1$, we obtain the desired representation of
$\widehat\sigma_1$.
\vskip0.2cm

$c)$ In view of $(a)$, we should only show that $r\ne p$.
Assume to the contrary that $r=p$; since $p$ is prime, 
$(b)$ implies that $\widehat\sigma_1=BC$, where $B$ and $C$
are disjoint $p$-cycles. 

Notice that $p\ne 2$. For otherwise $n=4$ and 
$[\widehat\sigma_1] = [\widehat\sigma_2]=[2,2]$;
however, in ${\mathbf S}(4)$ any two permutations of
cyclic type $[2,2]$ commute and, by
Lemma \ref{Lm: if psi(sigma(i)) and psi(sigma(i+1)) commute then psi is cyclic},
the homomorphism $\psi$ is cyclic.
Since $\psi$ is transitive, $\widehat\sigma_1$ must be
a $4$-cycle, and a contradiction ensues.

Thus $p\ge 3$. It follows from Lemma \ref{Lm: p cycles in S(2p)} that for each 
$i\ne 2$ there exist natural numbers 
$m_i, n_i$ \ ($1\le m_i,n_i<p$) such that 
$\widehat\sigma_i = B^{m_i}C^{n_i}$. (Cases $(ii)$, $(iii)$
described in Lemma \ref{Lm: p cycles in S(2p)} cannot occur here, since
$[\widehat\sigma_i] = [\widehat\sigma_1] = [p,p]$ and $p\ge 3$.)
Applying Lemma \ref{Lm: p cycles in S(2p)} to  
$\widehat\sigma_2$ and $\widehat\sigma_4$, we conclude that
$\widehat\sigma_2$ also is of the form $B^s C^t$, which contradicts the
transitivity of $\psi$. 
\end{proof}

\begin{Lemma}
\label{Lm: if [psi(sigma(1))]=[3] or |Fix(psi(sigma(1))|>n-4 then psi is intransitive} 
$a)$ If $k<n$ and $[\widehat\sigma_1]=[2]$,
then $\psi$ is intransitive.
\vskip0.2cm

$b)$ If $[\widehat\sigma_1]=[3]$ and $k>3$,
then $\widehat\sigma_1=\widehat\sigma_3$; consequently,
if in addition $k>4$, then $\psi$ is cyclic.
\vskip0.2cm
 
$c)$ If $4<k<n$ and $\#\Fix\,\widehat\sigma_1>n-4$,
then $\psi$ is intransitive.
\end{Lemma}

\begin{proof} 
$a)$ Clearly $[\widehat\sigma_i]=[2]$ for any $i$;
it is readily seen that in this case either $\psi$ is cyclic
and $\displaystyle\#(\cup_{i=1}^{k-1}\supp\,\widehat\sigma_i)=2<k<n$ or
$\psi$ is non-cyclic, $\widehat\sigma_i\infty\widehat\sigma_{i+1}$ and
$\displaystyle\#(\cup_{i=1}^{k-1}\supp\,\widehat\sigma_i)\le k<n$. Anyway
$\psi$ is intransitive.

$b)$ Since $[\widehat\sigma_i]=[3]$ for every $i$,
Lemma \ref{Lm: braid-like triple of 3 cycles} implies that 
$\widehat\sigma_1=\widehat\sigma_3$; if $k>4$, then $\psi$ is cyclic
by Lemma
\ref{Lm: if k ne 4 and psi(sigma(i))=psi(sigma(j)) for i ne j then psi is cyclic}.

$c)$ Let $m=\#\Fix\,\widehat\sigma_1$; so,
either $m=n$, or $m=n-2$, or $m=n-3$. If $m=n$, 
then $\widehat\sigma_1=\id$ and
$\psi$ is trivial. If $m=n-2$, then $[\widehat\sigma_1]=[2]$
and $\psi$ is intransitive by $(a)$. Finally,
if $m=n-3$, then $[\widehat\sigma_1]=[3]$
and, by $(b)$, $\psi$ is cyclic; in this case
all $\widehat\sigma_i$ coincide with the same $3$-cycle,
and $\psi$ is intransitive, since $n>k>4$.
\end{proof}

\noindent In the following lemma we show that the support of the permutation 
$\widehat\sigma_1=\psi(\sigma_1)$ must be relatively small provided
$\psi$ is transitive and all cycles that occur
in the cyclic decomposition of $\widehat\sigma_1$ have different
lengths. For a real $x\ge 0$, we denote the integral
part of $x$ by $E(x)$. 
\index{$E(x)$\hfill}

\begin{Lemma}
\label{Lm: condition no cycles of the same length restricts support} 
Let $k>4$ and let $\psi\colon B_k\to{\mathbf S}(n)$ be a transitive
homomorphism. Assume that $\widehat\sigma_1$ is a product of $\mu$
disjoint cycles of pairwise distinct lengths $r_1,...,r_\mu$
$(1<r_\nu<n$ for all $\nu=1,...,\mu)$.
Then each two permutations $\widehat\sigma_i$ and $\widehat\sigma_j$
with $|j-i|\ge 2$ are disjoint and hence
%&
\begin{equation}
\label{eq: upper bound for support of psi(sigma(1))}
\sum_{\nu=1}^\mu r_\nu \le n/E(k/2)\,.
\end{equation}
%&
\end{Lemma}

\begin{proof}
Each $\widehat\sigma_i$, \ $1\le i\le k-1$, also
is a disjoint product of $\mu$ cycles of pairwise distinct lengths
$r_1,...,r_\mu$. Let us fix some $i$ and consider the cyclic
decomposition $\widehat\sigma_i = C_1\cdots C_\mu$, \ $[C_\nu]=[r_\nu]$.
Let $\Sigma_\nu=\supp C_\nu$ and $\Sigma=\supp\widehat\sigma_i
=\underset{\nu=1}{\overset{\mu}\bigcup}\Sigma_\nu$;
each set $\Sigma_\nu$ is $\widehat\sigma_i$-invariant.
For every $j$ with $|j-i|\ge 2$
the permutation $\widehat\sigma_j$ commutes with $\widehat\sigma_i$;
since all $r_\nu$ are distinct, each $\Sigma_\nu$ is
$\widehat\sigma_j$-invariant and, by Lemma
\ref{Lm: invariant sets and components of commuting permutations}$(b)$,
the restriction $\widehat\sigma_j|\Sigma$ of $\widehat\sigma_j$
to $\Sigma$ is of the form
%&
\begin{equation}
\label{eq: restriction of hat(sigma(j)) to supp(hat(sigma(i)))}
\widehat\sigma_j|\Sigma = C_1^{q_{j,1}}\cdots C_\mu^{q_{j,\mu}}
\end{equation}
%&
with certain integers $q_{j,\nu}$, \ $0\le q_{j,\nu}<r_\nu$ .
\vskip0.2cm

\noindent Let us show first that $\widehat\sigma_j$ is disjoint with
$\widehat\sigma_i$ whenever $|j-i|>2$. To this end, it suffices to show
that $\widehat\sigma_j|\Sigma=\id_\Sigma$ for such $j$,
or, which is the same, that in 
(\ref{eq: restriction of hat(sigma(j)) to supp(hat(sigma(i)))})
all $q_{j,1}=...=q_{j,\mu}=0$.

Assume to the contrary that for some $j_\circ$ with
$|j_\circ-i|>2$ there is a non-zero $q_{j_\circ,\nu}$.
Then the permutation $D=C_\nu^{q_{j_\circ,\nu}}$ is an $r_\nu$-cycle
(for otherwise $D$ would be a product of a few cycles {\sl of the same length},
which is impossible since $D\preccurlyeq\widehat\sigma_{j_\circ}$).
Clearly $D$ is the only $r_\nu$-cycle in the cyclic decomposition
of $\widehat\sigma_{j_\circ}$. Since $|j_\circ-i|>2$,
the permutations $\widehat\sigma_s$
with $|s-i|=1$ commute with $\widehat\sigma_{j_\circ}$
so that the set $\supp D=\Sigma_\nu$ is $\widehat\sigma_s$-invariant.
However, $\Sigma_\nu$ is $\widehat\sigma_i$-invariant
and also $\widehat\sigma_j$-invariant for every
$j$ with $|j-i|\ge 2$. Thus, $\Sigma_\nu$ is invariant with
respect of all $\widehat\sigma_1,...,\widehat\sigma_{k-1}$ and hence
it is an $(\Img\psi)$-invariant set of cardinality $r_\nu$, $1<r_\nu<n$, 
which contradicts the transitivity of $\psi$.
\vskip0.2cm

\noindent We are left with the case when $|j-i|=2$.
Since $k>4$, there exists an index $t$,
\ $1\le t\le k-1$, neighboring to one of the indices $i,j$
and non-neighboring to another one. Because
the situation is symmetric with respect to $i$ and $j$,
we may assume that $|t-j|=1$ and $|t-i|\ge 2$.
Then the condition $|j-i|=2$ implies that in fact $|t-i|>2$ and,
as we have already proved, $\widehat\sigma_t$ is disjoint with
$\widehat\sigma_i$; hence $\widehat\sigma_t(m)=m$ for each $m\in\Sigma$.
Morever, for such $m$ we have also $\widehat\sigma_j(m)\in\Sigma$, for  
$\Sigma$ is $\widehat\sigma_j$-invariant.
In particular
$\widehat\sigma_t(\widehat\sigma_j(m))=\widehat\sigma_j(m)
=\widehat\sigma_j(\widehat\sigma_t(m))$.
Thus, taking into account that $|t-j|=1$, for each $m\in\Sigma$ we obtain
%&
\begin{equation}
\label{eq: widehat(sigma(j))=widehat(sigma(j))widehat(sigma(j)) on Sigma}
\aligned
\widehat\sigma_j(m)=\widehat\sigma_t(\widehat\sigma_j(m))
&=\widehat\sigma_t((\widehat\sigma_j(\widehat\sigma_t(m)))
=(\widehat\sigma_t\widehat\sigma_j\widehat\sigma_t)(m)\\
&=(\widehat\sigma_j\widehat\sigma_t\widehat\sigma_j)(m)
=\widehat\sigma_j(\widehat\sigma_t(\widehat\sigma_j(m)))
=\widehat\sigma_j(\widehat\sigma_j(m))\,.
\endaligned
\end{equation}
%&
Since $\widehat\sigma_j(\Sigma)=\Sigma$, 
(\ref{eq: widehat(sigma(j))=widehat(sigma(j))widehat(sigma(j)) on Sigma})
implies that $\widehat\sigma_j(m)=m$ for all $m\in\Sigma$ so that
$\widehat\sigma_j|\Sigma=\id_\Sigma$ and
the permutations $\widehat\sigma_j$ and $\widehat\sigma_i$ are disjoint.
\end{proof}

\subsection{Transitive homomorphisms $B_k\to{\mathbf S}(n)$
and prime numbers}
\label{Ss: 2.3. Transitive homomorphisms  B(k) to  S(n) and prime numbers}
The following lemma is the heart of Artin's methods developed in 
\cite{Art47b}.
It was not formulated explicitly but the proof
(for $k=n$) was given in the course of the proof
of Lemma 6 in the quoted paper. For completeness of the exposition
we present the proof of this very important lemma.

\begin{Lemma}[\bcaps Artin Fixed Point Lemma] 
\label{Lm: Artin Fixed Point Lemma} 
Let $k>4$ and $n$ be natural numbers. Suppose that there is
a prime $p>2$ such that
%&
\begin{equation}\label{2.8}
n/2 < p\le k-2\,.
\end{equation}
%&
Then for every non-cyclic transitive homomorphism
$\psi\colon B_k\to{\mathbf S}(n)$
the permutation $\widehat\sigma_1=\psi(\sigma_1)$
has at least $k-2$ fixed points; in particular, $n\ge k$.
\end{Lemma}

\begin{proof} 
The permutations $\widehat\alpha=\psi(\alpha)$ and
$\widehat\sigma_1$ generate the whole image
$\Img\psi\subseteq {\mathbf S}(n)$ of $\psi$.
For each $i$, \ $3\le i\le k-p+1$, set
$T_i = \widehat\alpha_{i,p+i-1}=\psi(\alpha_{i,p+i-1})$.
It follows from relations (\ref{alpha beta}) and (\ref{conjugation by alpha})
that $T_{i+1}=\widehat\alpha T_i \widehat\alpha^{-1}$. Thus, 
$T_3,...,T_{k-p+1}$ are conjugate to each other
and have the same cyclic type. 
\vskip0.2cm

Corollary
\ref{Crl: if psi is non-cyclic then psi(sigma(i)) infty psi(sigma(i+1))}
implies that $T_i\ne 1$ and $\ord T_i\equiv 0 \ (\modl p)$.
Since $p$ is prime and $p>n/2$, the cyclic decomposition of $T_i$
must contain a unique $p$-cycle $C_i$; let us show that its support
$\Sigma_i=\supp C_i$ is contained in $\Fix\widehat\sigma_1$.
Indeed, $\widehat\sigma_1$ commutes with $T_i$ and Lemma
\ref{Lm: invariant sets and components of commuting permutations}$(b)$
shows that $\widehat\sigma_1|\Sigma_i = C_i^{q_i}$
for a certain $q_i$, \ $0\le q_i < p$. In fact $q_i=0$;
for otherwise, $\GCD(q_i,p)=1$ and $C_i^{q_i}$ is 
a cycle of length $p>n/2$ that occurs in the cyclic decomposition of
$\widehat\sigma_1$, which contradicts
Lemma \ref{Lm: Artin Lemma on cyclic decomposition of hatsigma(1)}$(a)$.
Hence $\widehat\sigma_1|\Sigma_i=\id_{\Sigma_i}$.
\vskip0.2cm

\noindent For $3\le r\le k-p+1$, each union $S_r=\Sigma_3\cup...\cup\Sigma_r$
is contained in $\Fix\widehat\sigma_1$; in particular, $S_r$ is a non-trivial
$\widehat\sigma_1$-invariant set.
It suffices to show that $\#S_{k-p+1}\ge k-2$.
Since $\#S_3=\#\Sigma_3=p$, this will be achieved 
by showing that $\Sigma_{r+1}\not\subseteq S_r$ for all $r=3,...,k-p$
(if this is the case, then $\#S_{k-p+1}\ge p+(k-p-2)=k-2$). 

Suppose to the contrary that $\Sigma_{r+1}\subseteq S_r$ for
some $r$, \ $3\le r\le k-p$, so that $S_r = S_{r+1}$.
Since $C_j$ is the only $p$-cycle
in the cyclic decomposition of $T_j$ and
$T_{i+1}=\widehat\alpha T_i \widehat\alpha^{-1}$, we have
$\widehat\alpha C_i \widehat\alpha^{-1}=C_{i+1}$ and
$\widehat\alpha(\Sigma_i) = \Sigma_{i+1}$. 
Hence $\widehat\alpha(S_r)=\widehat\alpha(\Sigma_3\cup...\cup\Sigma_r)
=\Sigma_4\cup...\cup\Sigma_{r+1}\subseteq S_{r+1}=S_r$ and
$S_r$ is a non-trivial $(\Img\psi)$-invariant set,
which contradicts the transitivity of $\psi$.
\end{proof}

\begin{Remark}
\label{Rmk: bounds of Artin Fixed Point Lemma}
The mapping
$\sigma_1\mapsto (1,2)(3,4)(5,6)$, \
$\alpha\mapsto (1,2,3,4,5)$,
extends to a non-cyclic transitive homomorphism
$\psi_{5,6}\colon B_5\to{\mathbf S}(6)$.
This shows that the assertion of Lemma \ref{Lm: Artin Fixed Point Lemma}
becomes false if we replace the inequalities (\ref {2.8}) by the slightly
weaker inequalities $n/2\le p\le k-2$. Nevertheless
\vskip0.2cm

\noindent {\sl the conclusion of Artin Fixed Point Lemma holds true
whenever there is a prime $p>3$ that satisfies 
$n/2\le p\le k-3$.}
\vskip0.2cm

\noindent For $p>n/2$ this follows directly from 
Lemma \ref{Lm: Artin Fixed Point Lemma}.
Thus, to justify our assertion, we need only
to consider the case when $6<n = 2p\le 2(k-3)$.

Define $T_i$ as in the proof of Lemma \ref{Lm: Artin Fixed Point Lemma}.
As above, $\ord T_i\equiv 0 \ (\modl p)$ so that the cyclic decomposition
of $T_i$ must contain a cycle $C_i$ of length divisible by $p$;
but now we have $p=n/2$ instead of the strict inequality $p>n/2$
and hence either $(i)$ $T_i$ contains a single cycle of length $p$ or 
$(ii)$ $[T_i] = [p,p]$ or $(iii)$ $[T_i] = [2p]$.

In case $(i)$ the argument above has to be modified
just at one point: to show that all $q_i=0$, one should refer to
Lemma \ref{Lm: Artin Lemma on cyclic decomposition of hatsigma(1)}$(c)$ 
instead of
Lemma \ref{Lm: Artin Lemma on cyclic decomposition of hatsigma(1)}$(a)$.

Let us show that cases $(ii)$ and $(iii)$ cannot occur.
Notice that $k-p+1\ge 4$; therefore, we can deal with the permutation
$A=T_4$.

In case $(ii)$, $A=T_4=BC$, where $B$ and $C$ are disjoint
$p$-cycles. Let $B = (b_0,\ldots ,b_{p-1})$, \
$C = (c_0,\ldots ,c_{p-1})$. Since the permutations
$D=\widehat\sigma_1$ and $D'=\widehat\sigma_2$ commute with
$A=T_4=BC$, Lemma \ref{Lm: p cycles in S(2p)} applies to the couples
$A$, $D$ and $A$, $D'$ respectively.
Notice that $D$ and $D'$ are non-trivial, $D\sim D'$ and hence $[D]=[D']$;
so, if one of these permutations is of the form $(i)$ described in Lemma
\ref{Lm: p cycles in S(2p)}, then it contains a $p$-cycle;
however this contradicts
Lemma \ref{Lm: Artin Lemma on cyclic decomposition of hatsigma(1)}$(c)$.
The same lemma shows also that $[D]=[D']\ne [2p]$.
So $D$, $D'$ must be of the form $(ii)$ described in Lemma
\ref{Lm: p cycles in S(2p)}, that is,
$D=\widehat\sigma_1 = D_0\cdots D_{p-1}$ and
$D'=\widehat\sigma_2 = D_0'\cdots D_{p-1}'$, where
$D_i=(b_i,c_{\modo{i+r}_p})$, \
$D_j'=(b_j,\,c_{\modo{j+r'}_p})$, $r$ and $r'$ do not depend on $i,j$, and
$0\le r,r'<p$. Clearly
$(\widehat\sigma_1\widehat\sigma_2\widehat\sigma_1)(b_0) =
(\widehat\sigma_2\widehat\sigma_1\widehat\sigma_2)(b_0)$.
It is readily seen that the left hand side of the latter relation equals
$c_{\modo{2r-r'}_p}$, and the right hand side equals
$c_{\modo{2r'-r}_p}$. Consequently, $2r-r'\equiv 2r'-r \ (\modl \,p)$,
that is, $3(r-r')\equiv 0 \ (\modl p)$.
Since $p>3$ is prime, we have $r=r'$; therefore, $D_i = D_i'$ for all $i$.
Thus, $\widehat\sigma_1 =\widehat\sigma_2$ and $\psi$
is cyclic, which contradicts our assumption.

In case $(iii)$, $T_4$ is a $2p$-cycle.
Since $2p = n$ and both $\widehat\sigma_1$ and
$\widehat\sigma_2$ commute with $T_4$, they also commute with each other 
(Lemma \ref{Lm: invariant sets and components of commuting permutations}$(b)$),
and a contradiction ensues.
\vskip0.2cm

\noindent The homomorphism $\nu_6\colon B_6\to{\mathbf S}(6)$
(see Artin Theorem in Sec. \ref{Ss: Transitive homomorphisms B(k) to S(k)})
shows that the condition $p>3$ above is essential.
\hfill $\bigcirc$
\end{Remark}

\begin{Remark}
\label{Rmk: existence of desired primes} 
{\sl A prime $p\in ((k+l)/2,k-2]$
does certaily exist in each of the following three cases:
$(a)$ $6\ne k\ge 5$, \ $l=0$; \ $(b)$ $k\ge 7$, \ $l=1$; \
$(c)$ $8,12\ne k\ge 7$, \ $l=2$; \ $(d)$ $11,12\ne k\ge 9$, \ $l=3$.}
\vskip0.2cm

Case $(a)$ is known as ``Bertrand Postulate";
it was proven by P. L. Chebyshev in the 19th century.
In fact all cases may be treated on the ground of the following
inequality due to P. Finsler \cite{Fin45}
(see also \cite{Trost}, p. 60, Satz 32):
\ $\pi (2m) - \pi (m) > m/(3\log (2m))$ for all natural $m>1$,
where $\pi (x)$ is the number of all primes $p\le x$.
\end{Remark}

\vfill
\newpage

%Sec. 3
\section{Homomorphisms \ $B_k\to{\mathbf S}(n)$
\ and \ $B_k\to B_n$, \ \ $n\le k$}
\label{sec: B(k) to S(n), n<=k}

In this section we prove parts $(a)$ and $(b)$ of
Theorems A and Theorem B (see Theorem
\ref{Thm: homomorphisms B(k) to S(n), n<k} and Theorem
\ref{Thm: for k ne 4 psi(PB(k))<PB(k) for any non-integral psi} respectively). 
To prove Theorem \ref{Thm: homomorphisms B(k) to S(n), n<k},
we follow the methods of E. Artin; the main new point is that we
make extensive use of the fact that for $k>4$ the commutator subgroup
$B'_k$ is perfect (Sec. \ref{Ss: Commutator subgroup B'(k)}).
Theorem \ref{Thm: homomorphisms B(k) to S(n), n<k} applies in order to
improve Artin Theorem on homomorphisms $B_k\to{\mathbf S}(k)$
(see Theorem \ref{Thm: Improvement of Artin Theorem} below).

\begin{Theorem}\label{Thm: homomorphisms B(k) to S(n), n<k}
\index{Theorem!on homomorphisms $B_k\to{\mathbf S}(n)$, $n<k$\hfill}
\index{Homomorphisms $B_k\to{\mathbf S}(n)$, $n<k$\hfill}
Let $k\ne 4$ and $n<k$. Then

$a)$ every homomorphism $\psi\colon B_k\to{\mathbf S}(n)$ is cyclic;

$b)$ every homomorphism $\phi\colon B_k\to B_n$ is integral.
\end{Theorem}

\begin{proof} 
For $k<3$ both statements are trivial, since
$B_2\cong{\mathbb Z}$. Assume that $k>4$. 

$a)$ Suppose that there exists a non-cyclic homomorphism
$\psi\colon B_k\to{\mathbf S}(n)$ and
set $H = \Img\psi\subseteq {\mathbf S}(n)$.
Then there exists at least one $H$-orbit $Q\subseteq{\boldsymbol\Delta}_n$
of a certain length $m = \#Q\le n<k$ such that the reduction
$\psi_Q\colon B_k\to{\mathbf S}(Q) \cong {\mathbf S}(m)$ 
of the homomorphism $\psi$ to $Q$ is a non-cyclic transitive homomorphism
(see Observation at the end of Definition \ref{Def: symmetric groups}).
Since $k>4$ and $m<k$, Chebyshev Theorem
(Remark \ref{Rmk: existence of desired primes}$(a)$) supplies us with
prime $p>2$ such that $m/2<p\le k-2$ and 
Lemma \ref{Lm: Artin Fixed Point Lemma} implies that the permutation 
$\psi_Q(\sigma_1)\in{\mathbf S}(m)$ has at least $k-2>m-2$
fixed points. Hence $\psi_Q(\sigma_1)=\id_Q$
and $\psi_Q$ is trivial, which contradicts our choice of $Q$.
\vskip0.2cm

$b)$ Consider the composition  
$\psi = \mu\circ\phi\colon B_k\stackrel{\phi}{\longrightarrow}
B_n\stackrel{\mu}{\longrightarrow}{\mathbf S}(n)$ 
of $\phi$ with the canonical epimorphism $\mu\colon B_n\to{\mathbf S}(n)$.
By $(a)$, $\psi$ is cyclic and its restriction to $B'_k$ is trivial. 
Thus, $\phi(B'_k)\subseteq\Ker\mu=PB_n$. Being a perfect group $B'_k$
does not possess non-trivial homomorphisms to the pure braid group ${PB_n}$ 
(Corollary \ref{Crl: Perfect group has no homomorphisms to PB(k)}). 
Hence the restriction of $\phi$ to $B'_k$ is trivial and $\phi$ is integral. 
\end{proof}

\begin{Remark}
\label{Rmk: epimorphisms B(4) to B(3) and S(3)}
The condition $k\ne 4$ in Theorem \ref{Thm: homomorphisms B(k) to S(n), n<k}
is essential. To see this, take the canonical systems of generators
$\{\sigma_1,\sigma_2,\sigma_3\}$ in $B_4$ and $\{\sigma_1',\sigma_2'\}$
in $B_3$ and consider the surjective homomorphism $\pi\colon B_4\to B_3$
defined by 
$$
\pi(\sigma_1)=\pi(\sigma_3)=\sigma_1',\qquad
\pi(\sigma_2)=\sigma_2'.
$$
This example shows that for $k=4$ statement $(b)$
of Theorem \ref{Thm: homomorphisms B(k) to S(n), n<k} is false;
$(a)$ is also false, since the composition 
$\mu\circ\pi\colon B_4\to{\mathbf S}(3)$ of $\pi$ with the canonical
epimorphism $\mu\colon B_3\to{\mathbf S}(3)$ is surjective.
We call $\pi\colon B_4\to B_3$
the {\em canonical epimorphism}. It is easily seen that the kernel 
of $\pi$ coincides with the normal subgroup $\mathbf T\vartriangleleft B_4$
described in \ref{Ss: Commutator subgroup B'(k)}
(Gorin-Lin Theorem $(c)$). In terms of special generators
$\alpha,\beta\in B_4$ and $\alpha',\beta'\in B_3$,
the canonical epimorphism $\pi$ looks as follows:
\ $\pi\colon \alpha\mapsto\beta'$,
\ \ $\beta\mapsto(\beta')^{-1}(\alpha')^2\beta'$.
\hfill $\bigcirc$
\end{Remark}

\noindent Our next goal is to prove Theorem B. To this end,
we need some preparations. We start with some additional properties
of the pure braid group $PB_k$. In what follows, we use the notation
introduced in Sec. \ref{sec: Introduction}. In particular, 
$\sigma_1,\cdots \sigma_{k-1}$
are the canonical generators of $B_k$ and the elements
$\alpha_{ij}$, \ $1\le i<j\le k$, are defined by (\ref{alpha beta}).

\begin{Lemma} 
\label{Lm: alpha(1,t) pow t x sigma(t)x...x sigma(1)= alpha(1,t+1) pow t} 
Let $1<t\le k-1$. Then 
%&
\begin{equation}\label{3.1}
(\alpha_{1t})^t\cdot \sigma_t\sigma_{t-1}\cdots \sigma_1 
= (\alpha_{1,t+1})^t.
\end{equation}
%&
\end{Lemma}

\begin{proof} 
Let $q,r,s$ be integers such that
$0\le r\le s\le t$, \ $1\le q\le t$ and $q+r\le t$. 
We prove that
%&
\begin{equation}\label{3.2}
(\alpha_{1t})^s\cdot \sigma_t\sigma_{t-1}\cdots \sigma_q = 
(\alpha_{1t})^{s-r}\cdot \sigma_t\sigma_{t-1}
\cdots \sigma_{q+r}(\alpha_{1,t+1})^r.
\end{equation}
%&
The proof is by induction on $r$ with the case $r=0$ trivial.
Suppose that (\ref{3.2}) holds for some $r\ge 0$
such that $r+1\le s$ and $q+r+1\le t$. Notice that 
%&
\begin{equation}\label{3.3}
\alpha_{1t}\sigma_t = \sigma_1\cdots \sigma_{t-1}\cdot \sigma_t 
= \alpha_{1,t+1}\,.
\end{equation}
%&
Therefore, using relations (\ref{conjugation by alphaij}), we obtain
$$
\aligned
(\alpha_{1t})^s \cdot \sigma_t\sigma_{t-1}\cdots \sigma_q
&=(\alpha_{1t})^{s-r}\cdot \sigma_t\sigma_{t-1}
\cdots \sigma_{q+r}\cdot (\alpha_{1,t+1})^r\\
&=(\alpha_{1t})^{s-r-1}\alpha_{1t}\cdot \sigma_t\sigma_{t-1}
\cdots \sigma_{q+r}\cdot (\alpha_{1,t+1})^r\\
&=(\alpha_{1t})^{s-r-1}\alpha_{1,t+1}\cdot\sigma_{t-1}
\cdots \sigma_{q+r}\cdot(\alpha_{1,t+1})^r\\
&=(\alpha_{1t})^{s-r-1}\cdot\sigma_t
\cdots \sigma_{q+r+1}\alpha_{1,t+1}\cdot (\alpha_{1,t+1})^r\\
&=(\alpha_{1t})^{s-(r+1)}\cdot\sigma_t
\cdots \sigma_{q+(r+1)}\cdot (\alpha_{1,t+1})^{r+1}\,;
\endaligned
$$
this completes the inductive step and proves (\ref{3.2}). 
When $q=1$, \ $r=t-1$ and $s=t$ relation (\ref{3.2}) takes the form
$$
(\alpha_{1t})^t\cdot \sigma_t\sigma_{t-1}\cdots \sigma_1 = 
\alpha_{1t}\cdot \sigma_t\cdot (\alpha_{1,t+1})^{t-1}\,;
$$
in view of (\ref{3.3}), the latter relation coincides with (\ref{3.1}).
\end{proof}

\noindent Assume now that $1<r<k$. Let $s_{i,j}\in B_k$ ($1\le i<j\le k$) 
and $s'_{i,j}$ ($1\le i<j\le r$) be the canonical generators of the groups
$PB_k$ and $PB_r$ respectively. Consider the group epimorphism 
$\xi_{k,r}\colon PB_k\to PB_r$ defined by (\ref{homomorphism xi}) and set
$$
R_{t} = s_{1,t}s_{2,t}\cdots s_{t-1,t} \ \ (2\le t\le k),\qquad
R'_{t} = s'_{1,t}s'_{2,t}\cdots s'_{t-1,t}
\ \ (2\le t\le r)\,.
$$
 
\begin{Lemma}
\label{Lm: R(2)...R(t)=(sigma(1)...sigma(t-1)) pow t} 
For every integer $t$ that satisfies $2\le t\le k$
the following relations hold:
%&
\begin{equation}\label{3.4}
R_t = \sigma_{t-1}\sigma_{t-2}\cdots \sigma_2\sigma_1^2
\sigma_2\cdots \sigma_{t-2}\sigma_{t-1},\qquad
R_2 R_3\cdots R_t = (\sigma_1\sigma_2\cdots 
\sigma_{t-1})^t.
\end{equation}
%&
\end{Lemma}

\begin{proof} 
The proof is by induction on $t$ with the case
$t=2$ trivial ($R_2=s_{1,2}$ and $s_{1,2}= \sigma_1^2$).
Assume that relations (\ref{3.4}) hold for some $t$, \ $2\le t<k$.
Then 
$$
\aligned
R_{t+1} &= s_{1,t+1}s_{2,t+1}\cdots s_{t-1,t+1}s_{t,t+1}
= \sigma_t s_{1,t}\sigma_t^{-1}\cdot \sigma_t s_{2,t}\sigma_t^{-1}
\cdots \sigma_t s_{t-1,t}\sigma_t^{-1}\cdot \sigma_t^2\\
&=\sigma_t s_{1,t} s_{2,t}\cdots s_{t-1,t}\sigma_t
=\sigma_t R_t \sigma_t
= \sigma_t\sigma_{t-1}\sigma_{t-2}\cdots \sigma_2\sigma_1^2
\sigma_2\cdots \sigma_{t-2}\sigma_{t-1}\sigma_t
\endaligned
$$
and, according to (\ref{3.1}),
$$
\aligned
R_2 R_3\cdots R_tR_{t+1}&=(\sigma_1\sigma_2\cdots\sigma_{t-1})^t
\cdot\sigma_t\sigma_{t-1}\cdots \sigma_2\sigma_1^2
\sigma_2\cdots \sigma_{t-1}\sigma_t\\
&=\left[(\alpha_{1t})^t \cdot\sigma_t\sigma_{t-1}
\cdots\sigma_2\sigma_1\right]
\cdot\sigma_1\sigma_2\cdots \sigma_{t-1}\sigma_t\\
&=(\alpha_{1,t+1})^t\cdot \sigma_1\sigma_2\cdots \sigma_{t-1}\sigma_t
=(\alpha_{1,t+1})^{t+1}\,;
\endaligned
$$
this completes the inductive step and proves the lemma.
\end{proof}

\noindent Recall that $CB_m\subseteq PB_m$ is the cyclic subgroup of
$B_m$ generated by the element $A_m$.
This subgroup coincides with the center
of $B_m$ whenever $m>2$; clearly $CB_2=PB_2$.

\begin{Lemma}
\label{Lm: correspondence of generators of centers of braid groups} 
Let $1<r<k$. Then 
$$
\xi_{k,r}(A_k) = A_r\,;
$$
in particular, $\xi_{k,r}|CB_k\colon CB_k\to CB_r$ is an isomorphism.
\end{Lemma}

\begin{proof} 
According to Lemma \ref{Lm: R(2)...R(t)=(sigma(1)...sigma(t-1)) pow t},
%&
\begin{equation}\label{3.6}
\xi_{k,r}(A_k) = \xi_{k,r}\left((\sigma_1\cdots\sigma_{k-1})^k\right)
= \xi_{k,r}(R_2\cdots R_k)\,.
\end{equation}
%&
Clearly, $\xi_{k,r}(R_j) = R_j'$ for $j\le r$ and
$\xi_{k,r}(R_j)=1$ for $j>r$.
Using these relations and taking into account relations 
(\ref{3.6}), (\ref{3.4}) 
(the latter one for the elements $R'$'s and $\sigma'$'s), 
we obtain \
$\xi_{k,r}(A_k) = R'_2\cdots R'_r= (\sigma'_1\cdots \sigma'_{r-1})^r = A_r$.
\end{proof}

\begin{Corollary}
\label{Crl: PB(k)= PB(k)(2) x CB(k)}
For any $k\ge 3$, the pure braid group $PB_k$ is the direct product
of its subgroups $PB_k^{(2)}$ and $CB_k$.
\end{Corollary}

\begin{proof} 
We have already noted in Sec. \ref{Ss: Pure braid group} that the 
kernel of the epimorphism $\xi_{k,2}\colon PB_k\to PB_2=CB_2$
coincides with $PB_k^{(2)}$; so we have the exact sequence
$$
1\to PB_k^{(2)}\to {PB_k}\stackrel{\xi_{k,2}}{\longrightarrow} CB_2\to 1\,.
$$
By Lemma \ref{Lm: correspondence of generators of centers of braid groups},
$\xi_{k,2}$ carries the subgroup $CB_k$ onto $CB_2$ isomorphically; since
$CB_k$ is the center of $B_k$, this proves the lemma.
\end{proof}

\begin{Lemma}
\label{Lm: if psi: B(k) to S(k) is cyclic then every phi: Ker(psi) to PB(k) is integral}
Let $k\ne 4$ and let $G$ be the kernel of a cyclic homomorphism
$\psi\colon B_k\to{\mathbf S}(k)$. Then every homomorphism
$\phi\colon G\to{PB_k}$ is integral.
\end{Lemma}

\begin{proof} 
Assume first that $k>4$. Since $\psi$ is cyclic and $B'_k$ is perfect,
we have $G=\Ker\psi\supseteq B'_k$ and $B'_k\supseteq G'\supseteq (B'_k)'
=B'_k$; hence $G' = B'_k$ and $G'$ is perfect. By Corollary 
\ref{Crl: Perfect group has no homomorphisms to PB(k)}, 
$\phi(G') = \{1\}$ and the homomorphism $\phi$ is abelian.
Furthermore, $G/G' = G/B'_k$ is a subgroup of $B_k/B'_k\cong{\mathbb Z}$;
thus, $G/G'$ is cyclic and the homomorphism $\phi$ must be cyclic.
Since $PB_k$ is torsion free, $\phi$ is integral.
\vskip0.2cm

\noindent Now let $k=3$. By Corollary \ref{Crl: PB(k)= PB(k)(2) x CB(k)},
$PB_3\cong PB_3^{(2)}\times CB_3$. Let $\pi_1$ and $\pi_2$ be the projections
of $PB_3$ onto the first and the second factor
respectively. We know also that $CB_3\cong{\mathbb Z}$
and $PB_3^{(2)}\cong {\mathbb F}_2$ (Markov Theorem). 

If $\psi$ is trivial, then $G = B_3$ and the homomorphisms 
$\pi_1\circ \phi\colon B_3\to PB_3^{(2)}$ and
$\pi_2\circ \phi\colon B_3\to CB_3$
are integral (Remark \ref{Rmk: If G' is finitely generated then
any homomorphism G to F is integral}); consequently,
the homomorphism $\phi$ is abelian and therefore integral.

Now let $\psi$ be non-trivial. Then either $(a)$ $[\psi(\sigma_1)] = [2]$  or
$(b)$ $[\psi(\sigma_1)] = [3]$. Using Reidemeister-Schreier process, 
it is easy to show that in case $(a)$ the group $G = \Ker\psi$ is generated
by the elements $u = \sigma_2\sigma_1^{-1}$ and $v=\sigma_1^2$ 
that satisfy a single defining relation $(uvu)^2 = vuv$.
In case $(b)$ the group $G$ is generated by the elements
$u = \sigma_2\sigma_1^{-1}$, $x=\sigma_1\sigma_2\sigma_1^{-2}$,
and $y=\sigma_1^2\sigma_2$ that satisfy the defining relations
$yuy^{-1} = u^{-1}$, \ $yxy^{-1} = x^{-1}$.

In case $(a)$, applying 
Lemma \ref{Lm: if reduced f(u,uinv,v,vinv)=1 in a free group then uv=vu}
to the elements $u,v$ and to the homomorphisms
$\pi_1\circ\phi$, $\pi_2\circ\phi$, 
we obtain that the elements $\widehat{u}=\phi(u)$ and
$\widehat{v}=\phi(v)$ commute. Since they satisfy the relation
$(\widehat{u}\widehat{v}\widehat{u})^2
= \widehat{v}\widehat{u}\widehat{v}$,
we have $\widehat{u}^3=1$. But the group $PB_3$ is torsion free;
hence $\widehat{u} = 1$. Consequently, the subgroup
$\Img\phi\subset PB_3$ is generated by a single element $\widehat{v}$;
therefore $\phi$ is integral. 

Similarly, in case $(b)$ we obtain that
the element $\widehat{y} = \phi(y)$ commutes with the elements
$\widehat{u} = \phi(u)$ and $\widehat{x} = \phi(x)$.
Since $\widehat{y}\widehat{u}\widehat{y}^{-1} = \widehat{u}^{-1}$,
$\widehat{y}\widehat{x}\widehat{y}^{-1} = \widehat{x}^{-1}$ and the group
$PB_3$ is torsion free, we obtain
$\widehat{u} = \widehat{x} = 1$.  Thus, the group $\Img\phi$
is generated by a single element $\widehat{y}$ and $\phi$ is integral.
\end{proof}

\begin{Lemma}\label{Lm: non-surjective homomorphism B(k) to S(k) is cyclic} 
For $k\ne 4$, every non-surjective homomorphism
$\psi\colon B_k\to{\mathbf S}(k)$ is cyclic.
\end{Lemma}

\begin{proof} 
If $\psi$ is transitive, then the lemma follows directly from Artin Theorem.
So we may assume that $\psi$ is intransitive; then $\#Q<k$
for any $(\Img\psi)$-orbit $Q\subset{\boldsymbol\Delta}_k$.
Theorem \ref{Thm: homomorphisms B(k) to S(n), n<k} 
implies that the reduction of $\psi$ to any such orbit is cyclic;
hence $\psi$ itself is cyclic.
\end{proof}

\noindent Lemma \ref{Lm: non-surjective homomorphism B(k) to S(k) is cyclic}
shows that any non-cyclic homomorphism $\psi\colon B_k\to{\mathbf S}(k)$
is transitive provided $k\ne 4$. This implies the following useful improvement
of Artin Theorem:

\begin{Theorem}[\bcaps Improvement of Artin Theorem]
\label{Thm: Improvement of Artin Theorem}
\index{Artin Theorem!improvement of\hfill}
\index{Theorem!of Artin!improvement of\hfill}
\index{Improvement of Artin Theorem\hfill}
Statements $(a)$, $(b)$ $($and also $(d)$ whenever $k\ne 4)$
of Artin Theorem hold true for every non-cyclic homomorphism
$\psi\colon B_k\to{\mathbf S}(k)$ even without
the additional assumption that $\psi$ is transitive.
Thus, for $k\ne 4,6$ every non-cyclic homomorphism
$\psi\colon B_k\to{\mathbf S}(k)$ is conjugate
to $\mu_k$ and any non-cyclic homomorphism
$\psi\colon B_6\to{\mathbf S}(6)$ is conjugate
either to $\mu_6$ or to $\nu_6$.
\hfill $\square$
\end{Theorem}

\begin{Lemma}
\label{Lm: for k ne 4 every B(k) to B'(k) is integral}
For $k\ne 4$ every homomorphism $\phi\colon B_k\to B'_k$ is integral.
\end{Lemma}

\begin{proof} 
If $k = 3$, the lemma follows from Remark
\ref{Rmk: If G' is finitely generated then any homomorphism G to F is integral},
since $B'_3\cong {\mathbb F}_2$.
Let $k>4$ and let $\mu'\colon B'_k\to{\mathbf S}(k)$
be the canonical homomorphism. Consider the composition
$$
\psi=\mu'\circ\phi\colon B_k\stackrel{\phi}{\longrightarrow} B'_k
\stackrel{\mu'}{\longrightarrow} {\mathbf S}(k)\,.
$$
Clearly $\Img\psi\subseteq\Img\mu'={\mathbf A}(k)$; hence
$\psi$ is non-surjective and,
by Lemma \ref{Lm: non-surjective homomorphism B(k) to S(k) is cyclic},
cyclic. Consequently, $(\mu'\circ\phi)(B'_k) = \psi(B'_k) = \{1\}$
and $\phi (B'_k)\subseteq\Ker\mu'\subseteq PB_k$.
Since $B'_k$ is perfect, $\phi(B'_k) = \{1\}$ and $\phi$ is cyclic.
\end{proof}

\begin{Remark}
\label{Rmk 3.3} 
If $k>4$ and the restriction $\phi'$ of a homomorphism $\phi\colon B_k\to G$
to $B'_k$ is abelian, then $\phi(B'_k) =\{1\}$ (since $B'_k$ is perfect)
and $\phi$ is cyclic. For $k=3$ or $4$ this is not necessarily so;
for instance, the canonical epimorphism $\mu\colon B_3\to{\mathbf S}(3)$
carries $B'_3$ onto the cyclic group
$\mathbf A(3)\cong{\mathbb Z}/3{\mathbb Z}$. Moreover, the natural projection 
$\phi\colon B_3\to B_3/(B'_3)'$ is non-abelian and
$\phi(B'_3)\cong{\mathbb Z}\oplus{\mathbb Z}$; so, even if we assume that $G$
is torsion free, this will not save the situation. 
The following statement nevertheless holds.
\end{Remark}

\begin{Lemma} 
\label{Lm: k=3,4; if restriction phi' of phi: B(k)->G is integral then phi is integral} 
Let $k=3$ or $4$ and let $\phi'\colon B'_k\to G$ be the restriction of  
a homomorphism $\phi\colon B_k\to G$. 
If $\phi'$ is integral, then it is trivial and $\phi$ is integral.
\end{Lemma}

\begin{proof} 
Assume first that $k=3$. Let $u = \sigma_2\sigma_1^{-1}$ and  
$v = \sigma_1\sigma_2\sigma_1^{-2}$ be the canonical 
generators of the group $B'_3$ (\ref{generators of B'_3}).
It is easily seen that relations $\sigma_1 u\sigma_1^{-1} = v$ and 
$\sigma_1 v\sigma_1^{-1} = u^{-1}v$ hold true; so 
Lemma \ref{Lm: non-zero det(M2-I)->integral psi is trivial}
applies to the homomorphism 
$\psi=\phi'$ (the matrix $M=\begin{pmatrix} 0 & 1\\ -1 & 1 
                                      \end{pmatrix}$
satisfies the conditions $\det M=1$ and $\det (M^2-I)=3\ne 0$).
Therefore, $\phi' = 1$ and $\phi$ is integral.

Now let $k = 4$ and let ${\mathbf T} = \Ker\pi$ be the kernel
of the canonical epimorphism $\pi\colon B_4\to B_3$ (see Remark \ref{3.1}).
Since $\phi'$ is integral, $\Ker\phi\supseteq\Ker\phi'\supseteq (B'_4)'$
(the second commutator subgroup of the group $B_4$). 
By Gorin-Lin Theorem $(c)$, $\mathbf T$ coincides with the intersection
of the lower central series of the group $B'_4$.
Hence $\mathbf T\subset (B'_4)'\subseteq \Ker\phi$
and $\phi$ may be decomposed as
$$
\phi = \psi\circ\pi\colon B_4\stackrel{\pi}{\longrightarrow}
B_4/{\mathbf T}\cong B_3\stackrel{\psi}{\longrightarrow} G\,,
$$
with a certain homomorphism $\psi$. Let $\psi'$ be the restriction of $\psi$
to $B'_3$. Then $\psi'(B'_3)=\psi(\pi(B'_4))=\phi'(B'_4)$.   
Hence $\psi'$ is integral. Since for $k=3$ the lemma was already proved,
$\psi$ is integral and the original homomorphism 
$\phi=\psi\circ\pi$ is integral as well.
\end{proof}

\noindent The following simple lemma provides us with a very useful tool
for the study of endomorphisms of $B_k$.

\begin{Lemma}
\label{Lm: if k ne 4 and psi=mu circ phi: B(k) to B(k) to S(k) is cyclic then phi is integral}
Let $k\ne 4$ and let $\phi$ be an endomorphism of the group $B_k$ such that  
the composition 
$$
\psi=\mu\circ\phi\colon B_k\stackrel{\phi}{\longrightarrow} B_k
\stackrel{\mu}{\longrightarrow}{\mathbf S}(k)
$$
of $\phi$ with the canonical epimorphism $\mu$ is cyclic.
Then $\phi$ is integral.
\end{Lemma}

\begin{proof} 
Let $G = \Ker\psi$.
Then $(\mu\circ\phi)(G)=\psi(G) = \{1\}$ and
$\phi(G)\subseteq \Ker\mu = {PB_k}$.
By Lemma
\ref{Lm: if psi: B(k) to S(k) is cyclic then every phi: Ker(psi) to PB(k) is integral},
the homomorphism $\phi|G\colon G\to{PB_k}$ is integral. 
Since $\psi$ is cyclic, $G = \Ker\psi\supseteq B'_k$; hence the restriction
$\phi'$ of $\phi$ to $B'_k$ is integral. For $k>4$,
the group $B'_k$ is perfect, which implies that $\phi'$ is trivial and $\phi$
is integral. For $k=3$ the same conclusion follows from Lemma
\ref{Lm: k=3,4; if restriction phi' of phi: B(k)->G is integral then phi is integral}. 
\end{proof}

\noindent Let $\chi\colon B_k\to{\mathbb Z}$ be the canonical integral
projection (see Sec. \ref{Ss: Commutator subgroup B'(k)}({\bf a)}). 

\begin{Lemma} 
\label{Lm: Ker psi is contained in B'(k) for non-trivial psi: B(k) to B(k)} 
$\Ker\phi\subseteq B'_k$ for any non-trivial endomorphism
$\phi\colon B_k\to B_k$. Moreover, if $k\ne 4$ and $\phi$ is non-integral,
then $\phi^{-1}(B'_k) = B'_k$.
\end{Lemma}

\begin{proof} 
If $\phi$ is integral, then $\Ker\phi=B'_k$ (for $B_k$ is torsion free).
Assume that $\phi$ is non-integral. 
\vskip0.2cm

If $k=4$, then there is a non-trivial homomorphism 
$\eta\colon \Img\phi\to{\mathbb Z}$ 
(Corollary \ref{Crl: subgroups of B3 and B4}).
The composition 
$$
\xi=\eta\circ\phi\colon B_4\stackrel{\phi}{\longrightarrow} B_4
\stackrel{\eta}{\longrightarrow} {\mathbb Z}
$$
is also non-trivial; so $\Ker\xi=B'_4$ and $\Ker\phi\subseteq\Ker\xi=B'_4$.
\vskip0.2cm

Finally, if $k\ne 4$, then
Lemma \ref{Lm: for k ne 4 every B(k) to B'(k) is integral}
implies that the subgroup $\Img\phi\subseteq B_k$ is not contained
in $B'_k$; therefore, the composition
$$
\zeta=\chi\circ\phi\colon B_k\stackrel{\phi}{\longrightarrow} B_k
\stackrel{\chi}{\longrightarrow} {\mathbb Z}
$$
is non-trivial. Hence $B'_k=\Ker\zeta=\Ker(\chi\circ\phi)
=\phi^{-1}(\Ker\chi)=\phi^{-1}(B'_k)$ and $\Ker\phi\subseteq\Ker\zeta=B'_k$.
\end{proof}

\noindent Now we prove our Theorem B, which gives an essential
strengthening of E. Artin theorem on automorphisms of $B_k$,
\cite{Art47b}.

\begin{Theorem} 
\label{Thm: for k ne 4 psi(PB(k))<PB(k) for any non-integral psi}
If $k\ne 4$, then $\phi(PB_k)\subseteq PB_k$, \ 
$\phi^{-1}(PB_k) = PB_k$ \ and $\Ker\phi\subseteq J_k$
for every non-integral endomorphism $\phi\colon B_k\to B_k$.
\end{Theorem}

\begin{proof} 
According to Lemma
\ref{Lm: if k ne 4 and psi=mu circ phi: B(k) to B(k) to S(k) is cyclic then phi is integral},
the composition
$$
\psi=\mu\circ\phi\colon B_k\stackrel{\phi}{\longrightarrow} B_k
\stackrel{\mu}{\longrightarrow} {\mathbf S}(k)
$$ 
is non-cyclic.
By Lemma \ref{Lm: non-surjective homomorphism B(k) to S(k) is cyclic},
$\psi$ is surjective
and hence transitive. It follows from Artin Theorem
(Sec. \ref{Ss: Transitive homomorphisms B(k) to S(k)})
that $\Ker\psi = {PB_k}$. Thus,
$$
{PB_k} = \Ker\psi = \phi^{-1}(\Ker\mu) = \phi^{-1}({PB_k})
$$ 
and $\phi({PB_k}) = \phi(\phi^{-1}({PB_k}))
\subseteq {PB_k}$. Moreover, 
$\Ker\phi\subseteq \Ker\,(\mu\circ\phi)
= \Ker\psi = {PB_k}$. Finally,
by Lemma \ref{Lm: Ker psi is contained in B'(k) for non-trivial psi: B(k) to B(k)},
$\Ker\phi\subseteq B'_k$ and hence $\Ker\phi\subseteq {PB_k}\cap B'_k = J_k$.
\end{proof}

\begin{Remark}\label{Rmk 3.4} 
It follows from relations
(\ref{generators of B'_3}), (\ref{eq: commut}), (\ref{eq: braid like}) that
$$
\aligned
(\sigma_1\sigma_2)\cdot (\sigma_3\sigma_2)\cdot (\sigma_1\sigma_2)
&=(\sigma_1\sigma_2\sigma_3)\cdot \sigma_2\sigma_1\sigma_2
=\sigma_3\sigma_2\sigma_3\cdot (\sigma_1\sigma_2\sigma_3)\\
&=(\sigma_3\sigma_2)\cdot \sigma_1\cdot (\sigma_3\sigma_2\sigma_3)
=(\sigma_3\sigma_2)\cdot \sigma_1\cdot (\sigma_2\sigma_3\sigma_2)\\
&\qquad\qquad\qquad\qquad\qquad
=(\sigma_3\sigma_2)\cdot (\sigma_1\sigma_2)\cdot (\sigma_3\sigma_2).
\endaligned
$$
Therefore, we can define an endomorphism $\phi$ of the group
$B_4$ by
$$
\phi(\sigma_1)=\phi(\sigma_3) = \sigma_1\sigma_2, \qquad
\phi(\sigma_2) = \sigma_3\sigma_2.
$$
This endomorphism is non-abelian (for
$\phi(\sigma_1)\ne \phi(\sigma_2)$), but
$\phi(\sigma_1^2) = (\sigma_1\sigma_2)^2 \not\in PB_4$;
thus $\phi(PB_4)\not\subseteq PB_4$. Besides,
$\Ker\phi = \mathbf T\not\subseteq J_4$. Moreover,
$\phi(\sigma_1^3)=(\sigma_1\sigma_2)^3\in PB_4$, which
shows that $\sigma_1^3\in\phi^{-1}(PB_4)$; but
$\sigma_1^3\not\in PB_4$, and therefore
$\phi^{-1}(PB_4)\not\subseteq PB_4$.
This example shows that the condition $k\ne 4$ in Theorem
\ref{Thm: for k ne 4 psi(PB(k))<PB(k) for any non-integral psi}
is essential.

For any $k\ge 3$, there is an {\it integral} endomorphism
$\phi\colon B_k\to B_k$ with
$\phi({PB_k})\not\subseteq{PB_k}$.
Indeed, take %an element
$c\in B_k$ such that $c\not\in{PB_k}$
and define $\phi$ by $\phi(\sigma_1) =\ldots = \phi (\sigma_{k-1}) = c$.
\hfill $\bigcirc$
\end{Remark}
\vskip0.3cm

The rest of this section is devoted to some results on
endomorphisms of the groups $B_3$, \ $B_4$
and on homomorphisms from $B_4$ into
$B_3$ and ${\mathbf S}(3)$.

\begin{Theorem}\label{Thm 3.13} Any non-integral endomorphism $\phi$ of \
$B_3$ is an embedding.
\end{Theorem}

\begin{proof} 
Since $\phi(B'_3)\subseteq B'_3$,
the restriction $\phi'$ of $\phi$ to $B'_3$
may be regarded as an endomorphism of the group
$B'_3\cong{\mathbb F}_2$.
The image $G = \Img\phi'$ is a free group
of rank $r\le 2$. Since $\phi$ is non-integral, Lemma
\ref{Lm: k=3,4; if restriction phi' of phi: B(k)->G is integral then phi is integral}
implies that $\phi'$ is non-integral; hence $G\cong {\mathbb F}_2$
and $\phi'\colon B'_3\to G$
is an isomorphism. By Lemma
\ref{Lm: Ker psi is contained in B'(k) for non-trivial psi: B(k) to B(k)},
$\Ker\phi\subset B'_3$. Thus, $\Ker\phi = \Ker\phi' = \{1\}$.
\end{proof}

\begin{Remark}\label{Rmk 3.5} 
For any $k\ge 2$, there exist
{\sl proper embeddings} $B_k\hookrightarrow B_k$.
For $k = 2$ this is evident. If $k>2$, take any $m\in{\mathbb Z}$ and define
the endomorphism $\phi_{k,m}$ by
$$
\phi_{k,m}\colon B_k\ni g\mapsto
(A_k)^{m\chi (g)}\cdot g\in B_k,
$$
where $\chi\colon B_k\to{\mathbb Z}$
is the canonical integral projection.
For any $g\in B'_k$, we have $\chi(g)=0$ and
$\phi_{k,m}(g)=g$; in particular, $\phi_{k,m}$ is non-integral.
By Lemma
\ref{Lm: Ker psi is contained in B'(k) for non-trivial psi: B(k) to B(k)},
$\Ker\phi_{k,m}\subseteq B'_k$; hence $\phi_{k,m}$ is an embedding.
If $h=\phi_{k,m}(g)\in \Img\phi_{k,m}$, then
$$
\chi(h) = \chi\left(\phi_{k,m}(g)\right)
= \chi\left((A_k)^{m\chi(g)}\cdot g\right)
=(mk(k-1)+1)\chi(g);
$$
consequently, $\chi(h)$ is divisible by the number $s(m) = mk(k-1)+1$. 
On the other hand, if $h\in B_k$ and $\chi(h)=ts(m)$ for some
$t\in{\mathbb Z}$, take $g=(A_k)^{-tm}\cdot h$;
then 
$$
\chi(g)=-tmk(k-1)+t(mk(k-1)+1)=t,
$$ 
and hence
$$
\phi_{k,m}(g)=(A_k)^{m\chi(g)}\cdot g = 
(A_k)^{mt}\cdot (A_k)^{-tm}\cdot h =h.
$$
It follows that the image of $\phi_{k,m}$ coincides
with the normal subgroup $\chi^{-1}(s(m){\mathbb Z})\subseteq B_k$. 
If $m\ne 0$, then $s(m)\ne \pm 1$, 
$\Img\phi_{k,m} = \chi^{-1}(s(m){\mathbb Z})\ne B_k$,
and $\phi_{k,m}$ is a proper embedding. 
\end{Remark}
\vskip0.3cm

Let $\pi\colon B_4\to B_3$
and $\mu\colon B_3\to{\mathbf S}(3)$
be the canonical epimorphisms.

\begin{Theorem}\label{Thm 3.14} 
$a)$ Any non-cyclic homomorphism
$\psi\colon B_4\to{\mathbf S}(3)$
is conjugate to the composition 
$\mu\circ\pi\colon B_4\stackrel{\pi}{\longrightarrow}
B_3\stackrel{\mu}{\longrightarrow} {\mathbf S}(3)$.

$b)$ Let $\phi\colon B_4\to B_3$
be a non-integral homomorphism. Then there exists a 
monomorphism $\xi\colon B_3\to B_3$ such that
$$
\phi = \xi\circ\pi\colon B_4\stackrel{\pi}{\longrightarrow}
B_3\stackrel{\xi}{\longrightarrow} B_3.
$$
In particular $\Ker\phi = \mathbf T$. Moreover,
if $\phi$ is surjective, then $\xi$ is an automorphism of
the group $B_3$.
\end{Theorem}

\begin{proof} 
$a)$ Clearly, $\psi(B'_4)\subseteq
{\mathbf S}'(3) = \mathbf A(3)\cong{\mathbb Z}/3{\mathbb Z}$.
Consequently,
$$
\Ker\psi\supseteq (B'_4)'\supset \mathbf T
= \Ker\pi.
$$
Therefore, there exists a
homomorphism $\phi\colon B_3\to{\mathbf S}(3)$
such that $\psi = \phi\circ\pi$. Since $\psi$ is non-cyclic,
$\phi$ is non-cyclic too; by Artin Theorem,
$\phi$ is conjugate to $\mu$. Hence $\psi$ is conjugate
to the composition $\mu\circ\pi$.

$b)$ Since $\phi(B'_4)\subseteq B'_3$
and $\mathbf T = \Ker\pi$ is the intersection of the
lower central series of the group $B'_4$,
the image $\phi(\mathbf T)$ is contained in the
intersection $H$ of the lower central series of the group
$B'_3$. But $B'_3\cong {\mathbb F}_2$,
and thus $H = \{1\}$. Consequently,
$\phi(\mathbf T) = \{1\}$ and
$\Ker\pi = \mathbf T\subseteq \Ker\phi$.
Therefore, there exists an endomorphism $\xi$ of the group
$B_3$ such that $\phi = \xi\circ\pi$. Since $\phi$
is non-integral, $\xi$ is also non-integral;
according to Theorem \ref{Thm 3.13}, $\xi$ is injective. Hence,
$\Ker\phi = \pi^{-1}(\Ker\xi) =
\pi^{-1}(\{1\}) = \Ker\pi = \mathbf T$.

Finally, if $\phi$ is surjective, $\xi$ is surjective, too.
Hence $\xi$ is an automorphism of $B_3$.
\end{proof}

\begin{Theorem}\label{Thm 3.15} $\Ker\phi = \mathbf T$ for any
non-integral noninjective endomorphism $\phi$ 
of $B_4$.
\end{Theorem}

\begin{proof} 
Let $\psi = \pi\circ\phi$, where 
$\pi\colon B_4\to B_3$ is
the canonical epimorphism. The homomorphism $\psi$ is 
non-trivial (for otherwise $\Img\phi\subseteq \Ker\pi =
\mathbf T\cong {\mathbb F}_2$ and, 
by Remark
\ref{Rmk: If G' is finitely generated then any homomorphism G to F is integral},
$\phi$ is integral).
Consider the following two cases:
$a)$ $\psi$ is integral, and
$b)$ $\psi$ is non-integral.

$a)$ In this case $\phi(B'_4)\subseteq 
\Ker\pi = \mathbf T\cong {\mathbb F}_2$, and
(as in the proof of Theorem \ref{Thm 3.14}$(b)$) we obtain a 
homomorphism $\xi\colon B_3\to B_4$
such that $\phi = \xi\circ\pi$. Since $\phi$ is non-integral, 
$\xi$ is non-integral. By Lemma 
\ref{Lm: k=3,4; if restriction phi' of phi: B(k)->G is integral then phi is integral}, 
the restriction $\xi'$ of $\xi$ to $B'_3$ 
is non-integral. It is easily seen that 
$\xi'(B'_3)\subseteq \mathbf T$.
Consequently, $\Img\xi'\cong {\mathbb F}_r$,
where $r\le 2$. Since $\xi'$ is non-integral, $r = 2$.
Thus, we obtain a surjective endomorphism 
${\mathbb F}_2\cong B'_3\to 
\Img\xi'\cong {\mathbb F}_2$
of the Hopfian group ${\mathbb F}_2$; hence 
$\xi'$ must be injective. On the other
hand, it is easy to check that 
$\Ker\xi\subseteq B'_3$.
Thus, $\Ker\xi = \Ker\xi'$ and $\xi$ is injective. Therefore,
$$
\Ker\phi = \Ker\,(\xi\circ\pi)
= \pi^{-1}(\Ker\xi) = \pi^{-1}(\{1\})
= \Ker\pi = \mathbf T.
$$

$b)$ $\mathbf T$ is a completely characteristic subgroup 
of the group $B_4$. Hence $\phi(\mathbf T)\subseteq \mathbf T$. 
Let $\widetilde\phi\colon \mathbf T\to \mathbf T$
be the restriction of $\phi$ to $\mathbf T$. Since $\psi$
is non-integral, Theorem \ref{Thm 3.14}$(b)$ implies 
that $\Ker\psi = \mathbf T$; hence $\Ker\phi\subseteq \mathbf T$ and
$\Ker\phi = \Ker\widetilde\phi$. Clearly,
$\widetilde\phi(\mathbf T)\cong {\mathbb F}_r$ where $r\le 2$.
If $r = 2$, then $\widetilde\phi$ is injective
(since $\mathbf T\cong {\mathbb F}_2$ is Hopfian), and $\phi$  
is injective too. Finally, if $r<2$, then $\widetilde\phi$ is integral.
In this case it follows from relations (\ref{1.14}), (\ref{1.15})
and Lemma \ref{2.2} that the homomorphism 
$\widetilde\phi$ is trivial; hence $\Ker\phi
= \Ker\widetilde\phi = \mathbf T$. 
\end{proof}

\newpage

%Sec. 4
\section{Transitive homomorphisms
$B_k\to{\mathbf S}(n)$ for small $k$ and $n$}
\label{sec: Transitive homomorphisms of B(k) to S(n) for small k and n}

\noindent For $k$ large enough, transitive homomorphisms 
$B_k\to{\mathbf S}(n)$, \ $k<n\le 2k$, 
can be studied using some general methods based mainly on 
Lemma \ref{Lm: Artin Fixed Point Lemma},
Theorem \ref{Thm: homomorphisms B(k) to S(n), n<k}, Gorin-Lin Theorem and 
the techniques which will be developed in 
Sec. \ref{sec: Retraction of homomorphisms; homomorphisms and cohomology} and
Sec. \ref{sec: Omega-homomorphisms and cohomology: some computations}
However for small $k$ these methods do not work. For this reason 
we consider the case of small $k$ in this section.
\vskip0.2cm

\noindent Given a group homomorphism 
$\psi\colon B_k\to H$,
we denote the $\psi$-images of the canonical generators 
$\sigma_1,...,\sigma_{k-1}$ and of the corresponding special
generators 
$\alpha=\sigma_1\cdots\sigma_{k-1}, \ \beta=\alpha\sigma_1$
by $\widehat\sigma_1,...,\widehat\sigma_{k-1}$
and $\widehat\alpha, \widehat\beta$,
respectively. Assume that the group $H$ is finite and the homomorphism
$\psi$ is non-cyclic. Then it follows from Lemma
\ref{Lm: k|m if psi(alpha m)<->psi(beta),(k-1)|m if psi(beta m)<->psi(alpha)}
that $\ord\widehat\beta$ is divisible by $k-1$.
Moreover, $\ord \widehat\alpha$
is divisible by $k$ whenever $k\ne 4$; if $k=4$, then either 
$\ord\widehat\alpha$ is divisible by $4$ or
$\widehat\sigma_1 = \widehat\sigma_3$ and $\ord\widehat\alpha$
is divisible by $2$ (but not by $4$). 
The following proposition follows immediately from
these remarks.
 
\begin{Proposition}\label{Prp 4.1} 
Let $4\le n\le 7$. Then any non-cyclic transitive homomorphism
$\psi\colon B_3\to{\mathbf S}(n)$
is conjugate to one of the following homomorphisms
$\psi^{(i)}_{3,n}$:

\

$a)$ $n=4:$ \ $\psi^{(1)}_{3,4}\colon 
\alpha\mapsto (1,2,3), \ \beta\mapsto (1,4)$; \ \
$\psi^{(2)}_{3,4}\colon 
\alpha\mapsto (1,2,3), \ \beta\mapsto (1,2)(3,4)$

\ \ \ \ $(\Img \psi^{(1)}_{3,4} = {\mathbf S}(4), \
\Img \psi^{(2)}_{3,4} = \mathbf A(4))$.

\

$b)$ $n=5:$ \
$\psi_{3,5}\colon \ \,\alpha \mapsto (1,2,3), \ \
\beta \mapsto (1,4)(2,5)$ \ \ $(\Img \psi_{3,5} = \mathbf A(5))$.

\

$c)$ $n=6:$ 
$$
\aligned
&\ \ \ \ \psi^{(1)}_{3,6}\colon 
\alpha\mapsto (1,2,3)(4,5,6), \ \ \beta\mapsto (1,2)(3,4)(5,6), \ \
\sigma_1\mapsto (2,3,6,4);\\
&\ \ \ \ \psi^{(2)}_{3,6}\colon 
\alpha\mapsto (1,2,3)(4,5,6),\ \
\beta\mapsto (1,4)(2,6)(3,5),\ \
\sigma_1\mapsto (1,6)(2,5)(3,4);\\
&\ \ \ \ \psi^{(3)}_{3,6}\colon 
\alpha\mapsto (1,2,3)(4,5,6),\ \
\beta\mapsto (1,2)(3,4),\ \ \ \qquad
\sigma_1\mapsto (2,3,6,5,4);\\
&\ \ \ \ \psi^{(4)}_{3,6}\colon 
\alpha\mapsto (1,2,3)(4,5,6),\ \
\beta\mapsto (1,4)(2,5),\ \ \ \qquad
\sigma_1\mapsto (1,6,5)(2,4,3);\\
&\ \ \ \ \psi^{(5)}_{3,6}\colon 
\alpha\mapsto (1,2,3)(4,5,6),\ \
\beta\mapsto (1,2),\ \ \ \ \qquad\qquad
\sigma_1\mapsto (2,3)(4,6,5);\\
&\ \ \ \ \psi^{(6)}_{3,6}\colon 
\alpha\mapsto (1,2,3)(4,5,6),\ \
\beta\mapsto (1,4),\ \ \ \ \qquad\qquad
\sigma_1\mapsto (1,6,5,4,3,2);\\
&\ \ \ \ \psi^{(7)}_{3,6}\colon 
\alpha\mapsto (1,2,3),\qquad\qquad
\beta\mapsto (1,4)(2,5)(3,6),\ \ 
\sigma_1\mapsto (1,4,3,6,2,5).\\
&\\
&d) \ n=7:\\
&\ \ \ \ \ \psi^{(1)}_{3,7}\colon 
\alpha\mapsto (1,2,3)(4,5,6),\ \
\beta\mapsto  (1,4)(2,7),\ \ \ \qquad
\sigma_1\mapsto (1,6,5,4,3,2,7);\\
&\ \ \ \ \ \psi^{(2)}_{3,7}\colon 
\alpha\mapsto (1,2,3)(4,5,6),\ \
\beta\mapsto (1,2)(3,4)(5,7),\ \ \
\sigma_1\mapsto (2,3,6,5,7,4);\\
& \ \ \ \ \psi^{(3)}_{3,7}\colon 
\alpha\mapsto (1,2,3)(4,5,6),\ \
\beta\mapsto (1,4)(2,5)(3,7),\ \ \
\sigma_1\mapsto (1,6,5)(2,4,3,7).
\endaligned
$$
\end{Proposition}

\begin{Remark}\label{Rmk 4.1} 
One of the $7$ homomorphisms
$\psi^{(i)}_{3,6}\colon B_3\to{\mathbf S}(6)$
listed above, namely, $\psi^{(1)}_{3,6}$,
appears in a way, which deserves some comments.

Take any $z=(z_1,z_2,z_3)\in{\mathbb C}^3$; let
$\lambda_{1},\lambda_{2},\lambda_{3}$ be the roots of the
polynomial $p_3(t,z)=t^3 + z_1 t^2 + z_2 t + z_3$.
Let 
$P(t,w)=t^6 + w_1 t^5 + w_2 t^4 
+ w_3 t^3 + w_4 t^2 + w_5 t + w_6$ be the monic polynomial in $t$
of degree $6$ with the roots
$\mu_{1}^\pm,\mu_{2}^\pm,\mu_{3}^\pm$ 
defined by the quadratic equations
%&
\begin{equation}\label{4.1}
\aligned
\ &(\mu_{1}^\pm-\lambda_{1})^2 =
(\lambda_{1}-\lambda_{2})(\lambda_{1}-\lambda_{3}), \\
&(\mu_{2}^\pm-\lambda_{2})^2 =
(\lambda_{2}-\lambda_{3})(\lambda_{2}-\lambda_{1}), \\
&(\mu_{3}^\pm-\lambda_{3})^2 =
(\lambda_{3}-\lambda_{1})(\lambda_{3}-\lambda_{2}).
\endaligned
\end{equation}
%&

The set of the $6$ numbers $\mu$'s (taking into account possible
multiplicities) is invariant under any permutation of
the roots $\lambda_i$. Since the coefficients $w_i$ are the elementary 
symmetric polynomials in $\mu$'s, they are polynomials in $z_1,z_2,z_3$.
Thereby, we obtain a polynomial mapping 
$f\colon{\mathbb C}^3\ni z\mapsto w\in{\mathbb C}^6$. It is easy to compute
the coordinate functions of this mapping:
%&
\begin{equation}\label{4.2}
\aligned
f_{1}(z) &= 2z_{1};\\
f_{2}(z) &= 5z_{2};\\
f_{3}(z) &= 20z_{3};
\endaligned
\qquad
\aligned
f_{4}(z) &= 20z_{1}z_{3}-5z^{2}_{2};\\
f_{5}(z) &= 8z^{2}_{1}z_{3} - 2z_{1}z^{2}_{2} - 4z_{2}z_{3}; \\
f_{6}(z) &= 4z_{1}z_{2}z_{3} - z^{3}_{2} - 8z^{2}_{3}.
\endaligned
\end{equation}
%&
The formulae for $f_1,f_2,f_3$ show that $f$ is an embedding.
Computing the discriminants
$D_P(w)$ and $D_{p_3}(z)=d_3(z)$ of the polynomials $P(t,w)$ and $p_3(t,z)$,
respectively, we obtain the relation
$D_P(w) = -4^9\cdot [d_3(z)]^5$.
In particular, {\sl if the polynomial $p_3(t,z)$ has no multiple roots,
then the polynomial $P(t,w)=p_6(t,f(z))$ has no multiple roots}. Hence,
the restriction of $f$ to the domain 
${\mathbf G}_3 = \{z\in{\mathbb C}^3| \,d_3(z)\ne 0\}$
(see Sec. \ref{sec: Introduction}) defines the polynomial mapping 
$$
f\colon{\mathbf G}_3\ni z\mapsto w=f(z)\in{\mathbf G}_6
=\{w\in{\mathbb C}^6| \,d_6(w)\ne 0\}.
$$
Moreover, formulae (\ref{4.1}) show that {\sl for any
$z\in{\mathbf G}_3$ the polynomials $p_3(t,z)$ and 
$P(t,w)=p_6(t,f(z))$ have no common roots}.

On the other hand, it was proven in \cite{Lin96a} that
{\sl for any $k>3$, any natural $n$, and any holomorphic mapping
$F\colon{\mathbf G}_k\to{\mathbf G}_n$
there must be a point $z^\circ\in{\mathbf G}_k$ such that the
polynomials $p_k(t,z^\circ)$ and $p_n(t,F(z^\circ))$
have common roots}.
This means that the mapping $f\colon{\mathbf G}_3\to{\mathbf G}_6$
constructed above is very exceptional. 

Take $z^\circ\in{\mathbf G}_3$ and fix an isomorphism 
$\tau\colon B_3\stackrel{\cong}{\longrightarrow}\pi_1({\mathbf G}_3,z^\circ)$.
Any element $s\in B_3$ produces the permutation $\widehat s$ of the
$6$ roots of the polynomial $P(t,f(z))$ along the loop in ${\mathbf G}_3$
(based at $z^\circ$) representing the $3$-braid $s$. This gives rise to a
homomorphism $B_3\to{\mathbf S}(6)$; up to conjugation, this
homomorphism does not depend on $z_\circ$, $\tau$ 
and coincides with $\psi^{(1)}_{3,6}$. Since $\psi^{(1)}_{3,6}$ 
is non-cyclic, {\sl the mapping $f$ is unsplittable}.
This mapping is related to a holomorphic section of the universal
Teichm\"uller family $\mathbf V'(0,4)\to \mathbf T(0,4)$ 
over the Teichm\"uller space $\mathbf T(0,4)$ and to elliptic functions. 
In fact, this is the way how $f$ was found; however, now it is written
down explicitly, and one can ask whether it may be found in a shorter way 
(say in some paper of the last century!). 
Let me also mention that the points
$\mu$'s lie on the bisectors of the triangle $\Delta$
with the vertices $\lambda$'s (these bisectors are well defined, even if
$\Delta$ degenerates to a segment with a marked interior point),
and the distance between $\lambda_i$ and $\mu_i^\pm$ is the geometric
mean of the corresponding legs of the triangle $\Delta$ 
(the latter observation is due to E. Gorin).
\end{Remark}
\vskip0.3cm

\begin{Lemma}\label{Lm 4.2} 
$\psi(\sigma_1) = \psi(\sigma_3)$ for any
transitive homomorphism $\psi\colon B_4\to{\mathbf S}(5)$.
\end{Lemma}

\begin{proof} 
Suppose that $\psi(\sigma_1)\ne\psi(\sigma_3)$. 
Then $\psi$ is non-cyclic, $4$ divides
$\ord\widehat\alpha$, and $3$ divides $\ord\widehat\beta$. Hence,
$[\widehat\alpha] = [4]$ and $[\widehat\beta] = [3]$.
Regarding ${\mathbf S}(5)$ as ${\mathbf S}(\{0,1,2,3,4\})$, 
we may assume that $\widehat\alpha = (0,1,2,3)$.
Since $\psi$ is transitive, $4\in \supp \widehat\beta$; consequently,
$\widehat\beta = (p,q,4)$, where $p,q\in \{0,1,2,3\}$
and $p\ne q$. Put $A = \widehat\beta\widehat\alpha\widehat\beta$ and
$B =\widehat\alpha^2\widehat\beta\widehat\alpha^{-3}\widehat\beta
\widehat\alpha^2$.
It follows from (\ref{special presentation1}) (with $i = 2$)
that $A = B$; in particular, $A(4) = B(4)$ and
$A(p) = B(p)$. \ $A(4) = B(4)$ implies
that $q = |p+1|_4$ ($|N|_4$ denotes the
residue of $N\in\mathbb N$ modulo $4$, \ $0\le|N|_4\le 3$).
Combining this with $A(p) = B(p)$,
we obtain $|p+2|_4 = |p+3|_4$, which is impossible.
\end{proof}

\begin{Proposition}\label{Prp 4.3} 
Any non-cyclic transitive
homomorphism
$\psi\colon B_4\to{\mathbf S}(5)$
is conjugate to the homomorphism \ \
$\psi_{4,5}\colon \alpha \mapsto (1,4)(2,5), \ \
\beta \mapsto (3,5,4)$ \ \ $(\Img \psi_{4,5} = \mathbf A(5))$.
Moreover, $\psi(\sigma_1)=\psi(\sigma_3)$.
\end{Proposition}

\begin{proof} 
It follows from Lemma \ref{Lm 4.2} that the
homomorphism $\psi$ may be represented as a composition 
of the canonical epimorphism
$\pi\colon B_4\to B_3$ (Remark \ref{Rmk: epimorphisms B(4) to B(3) and S(3)})
with a non-cyclic transitive homomorphism
$B_3\to{\mathbf S}(5)$; Proposition \ref{Prp 4.1}$(b)$
completes the proof. 
\end{proof}

Our next goal is to describe all transitive homomorphisms 
$\psi\colon B_4\to{\mathbf S}(6)$ that satisfy
the condition $\psi(\sigma_1)\ne\psi(\sigma_3)$. 
We start with some examples of such homomorphisms.

Recall that ${\mathbf S}(6)$ is the only symmetric group 
admitting outer automorphisms. Any such automorphism is conjugate 
to the automorphism $\varkappa$ defined by
%&$$
\begin{equation}\label{4.3}
\varkappa\colon \widetilde\sigma_1\mapsto (1,2)(3,4)(5,6),\qquad
\widetilde\alpha\mapsto (1,2,3)(4,5),
\end{equation}
%&$$
where $\widetilde\sigma_1 = (1,2)$ and $\widetilde{\alpha} = (1,2,3,4,5,6)$
(for instance, this can be proven using Artin Theorem).

Define the following two embeddings
$\xi,\eta\colon {\mathbf S}(4)\hookrightarrow {\mathbf S}(6)$.
The embedding $\xi$ is just induced by the natural inclusion
${\boldsymbol\Delta}_4 = \{1,2,3,4\}\hookrightarrow \{1,2,3,4,5,6\} 
= {\boldsymbol\Delta}_6$. Further, ${\mathbf S}(4)$ 
may be regarded as the group of all isometries of the tetrahedron; 
thereby, ${\mathbf S}(4)$ acts naturally on the set 
$E\cong {\boldsymbol\Delta}_6$ consisting of
the $6$ edges of the tetrahedron, which defines the
embedding $\eta\colon {\mathbf S}(4)\hookrightarrow {\mathbf S}(6)$.
With the canonical generators 
$\widetilde\sigma_i = (i,i+1)\in{\mathbf S}(4)$, \ $1\le i\le 3$,
the embedding $\eta$ looks as follows:
$$
\eta(\widetilde\sigma_1) = (1,2)(3,4),\qquad
\eta(\widetilde\sigma_2) = (2,5)(4,6),\qquad
\eta(\widetilde\sigma_3) = (1,4)(2,3).
$$
Let $\mu_4\colon B_4\to{\mathbf S}(4)$ 
be the canonical epimorphism and 
$\nu_{4,2}\colon B_4\to{\mathbf S}(4)$
be the homomorphism described in Artin Theorem.
It is easy to check that each of the compositions
$$
\psi^{(1)}_{4,6} = \varkappa\circ\xi\circ\mu_4,\qquad
\psi^{(2)}_{4,6} = \varkappa\circ\xi\circ\nu_{4,2},\qquad
\psi^{(3)}_{4,6} = \eta\circ\mu_4,\qquad
\psi^{(4)}_{4,6} = \eta\circ\nu_{4,2}
$$
defines a non-cyclic transitive homomorphism
$B_4\to{\mathbf S}(6)$ 
such that $\psi^{(i)}_{4,6}(\sigma_1)\ne
\psi^{(i)}_{4,6}(\sigma_3)$ for each $i=1,2,3,4$. These
homomorphisms act on the canonical generators 
$\sigma_1,\sigma_2,\sigma_3$ of $B_4$ as follows:
%&
\begin{equation}\label{4.4}
\aligned
\ &\psi^{(1)}_{4,6}\colon \sigma_1\mapsto (1,2)(3,4)(5,6), \
\sigma_2\mapsto (1,5)(2,3)(4,6), \ 
\sigma_3\mapsto (1,3)(2,4)(5,6); \\
%\\
&\psi^{(2)}_{4,6}\colon \sigma_1\mapsto (1,2,4,3), \ \ \ \ \ \ \ \ \
\sigma_2\mapsto (1,5,4,6), \ \ \ \ \ \ \ \ 
\sigma_3\mapsto (3,4,2,1);\\
%\\
&\psi^{(3)}_{4,6}\colon \sigma_1\mapsto (1,2)(3,4), \ \  \qquad
\sigma_2\mapsto (2,5)(4,6), \ \qquad
\sigma_3\mapsto (1,4)(2,3);\\
%\\
&\psi ^{(4)}_{4,6}\colon \sigma_1\mapsto (4,3,2,1)(5,6), \ \
\sigma_2\mapsto (4,6,2,5)(1,3), \ \ 
\sigma_3\mapsto (1,2,3,4)(5,6).\\
\endaligned
\end{equation}
%&$$
We shall show that any transitive homomorphism 
$\psi\colon B_4\to{\mathbf S}(6)$
that satisfies the condition 
$\psi(\sigma_1)\ne\psi(\sigma_3)$ is conjugate 
to one of the homomorphisms 
$\psi^{(i)}_{4,6}$, $1\le i\le 4$.

\begin{Lemma}\label{Lm 4.4} 
Let $\psi\colon B_4\to{\mathbf S}(6)$
be a homomorphism that satisfies $\widehat\sigma_1\ne\widehat\sigma_3$.
Then $a)$ $[\widehat\sigma_1]\ne [2,3]$; \
$b)$ $[\widehat\sigma_1]\ne [5]$; \
$c)$ $[\widehat\sigma_1]\ne [6]$; \
$d)$ $[\widehat\sigma_1]\ne [3,3]$.
\end{Lemma}

\begin{proof} 
$a)$ Assume that $\widehat\sigma_1=C_2 C_3$, where
$\supp C_2\cap\supp C_3=\varnothing$, $[C_2]=[2]$, and $[C_3]=[3]$.
Since $[\widehat\sigma_3] = [2,3]$ and
$\widehat\sigma_1\widehat\sigma_3=\widehat\sigma_3\widehat\sigma_1$,
Lemma \ref{Lm: invariant sets and components of commuting permutations}
implies that
$$
\widehat\sigma_3 = C_2C_3^2 = (C_2C_3)^5 = \widehat\sigma_1^5.
$$ 
Since $\psi$ is non-cyclic, the latter relation shows that
$\widehat\sigma_2$ forms braid-like couples with $\widehat\sigma_1$ and
$\widehat\sigma_1^5$. 
By Lemma \ref{Lm: braid-like couples (a,b) and (aq,b)} 
(with $q=5$, \ $\nu =\GCD(q+1,4)=2$ 
and $\nu (q-1)=2\cdot 4=8$), we have $\widehat\sigma_1^8=1$,
which contradicts the property $C_3\preccurlyeq \widehat\sigma_1$.

$b)$ For $[\widehat\sigma_1]=[5]$
Lemma \ref{Lm: invariant sets and components of commuting permutations}
shows that
$\widehat\sigma_3=\widehat\sigma_1^q$, where $q=3$ or $q=4$. Hence,
$\widehat\sigma_2$ forms braid-like couples with $\widehat\sigma_1$ and
$\widehat\sigma_1^q$. By Lemma
\ref{Lm: braid-like couples (a,b) and (aq,b)}, 
either $\widehat\sigma_1^8=1$ or
$\widehat\sigma_1^3=1$ respectively; but this is impossible.

$c)$ Similarly, for $[\widehat\sigma_1]=[6]$ we obtain
$\widehat\sigma_3=\widehat\sigma_1^5$ and $\widehat\sigma_1^8=1$,
which is impossible.

$d)$ Assume that $\widehat\sigma_1 = BC$,
where $B,C$ are disjoint 3-cycles. Then Lemma \ref{Lm: p cycles in S(2p)}
implies that $\widehat\sigma_3 = BC^2$ (cases $(ii)$,
$(iii)$ described in this lemma cannot occur here, because of
$[\widehat\sigma_3] = [3,3]$).
Since $\psi(\sigma_1)\ne \psi(\sigma_3)$, Lemma
\ref{Lm: k|m if psi(alpha m)<->psi(beta),(k-1)|m if psi(beta m)<->psi(alpha)}$(a)$
shows that $4$ divides $\ord \widehat\alpha$. Therefore, either
$[\widehat\alpha] = [4]$ or $[\widehat\alpha] = [4,2]$. In any case,
$[\widehat\alpha^2] = [2,2]$. 
It follows from relation (\ref{conjugation by alpha power}) that
$\widehat\sigma_3 = \widehat\alpha^2\widehat\sigma_1\widehat\alpha^{-2}$.
So, either 
$$
B = \widehat\alpha^2 B\widehat\alpha^{-2} \ \
{\text {and}} \ \ C^2 = \widehat\alpha^2 C\widehat\alpha^{-2}, \ \ \
\text{or} \ \ \
B = \widehat\alpha^2 C\widehat\alpha^{-2} \ \
{\text {and}} \ \ C^2 = \widehat\alpha^2 B\widehat\alpha^{-2}.
$$
However, it easy to see that this contradicts the condition 
$[\widehat\alpha^2] = [2,2]$.
\end{proof}

\begin{Proposition}\label{Prp 4.5} 
Any transitive homomorphism
$\psi\colon B_4\to{\mathbf S}(6)$
that satisfies $\psi(\sigma_1)\ne\psi(\sigma_3)$
is conjugate to one of the homomorphisms 
$\psi^{(i)}_{4,6}$ \ $(1\le i\le 4)$ defined by (\ref{4.4}).
\end{Proposition}

\begin{proof} 
Lemma \ref{Lm: if [psi(sigma(1))]=[3] or |Fix(psi(sigma(1))|>n-4 then psi is intransitive}
and Lemma \ref{Lm 4.4} show that $\widehat\sigma_1$ has one of the
following cyclic types:
$a)$ $[\widehat\sigma_1] = [2,2,2]$; \
$b)$ $[\widehat\sigma_1] = [4]$; \
$c)$ $[\widehat\sigma_1] = [2,2]$; \
$d)$ $[\widehat\sigma_1] = [4,2]$.

$a)$ We may assume that 
$\widehat\sigma_1 = (1,2)(3,4)(5,6)$. This
permutation is odd; hence 
$\widehat\alpha = \widehat\sigma_1\widehat\sigma_2\widehat\sigma_3$
is also odd. Since $4$ divides $\ord \widehat\alpha$,
we have $[\widehat\alpha] = [4]$; so $\widehat\alpha$
has precisely two fixed points, which cannot be
in the support of a transposition that occurs above in
$\widehat\sigma_1$ (for $\widehat\sigma_1$ and $\widehat\alpha$
generate $\Img\psi$, whereas $\psi$ is transitive). 
Thus, up to a $\widehat\sigma_1$-admissible conjugation
(i. e., a conjugation that does not change
the above form of $\widehat\sigma_1$), we have
$\Fix\,\widehat\alpha= \{1,4\}$. Since
$\widehat\sigma_3 = \widehat\alpha^2\widehat\sigma_1\widehat\alpha^{-2}$
commutes with $\widehat\sigma_1$ and
$\widehat\sigma_3\ne\widehat\sigma_1$, we obtain
that $\widehat\alpha^2=(2,3)(5,6)$. It follows
that (up to a $\widehat\sigma_1$-admissible conjugation)
$\widehat\alpha = (2,5,3,6)$; so,
$\widehat\sigma_2 = \widehat\alpha\widehat\sigma_1\widehat\alpha^{-1}
= (1,5)(2,3)(4,6)$, \ 
$\widehat\sigma_3 = \widehat\alpha^2\widehat\sigma_1\widehat\alpha^{-2}
= (1,3)(2,4)(5,6)$,
and $\psi\sim\psi^{(1)}_{4,6}$.

$b)$ Since $\widehat\sigma_1,\widehat\sigma_3$ commute but
do not coincide, it follows from
Lemma \ref{Lm: invariant sets and components of commuting permutations} that 
$\widehat\sigma_3 = (\widehat\sigma_1)^3$. Therefore, up
to conjugation,
$\widehat\sigma_1 = (1,2,4,3)$ and $\widehat\sigma_3 = (3,4,2,1)$.
Since $\psi$ is transitive, we have
$\{5,6\}\subset \supp \widehat\sigma_2$. It
follows from Lemma \ref{Lm: braid-like couples of 4-cycles} that (up to a 
$\widehat\sigma_1,\widehat\sigma_3$-admissible conjugation)
$\widehat\sigma_2 = (1,5,4,6)$. Thus, $\psi\sim\psi^{(2)}_{4,6}$.

$c)$ We may assume that $\widehat\sigma_1= (1,2)(3,4)$.
It follows from Lemma \ref{Lm: braid-like couples of type [2,2]} that 
(up to a $\widehat\sigma_1$-admissible conjugation) either
$\widehat\sigma_2 = (1,2)(4,5)$ or $\widehat\sigma_2 = (2,5)(4,6)$.
Using
Lemma \ref{Lm: commuting permutations of type [2,2]},
Lemma \ref{Lm: braid-like couples of type [2,2]} 
and taking into account 
that $\widehat\sigma_3\ne\widehat\sigma_1$,
we obtain that in the first case
$\widehat\sigma_3 = (1,2)(5,6)$ (which contradicts
the transitivity of $\psi$), and in the second case
$\widehat\sigma_3 = (1,4)(2,3)$. Hence $\psi\sim\psi^{(3)}_{4,6}$.

$d)$ We may assume that $\widehat\sigma_1 = (4,3,2,1)(5,6)$.
Since $\widehat\sigma_1$, $\widehat\sigma_3$ commute but
do not coincide, Lemma \ref{Lm 4.4} 
implies that $\widehat\sigma_3 = (1,2,3,4)(5,6)$.
All $\widehat\sigma_i$ \ $(1\le i\le 3)$
are even; so $\widehat\alpha$ is even too;
since $4$ divides $\ord\widehat\alpha$,
we see that $[\widehat\alpha] = [4,2]$; hence
$[\widehat\alpha^2] = [2,2]$. Since $\psi$ is transitive,
the transposition $T\preccurlyeq\widehat\alpha$
cannot coincide with $(5,6)$. Therefore, it follows from 
the relation $\widehat\sigma_3
= \widehat\alpha^2\widehat\sigma_1\widehat\alpha^{-2}$
that $\widehat\alpha^2| {\{5,6\}} = (5,6)$ and
(up to a $\widehat\sigma_1,\widehat\sigma_3$-admissible conjugation)
$\widehat\alpha^2| {\{1,2,3,4\}} = (1,3)$.
So, $\widehat\alpha^2 = (1,3)(5,6)$.
Since $[\widehat\alpha] = [4,2]$, it follows that
$\widehat\alpha = (2,4)(5,1,6,3)$ (up to conjugation of the above type).
Thereby,
$\widehat\sigma_2 = \widehat\alpha\widehat\sigma_1\widehat\alpha^{-1}
= (2,5,4,6)(1,3)$ and $\psi\sim\psi^{(4)}_{4,6}$. 
\end{proof}

\begin{Proposition}\label{Prp 4.6} 
Any non-cyclic transitive homomorphism
$\psi\colon B_4\to{\mathbf S}(6)$
is either conjugate to one of the homomorphisms
$\psi^{(i)}_{4,6}$ \ $(1\le i\le 4)$ defined by {\rm {(\ref{3.4})}} 
or conjugate to one of the compositions 
$$
\psi^{(i)}_{3,6}\circ\pi\colon B_4
\stackrel{\pi}{\longrightarrow} B_3
\stackrel{\psi^{(i)}_{3,6}}{\longrightarrow} {\mathbf S}(6),
$$
where $\pi\colon B_4\to B_3$ is the canonical
epimorphism and $\psi^{(i)}_{3,6}$ \ $(1\le i\le 7)$ are the homomorphisms
exhibited in {\rm {Proposition \ref{Prp 4.1}}}$(c)$.
\end{Proposition}

\begin{proof} 
In what follows we use the generators $u$, $v$, $w$, and $c_1$ 
of the commutator subgroup $B'_4$ described in
(\ref{generators of commutator}).
Proposition \ref{Prp 4.5} covers the case when
$\psi(\sigma_1)\ne\psi(\sigma_3)$. If $\psi(\sigma_1) = \psi(\sigma_3)$,
then $hat\sigma_3\widehat\sigma_1^{-1}=1$, so that the element
$c_1 = \sigma_3\sigma_1^{-1}\in \Ker\, \psi$ 
Hence the element $w = uc_1u^{-1}$ (see \ref{1.14}) also is
in $\Ker\psi$. Since the kernel $\mathbf T$ of
$\pi$ is generated by $c_1$ and $w$ (see Gorin-Lin
Theorem $(c)$ and Remark \ref{Rmk: epimorphisms B(4) to B(3) and S(3)}),
it follows that $\Ker\pi = \mathbf T\subseteq \Ker\psi$.
Therefore, there exists a homomorphism 
$\psi_{3,6}\colon B_3\to{\mathbf S}(6)$ such
that $\psi = \psi_{3,6}\circ\pi$. Clearly, 
$\psi_{3,6}$ must be non-cyclic and transitive; Proposition \ref{Prp 4.1}$(c)$
completes the proof.
\end{proof}

\begin{Remark}\label{Rmk 4.2} 
The homomorphism $\psi_{4,6}^{(3)}$ is conjugate
(by $(1,3,2)(5,6)$) to the homomorphism 
$\widetilde\nu_6'$ that is defined as follows. Let $\nu_6'$
be the restriction of Artin's homomorphism 
$\nu_6\colon B_6\to{\mathbf S}(6)$ to
$B'_6$. 
The mapping of the generators
$\sigma_i\mapsto c_i=\sigma_{i+2}\sigma_1^{-1}\in B'_6$, \
$i=1,2,3$, extends to an embedding 
$\lambda_6'\colon B_4\hookrightarrow B'_6$
(Remark \ref{Rmk: embeddings lambda(k,m) of B(k-2) to B(k)}). 
The homomorphism $\widetilde\nu_6'$ is
the composition of $\lambda_6'$ with $\nu_6'$.
\end{Remark}
\vskip0.3cm
\begin{Remark}\label{Rmk 4.3} 
The trivial embedding
${\mathbf S}(5)\hookrightarrow {\mathbf S}(6)$ 
is of little moment. However, its composition with the outer 
automorphism $\varkappa$ of ${\mathbf S}(6)$ is more interesting.
This composition can be also described as follows.
It is well known that $\mathbf A(5)$ may be regarded as 
the group of all rotations of the icosahedron. 
In particular $\mathbf A(5)$ acts on the set 
$LD\cong{\boldsymbol\Delta}_6$ of all the 6 "long diagonals" 
of the icosahedron. 
This action of $\mathbf A(5)$ on ${\boldsymbol\Delta}_6$ 
extends to an action of ${\mathbf S}(5)$, which leads to an embedding 
$\phi_{5,6}\colon {\mathbf S}(5)\to{\mathbf S}(6)$.
In terms of the generators
$\widetilde\sigma_i = (i,i+1)\in{\mathbf S}(5)$, \ $1\le i\le 4$,
it looks as follows:
%&
\begin{equation}\label{4.5}
\phi_{5,6}\colon \left\{
\aligned
\widetilde\sigma_1&\mapsto (1,2)(3,4)(5,6),\qquad
\widetilde\sigma_2\mapsto (1,5)(2,3)(4,6),\\
\widetilde\sigma_3&\mapsto (1,3)(2,4)(5,6),\qquad
\widetilde\sigma_4\mapsto (1,2)(3,5)(4,6).
\endaligned
\right .
\end{equation}
%&
Of course, it is easy to check directly that these formulae
indeed define a group homomorphism, which
is certainly transitive and non-abelian,
and therefore faithful (since ${\mathbf S}(5)$ does not possesses
proper non-abelian quotient groups). Hence the composition
%&
\begin{equation}\label{4.6}
\psi_{5,6} = \phi_{5,6}\circ\mu_5
\colon B_5\stackrel{\mu_5}{\longrightarrow} {\mathbf S}(5)
\stackrel{\phi_{5,6}}{\longrightarrow} {\mathbf S}(6)
\end{equation}
%&
is a non-cyclic transitive homomorphism with
$\Ker\psi_{5,6} = PB_5$ and $\Img \psi_{5,6}\cong {\mathbf S}(5)$.
It is easily seen that $\psi_{5,6}$ coincides with the composition 
%&
\begin{equation}\label{4.6'}
\psi_{5,6}=\nu_6\circ j^5_6\colon B_5\stackrel{j^5_6}{\hookrightarrow} B_6
\stackrel{\nu_6}{\longrightarrow} {\mathbf S}(6), 
\end{equation}
%&
where $j^5_6\colon B_5\ni\sigma_i\mapsto \sigma_i\in B_6$,
\ $1\le i\le 4$, and $\nu_6$ is Artin's homomorphism.
\end{Remark}
\vskip0.3cm

\begin{Proposition}\label{Prp 4.9}  
Any non-cyclic transitive homomorphism
$\psi\colon B_5\to{\mathbf S}(6)$ is conjugate
to the homomorphism $\psi_{5,6}$ defined by (\ref{4.5}), (\ref{4.6}) 
(or by (\ref{4.6'}), which is the same). 
In particular, $\Ker\psi = PB_5$ 
and $\Img\psi\cong{\mathbf S}(5)$.
\end{Proposition}

\begin{proof} 
Since $3=6/2$ is prime, 
Lemma \ref{Lm: Artin Lemma on cyclic decomposition of hatsigma(1)}$(c)$
implies that the cyclic decomposition
of $\widehat\sigma_1$ cannot contain a cycle
of length $\ge 3$. Lemma
\ref{Lm: if [psi(sigma(1))]=[3] or |Fix(psi(sigma(1))|>n-4 then psi is intransitive}$(a)$
eliminates the case $[\widehat\sigma_1]=[2]$. The non-commuting permutations
$\widehat\sigma_3$, $\widehat\sigma_4$ commute
with $\widehat\sigma_1$. On the other hand, any two permutations
of cyclic type $[2,2]$ supported on the same $4$ points commute;
therefore, Lemma \ref{Lm: p cycles in S(2p)} eliminates the case
$[\widehat\sigma_1]=[2,2]$. Thus, the only possible case is
$[\widehat\sigma_1]=[2,2,2]$. Two distinct permutations of this cyclic
type in ${\mathbf S}(6)$ commute if and only if they have precisely one
common transposition. Therefore, without loss of generality,
we may assume that
%&
\begin{equation}\label{4.7}
\widehat\sigma_1=(1,2)(3,4)(5,6),\qquad
\widehat\sigma_3=(1,3)(2,4)(5,6).
\end{equation}
%&
By the same reason, $\widehat\sigma_4$ has one common transposition
with $\widehat\sigma_1$ but not with $\widehat\sigma_3$; this common
transposition may be either $(1,2)$ or $(3,4)$.
The renumbering of the symbols
$1\rightleftarrows 3$, \ $2\rightleftarrows 4$
takes each of these cases into another one and does not change
the forms (\ref{4.7}); so we may assume
that this common transposition is $(1,2)$. Since $\widehat\sigma_4$
has no common transpositions with $\widehat\sigma_3$, we have
either $\widehat\sigma_4=(1,2)(3,5)(4,6)$
or $\widehat\sigma_4=(1,2)(3,6)(4,5)$. The second case can be
obtain from the first one by 
$5\rightleftarrows 6$, which does not change
the forms (\ref{4.7}); hence we may assume that
%&
\begin{equation}\label{4.8}
\widehat\sigma_4=(1,2)(3,5)(4,6).
\end{equation}
%&
An argument of the same kind shows that
$\widehat\sigma_2$ must contain a single common transposition with
$\widehat\sigma_4$, but not with $\widehat\sigma_1$ and $\widehat\sigma_3$.
Each of the transpositions $(1,2)$, $(3,4)$, $(5,6)$,
$(1,3)$, $(2,4)$ is contained in 
$\widehat\sigma_1$ or in $\widehat\sigma_3$;
hence either $\widehat\sigma_2=(1,5)(2,3)(4,6)$
or $\widehat\sigma_2=(1,4)(2,6)(3,5)$. However, the second case
may be obtained from the first one by
$1\rightleftarrows 2$, \ $3\rightleftarrows 4$, \ $5\rightleftarrows 6$,
which does not change the forms (\ref{4.7}),(\ref{4.8}).
This shows that $\psi\sim \psi_{5,6}$.
\end{proof}
\vfill
\newpage

%Sec. 5
\section{Retractions of homomorphisms
and cohomology}
\label{sec: Retraction of homomorphisms; homomorphisms and cohomology}

Our aim is to study homomorphisms
$\psi\colon B_k\to{\mathbf S}(n)$ up to conjugation.
In this section we develop an approach to this problem. In general
terms, this approach may be described as follows.
\vskip0.2cm

Let $2\le r\le n$ and let $\mathfrak C_r$ be the $r$-component
of $\widehat\sigma_1=\psi(\sigma_1)$, i. e., the set of all $r$-cycles
that occur in the cyclic decomposition of $\widehat\sigma_1$ 
(see Definition \ref{Def: symmetric groups}{\bf (b)}).
Since $\widehat\sigma_i=\psi(\sigma_i)$ with $3\le i\le k-1$ commute
with $\widehat\sigma_1$, all $r$-cycles
$C_i'=\widehat\sigma_i C \widehat\sigma_i^{-1}$ ($3\le i\le k-1$) 
also belong to $\mathfrak C_r$.
Thereby we obtain an action $\Omega_\psi$
of $B_{k-2}$ on $\mathfrak C_r$, that is,
a homomorphism $\Omega_\psi\colon B_{k-2}\to{\mathbf S}(\mathfrak C_r)$.
Since $\#\mathfrak C_r\le n/r\le n/2$, the homomorphism $\Omega_\psi$ 
is ``simpler" than $\psi$ and we may hope
to ``recognize" first $\Omega_\psi$ and then handle $\psi$ itself.

\subsection{Components and corresponding exact sequences}
\label{Ss: 5.1. Components and corresponding exact sequences}
We denote by $\{\psi\}$ the conjugacy class of a homomorphism 
$\psi\in\Hom(B_k, {\mathbf S}(n))$; that is, 
$\psi'\in \{\psi\}$ if and only if $\psi'\sim\psi$.
Recall that for any $r\ge 2$ the $r$-component
$\mathfrak C = \mathfrak C_r(A)$ of a permutation $A\in{\mathbf S}(n)$ 
is the set of all the $r$-cycles that occur in the cyclic
decomposition of $A$ (see Definition \ref{Def: symmetric groups}). 
For natural numbers $r,t$ \ ($2\le r\le n$, \ $t\le n/r$), 
we denote by $\Hom_{r,t}(B_k, {\mathbf S}(n))$ the subset of
$\Hom(B_k, {\mathbf S}(n))$ consisting of all homomorphisms $\psi$ 
that satisfy the following condition:
\vskip0.3cm

(!) {\sl the permutation 
$\widehat\sigma_1 = \psi(\sigma_1)\in{\mathbf S}(n)$ has an $r$-component
$\mathfrak C$ of length $t$}.
\vskip0.3cm

Let $\psi\in\Hom_{r,t}(B_k, {\mathbf S}(n))$ and let
$\mathfrak C = \{C_1,...,C_t\}$ be the $r$-component of the
permutation $\widehat\sigma_1 = \psi(\sigma_1)$ (so,
$C_1,...,C_t$ are disjoint $r$-cycles). The union 
$$
\Sigma = \supp{\mathfrak C} = \bigcup_{m=1}^{t} \supp C_m
\subseteq \supp \widehat\sigma_1\subseteq {\boldsymbol\Delta}_n
= \{1,\ldots,n\}
$$
of the supports $\Sigma_m = \supp C_m$ of all the
cycles $C_m$ is called the {\it support of the $r$-component}
$\mathfrak C$; we denote this set $\Sigma$ also by $\supp \mathfrak C$.
\vskip0.2cm

A homomorphism $\psi\in\Hom_{r,t}(B_k, {\mathbf S}(n))$
is said to be {\it normalized} if
$$
\ \ \ (!!) \qquad
\aligned
\ &\Sigma = \supp \mathfrak C = \{1,2,\ldots,tr\}, \ \ \
\mathfrak C = \{C_1,...,C_t\},  \\
&C_m = ((m-1)r+1, (m-1)r+2,\ldots , mr),\ \ \ m = 1,\ldots ,t,
\ \ {\text {and}}\\
&\psi(\sigma_1)| \Sigma = C_1\cdots C_t.
\endaligned 
\qquad 
$$
The following two statements are evident:
\vskip0.3cm

\noindent{\bfit Claim 1.}
{\sl If $\psi\in\Hom_{r,t}(B_k, {\mathbf S}(n))$, then
$\{\psi\}\subseteq\Hom_{r,t}(B_k, {\mathbf S}(n))$ and the class
$\{\psi\}$ contains at least one normalized homomorphism. 
Two normalized homomorphisms 
$\psi,\widetilde\psi\in\Hom_{r,t}(B_k, {\mathbf S}(n))$
are conjugate if and only if there is a permutation 
$\widetilde s\in{\mathbf S}(n)$
such that the set $\Sigma = \supp \mathfrak C = \{1,2,\ldots,tr\}$
is $\widetilde s$-invariant and}
%&
\begin{equation}\label{5.1}
\aligned
\ &\widetilde\psi(b) = \widetilde s \psi(b) \widetilde s^{-1}
\ \ \ {\text {for all}} \ b\in B_k
\qquad {\text {and}} \\
&\widetilde s\cdot C_1\cdots C_t \cdot \widetilde s^{-1}
= \widetilde s\cdot \psi(\sigma_1)| \Sigma \cdot \widetilde s^{-1}
= \widetilde\psi(\sigma_1)| \Sigma
= C_1\cdots C_t. \hskip20pt \square
\endaligned
\end{equation}
%&
\vskip0.2cm

\noindent{\bfit Claim 2.}
{\sl To classify the homomorphisms of the class
$\Hom_{r,t}(B_k, {\mathbf S}(n))$ up to conjugation,
it is sufficient to classify the normalized homomorphisms
of this class up to conjugation by permutations $\widetilde s$
that satisfy $(\ref{5.1})$. 
Such permutations $\widetilde s$ form a subgroup
$\widetilde G\subset{\mathbf S}(n)$ that
coincides with the centralizer $\mathbf C(\mathcal C,{\mathbf S}(n))$
of the element $\mathcal C = C_1\cdots C_t$ in ${\mathbf S}(n)$.
This subgroup $\widetilde G$ is naturally isomorphic to
the direct product $G\times {\mathbf S}(\Sigma')$,
where $G = \mathbf C(\mathcal C,{\mathbf S}(\Sigma))$ is the
centralizer of the element $\mathcal C$ in the symmetric
group ${\mathbf S}(\Sigma)\cong {\mathbf S}(rt)$ and ${\mathbf S}(\Sigma')$
is the symmetric group of the complement
$\Sigma'={\boldsymbol\Delta}_n\setminus\Sigma$.}
\hfill $\square$
\vskip0.3cm

\noindent We denote by $H\cong ({\mathbb Z}/r{\mathbb Z})^t$
the abelian subgroup of the symmetric group ${\mathbf S}(\Sigma)$
generated by the $r$-cycles $C_1,...,C_t$ defined in (!!); in particular,
$H$ contains the product $\mathcal C = C_1\cdots C_t$ and therefore
$H\subset G$. Clearly, $H$ {\sl is an abelian normal subgroup in $G$,
and the quotient group $G/H$ is isomorphic to the symmetric group
${\mathbf S}(\mathfrak C)\cong{\mathbf S}(t)$ permuting
the cycles $C_1,...,C_t$.} Thereby, we obtain the exact sequence
%&
\begin{equation}\label{5.2}
1\to H\to G\stackrel{\pi}{\longrightarrow}{\mathbf S}(t)\to 1,
\end{equation}
%&$$
where $\pi$ is the natural projection onto the quotient group;
in fact, this projection $\pi$ may be described explicitly as follows.
Since any element $g\in G$ commutes with the element
$\mathcal C = C_1\cdots C_t$, we have
$$
C_1\cdots C_t = g\cdot C_1\cdots C_t \cdot g^{-1}
= gC_1g^{-1} \cdots gC_tg^{-1}\qquad\text{for any} \ \ g\in G,
$$
where $C_1,\ldots ,C_t$ are disjoint $r$-cycles
as well as $gC_1g^{-1},\ldots , gC_tg^{-1}$.
Hence there is a unique permutation $s=s_g\in{\mathbf S}(t)$ such that
$g C_m g^{-1} = C_{s(m)}$ for all $m$; clearly, $\pi(g) = s_g$.

\begin{Lemma}\label{Lm 5.1} 
The exact sequence $(\ref{5.2})$ splits, 
that is, there exists a homomorphism
$\rho\colon {\mathbf S}(t)\to G$ such that
$\pi\circ\rho =\id_{{\mathbf S}(t)}$.
In fact, the group $G$ is
the semi-direct product of the groups $H$ and ${\mathbf S}(t)$ 
defined by the natural action of the symmetric group 
${\mathbf S}(t)$ on the direct product $({\mathbb Z}/r{\mathbb Z})^t$.
\end{Lemma}

\begin{proof} 
For any element $s\in{\mathbf S}(t)$,
define the permutation $\rho (s)\in{\mathbf S}(\Sigma)$ by
%&$$
\begin{equation}\label{5.3}
\rho (s)((m-1)r+q) = (s(m)-1)r+q\qquad
(1\le m\le t, \ \ 1\le q\le r);
\end{equation}
%&$$
thereby, we obtain the homomorphism
$$
\rho\colon {\mathbf S}(t)\ni s\mapsto \rho (s)\in{\mathbf S}(\Sigma).
$$
Using the explicit forms of the cycles $C_1,...,C_t$
(see (!!)), it is easy to check that for any $m = 1,\ldots ,t$ 
and any $s\in{\mathbf S}(t)$
%&$$
\begin{equation}\label{5.4}
\rho (s) C_m \rho (s)^{-1} = C_{s(m)}, 
\end{equation}
%&$$
which implies $\rho (s)\mathcal C \rho (s)^{-1} = \mathcal C$. So,
$\rho (s)\in G$ and hence $\rho$ may be considered as a
homomorphism from ${\mathbf S}(t)$ to $G$. 
It follows from (\ref{5.4}) and the above description of $\pi$ that 
$\pi\circ\rho = \id_{{\mathbf S}(t)}$. 
\end{proof}

From now on, {\sl we fix the splitting
$\rho\colon {\mathbf S}(t)\to G\subset {\mathbf S}(\Sigma)$
defined by} \ref{5.3}.

\subsection{Retractions of homomorphisms to components}
\label{Ss: 5.2. Retractions of homomorphisms to components}
In what follows, we fix a normalized homomorphism
$\psi\in \Hom_{r,t}(B_k,{\mathbf S}(n))$.
and put $\widehat\sigma_i = \psi(\sigma_i)$,\ $1\le i\le k-1$.
We work with the $r$-component $\mathfrak C = \{C_1,...,C_t\}$
of the permutation $\widehat\sigma_1$ keeping in mind the particular
forms of the $r$-cycles $C_1,...,C_t$ exhibited in condition (!!).
\vskip0.2cm

\noindent Relations (\ref{conjugation by alpha power}) 
show that
$\widehat\sigma_j 
= \widehat\alpha^{j-1}\widehat\sigma_1\widehat\alpha^{-(j-1)}$
for any $j = 1,\ldots,k-1$ (as usual,
$\widehat\alpha = \widehat\sigma_1\cdots \widehat\sigma_{k-1}$).
For any natural numbers $q,m$ such that 
$1\le q\le k-1$ and $1\le m \le t$, we put
%&
\begin{equation}\label{5.5}
C^{(q)}_m = \widehat\alpha^{q-1} C_m \alpha^{-(q-1)},
\qquad {\mathfrak C}^{(q)} = \{C^{(q)}_1,\ldots,C^{(q)}_t\}.
\end{equation}
%&
Clearly, $C^{(1)}_m = C_m$ for each $m=1,\ldots t$, 
and ${\mathfrak C}^{(1)} = \mathfrak C$; moreover, the set 
${\mathfrak C}^{(q)} = \{C^{(q)}_1,\ldots,C^{(q)}_t\}$
coincides with the $r$-component $\mathfrak C_r(\widehat\sigma_q)$ of
the permutation $\widehat\sigma_q$, and formulae (\ref{5.5}) provide
the {\sl marked identifications} 
%&
\begin{equation}\label{5.6}
\mathfrak I_q\colon \mathfrak C_r(\widehat\sigma_q) =
{\mathfrak C}^{(q)}\cong {\boldsymbol\Delta}_t \ \
{\text {and}} \ \ 
\mathfrak J_q\colon {\mathbf S}(\mathfrak C_r(\widehat\sigma_q)) =
{\mathbf S}({\mathfrak C}^{(q)})\cong {\mathbf S}(t).
\end{equation}
%&
In what follows, {\sl we always have in mind these
identifications}.

Put
$$
\Sigma^{(q)}_{m} = \supp C^{(q)}_{m} =
\widehat\alpha^{q-1}(\supp C_m);
$$
clearly,
$$
\supp \mathfrak C_r(\widehat\sigma_q) = \supp {\mathfrak C}^{(q)} = \Sigma^{(q)} =
\Sigma ({\mathfrak C}^{(q)}) =
\widehat\alpha^{q-1}(\supp {\mathfrak C}) =
\bigcup_{m=1}^{t} \Sigma^{(q)}_m \subseteq{\boldsymbol\Delta}_t.
$$
The set $\Sigma^{(q)}$ is invariant under all the
permutations $\widehat\sigma_j$, \ $j\ne q-1,q+1$ (for
each of them commutes with $\widehat\sigma_q$).
By the same reason, for any $j\ne q-1,q+1$ and 
any $m = 1,\ldots, t$, the $r$-cycle
$$
\widetilde C^{(q)}_m = \widehat\sigma_j
\cdot C^{(q)}_m\cdot \widehat\sigma_j^{-1}
$$
occurs in the cyclic decomposition of $\widehat\sigma_q$,
and therefore $\widetilde C^{(q)}_m$ coincides with one of the cycles 
$C^{(q)}_1,\ldots,C^{(q)}_t$. 
Thereby, for each $j\in \{1,\ldots,k-1\}$, $j\ne q-1,q+1$, 
the correspondence
$$
C^{(q)}_{m}\mapsto\widetilde {C}^{(q)}_{m} = \widehat\sigma_j
\cdot C^{(q)}_m\cdot \widehat\sigma_j^{-1},
\qquad m=1,\ldots,t,
$$
gives rise to the permutation 
$\mathfrak g^{(q)}_j\in{\mathbf S}(\mathfrak C^{(q)})\cong 
{\mathbf S}(t)$, and we obtain the correspondence
%&
\begin{equation}\label{5.7}
\widehat\sigma_j\mapsto {\mathfrak g}^{(q)}_j
\in{\mathbf S}({\mathfrak C}^{(q)})
\stackrel{\mathfrak J_q}{\cong}
{\mathbf S}(t),\qquad j\in \{1,\ldots,k-1\},\qquad j\ne q-1,q+1,
\end{equation}
%&
such that
%&
\begin{equation}\label{5.8}
\widehat\sigma_j\cdot C^{(q)}_m\cdot\widehat\sigma_j^{-1}
= {\mathfrak g}^q_j (C^{(q)}_m) =
C^{(q)}_s,\qquad s={\mathfrak g}^q_j(m).
\end{equation}
%&
It is convenient to introduce some special notations for 
some of the above objects corresponding to the values
$q=1$ and $q=k-1$. Namely, we put 
%&
\begin{equation}\label{5.9}
\aligned
\ &g_i = {\mathfrak g}^{(1)}_{i+2}\qquad\qquad {\text {and}}
\qquad g^*_i = {\mathfrak g}^{(k-1)}_i\qquad \ \ 
{\text {for}} \ \ \ i = 1,\ldots,k-3,\\
%\\
&C^*_m = C^{(k-1)}_m\qquad \ {\text {and}}
\qquad \Sigma^{(k-1)}_m = \Sigma^*_m\qquad {\text {for}}
\ \ \ m = 1,\ldots,t,
\endaligned
\end{equation}
%&
and
%&
\begin{equation}\label{5.10}
{\mathfrak C}^* = {\mathfrak C}^{(k-1)} = 
\{C^*_1,...,C^*_t\}, \ \ \ \mathcal C^* = C^*_1\cdots C^*_t, \ \ \
\Sigma^{*} = \Sigma^{(k-1)} = \supp {\mathfrak C}^{*}.
\end{equation}
%&$$
We should also keep in mind that
$$
\aligned
\ &C_m = C^{(1)}_m, \qquad \mathcal C = C_1\cdots C_t,\qquad
\ \ \ \mathfrak C = \{C_1,...,C_t\} = {\mathfrak C}^{(1)}, \\
&\Sigma_m = \Sigma^{(1)}_m = \supp C_m,\qquad
\Sigma = \Sigma^{(1)} = \supp \mathfrak C.
\endaligned
$$
\vskip0.2cm

\noindent The construction of 
Sec. \ref{Ss: 5.1. Components and corresponding exact sequences} 
applies also to the 
$r$-cycles $C^*_1,\ldots , C^*_t$. Namely, we denote by $G^*$
the centralizer of the element $\mathcal C^* = C^*_1\cdots C^*_t$ 
in ${\mathbf S}(\Sigma^*)$, and denote by
$H^*\cong ({\mathbb Z}/r{\mathbb Z})^t$ 
the abelian normal subgroup in $G^*$ generated by all
the $r$-cycles $C^*_1,...,C^*_t$. Then
$G^*/H^*\cong {\mathbf S}(\mathfrak C^*)\cong{\mathbf S}(t)$, 
and we obtain the exact sequence
%&
\begin{equation}\label{5.2*}
1\to H^*\to G^*\stackrel{\pi^*}{\longrightarrow} {\mathbf S}(t)\to 1.
\end{equation}
%&
The projection $\pi^*$ may be described as follows.
Any element $g\in G^*$ commutes with the product 
$\mathcal C^* = C^*_1\cdots C^*_t$, and thus
$$
C^*_1\cdots C^*_t = g\cdot C^*_1\cdots C^*_t \cdot g^{-1}
= gC^*_1g^{-1} \cdots gC^*_tg^{-1}.
$$
Since $C^*_1,\ldots ,C^*_t$ and also $gC^*_1g^{-1},\ldots , gC^*_tg^{-1}$ 
are disjoint $r$-cycles, 
there is a unique permutation $s^*=s^*_g\in{\mathbf S}(t)$ such that
$g C^*_m g^{-1} = C^*_{s^*(m)}$ for all $m$; we put $\pi^*(g) = s^*_g$.
The following statement follows immediately from our definitions:
\vskip0.2cm

\noindent{\bfit Claim 3.} {\sl The conjugation by the element
$\widehat\alpha^{k-2} = \psi(\alpha^{k-2})$,
$$
c_\psi\colon {\mathbf S}(n)\ni A\mapsto \widehat\alpha^{k-2}\cdot A
\cdot \widehat\alpha^{-(k-2)}\in{\mathbf S}(n),
$$
provides the commutative diagram

$$
\CD
1@> >> H @> >> G @> \pi >> {\mathbf S}(t) @> >> 1\\
@.      @V\cong Vc_\psi V   @V\cong Vc_\psi V  @VV\id V\\
1@> >> H^* @> >> G^* @> \pi^* >> {\mathbf S}(t) @> >> 1
\endCD
\eqnum{CD[r,t;\psi]}
$$
\vskip0.2cm

\noindent The upper line of this diagram {\rm (the exact sequence (\ref{5.2}))}
is universal for all normalized homomorphisms
$\psi\in \Hom_{r,t}(B_k,{\mathbf S}(n))$, while the lower line
{\rm (the exact sequence (\ref{5.2*}))} and the vertical isomorphisms
$c_{\psi}$ may depend on $\psi$}.
\hfill $\square$
\vskip0.3cm

\noindent As we know, $\alpha^k$ is a central element in $B_k$
(see (\ref{alpha pow k=beta pow(k-1)}) 
or Sec. \ref{Ss: Center}); this implies some useful relations
between the permutations $\mathfrak g^{(q)}_j$ (with various $q,j$)
defined by (\ref{5.7}),(\ref{5.8}).

\begin{Lemma}\label{Lm 5.2} 
$a)$ If \ $1\le q\le k-3$ then
%&
\begin{equation}\label{5.11}
{\mathfrak g}^{(q)}_j = {\mathfrak g}^{(k-1)}_{j-q-1} = g^*_{j-q-1}
\qquad {\text {for}} \ \ q+2\le j\le k-1.
\end{equation}
%&

$b)$ If \ $3\le q\le k-1$ then
%&
\begin{equation}\label{5.11*}
{\mathfrak g}^{(q)}_j = {\mathfrak g}^{(1)}_{j+k-q+1} = g_{j+k-q-1}
\qquad {\text {for}} \ \ 1\le j\le q-2.
\end{equation}
%&
In particular
%&
\begin{equation}\label{5.12}
g^*_j = {\mathfrak g}^{(k-1)}_j 
= {\mathfrak g}^{(1)}_{j+2} = g_j\qquad {\text {for}} \ \ 
1\le j\le k-3.
\end{equation}
%&
\end{Lemma}

\begin{proof} 
$a)$ Take any $m\in \{1,...,t\}$ and any $q,j$
such that $1\le q\le k-3$, \  $q+2\le j\le k-1$, and put
$s={\mathfrak g}^{(q)}_j(m)$. It follows from (\ref{5.5}) and 
(\ref{5.8}) that 
%&
\begin{equation}\label{5.13}
\aligned
C^{(q)}_s = {\mathfrak g}^q_j (C^{(q)}_m) &=
\widehat\sigma_j\cdot C^{(q)}_m\cdot\widehat\sigma_j^{-1}\\
&= \widehat\sigma_j\cdot \widehat\alpha^{q-1} C_m \widehat\alpha^{-(q-1)}
\cdot\widehat\sigma_j^{-1} =
\widehat\sigma_j \widehat\alpha^{q-1}\cdot C_m \cdot
(\widehat\sigma_j \widehat\alpha^{q-1})^{-1}.
\endaligned
\end{equation}
%&
Since $\widehat\alpha^k$ commutes with any element in $\Img \psi$, we have
$$
\widehat\sigma_j = \widehat\alpha^{q+1}\widehat\sigma_{j-q-1}\widehat\alpha^{-(q+1)}
= \widehat\alpha^{-(k-q-1)}\widehat\sigma_{j-q-1}\widehat\alpha^{k-q-1},
$$
and thus
%&
\begin{equation}\label{5.14}
\widehat\sigma_j \widehat\alpha^{q-1} =
\widehat\alpha^{-(k-q-1)}\widehat\sigma_{j-q-1}\cdot \widehat\alpha^{k-2}.
\end{equation}
%&
Relations (\ref{5.13}),(\ref{5.14}) and (\ref{5.5}) 
(the latter one with $q=k-1$)
show that
%&
\begin{equation}\label{5.15}
\aligned
C^{(q)}_s &= 
\widehat\alpha^{-(k-q-1)}\widehat\sigma_{j-q-1}\widehat\alpha^{k-2}\cdot C_m
\cdot
\left(\widehat\alpha^{-(k-q-1)}\widehat\sigma_{j-q-1}
                           \widehat\alpha^{k-2}\right)^{-1}\\
&= \widehat\alpha^{-(k-q-1)}\widehat\sigma_{j-q-1}\cdot
\widehat\alpha^{k-2}\cdot C_m \cdot \widehat\alpha^{-(k-2)}\cdot
\left(\widehat\alpha^{-(k-q-1)}\widehat\sigma_{j-q-1}\right)^{-1}\\
&= \widehat\alpha^{-(k-q-1)}\cdot \widehat\sigma_{j-q-1}\cdot C^{(k-1)}_m
\cdot \widehat\sigma_{j-q-1}^{-1}\cdot \widehat\alpha^{k-q-1}
\endaligned
\end{equation}
%&
According to (\ref{5.8}) (with $q=k-1$ and $j-q-1$ instead of $j$), we have
$$
\widehat\sigma_{j-q-1}\cdot C^{(k-1)}_m \cdot \widehat\sigma_{j-q-1}^{-1}
=\mathfrak g^{(k-1)}_{j-q-1} \left(C^{(k-1)}_m\right)
= C^{(k-1)}_{s'}, \ \ \
{\text {where}} \ \ s' = \mathfrak g^{(k-1)}_{j-q-1} (m);
$$
thus, (\ref{5.15}) can be written as
%&
\begin{equation}\label{5.16}
\aligned
C^{(q)}_s &=
\widehat\alpha^{-(k-q-1)}\cdot C^{(k-1)}_{s'} \cdot \widehat\alpha^{k-q-1}\\
&= \widehat\alpha^{-(k-q-1)}\cdot\widehat\alpha^{k-2} C_{s'}
\widehat\alpha^{-(k-2)} \cdot \widehat\alpha^{k-q-1}
=\widehat\alpha^{q-1}\cdot C_{s'}
\cdot \widehat\alpha^{-(q-1)}
= C^{(q)}_{s'},
\endaligned
\end{equation}
%&
where
$$
s={\mathfrak g}^{(q)}_j(m),\qquad
s' = \mathfrak g^{(k-1)}_{j-q-1}(m).
$$
It follows from (\ref{5.16}) that $s=s'$, and thus
${\mathfrak g}^{(q)}_j(m) = \mathfrak g^{(k-1)}_{j-q-1}(m)$.
The latter relation holds for any $m\in \{1,...,t\}$,
which means that ${\mathfrak g}^{(q)}_j = \mathfrak g^{(k-1)}_{j-q-1}
= g^*_{j-q-1}$.

$b)$ For any $q=3,\ldots,k-1$, any $j=1,\ldots,q-2$ and every
$m=1,\ldots,t$ we have:
$$
\aligned
C^{(q)}_{s} &= \widehat\sigma_j\cdot C^{(q)}_m
\cdot \widehat\sigma_j^{-1}
%\\
= \widehat\sigma_j \widehat\alpha^{q-1}\cdot C_m
\cdot \widehat\alpha^{-(q-1)}\widehat\sigma_j^{-1}\\
%\\
&= \widehat\alpha^{q-1}\widehat\alpha^{k-q+1}\widehat\sigma_j\widehat\alpha^{-(k-q+1)}
\cdot C_m\cdot
\widehat\alpha^{k-q+1}\widehat\sigma_j^{-1}\widehat\alpha^{-(k-q+1)}
\widehat\alpha^{-(q-1)}\\
%\\
&= \widehat\alpha^{q-1}\widehat\sigma_{j+k-q+1}
\cdot C_m\cdot \widehat\sigma_{j+k-q+1}^{-1}
\widehat\alpha^{-(q-1)}
%\\
= \widehat\alpha^{q-1}\cdot C_{s'}\cdot \widehat\alpha^{-(q-1)}
= C^{(q)}_{s'},
\endaligned
$$
where $s = {\mathfrak g}^{(q)}_j(m)$ and
$s' = {\mathfrak g}^{(1)}_{j+k-q+1}(m)$. Consequently,
${\mathfrak g}^{(q)}_{j} = {\mathfrak g}^{(1)}_{j+k-q+1}
= g_{j+k-q+1}$.
Using the latter relations for $q=k-1$, we obtain (\ref{5.12}).
\end{proof}

\subsection{Homomorphism $\Omega$ relatad to an $r$-component
$\mathfrak C=\{C_1,...,C_t\}$}
\label{Ss: 5.3. Homomorphism Omega relatad to a component}
As above, we deal with a fixed homomorphism 
$\psi\colon B_k\to{\mathbf S}(n)$ and a certain
$r$-component $\mathfrak C=\{C_1,...,C_t\}$ of the permutation
$\widehat \sigma_1=\psi(\sigma_1)$; we keep the notation of the
previous subsection. Moreover, we assume that $k>3$ and denote by
$s_1,...,s_{k-3}$ the canonical generators of
the braid group $B_{k-2}$.
\vskip0.2cm

\noindent Consider the homomorphisms 
%&
\begin{eqnarray}
%\aligned
&&\phi\colon \ B_{k-2}\to{\mathbf S}(n),\quad
\phi(s_i) \ \ = \psi(\sigma_{i+2}) = \widehat\sigma_{i+2}\,,
\label{5.17}\\
&&\phi^*\colon B_{k-2}\to{\mathbf S}(n),\quad
\phi^*(s_i) \,= \psi(\sigma_i) \ \ \ = \widehat\sigma_i\,,
\label{5.17*}
%\endaligned
\end{eqnarray}
%&
where $i = 1,...,k-3$.
According to (\ref{5.7}), (\ref{5.8}), we have:
%&
\begin{equation}\label{5.18}
\aligned
\phi(s_i)\cdot C_m\cdot \phi(s_i)^{-1} &=
\widehat\sigma_{i+2}\cdot C_m\cdot \widehat\sigma_{i+2}^{-1} =
\mathfrak g^{(1)}_{i+2}\big(C^{(1)}_m\big)\\
&= C^{(1)}_s = C_{g_i(m)},
\qquad s=\mathfrak g^{(1)}_{i+2}(m)=g_i(m),
\endaligned
\end{equation}
%&
%&
\begin{equation}\label{5.18*}
\aligned
\phi^*(s_i)\cdot C^*_m\cdot \phi^*(s_i)^{-1} &=
\widehat\sigma_i\cdot C^*_m\cdot \widehat\sigma_i^{-1} =
\mathfrak g^{(k-1)}_i \big(C^{(k-1)}_m\big)\\
&= C^{(k-1)}_s = C^*_{g^*_i(m)},
\qquad s={\mathfrak g}^{(k-1)}_i(m) = g^*_i(m),
\endaligned
\end{equation}
%&
for any $i = 1,\ldots ,k-3$.
\vskip0.2cm

\noindent The image of $\phi$ is generated by the permutations
$\widehat\sigma_j$, \ $3\le j\le k-1$; since any such $\widehat\sigma_j$
commutes with $\widehat\sigma_1$,
Lemma \ref{Lm: invariant sets and components of commuting permutations}
implies that the set
$\Sigma = \supp \mathfrak C$ is invariant under the subgroup
$\Img \phi\subseteq {\mathbf S}(n)$. Similarly,
$\Img \phi^*$ is generated by the permutations
$\widehat\sigma_j$, \ $1\le j\le k-3$, and
the set $\Sigma^* = \supp {\mathfrak C}^*$ is invariant under
the subgroup $\Img \phi^*\subseteq {\mathbf S}(n)$.
\vskip0.2cm

\noindent Let
%&
\begin{equation}\label{5.19}
\varphi_{_{\Sigma}}\colon B_{k-2}\to
{\mathbf S}(\Sigma)\qquad\text{and}\qquad
\varphi^*_{_{\Sigma^{*}}}\colon B_{k-2}\to
{\mathbf S}(\Sigma^*)
\end{equation}
be the reductions of the homomorphisms $\phi$
and $\phi^*$ to the invariant sets $\Sigma$ and $\Sigma^*$,
respectively (see Definition \ref{Def: symmetric groups}$(d)$). That is,
%&
\begin{eqnarray}
&&\varphi_{_{\Sigma}}(s_i) \ = \phi(s_i)|_\Sigma
\ \ \ = \psi(\sigma_{i+2})|_\Sigma \ = \widehat\sigma_{i+2}|_\Sigma,
\label{5.20}\\
%\\
&&\varphi^*_{_{\Sigma^*}}(s_i) = \phi^*(s_i)|_{\Sigma^*}
= \psi(\sigma_i)|_{\Sigma^*} \ \ \ = \widehat\sigma_i|_{\Sigma^*}
\label{5.20*}
\end{eqnarray}
%&
($i = 1,...,k-3$). It follows from (\ref{5.18}),(\ref{5.18*}) that
%&
\begin{equation}\label{5.21}
\varphi_{_{\Sigma}}(s_i)\cdot C_m\cdot
\varphi_{_{\Sigma}}(s_i)^{-1} = C_{g_i(m)}
\end{equation}
%&
and
%&
\begin{equation}\label{5.21*}
\varphi^*_{_{\Sigma^*}}(s_i)\cdot C^*_m\cdot
\varphi^*_{_{\Sigma^*}}(s_i)^{-1} = C^*_{g^*_i(m)}
\end{equation}
%&
for all $m = 1,\ldots ,t$ and all $i = 1,\ldots ,k-3$, and thus
$$
\varphi_{_{\Sigma}}(s_i)\cdot \mathcal C\cdot
\varphi_{_{\Sigma}}(s_i)^{-1} = \mathcal C,\qquad
\varphi^*_{_{\Sigma^*}}(s_i)\cdot \mathcal C^*\cdot
\varphi^*_{_{\Sigma^*}}(s_i)^{-1} = \mathcal C^*.
$$
Therefore, 
$$
\Img\varphi_{_\Sigma}\subseteq G\,,
\qquad \Img \varphi^*_{_{\Sigma^*}}\subseteq G^*,
$$ 
which means that $\varphi_{_\Sigma}$ {\sl and $\varphi^*_{_{\Sigma^*}}$
may be regarded as homomorphisms from $B_{k-2}$ into
the groups $G$ and $G^*$ respectively.}
\vskip0.2cm

\noindent Relations (\ref{5.21}), $(\ref{5.21*})$ and the definitions of the
projection $\pi, \pi^*$ show that 
%&
\begin{equation}\label{5.22}
\pi(\varphi_{_{\Sigma}}(s_i)) = g_i\qquad {\text {and}}\qquad
\pi^*(\varphi^*_{_{\Sigma^*}}(s_i)) = g^*_i
 \ \ \ {\text {for all}} \ \
i=1,...,k-3.
\end{equation}
%&
Consider the compositions
%&
\begin{equation}\label{5.26}
\Omega = \pi\circ\varphi_{_{\Sigma}}
\colon B_{k-2}\stackrel{\varphi_{_{\Sigma}}}{-\!\!\!\!\longrightarrow} G
\stackrel{\pi}{\longrightarrow} {\mathbf S}(t)
\end{equation}
%&
and
%&
\begin{equation}\label{5.26*}
\Omega^* = \pi^*\circ\varphi^*_{_{\Sigma^*}}
\colon B_{k-2}\stackrel{\varphi^*_{_{\Sigma^*}}}
{-\!\!\!\!-\!\!\!\!-\!\!\!\!\longrightarrow} G^*
\stackrel{\pi^*}{\longrightarrow} {\mathbf S}(t).
\end{equation}
%&
The following simple lemma is in fact very important for us.

\begin{Lemma}\label{Lm: Omega=Omega*} 
\noindent $a)$ The homomorphisms
$$
\Omega\colon B_{k-2}\to{\mathbf S}(t) \ \ \text{and} \ \
\Omega^*\colon B_{k-2}\to{\mathbf S}(t)
$$
coincide.
\vskip0.2cm

\noindent $b)$ All the permutations $\mathfrak g^{(q)}_j$  
$(1\le q\le k-1; \ \ j\ne q-1,\,q,\,q+1)$ are conjugate to each other.
\end{Lemma}

\begin{proof} 
\noindent $a)$ Formulae (\ref{5.22}), (\ref{5.26}), 
and (\ref{5.26*}) show that
$\Omega (s_i) = g_i$ and $\Omega^*(s_i) = g^*_i$
for all $i = 1,\ldots ,k-3$. According to Lemma \ref{Lm 5.2}$(b)$ 
(see (\ref{5.12})), 
$g_i = g^*_i$ for all such $i$. Consequently, 
$\Omega  = \Omega^*$.
\vskip0.2cm

$b)$ Lemma \ref{5.2} implies that for
$1\le q\le k-1$ and $j\ne q-1,q,q+1$
the permutation $\mathfrak g^{(q)}_j$
coincides either with some $g_i$ or with some $g^*_i$.
Since the latter permutations coincide
with $\Omega (s_i)$ and the canonical generators $s_i$ 
are conjugate to each other, all the permutations 
$\mathfrak g^{(q)}_j$ are pairwise conjugate.
\end{proof} 

\begin{Definition}\label{Def 5.1} 
\noindent The homomorphism
$\Omega\colon B_{k-2}\to{\mathbf S}(t)$
defined by formulae (\ref{5.17}), (\ref{5.19}), 
(\ref{5.20}) and (\ref{5.26}) is called the 
{\em retraction} of the original normalized homomorphism
$\psi\in \Hom_{r,t}(B_k,{\mathbf S}(n))$
(to an $r$-component $\mathfrak C$ of $\widehat\sigma_1$).
According to Lemma \ref{5.3}, $\Omega$
coincides with the {\em co-retraction} $\Omega^*$ of $\psi$
defined by formulae (\ref{5.17*}), (\ref{5.19}), 
(\ref{5.20*}) and (\ref{5.26*}).
\end{Definition}

\noindent Since the set $\Sigma=\supp \mathfrak C$ is 
$(\Img \phi)$-invariant,
its complement $\Sigma'={\boldsymbol\Delta}_n\setminus\Sigma$ is also
$(\Img \phi)$-invariant, and we can consider the reduction
$\varphi_{_{\Sigma'}}\colon B_{k-2}
\to{\mathbf S}(\Sigma')$
of the homomorphism $\phi$ to $\Sigma'$:
%&
\begin{equation}\label{5.20'}
\varphi_{_{\Sigma'}}(s_i) = \phi(s_i)|_{\Sigma'}
= \psi(\sigma_{i+2})| {\Sigma'} \
= \widehat\sigma_{i+2}| {\Sigma'},\qquad i = 1,...,k-3.
\end{equation}

\begin{Lemma}
\label{Lm: if psi: B_k to S(n) is non-cyclic and varphi(Sigma') is abelian then phi, varphi(Sigma) and Omega are non-abelian} 
Assume that $k>6$ and that the homomorphism 
$\psi$ is non-cyclic. If the homomorphism
$\varphi_{_{\Sigma'}}$ is abelian, then the
homomorphisms $\phi$,\, $\varphi_{_{\Sigma}}$ and $\Omega$
are non-abelian.
\end{Lemma}

\begin{proof} 
If $\phi$ is abelian, then it is cyclic and
$\phi(s_3) = \phi(s_4)$ ($k-2>4$). So
$\psi(\sigma_5) = \psi(\sigma_6)$, contradicting
the assumption that $\psi$ is non-cyclic.
\vskip0.2cm

Since $\Sigma$ and $\Sigma'$ are disjoint, 
$\phi$ is the disjoint product of
the reductions $\varphi_{_{\Sigma}}$ and $\varphi_{_{\Sigma'}}$.
Since $\varphi_{_{\Sigma'}}$ is abelian and we have already proved
that $\phi$ is non-abelian, $\varphi_{_{\Sigma}}$ must be non-abelian.
\vskip0.2cm

Finally, assume that the homomorphism
$\Omega = \pi\circ\varphi_{_{\Sigma}}$ is abelian. Then
%&
\begin{equation}\label{5.27}
(\pi\circ\varphi_{_{\Sigma}})(B'_{k-2}) = \{1\}, \qquad
{\text{that is,}}\qquad
\varphi_{_{\Sigma}}(B'_{k-2})
\subseteq \Ker\pi = H.
\end{equation}
%&
Since $k-2>4$, the group $B'_{k-2}$ is perfect.
On the other hand, the group $H$ is abelian.
Hence (\ref{5.27}) implies that
$\varphi_{_{\Sigma}}(B'_{k-2}) = \{1\}$;
this means that the homomorphism $\varphi_{_{\Sigma}}$ is abelian,
which contradicts the statement proven above.
\end{proof}

\noindent The construction described above provides us with
the {\sl universal} exact sequence (\ref{5.2}) with the fixed
splitting $\rho$. This sequence and the homomorphisms 
$\varphi_{_\Sigma}$ and $\Omega$ defined by (\ref{5.20}),(\ref{5.27})
form the commutative diagram
%&
\begin{equation}\label{5.28}
\CD
@. @.  B_{k-2}@. @. \\
@. @.  @V{\varphi_{_{\Sigma}}}VV @/SE//{\,\Omega}/ @. \\ %{\Omega}
1@>>>H@>>>G@>>{\pi}>{\mathbf S}(t)@>>>1\\
\endCD
\end{equation}
%&
The homomorphism $\varphi_{_{\Sigma}}$ defined
by $\varphi_{_{\Sigma}}(s_i) = \psi(\sigma_{i+2})| \Sigma$, \
$1\le i\le k-3$  (see (\ref{5.17}) and (\ref{5.20})),
keeps a lot of information 
on the original normalized homomorphism $\psi$. Hence it seems reasonable 
to find out to which extent we can recover the homomorphism 
$\varphi_{_{\Sigma}}$ if we know the homomorphism $\Omega$.

\begin{Remark}\label{Rmk 5.1} 
Let us clarify the actual nature of this problem.
\vskip0.2cm

\noindent We are interested to classify (as far as possible) homomorphisms
$B_k\to{\mathbf S}(n)$ up to conjugation. If the permutation
$\widehat\sigma_1$ corresponding to such a
homomorphism $\psi$ has an $r$-component of length $t$, then,
without loss of generality, we may assume that $\psi$ is normalized.
So, we have diagram (\ref{5.28}) corresponding to this $\psi$.
Suppose that we can somehow find out what is the
homomorphism $\varphi_{_{\Sigma}}$. Then we know all the restrictions
$\psi(\sigma_3)| \Sigma,\ldots ,\psi(\sigma_{k-1})| \Sigma$.
This would provide us with an essential (and in some cases even sufficient)
information to determine the homomorphism $\psi$ itself.
The knowledge of all these restrictions
is certainly the best possible result, which we may hope to get
by studying diagram (\ref{5.28}).
\vskip0.2cm

\noindent Unfortunately, if $\psi$ is unknown to us, then
we know neither $\Omega$ nor $\varphi_{_{\Sigma}}$ in diagram (\ref{5.28}).
\vskip0.2cm

\noindent A reassuring circumstance is, however, 
that $k-2<k$ and $t\le n/r<n$.
Hence we may suppose that we succeeded in classifying
the homomorphisms $B_{k-2}\to{\mathbf S}(t)$ up
to conjugation, meaning that we have a finite list
of pairwise nonconjugate homomorphisms
$\Omega_p\colon B_{k-2}\to{\mathbf S}(t)$ \
($p=1,...,N$) such that any $\Omega\in\Hom(B_{k-2}, {\mathbf S}(t))$
is conjugate to one of $\Omega_p$'s. Moreover, suppose that for each
$\Omega_p$ we have classified up to conjugation
the homomorphisms $\varPhi\colon B_{k-2}\to G$
satisfying the commutativity condition $\pi\circ\varPhi = \Omega_p$.
If so, then for any $p=1,...,N$ we have a finite list
$\{\varPhi_{p,q_p}|  \ 1\le q_p\le M_p\}$ of
the pairwise nonconjugate representatives, and any
$\varPhi\in\Hom(B_{k-2},G)$ that satisfies
$\pi\circ\varPhi = \Omega_p$ is conjugate to one
of $\varPhi_{p,q_p}$.
\vskip0.2cm

\noindent Furthermore, let 
$\varphi_{_{\Sigma}}$ and $\Omega$ be the homomorphisms
related to our (unknown) normalized homomorphism
$\psi\in \Hom_{r,t}(B_k,{\mathbf S}(n))$. Then
$\Omega = s\Omega_p s^{-1}$ for some $p$ and some $s\in{\mathbf S}(t)$.
Using the splitting $\rho$, define the homomorphism
$\varPhi\colon B_{k-2}\to G$ by
$\varPhi = \rho(s^{-1})\varphi_{{\Sigma}}\rho(s)$.    %@
%Since $\pi\circ\rho = \id_{{\mathbf S}(t)}$ and
%$\pi\circ\varphi_{_{\Sigma}} = \Omega$, we have
%$\pi\circ\varPhi = s^{-1}\Omega s
%= s^{-1}\cdot s\Omega_p s^{-1}\cdot s = \Omega_p$.
%Hence,
It is easily seen that there are an element $g\in G$ and an index $q$
\ ($1\le q\le M_p$) such that
$\varPhi = g\varPhi_{p,q} g^{-1}$.
Since the element
$\widetilde g = \rho(s)g\in G\subset {\mathbf S}(\Sigma)
\subseteq {\mathbf S}(n)$,
we can define the homomorphism
$$
\widetilde\psi\colon B_k\to{\mathbf S}(n),
\qquad \widetilde\psi = \widetilde g^{-1}\psi \widetilde g
= g^{-1}\rho(s^{-1})\cdot\psi\cdot\rho(s)g.
$$
The condition $\widetilde g\in G$ means that
$\widetilde g\mathcal C \widetilde g^{-1} = \mathcal C$;
therefore,
\vskip0.3cm

\noindent $\widetilde\psi$ {\sl is a normalized homomorphism
in $\Hom_{r,t}(B_k,{\mathbf S}(n))$ conjugate to our original
homomorphism $\psi$}.
\vskip0.3cm

\noindent Let $\widetilde \varphi_{_\Sigma}$ and
$\widetilde \Omega$ be the homomorphisms related
to this homomorphism $\widetilde\psi$; then
$\pi\circ\widetilde \varphi_{_\Sigma} = \widetilde \Omega$.
The set $\Sigma = \supp \mathcal C$ is
$\widetilde\psi(\sigma_{i+2})$-invariant (for any $i$, \ $1\le i\le k-2$),
and (by definition) the permutation
$\widetilde \varphi_{_\Sigma}(s_i)$ coincides with the permutation
$$
\aligned
\widetilde\psi(\sigma_{i+2})| \Sigma =
\widetilde g^{-1}\cdot (\psi(\sigma_{i+2})| \Sigma)\cdot \widetilde g
&= g^{-1}\rho(s^{-1})\cdot \varphi_{_\Sigma}(s_i)\cdot \rho(s)g\\
= g^{-1}&\cdot \varPhi(s_i)\cdot g
= g^{-1}g\cdot \varPhi_{p,q}(s_i)\cdot g^{-1}g
= \varPhi_{p,q}(s_i),
\endaligned
$$
which means that \
$\widetilde\psi(\sigma_{i+2})| \Sigma = \varPhi_{p,q}(s_i)$ \
and \ $\widetilde \varphi_{_\Sigma} = \varPhi_{p,q}$.
These observations lead to the following
%\hfill $\bigcirc$
\end{Remark}

\noindent{\bfit Declaration.} 
Suppose that we solved the above mentioned
classification problems for homomorphisms
$B_{k-2}\to{\mathbf S}(t)$ and $B_{k-2}\to G$.
Hence we have the list of representatives $\{\varPhi_{p,q}\}$.
Then, without loss of generality, we may assume that the homomorphism
$\psi$ $($which we want to identify up to conjugation$)$,
besides the normalization condition $(!!)$, satisfies for some $p,q$
the condition
$$
\ \ (!!!) \qquad\qquad
\psi(\sigma_{i+2})| \Sigma = \varPhi_{p,q}(s_i)\qquad
{\text {for all}} \ \ i, \ 1\le i\le k-3.
\qquad\qquad\qquad \ \ \ \bigcirc
$$

We have almost nothing to say about the first classification problem.
If fact, this is the same problem which we started with, but
rather easier (since $t\le n/2$); in some cases it can be
solved, indeed. For instance, if $k\ne 6$ and $n<2k-4$ then $t<k-2$
and any homomorphism $\Omega\colon B_{k-2}\to{\mathbf S}(t)$
is cyclic (Theorem \ref{Thm: homomorphisms B(k) to S(n), n<k}$(a)$);
this puts a strict restriction to the original homomorphism $\psi$.
\vskip0.2cm

As to the second problem, it is as follows:

\begin{Problem} Given exact sequence (\ref{5.2}) and a homomorphism
$$
\Omega\colon B_{k-2}\to{\mathbf S}(t),
$$
find $($up to conjugation$)$ all the homomorphisms
$\varPhi\colon B_{k-2}\to G$ that
satisfy the commutativity relation $\pi\circ\varPhi = \Omega$.
\end{Problem}

We postpone the study of this problem to 
\ref{sec: Omega-homomorphisms and cohomology: some computations},
since we need first to develop an adequate tool; the next
subsection is devoted to this task.

\subsection{Homomorphisms and cohomology}
\label{Ss: 5.4. Homomorphisms and cohomology}
In this section we consider a diagram of the form
%&
\begin{equation}\label{5.29}
\CD
@. @.  B @. @. \\
@. @.  @. @/SE//{\,\Omega}/ @. \\ 
1@>>>H@>>>G@>>{\pi}> S @>>>1\,,\\
\endCD
\end{equation}
%&
where all the groups and all the homomorphisms are given,
and the horizontal line is an exact sequence
with some fixed splitting homomorphism
%&
\begin{equation}\label{5.30}
\rho\colon S\to G,\qquad \pi\circ\rho =\id_{S}.
\end{equation}
%&
Moreover, we assume that $H$ {\it is an abelian group}
and identify this group with its image under the given 
embedding $H\hookrightarrow G$.

\begin{Definition}\label{Def: Omega-homomorphisms} 
$a)$ {\bfit $\Omega$-homomorphisms.} A homomorphism %@
$\varPhi\colon B\to G$ is said to be an
{\em $\Omega$-homomorphism}, if $\pi\circ\varPhi = \Omega$.
The set of all $\Omega$-homomorphisms is
denoted by $\Hom_\Omega(B,G)$.
%\vskip0.2cm

$b)$ {\bfit The $S$- and $B$-actions on $H$ 
and the corresponding cohomology} 
\footnote{Whenever the exact sequence
$
1\to H\to G\stackrel{\pi}{\longrightarrow} S\to 1
$
with the splitting $\rho\colon S\to G$ in
diagram (\ref{5.29}) is given, the homomorphism
$\varepsilon=\varepsilon_\Omega$ and the action $T=T_\Omega$
defined by (\ref{5.31}),(\ref{5.32}) are determined 
by the homomorphism $\Omega$.
Therefore in our notation of the groups and homomorphisms
related to the corresponding cohomology
we use the sign of the homomorphism $\Omega$
instead of the traditional usage of the sign of an action.}.
We consider the composition
%&
\begin{equation}\label{5.31}
\varepsilon=\varepsilon_\Omega=
\rho\circ\Omega\colon B\stackrel{\Omega}{\longrightarrow}
S\stackrel{\rho}{\longrightarrow} G,\qquad
\pi\circ\varepsilon=\Omega,
\end{equation}
%&   
and define the left actions $\tau$ and $T=T_\Omega$ of the groups
$S$ and $B$ respectively, on the group $H$ by
%&
\begin{equation}\label{5.32}
\tau_s (h) = \rho (s)\cdot h\cdot\rho (s),\qquad 
T_b (h) = \tau_{\Omega(b)} (h) 
= \varepsilon(b)\cdot h\cdot\varepsilon(b^{-1}).
\end{equation}
%&
A mapping $z\colon B\to H$ with $z(1)=1$ is called
a $1$-{\it cochain on $B$ with values in $H$}. A $1$-cochain $z$ is
a $1$-{\it cocycle} if its {\it $1$-coboundary}
$\delta_\Omega^1 z\colon B\times B\to H$ is trivial, that is, if
%&
\begin{equation}\label{5.33}
\left(\delta_\Omega^1 z\right)(b_1,b_2)\Def
\left[T_{b_1}z(b_2)\right]\cdot \left[z(b_1b_2)\right]^{-1}
\cdot z(b_1)=1\qquad {\text{for all}}
\ \ b_1,b_2\in B.
\end{equation}
%&
The group of all $1$-cocycles is denoted by $\mathcal Z^1_\Omega(B,H)$.
The subgroup ${\mathcal B}^1_\Omega(B,H)\subseteq \mathcal Z^1_\Omega(B,H)$
consists of all {\em $0$-coboundaries}, that is, a $1$-cocycle
$z\colon B\to H$ belongs to ${\mathcal B}^1_\Omega(B,H)$ if
and only if
\vskip0.2cm

\begin{itemize} 

\item{\sl there is an element $h\in H$ such that} 

\end{itemize} 
%&
\begin{equation}\label{5.34}
z(b) = \left(\delta^0_\Omega h\right)(b)
\Def \left(T_b h\right)\cdot h^{-1}
\ \ \text{\sl for all} \ \, b\in B.
\end{equation}
%&
The cohomology group $H^1_\Omega(B,H)$ is defined by
$$
H^1_\Omega(B,H) = \mathcal Z^1_\Omega(B,H)/{\mathcal B}^1_\Omega(B,H).
$$
For any $\Omega$-homomorphism
$\varPhi\colon B\to G$, define the mapping
%&
\begin{equation}\label{5.35}
z_\varPhi\colon B\to G,\qquad
z_\varPhi(b) = \varPhi(b)\varepsilon(b^{-1}),
\end{equation}
%&
and vice versa, for any $1$-cocycle $z\in \mathcal Z^1_\Omega(B,H)$, define
the mapping
%&
\begin{equation}\label{5.36}
\varPhi_z\colon B\to G,\qquad
\varphi_z(b) = z(b)\varepsilon(b).
\end{equation}
%&
\hfill $\bigcirc$
\end{Definition}

The following simple lemma seems very well known; however, I could not
find it in standard textbooks in homological algebra.

\begin{Lemma}\label{Lm 5.5} 
$a)$ The mapping $z_\varPhi\colon B\to G$
defined by (\ref{5.35}) is, in fact, a $1$-cocycle of the group $B$
with values in $H$.

$b)$ The mapping $\varPhi_z\colon B\to G$
defined by {\rm (\ref{5.36})} is an $\Omega$-homomorphism, and, besides,
the $1$-cocycle $z_\varPhi$ corresponding to this
$\Omega$-homomorphism $\varPhi = \varPhi_z$ $($via statement
$(a))$ coincides with the 
original $1$-cocycle $z$.

Thereby, formulae (\ref{5.35}),(\ref{5.36}) define
the two (mutually inverse) one-to-one correspondence
%&
\begin{equation}\label{5.37}
\mathcal Z^1_\Omega(B,H)\ni z\mapsto\varPhi_z\in
\Hom_\Omega(B,G), \ \ \ \Hom_\Omega(B,G)\ni\varPhi\mapsto
z_\varPhi\in \mathcal Z^1_\Omega(B,H).
\end{equation}
%&
\end{Lemma}

\begin{proof} 
$a)$ Since $\pi\circ\varPhi = \Omega = \pi\circ\varepsilon$,
we have
$\pi\left(z_\varPhi(b)\right) =
\pi(\varPhi(b))\left(\pi\circ\varepsilon\right)(b^{-1}) =
\Omega(b)\Omega\left(b^{-1}\right) = 1$,
and thus $z_\varPhi(b)\in \Ker\pi = H$; moreover,
$z_\varPhi(1) = \varPhi(1)\varepsilon(1)=1$.
So, we can regard $z_\varPhi$ as a $1$-cochain of the group
$B$ with values in $H$. Further,
$$
\aligned
%1
\ &\left(\delta^1_\Omega z_\varPhi\right)(b_1,\,b_2) =
\left[ T_{b_1} z_{\varPhi}\left(b_2\right) \right]
\cdot \left[z_\varPhi\left(b_1 b_2\right)\right]^{-1}
\cdot z_\varPhi(b_1)\\
%2
&\\
&=T_{b_1}\left[\varPhi(b_2)\cdot \varepsilon\left(b_2^{-1}\right)\right]
\times
\left[\varPhi\left(b_1 b_2\right)
\cdot \varepsilon\left(\left(b_1 b_2\right)^{-1}\right)\right]^{-1}
\times \left[\varPhi\left(b_1\right)
\cdot \varepsilon\left(b_1^{-1}\right)\right]\\
%3
&\\
&=\varepsilon\left(b_1\right)\cdot \varPhi\left(b_2\right)
\cdot \varepsilon\left(b_2^{-1}\right)
\cdot \varepsilon\left(b_1^{-1}\right)
\times
\left[\varPhi\left(b_1 b_2\right)
\cdot \varepsilon\left(\left(b_1 b_2\right)^{-1}\right)\right]^{-1}
\times
\left[\varPhi(b_1)\cdot \varepsilon\left(b_1^{-1}\right)\right]\\
%4
&\\
&=\varepsilon(b_1) \varPhi(b_2) \varepsilon\left(b_2^{-1}\right)
 \varepsilon\left(b_1^{-1}\right)
 \varepsilon(b_1)
 \varepsilon(b_2)
 \varPhi\left(b_2^{-1}\right)
 \varPhi\left(b_1^{-1}\right)
 \varPhi(b_1) \varepsilon\left(b_1^{-1}\right)=1,
\endaligned
$$
which shows that $z_\varPhi$ is a $1$-cocycle.

$b)$ Since $z$ is a $1$-cocycle,
$\left(\delta^1_{\Omega} z\right)(b_1,\,b_2)=1$
for all $b_1,\,b_2\in B$, which means that
$\varepsilon(b_1)z(b_2)\varepsilon\left(b_1^{-1}\right)
\cdot \left[z(b_1 b_2)\right]^{-1}\cdot z(b_1)=1$. Since $H$ is abelian,
the latter relation may be written as
$$
z(b_1b_2)=z(b_1)\varepsilon(b_1)z(b_2)\varepsilon\left(b_1^{-1}\right);
$$
hence,
$$
\aligned
\varPhi_z(b_1b_2) &= z(b_1b_2)\varepsilon(b_1b_2)\\
&=z(b_1)\varepsilon(b_1)z(b_2)\varepsilon\left(b_1^{-1}\right)
\varepsilon(b_1b_2)
=z(b_1)\varepsilon(b_1)z(b_2)\varepsilon(b_2)
=\varPhi_z(b_1)\varPhi_z(b_2),
\endaligned
$$
which shows that $\varPhi_z\colon B\to G$ is a group
homomorphism.
Moreover, $z(b)\in H=\Ker\pi$ for any $b\in B$,
and $\pi\circ\varepsilon=\Omega$; thus,
$\pi(\varPhi_z(b))=\pi(z(b)\varepsilon(b))=\pi(z(b))\pi(\varepsilon(b))
=\Omega(b)$
and $\varPhi_z$ is an $\Omega$-homomorphism. Finally,
applying (\ref{5.35}) to the $\Omega$-homomorphism
$\varPhi = \varPhi_z$ and using (\ref{5.36}), we have
$$
z_\varPhi(b) = \varPhi(b)\varepsilon\left(b^{-1}\right)=
\varPhi_z(b)\varepsilon\left(b^{-1}\right)=
z(b)\varepsilon(b)\cdot \varepsilon\left(b^{-1}\right) = z(b),
$$
which concludes the proof.
\end{proof}

Our immediate goal is to study $\Omega$-homomorphisms
$B\to G$ up to conjugation. In view of the previous
lemma, it is useful to find out the binary relation
in $\mathcal Z_\Omega^1(B,H)$ corresponding to the conjugacy
relation ''$\sim$" for $\Omega$-homomorphisms.
The optimistic expectation that the equivalent cycles must be in the same
cohomology class is not very far from the truth.
Actually, it is so under some simple and soft
additional restriction on $\Omega$.

\begin{Definition}\label{Def 5.3} 
Two $\Omega$-homomorphisms
$\varPhi_1,\varPhi_2\colon B\to G$ are called
$H$-{\em conjugate}, if there exists an element $h\in H$
such that
%&
\begin{equation}\label{5.38}
\varPhi_2(b) = h\cdot\varPhi_1(b)\cdot h^{-1}
\end{equation}
%&
for all $b\in B$. If the latter condition holds, we write
$\varPhi_1\approx\varPhi_2$.
\end{Definition}

Clearly, $\approx$ {\sl is an equivalence relation on the set 
$\Hom_\Omega(B,G)$, which is stronger than the usual 
conjugacy relation} $\sim$ (that is, $\varPhi_1\approx\varPhi_2$
implies $\varPhi_1\sim\varPhi_2$).
%We denote by $\mathbf C(K_1,K)$ the centralizer of a subgroup $K_1$ 
%of a group $K$. Clearly,
There is the following evident inclusion involving the centralizers: \
$\pi\{\mathbf C(\pi^{-1}(\Img\Omega),G)\}\subseteq \mathbf C(\Img\Omega,S)$.
In general, this inclusion may be strict; however, if
\vskip0.3cm

$(i)$ $\mathbf C(\Img\Omega,S) = \{1\}$, that is,
{\sl the centralizer of the subgroup $\Omega(B)\subseteq S$ 
in $S$ is trivial}
(for instance, {\sl $\Omega$ is surjective, and the center of $S$
is trivial})
\vskip0.3cm

\noindent or the group $G$ is abelian, then we have
%&
\begin{equation}\label{5.39}
\pi\{\mathbf C(\pi^{-1}(\Img\Omega),G)\} = \mathbf C(\Im\Omega,S).
\end{equation}
%&

\begin{Proposition}
\label{Prp: Omega-homomorphisms and cocycles} 
Let
$\varPhi_1,\varPhi_2\colon B\to G$ be two
$\Omega$-homomorphisms, and let $z_1 = z_{\varPhi_1}$,
$z_2 = z_{\varPhi_2}$ be the corresponding $1$-cocycles.

$a)$ The relation $\varPhi_1\approx\varPhi_2$ holds if and only if
$z_1 z_2^{-1}\in {\mathcal B}^1_\Omega(B,H)$. Thus, the set of the
$\approx$-equivalence classes of $\Omega$-homomorphisms is 
in natural one-to-one correspondence with the cohomology group 
$H^1_\Omega(B,H)$.

$b)$ Assume that condition (\ref{5.39}) is held. Then the relations
$\varPhi_1\sim\varPhi_2$ and $\varPhi_1\approx\varPhi_2$ are
equivalent, and the set of the classes of conjugate 
$\Omega$-homomorphisms $B\to G$ is in natural 
one-to-one correspondence with the cohomology group $H^1_\Omega(B,H)$.
\end{Proposition}

%centralizer of the subgroup
%$\Img \Omega\subseteq S$ in the group $S$ is trivial
%$($that is, if $s\in S$ and $s\cdot\Omega(g)\cdot s^{-1} = 
%\Omega(g)$ for all $g\in G$, then $s = 1)$.

\begin{proof} 
$a)$ First, assume that
$\varPhi_1\approx\varPhi_2$; so, for some $h\in H$,
we have
$$
\varPhi_2(b) = h\cdot\varPhi_1(b)\cdot h^{-1}\qquad
{\text {for all}} \ b\in B.
$$
According to (\ref{5.35}), \
$
%\aligned
z_1(b) = \varPhi_1(b)\varepsilon(b^{-1}), \ \
%\\
z_2(b) = h\cdot\varPhi_1(b)\cdot h^{-1}\cdot \varepsilon(b^{-1}).
%\endaligned
$
Consequently,
%&
\begin{equation}\label{5.40}
\aligned
z_1(b)(z_2(b))^{-1} &=
\left[\varPhi_1(b)\cdot \varepsilon\left(b^{-1}\right)\right]\cdot
\left[h\cdot\varPhi_1(b)\cdot h^{-1}
\cdot \varepsilon\left(b^{-1}\right)\right]^{-1}\\
&\\
&= \left[\varPhi_1(b)\cdot \varepsilon\left(b^{-1}\right)\right]\cdot
\left[\varepsilon(b)\cdot h
\cdot \varPhi_{1}\left(b^{-1}\right)\cdot h^{-1}\right]\\
&\\
&= \left[\varPhi_1(b)\cdot \varepsilon\left(b^{-1}\right)\right]\cdot
\left[\varepsilon(b)
\cdot h\cdot \varepsilon\left(b^{-1}\right)\right]\cdot
\left[\varepsilon(b)\cdot\varPhi_{1}\left(b^{-1}\right)\right]\cdot [h^{-1}].
\endaligned
\end{equation}
%&
The four expressions in the brackets in the third line 
of (\ref{5.40}) belong to the {\it abelian} normal subgroup
$H\vartriangleleft G$, and
the first and the third of them are mutually
inverse; hence,
$$
z_1(b)(z_2(b))^{-1} = \left[\varepsilon(b)\cdot h\cdot
\varepsilon\left(b^{-1}\right)\right]\cdot h^{-1}
= \left(T_b h\right)\cdot h^{-1} = \left(\delta^0_\Omega h\right)(b)\qquad
{\text {for all}} \ b\in B,
$$
and $z_1z_2^{-1}\in {\mathcal B}^1_\Omega(B,H)$.

Now, let $z_1 z_2^{-1}\in {\mathcal B}^1_\Omega(B,H)$.
Then there is an element $h\in H$ such that for all $b\in B$
%&
\begin{equation}\label{5.41}
z_1 (b)\cdot \left[z_2 (b)\right]^{-1}
= \left(\delta^0_\Omega h\right)(b) = (T_b h)\cdot h^{-1}
= \left[\varepsilon(b)\cdot h
\cdot \varepsilon\left(b^{-1}\right)\right]\cdot h^{-1}.
\end{equation}
%&
By (\ref{5.36}), $\varPhi_j (b) = z_j (b)\varepsilon(b),
\ \ j = 1,2$. Using (\ref{5.41}) and commutativity of $H$, we have
$$
\aligned
\varPhi_2 (b)\cdot
\left[h\cdot\varPhi_1 \left(b^{-1}\right)\cdot h^{-1}\right] &=
\left[z_2 (b)\varepsilon(b)\right]\cdot
\left[h\cdot \varepsilon\left(b^{-1}\right)
\left(z_1 (b)\right)^{-1}\cdot h^{-1}\right] \\
&\\
&= z_2 (b)\cdot \left[\varepsilon(b)\cdot h\cdot
\varepsilon\left(b^{-1}\right)\right]
\cdot \left(z_1 (b)\right)^{-1}\cdot h^{-1} \\
&\\
&= \left[z_2 (b)\left(z_1 (b)\right)^{-1}\right]
\cdot \left[\varepsilon(b)\cdot h\cdot
\varepsilon\left(b^{-1}\right)\cdot h^{-1}\right]\\
&\\
&= z_2 (b)(z_1 (b))^{-1}\cdot
z_1 (b)(z_2 (b))^{-1} = 1;
\endaligned
$$
so, $\varPhi_2 (b)= h\cdot\varPhi_1 (b)\cdot h^{-1}$
for all $b\in B$ and $\varPhi_1\approx\varPhi_2$.

$b)$ In view of $(a)$, we should only prove that (under
condition \ref{5.39})
$\varPhi_1\sim\varPhi_2$ implies $\varPhi_1\approx\varPhi_2$.
The relation $\varPhi_1\sim\varPhi_2$ means that there exists
an element $g\in G$ such that
$\varPhi_2 (b) = g\cdot\varPhi_1 (b)\cdot g^{-1}$ for all
$b\in B$. Since $\varPhi_1$ and $\varPhi_2$ are
$\Omega$-homomorphisms, we have
$$
\Omega (b)
= (\pi\circ\varPhi_2)(b) =
\pi \left[g\cdot\varPhi_1 (b)\cdot g^{-1}\right]
= \pi (g)\cdot (\pi\circ\varPhi_1)(b)\cdot
\pi (g^{-1}) = \pi (g)\cdot\Omega (b)\cdot\pi (g^{-1})
$$
for all $b\in B$; thus, $\pi (g)\in \mathbf C(\Img\Omega,S)$.
It follows from (\ref{5.39}) that there is an element
$\widetilde g\in \mathbf C(\pi^{-1}(\Img\Omega),G)$ such that
$\pi(\widetilde g)=\pi(g)$; clearly, the element
$h=g \widetilde g^{-1}$ is in $H$. The element $\widetilde g$
commutes with any element of the subgroup $\pi^{-1}(\Img\Omega)$.
This subgroup contains the image of any $\Omega$-homomorphism
$B\to G$; hence $\widetilde g$ commutes with all the elements
$\varPhi_1(b)$, \ $b\in B$, and
$$
\varPhi_2(b) =  g\cdot \varPhi_1(b)\cdot g^{-1}
= h\widetilde g\cdot \varPhi_1(b)\cdot (h\widetilde g)^{-1}
= h\cdot \varPhi_1(b)\cdot h^{-1}.
$$
This shows that $\varPhi_1\approx\varPhi_2$.
\end{proof}

\begin{Remark}\label{Rmk 5.2} 
If we replace $\Omega$ by a conjugate homomorphism
$\Omega'$,
$$
\Omega'(b) = s\Omega(b)s^{-1}\,,
$$
and define the corresponding
$\varepsilon'=\varepsilon_{\Omega'}$
and $T'=T_{\Omega'}$ according to (\ref{5.31}),(\ref{5.32}),
then we have the bijection
$$
\mathcal Z^1_{\Omega}(B,H)\ni z\mapsto z'\in 
\mathcal Z^1_{\Omega'}(B,H),\qquad
z'(b)=\rho(s)z(b)\rho(s^{-1}),
$$
which induces an isomorphism of the cohomology groups
$H^1_{\Omega}(B,H)\stackrel{\cong}{\longrightarrow}
H^1_{\Omega'}(B,H)$. We have also the bijection
$$
\Hom_\Omega(B,G)\ni\varPhi\mapsto\varPhi'\in\Hom_{\Omega'}(B,G),
\qquad \varPhi'(b)=\rho(s)\varPhi(b)\rho(s^{-1}),
$$
which is compatible with the equivalence relations
$\approx$, \ $\approx'$. 
The matching between cocycles and ($\Omega$- or $\Omega'$-) 
homomorphisms defined by (\ref{5.35}),(\ref{5.36}) is also compatible 
with the above bijections.
Moreover, if $\Omega$ satisfies (\ref{5.39}), 
then $\Omega'$ does as well.
Combined with Remark \ref{Rmk 5.1}, 
this shows that in our problem we may freely pass
from a homomorphism $\Omega$ to a conjugate one.
\end{Remark}
\vskip0.3cm

\begin{Remark}\label{Rm 5.3} 
Even if condition (\ref{5.39}) is not held,
we may compute the cohomology group $H^1_\Omega(B,H)$, 
choose some $1$-cocycle $z_{\mathcal H}$ in each cohomology class
${\mathcal H}$, and then take the corresponding
$\Omega$-homomorphisms $\varPhi_{\mathcal H} 
= \varPhi_{z_{\mathcal H}}$.
The homomorphisms $\varPhi_{\mathcal H}$ corresponding to
distinct cohomology classes ${\mathcal H}$ cannot be $H$-conjugate;
but some of them can be conjugate (by means of an element in $G$).
Even if this happens, the homomorphisms $\varPhi_{\mathcal H}$,
${\mathcal H}\in H^1_\Omega(B,H)$, form a {\it complete system} of
$\Omega$-homomorphisms $B\to G$, meaning that 
any $\Omega$-homomorphisms $\varPhi\colon B\to G$
is conjugate to some of $\varPhi_{\mathcal H}$.
Such a system $\{\varPhi_{\mathcal H}\}$ provides us with a solution
of our problem (maybe, not with ``the best" one; 
see Remark \ref{Rmk 5.1}). 
Anyway, this procedure reduces our nonlinear classification problem  
to computation of cohomology, which is, in a sense, a linear algebra
problem.
\end{Remark}
\vskip0.3cm
\vfill

\newpage

%Sec. 6
\section{$\Omega$-homomorphisms and cohomology: some computations}
\label{sec: Omega-homomorphisms and cohomology: some computations}

\noindent In this section we study some particular diagrams 
of the form (\ref{5.29}) and compute the corresponding 
cohomology and $\Omega$-homomorphisms.
Namely, we consider a diagram
%&
\begin{equation}\label{6.1}
\CD
@. @.  B_n @. @. \\
@. @.  @. @/SE//{\,\Omega}/ @. \\ 
1@>>>H@>>>G@>>{\pi}> {\mathbf S}(t)@>>>1\,,\\
\endCD
\end{equation}
%&
where $\Omega\colon B_n\to{\mathbf S}(t)$
is some given homomorphism of the braid group $B_n$
to the symmetric group ${\mathbf S}(t)$ \ ($n,t\ge 2$).
We fix some abelian group $A$, which is written as
{\sl additive}, and assume that $H$ is the direct sum of $t$ copies of $A$:
%&
\begin{equation}\label{6.2}
H = A^{\oplus t} = \bigoplus_{j=1}^t A.
\end{equation}
%&
We denote elements of the group $H$ by bold letters (say $\mathbf h$)
and regard them as ``vectors" with
$t$ ``coordinates" in $A$: \
$\mathbf h = (a^1,...,a^t)\in H$, \ \ $a^1,...,a^t\in A$.
We consider the {\sl standard left action} $\tau$ of
the symmetric group ${\mathbf S}(t)$ on this group $H$.
Namely, for any element $\mathbf h = (a^1,...,a^t)\in H$
and any $s\in{\mathbf S}(t)$, we put
%&
\begin{equation}\label{6.3}
\tau_s \mathbf h = (a^{s^{-1}(1)},\ldots ,a^{s^{-1}(t)}).
\end{equation}
%&
The group $G$ is assumed to be the {\sl semidirect product}
$H\leftthreetimes_\tau{\mathbf S}(t)$ of the groups $H$ and ${\mathbf S}(t)$
corresponding to the action $\tau$. That is,
$G$ is the set of all pairs $(\mathbf h,s)$, \
$\mathbf h\in H$, \ $s\in{\mathbf S}(t)$, with the multiplication
%&
\begin{equation}\label{6.4}
(\mathbf h,s)\cdot (\mathbf h' ,s') =
(\mathbf h + \tau_s \mathbf h', s\cdot s').
\end{equation}
%&
The injection \ $\mathbf j\colon H\hookrightarrow G$, \ the projection \
$\pi\colon G\to{\mathbf S}(t)$, \ the splitting
homomorphism \ $\rho\colon {\mathbf S}(t)\to G$, \
and the homomorphism \ $\varepsilon\colon B_n\to G$ \ are
defined as follows:
%&
\begin{equation}\label{6.5}
\mathbf j(\mathbf h) = (\mathbf h,1), \ \ \pi (\mathbf h,s) = s,
\ \ \rho (s) = (\mathbf 0,s),
\ \ \varepsilon(b)=(\rho\circ\Omega)(b)=(\mathbf 0,\Omega(b)).
\end{equation}
%&
We identify any element $\mathbf h\in H$ with its image
$\mathbf j(\mathbf h) = (\mathbf h,1)\in G$ (however, we must remember
that the group $H$ is additive, and $G$ is multiplicative).
The left action $T$ of the group $B_n$ on the group $H$
is defined by the given homomorphism $\Omega$ and the action $\tau$:
%&
\begin{equation}\label{6.6}
T_b \mathbf h=\varepsilon(b) (\mathbf h,1) \varepsilon(b^{-1})
=(\mathbf 0,\Omega(b)) (\mathbf h,1) (\mathbf 0,\Omega(b^{-1}))
=(\tau_{\Omega(b)} \mathbf h,1) = \tau_{\Omega(b)} \mathbf h.
\end{equation}
%&
For a cocycle $\mathbf z\in \mathcal Z^1_T(B_n,H)$,
we have $\mathbf z(1) = \mathbf 0$ and
$$
(\delta_T^1 \mathbf z)(b_1,b_2) = T_{b_1}\mathbf z(b_2)
- \mathbf z(b_1b_2) + \mathbf z(b_1) = \mathbf 0
\qquad {\text{for all}} \ \ b_1,b_2\in B_n
$$
(this is the additive version of (\ref{5.33})); thus,
%&
\begin{equation}\label{6.7}
\mathbf z(b_1b_2) = \mathbf z(b_1) + T_{b_1}\mathbf z(b_2).
\end{equation}
%&
In particular, setting $b_1 = b$ and $b_2 = b^{-1}$, we obtain
%&
\begin{equation}\label{6.8}
\mathbf z(b^{-1}) = -T_{b^{-1}}\mathbf z(b).
\end{equation}
%&
It follows from (\ref{6.7}),(\ref{6.8}) that {\sl any $1$-cocycle $\mathbf z$
is completely determined by its values}
%&
\begin{equation}\label{6.9}
\mathbf h_i = \mathbf z(s_i) =
(z_i^1,...,z_i^t)\in H,\qquad z_i^j\in A,
\end{equation}
%&
{\sl on the canonical generators $s_1,...,s_{n-1}$ of the group}
$B_n$.

Assume that some elements
%&
\begin{equation}\label{6.10}
\mathbf h_i = (a_i^1,...,a_i^t)\in H,\qquad
a_i^j\in A,\qquad 1\le i\le n-1,
\end{equation}
%&
are given, and we are looking for a cocycle
$\mathbf z\in \mathcal Z^1_T(B_n,H)$ with the values
%&
\begin{equation}\label{6.11}
\mathbf z(s_i) = \mathbf h_i,\qquad 1\le i\le n-1.
\end{equation}
%&
Since $s_p s_q = s_q s_p$ whenever $1\le p,q\le n-1$ and
$\modo{p-q}\ge 2$, we have for any such $p,q$ the relation
$\mathbf z(s_p s_q) = \mathbf z(s_q s_p)$. In view of (\ref{6.7}),
the latter relations may be written as
%&
\begin{equation}\label{6.12}
T_{s_p}\mathbf z(s_q) + \mathbf z(s_p)
= T_{s_q}\mathbf z(s_p) + \mathbf z(s_q),
\end{equation}
%&
which shows that the elements {\ref{6.10}) must satisfy the relations
%&
\begin{equation}\label{6.13}
\mathbf h_q - T_{s_p}\mathbf h_q =
\mathbf h_p - T_{s_q}\mathbf h_p ,\qquad \ 1\le p,q\le n-1, \ \
|p-q|\ge 2.
\end{equation}
%&
We must also take into account the relations
$s_p s_{p+1} s_p = s_{p+1} s_p s_{p+1}$, \ $1\le p<n-1$, which
leads to the conditions
%&
\begin{equation}\label{6.12'}
T_{s_p s_{p+1}}\mathbf z(s_p) + T_{s_p}\mathbf z(s_{p+1}) + \mathbf z(s_p) =
T_{s_{p+1} s_p}\mathbf z(s_{p+1}) + T_{s_{p+1}}\mathbf z(s_p) 
+ \mathbf z(s_{p+1})
\end{equation}
%&
and
%&
\begin{equation}\label{6.13'}
\mathbf h_p - T_{s_{p+1}}\mathbf h_p + T_{s_p s_{p+1}}\mathbf h_p =
\mathbf h_{p+1} - T_{s_p}\mathbf h_{p+1} + T_{s_{p+1} s_p}\mathbf h_{p+1},
\ \ 1\le p\le n-2.
\end{equation}
%&
The following lemma is evident.

\begin{Lemma}\label{Lm 6.1} 
\noindent A cocycle $\mathbf z\in \mathcal Z^1_T(B_n,H)$
with the values $\mathbf z(s_i) = \mathbf h_i$ \ $(1\le i\le n-1)$
does exist if and only if the elements $\mathbf h_i$ satisfy
relations (\ref{6.13}), (\ref{6.13'}). If these relations hold, then
the corresponding cocycle $\mathbf z$ is uniquely determined by
the elements $\mathbf h_i$. Moreover, this cocycle $\mathbf z$
is a coboundary if and only if there exists an element
$\mathbf h\in H$ such that
%&
\begin{equation}\label{6.14}
\mathbf h_p = T_{s_p}\mathbf h - \mathbf h\qquad
{\text {for all}} \ \ p = 1,...,n-1.
\end{equation}
%&
The cohomology group $H^1_T(B_n,H)$
is isomorphic to the quotient group $\mathcal Z/{\mathcal B}$,
where $\mathcal Z$ consists of all the solutions
$[\mathbf h_1,...,\mathbf h_{n-1}]$
of the linear system (\ref{6.13}),(\ref{6.13'}), and
${\mathcal B}\subseteq \mathcal Z$ is the subgroup consisting of all
the solutions $[\mathbf h_1,...,\mathbf h_{n-1}]$
such that there is an element $\mathbf h\in H$ that satisfies (\ref{6.14})
\end{Lemma}

If the action $T$ is given explicitly, the computation of
the quotient group $\mathcal Z/{\mathcal B}$ is a routine (however,
it can be very long).
\vskip0.2cm

In the sequel we use the coordinate representation
(\ref{6.10}) of the vectors $\mathbf h_i\in H$.

\begin{Lemma}\label{Lm 6.2} 
Let $t = n$ and let $\Omega=\mu\colon B_n\to{\mathbf S}(n)={\mathbf S}(t)$
be the canonical epimorphism. Then the system
of equations $(\ref{6.13})$ is equivalent to the system
%&
\begin{equation}\label{6.15}
a_q^p = a_q^{p+1},\qquad 1\le p,q\le n-1, \ \ |p-q|\ge 2,
\end{equation}
%&
and $(\ref{6.13'})$ is equivalent to
%&
\begin{equation}\label{6.15'}
\aligned
\ &a_q^p = a_{q+1}^p,\qquad\qquad\qquad \ \ 1\le q \le n-2,
\ \ 1\le p\le n, \ \ p\ne q,q+1,q+2;\\
&a_q^{q+2} = a_{q+1}^q,\qquad\qquad \ \ \ \ \ 1\le q\le n-2;\\
&a_q^q + a_q^{q+1} = a_{q+1}^{q+1} + a_{q+1}^{q+2}, \ \ 1\le q\le n-2.
\endaligned
\end{equation}
%&
\end{Lemma}

\begin{proof} 
Take any element $\mathbf h = (a^1,\ldots ,a^n)\in H$.
Using Definition \ref{Def: Omega-homomorphisms}$b)$ of the actions $\tau$
and $T$ (see formulae (\ref{5.31}) and (\ref{5.32})) and
taking into account that for $\Omega=\mu$ we have
$$
\Omega(s_p)=(p,p+1), \ \ \Omega(s_p s_{p+1})=(p,p+1,p+2),\ \
\text{and} \ \
\Omega(s_{p+1} s_p)=(p+2,p+1,p)\,,
$$
we can readily compute that
%&
\begin{equation}\label{6.16}
\aligned
T_{s_p} \mathbf h \ \ \ \ \ &=(a^1,\ldots ,a^{p-1},
\underbrace {a^{p+1},a^p},a^{p+2},a^{p+3},\ldots ,a^n),\\
T_{s_p s_{p+1}} \mathbf h &=(a^1,\ldots ,a^{p-1}, 
\underbrace{a^{p+2},a^p,a^{p+1}},a^{p+3},\ldots ,a^n),\\
T_{s_{p+1} s_p} \mathbf h&=(a^1,\ldots ,a^{p-1}, 
\underbrace{a^{p+1},a^{p+2},a^p}, a^{p+3},\ldots ,a^n)
\endaligned
\end{equation}
%&
(we underbrace the ``nonregular" permuted parts).
Using these formulae, we can write
relations (\ref{6.13}), (\ref{6.13'}) in the coordinates; after
evident cancellations, this leads to (\ref{6.15}) and (\ref{6.15'}),
respectively.
\end{proof}

Now we can compute certain cohomology and homomorphisms.

\begin{Remark}\label{Rmk 6.1} 
Assume that $A ={\mathbb Z}/r{\mathbb Z}$ and
take the following $t$ disjoint $r$-cycles:
$$
C_m = ((m-1)r+1, (m-1)r+2,\ldots , mr)\in{\mathbf S}(rt),
\ \ \ m = 1,\ldots ,t.
$$
Identify any
$\mathbf h = (a^1,...,a^t)\in ({\mathbb Z}/r{\mathbb Z})^{\oplus t}$
with the product $C_1^{a^1}\cdots C_t^{a^t}\in{\mathbf S}(rt)$.
Thereby, we obtain an embedding
$({\mathbb Z}/r{\mathbb Z})^{\oplus t}\hookrightarrow {\mathbf S}(rt)$.
Using this embedding and Lemma \ref{Lm 5.1}, we may identify the second
horizontal line of diagram (\ref{6.1}) 
with the exact sequence (\ref{5.2}). 
This identification is compatible with the actions, splittings, etc.
This means that for the group $A ={\mathbb Z}/r{\mathbb Z}$
any $\Omega$-homomorphism in diagram (\ref{6.1}) may be regarded
as a homomorphism $B_n\to G\subset {\mathbf S}(rt)$.
\vskip0.2cm

\noindent When $A$ is a commutative ring with unity $1$
(say $A ={\mathbb Z}/r{\mathbb Z}$), we set
%&
\begin{equation}\label{6.17}
\mathbf e_i = (\underbrace{0,\ldots ,0}_{i-1 \ \text{times}} ,
1, \underbrace{0,\ldots ,0}_{t-i \ \text{times}} )\in H, \ \
1\le i\le t.
\end{equation}
%&
Clearly, in this case $H$ is a free $A$-module with the free base
$\mathbf e_1,\ldots ,\mathbf e_t$, and the action $T$ on $H$ is compatible
with the $A$-module structure of $H$.
Hence the cohomology group is also an $A$-module.
\end{Remark}
\vskip0.3cm

\begin{Theorem}\label{Thm 6.3} 
Let $t = n$ and let $\Omega=\mu\colon B_n\to{\mathbf S}(n)
= {\mathbf S}(t)$ be the canonical epimorphism.
Then
$$
H^1_{\mu}(B_n,A^{\oplus n})\cong A\oplus A.
$$
If $A$ is a ring with unit, then the $A$-module
$H^1_{\mu}(B_n,A^{\oplus n})$ is generated
by the cohomology classes of the following two cocycles
$\mathbf z_1$, $\mathbf z_2$:
%&
\begin{equation}\label{6.18}
\aligned
\ &\mathbf z_1(s_i) = \mathbf e_{i+1},\\
&\mathbf z_2(s_i) = \mathbf e_1 + \ldots +
\mathbf e_{i-1} + \mathbf e_{i+2} + \ldots + \mathbf e_n.
\endaligned
\qquad (1\le i\le n-1)
\end{equation}
%&
\end{Theorem}

\begin{proof} 
$a)$ Any solution $\mathbf h_1,...,\mathbf h_{n-1}$
of the system of equations (\ref{6.15}),(\ref{6.15'}) is of the form
$$
\mathbf h_i 
=(\underbrace {b,...,b}_{i-1 \ \text{times}} , \
c_i, \ a-c_i, \ 
\underbrace {b,...,b}_{n-i-1 \ \text{times}}),
\qquad i = 1,...,n-1,
$$
where $a, b$ and $c_1,...,c_{n-1}$ are arbitrary
elements of the group $A$ (the elements $a$ and $b$ do not depend
on $i$). Hence
$$
\mathcal Z = \mathcal Z^1_{\mu}(B_n,H)\cong A^{\oplus n+1}
= \bigoplus_{i=1}^{n+1} A\,.
$$
Such a solution satisfies the system of equations (\ref{6.14}) for
some $\mathbf h\in H$ {\sl if and only if} $a = b = 0$. That is,
the subgroup
$$
{\mathcal B} = {\mathcal B}^1_{\mu}(B_n,H)\subseteq
\mathcal Z^1_{\mu}(B_n,H)
$$
consists of all $[\mathbf h_1,...,\mathbf h_{n-1}]$ of the form
$$
\mathbf h_i=(\underbrace {0,...,0}_{\text {$i-1$ times}}, \
c_i, \ -c_i, \ 
\underbrace {0,...,0}_{\text {$n-i-1$ times}})\,,
\ \ \ i = 1,...,n-1.
$$
This implies that any cohomology class in $H^1_{\mu}(B_n,H)$
contains a unique cocycle $\mathbf z\in \mathcal Z^1_{\mu}(B_n,H)$
that takes on the generators $s_1,...,s_{n-1}$ the values of the form
%&
\begin{equation}\label{6.19}
\mathbf z(s_i) = \mathbf h_i
= (\underbrace {b,...,b}_{\text {$i-1$ times}}, \ 0, \ a, \
\underbrace{b,...,b}_{\text {$n-i-1$ times}})\,,
\ \ \ i = 1,...,n-1,
\end{equation}
%&
where elements $a,b\in A$ do not depend on $i$.
Therefore, $H^1_{\mu}(B_n,A^{\oplus n})\cong A\oplus A$.
\vskip0.2cm

\noindent Now, if $A$ is a ring (with unity), 
any cocycle of the form (\ref{6.19}) may
be represented as $\mathbf z = a\cdot \mathbf z_1 + b\cdot \mathbf z_2$,
where $\mathbf z_1, \mathbf z_2$ are defined by (\ref{6.18}).
\end{proof}

\begin{Remark}\label{Rmk 6.2} 
Remark \ref{Rmk 5.2} and Theorem \ref{Thm: Improvement of Artin Theorem}
show that for $n\ne 4,6$
Theorem \ref{Thm 6.3} applies in fact to any {\sl non-cyclic} homomorphism
$\Omega\colon B_n\to{\mathbf S}(n)$. Namely, for such
a homomorphism $\Omega$, we have
$H^1_\Omega(B_n, A^{\oplus n})
\cong H^1_\mu(B_n, A^{\oplus n})$.
\end{Remark}
\vskip0.3cm

\noindent Now we should compute the $A^{\oplus 6}$-cohomology for the
two exceptional homomorphisms
$\Omega=\psi_{5,6}\colon B_5\to{\mathbf S}(6)$ and
$\Omega=\nu_6\colon B_6\to{\mathbf S}(6)$,
where  $\psi_{5,6}$ is defined by (\ref{4.6'})
and $\nu_6$ is Artin's homomorphism 
(see Sec. \ref{Ss: Transitive homomorphisms B(k) to S(k)},
Artin Theorem $(b)$). We
skip some completely elementary details, since the computations are
very long.
\vskip0.3cm

\noindent {\bfit Notation.} For an additive abelian group $A$,
we denote by $A_2$ the subgroup
in $A$ consisting of the zero element and all elements of order $2$: 
\ $A_2 = \{a\in A| \, 2a=0\}$.
\hfill $\bigcirc$

\begin{Theorem}\label{Thm 6.4} Let $n = 5$, $t=6$ and let
$\Omega=\psi_{5,6}\colon B_5\to{\mathbf S}(6)$ be the
homomorphism defined by (\ref{4.6}),(\ref{4.6'}). Then
$$
H^1_{\psi_{5,6}}(B_5,A^{\oplus 6}) \cong A_2\oplus A.
$$
Moreover, any cohomology class contains a cocycle $\mathbf z$
that takes on the canonical generators
$s_i\in B_5$ the values $\mathbf z(s_i) = \mathbf h_i$ of the form
%&
\begin{equation}\label{6.20}
\aligned
\ &\mathbf h_1 = ( \ \ 0, \ \ x+2y, \qquad 0,
\qquad x+2y, \qquad 0, \ \ \ \ x+2y)\\
&\mathbf h_2 = ( \ \ 0, \qquad 0, \qquad x+2y,
\ \ x+2y, \ \ x+2y, \ \ \ \ 0 \ \ \ \ )\\
&\mathbf h_3 = (-y, \ \ \ -y, \qquad x+3y,
\ \ x+3y, \qquad x, \qquad 2y \ \ \ )\\
&\mathbf h_4 = ( \ \ x, \ \ \ \ 2y, \qquad \ \ 2y,
\qquad \ 2y, \qquad \ \ \ \ x, \qquad \ x \ \ \ )
\endaligned
%\eqnum{6.20}
\end{equation}
%&
where $x\in A_2$ and $y\in A$.
\end{Theorem}

\begin{proof}
The action $T=T_{\psi_{5,6}}$ of $B_5$ on
the group $A^{\oplus 6}$ corresponding to the homomorphism
$\psi_{5,6}$ is given by
$$
\aligned
\ &T_{s_1}(a^1,a^2,a^3,a^4,a^5,a^6) = (a^2,a^1,a^4,a^3,a^6,a^5), \\
&T_{s_2}(a^1,a^2,a^3,a^4,a^5,a^6) = (a^5,a^3,a^2,a^6,a^1,a^4), \\
&T_{s_3}(a^1,a^2,a^3,a^4,a^5,a^6) = (a^3,a^4,a^1,a^2,a^6,a^5), \\
&T_{s_4}(a^1,a^2,a^3,a^4,a^5,a^6) = (a^2,a^1,a^5,a^6,a^3,a^4).
\endaligned
$$
Let $\mathbf z$ be a cocycle with the values
$\mathbf z(s_i) = \mathbf h_i = (a_i^1,...,a_i^6)\in A^{\oplus 6}$,
\ $i=1,2,3,4$. Then the system of equations
corresponding to the commutativity relations
$s_i\rightleftarrows s_j \ (\modo{i-j}>1)$ and the braid relations
$s_i\infty s_{i+1}$ between the generators $s_i$ looks as follows:
$$
\aligned
&\qquad a_1^1 - a_1^3 = a_3^1 - a_3^2;
\qquad a_1^2 - a_1^4 = a_3^2 - a_3^1;
\qquad a_1^3 - a_1^1 = a_3^3 - a_3^4;\\
&\qquad a_1^4 - a_1^2 = a_3^4 - a_3^3;
\qquad a_1^5 - a_1^6 = a_3^5 - a_3^6;
\qquad a_1^6 - a_1^5 = a_3^6 - a_3^5;
\endaligned
\eqnum{s_1\rightleftarrows s_3}
$$

$$
\aligned
&\qquad a_1^1 - a_1^2 = a_4^1 - a_4^2;
\qquad a_1^2 - a_1^1 = a_4^2 - a_4^1;
\qquad a_1^3 - a_1^5 = a_4^3 - a_4^4;\\
&\qquad a_1^4 - a_1^6 = a_4^4 - a_4^3;
\qquad a_1^5 - a_1^3 = a_4^5 - a_4^6;
\qquad a_1^6 - a_1^4 = a_4^6 - a_4^5;
\endaligned
\eqnum{s_1\rightleftarrows s_4}
$$

$$
\aligned
&\qquad a_2^1 - a_2^2 = a_4^1 - a_4^5;
\qquad a_2^2 - a_2^1 = a_4^2 - a_4^3;
\qquad a_2^3 - a_2^5 = a_4^3 - a_4^2;\\
&\qquad a_2^4 - a_2^6 = a_4^4 - a_4^6;
\qquad a_2^5 - a_2^3 = a_4^5 - a_4^1;
\qquad a_2^6 - a_2^4 = a_4^6 - a_4^4;
\endaligned
\eqnum{s_2\rightleftarrows s_4}
$$

$$
\aligned
&a_1^1 - a_1^5 + a_1^3 = a_2^1 - a_2^2 + a_2^6; \ \
a_1^2 - a_1^3 + a_1^3 = a_2^2 - a_2^1 + a_2^4; \\
&a_1^3 - a_1^2 + a_1^6 = a_2^3 - a_2^4 + a_2^1;\ \
a_1^4 - a_1^6 + a_1^2 = a_2^4 - a_2^3 + a_2^5; \\
&a_1^5 - a_1^1 + a_1^4 = a_2^5 - a_2^6 + a_2^2; \ \
a_1^6 - a_1^4 + a_1^1 = a_2^6 - a_2^5 + a_2^3;
\endaligned
\eqnum{s_1\infty s_2}
$$

$$
\aligned
&a_2^1 - a_2^3 + a_2^6 = a_3^1 - a_3^5 + a_3^2; \ \
a_2^2 - a_2^4 + a_2^1 = a_3^2 - a_3^3 + a_3^6; \\
&a_2^3 - a_2^1 + a_2^4 = a_3^3 - a_3^2 + a_3^5;\ \
a_2^4 - a_2^2 + a_2^5 = a_3^4 - a_3^6 + a_3^3; \\
&a_2^5 - a_2^6 + a_2^3 = a_3^5 - a_3^1 + a_3^4; \ \
a_2^6 - a_2^5 + a_2^2 = a_3^6 - a_3^4 + a_3^1;
\endaligned
\eqnum{s_2\infty s_3}
$$

$$
\aligned
&a_3^1 - a_3^2 + a_3^5 = a_4^1 - a_4^3 + a_4^4; \ \
a_3^2 - a_3^1 + a_3^6 = a_4^2 - a_4^4 + a_4^3; \\
&a_3^3 - a_3^5 + a_3^2 = a_4^3 - a_4^1 + a_4^6; \ \
a_3^4 - a_3^6 + a_3^1 = a_4^4 - a_4^2 + a_4^5; \\
&a_3^5 - a_3^3 + a_3^4 = a_4^5 - a_4^6 + a_4^1; \ \
a_3^6 - a_3^4 + a_3^3 = a_4^6 - a_4^5 + a_4^2.
\endaligned
\eqnum{s_3\infty s_4}
$$
Straightforward computations show that the general solution
$\mathbf h_i = (a_i^1,...,a_i^6)$ \ $(1\le i\le 4)$ of this system
depends linearly (over ${\mathbb Z}$) on seven parameters
$x_1,\ldots ,x_7\in A$ that must satisfy the only relation
%&
\begin{equation}\label{6.21}
2x_1 = 0.
%\eqnum{6.21}
\end{equation}
%&
Explicitly, the general solution is of the form
$$
\aligned
\ &a_1^1 = x_3 \\
&a_1^2 = x_1+2x_2-x_3\\
&a_1^3 = x_4
\endaligned
\qquad\qquad\qquad
\aligned
\ &a_1^4 = x_1+2x_2-x_4\\
&a_1^5 = x_5 \\
&a_1^6 = x_1+2x_2-x_5
\endaligned\qquad\qquad\quad
\eqnum{\mathbf h_1}
$$
\vskip0.3cm
$$
\aligned
\ &a_2^1 = x_6 \\
&a_2^2 = x_7 \\
&a_2^3 = x_1+2x_2-x_7
\endaligned
\ \qquad
\aligned
\ &a_2^4 = x_1+2x_2-x_3-x_4+x_5+x_6-x_7\\
&a_2^5 = x_1+2x_2-x_6\\
&a_2^6 = x_3+x_4-x_5-x_6+x_7
\endaligned
\eqnum{\mathbf h_2}
$$
\vskip0.3cm
$$
\aligned
\ &a_3^1 = -x_2+x_3+x_7 \\
&a_3^2 = -x_2+x_4+x_7\\
&a_3^3 = x_1+3x_2-x_3-x_7
\endaligned
\ \qquad
\aligned
\ &a_3^4 = x_1+3x_2-x_4-x_7 \\
&a_3^5 = x_1+x_5\\
&a_3^6 = 2x_2-x_5
\endaligned\qquad\quad\quad \
\eqnum{\mathbf h_3}
$$
\vskip0.3cm
$$
\aligned
\ &a_4^1 = x_1+x_3\\
&a_4^2 = 2x_2-x_3\\
&a_4^3 = 2x_2-x_3+x_6-x_7
\endaligned
\ \qquad
\aligned
\ &a_4^4 = 2x_2-x_3-x_4+x_5+x_6-x_7\\
&a_4^5 = x_1+x_3-x_6+x_7\\
&a_4^6 = x_1+x_3+x_4-x_5-x_6+x_7
\endaligned
\eqnum{\mathbf h_4}
$$
To select the solutions corresponding to
coboundaries, we must find out all the solutions
$\mathbf h_i = (a_i^1,...,a_i^6)$ \ $(1\le i\le 4)$ such that
the system of equations
$$
\mathbf h_i = T_{s_i}\mathbf h - \mathbf h, \qquad  1\le i\le 4,
$$
has a solution $\mathbf h = (u^1,...,u^6)\in A^{\oplus 6}$.
In coordinates, this system of equations looks as follows:
$$
\aligned
\ &a_1^1 = x_3 = u^2 - u^1;  \\
&a_1^2 = x_1+2x_2-x_3 = u^1 - u^2;  \\
&a_1^3 = x_4 = u^4 - u^3;
\endaligned
\ \ \
\aligned
\ & \ a_1^4 = x_1+2x_2-x_4 = u^3 - u^4;  \\
& \ a_1^5 = x_5 = u^6 - u^5; \\
& \ a_1^6 = x_1+2x_2-x_5 = u^5 - u^6;
\endaligned
\eqnum{\mathbf h_1'}
$$
\vskip0.3cm
$$
\aligned
\ &a_2^1 = x_6 = u^5 - u^1;\\
&a_2^2 = x_7 = u^3 - u^2;\\
&a_2^3 = x_1+2x_2-x_7 = u^2 - u^3;\\
&a_2^4 = x_1+2x_2-x_3-x_4+x_5+x_6-x_7 = u^6 - u^4;\\
&a_2^5 = x_1+2x_2-x_6 = u^1 - u^5;\\
&a_2^6 = x_3+x_4-x_5-x_6+x_7 = u^4 - u^6;
\endaligned\qquad\qquad\qquad\qquad\quad \ \
\eqnum{\mathbf h_2'}
$$
\vskip0.3cm
$$
\aligned
\ &a_3^1 = -x_2+x_3+x_7 = u^3 - u^1;\\
&a_3^2 = -x_2+x_4+x_7 = u^4 - u^2;\\
&a_3^3 = x_1+3x_2-x_3-x_7 = u^1 - u^3;
\endaligned
\ \ \
\aligned
\ &a_3^4 = x_1+3x_2-x_4-x_7 = u^2 - u^4;\\
&a_3^5 = x_1+x_5 = u^6 - u^5;\\
&a_3^6 = 2x_2-x_5 = u^5 - u^6;
\endaligned
\eqnum{\mathbf h_3'}
$$
\vskip0.3cm
$$
\aligned
\ &a_4^1 = x_1+x_3 = u^2 - u^1;\\
&a_4^2 = 2x_2-x_3 = u^1 - u^2;\\
&a_4^3 = 2x_2-x_3+x_6-x_7 = u^5 - u^3;\\
&a_4^4 = 2x_2-x_3-x_4+x_5+x_6-x_7 = u^6 - u^4;\\
&a_4^5 = x_1+x_3-x_6+x_7 = u^3 - u^5;\\
&a_4^6 = x_1+x_3+x_4-x_5-x_6+x_7 = u^4 - u^6.
\endaligned\qquad\qquad\qquad\qquad\qquad\qquad
\eqnum{\mathbf h_4'}
$$
It has a solution $(u^1,...,u^6)\in A^{\oplus 6}$ if and only if the
parameters $x_i$ satisfy the relations
%&
\begin{equation}\label{6.22}
x_1 = x_2 = 0.
%\eqnum{6.22}
\end{equation}
%&
Hence the group $\mathcal Z_{\psi_{5,6}}^{1}(B_5,A^{\oplus 6})$ of all
cocycles is isomorphic to the direct sum
$$
\mathcal Z = A_2\oplus A^{\oplus 6}
= \{\left.(x_1,...,x_7)\in A^{\oplus 7}\right|\, 2x_1=0\},
$$
and the group ${\mathcal B}_{\psi_{5,6}}^{1}(B_5,A^{\oplus 6})$ of all
coboundaries is isomorphic to the subgroup ${\mathcal B}\subset \mathcal Z$
$$
{\mathcal B} = A^{\oplus 5}
= \{\left.(x_1,...,x_7)\in A^{\oplus 7}\right|\, x_1=x_2=0\}.
$$
This means that
$$
H_{\psi_{5,6}}^1(B_5,A^{\oplus 6})\cong \mathcal Z/{\mathcal B}\cong
A_2\oplus A.
$$
Clearly, any cohomology class in $H_{\psi_{5,6}}^1(B_5,A^{\oplus 6})$
contains a cocycle of the form $(\mathbf h_1)-(\mathbf h_4)$ with
$x_3=x_4=x_5=x_6=x_7=0$; this proves (\ref{6.20})
(with $x=x_1\in A_2, \ y = x_2\in A$).
\end{proof}

\begin{Remark}\label{Rmk 6.3} In the next theorem we
use some details of the proof of Theorem \ref{Thm 6.4}. 
To simplify our notations,
we denote the system of equations
$(s_1\rightleftarrows s_3)-(s_3\infty s_4)$ by $(\mathcal S)$, and
formulae $(\mathbf h_1)-(\mathbf h_4)$ by $({\mathcal H})$.
We denote by $({\mathcal H}_0)$ the formulae for $a_i^j$ given by
$(\mathbf h_1)-(\mathbf h_4)$ {\it with the
particular value} $x_1 = 0$. Finally, we denote by $(\mathcal U_0)$
the system of equations $(\mathbf h_1')-(\mathbf h_4')$
with the same particular value $x_1 = 0$. It follows from the
proof of Theorem \ref{Thm 6.4} that {\sl system
$(\mathcal U_0)$ has a solution $(u^1,...,u^6)\in A^{\oplus 6}$
if and only if} $x_2=0$. \hfill $\bigcirc$
\end{Remark}
\vskip0.3cm

\begin{Theorem}\label{Thm 6.5} 
Let $t = n = 6$, and let
$\Omega=\nu_6\colon B_6\to{\mathbf S}(6)$ be
Artin's homomorphism. Then
$$
H^1_{\nu_6}(B_6, A^{\oplus 6})\cong A.
$$
Any cohomology class may be represented by a cocycle
that takes on the canonical generators $s_i\in B_6$
the values $\mathbf z(s_i) = \mathbf h_i$ of the form
%&
\begin{equation}\label{6.23}
\aligned
\ &\mathbf z(s_1) = \mathbf h_1 = (\ \ \ 0, \ \ 2y, \ \ \ 0,
\ \ \ 2y, \ \ \ 0, \ \ 2y)\\
&\mathbf z(s_2) = \mathbf h_2 = (\ \ \ 0, \ \ \ 0, \ \ \ 2y,
\ \ 2y, \ \ 2y, \ \ 0 \ )\\
&\mathbf z(s_3) = \mathbf h_3 = (\ -y, \  -y, \ \ 3y,
\ \ 3y, \ \ \ 0, \ \ 2y )\\
&\mathbf z(s_4) = \mathbf h_4 = ( \ \ \ 0, \ \ 2y, \ \ 2y,
\ \ 2y, \ \ \ 0, \ \ \ 0 \ )\\
&\mathbf z(s_5) = \mathbf h_5 = (-2y, \ \ 0, \ \ \ 2y,
\ \ 4y, \ \ \ 0, \ \ 2y )
\endaligned
%\eqnum{6.23}
\end{equation}
%&
\end{Theorem}

\begin{proof}
According to (\ref{4.6'}), we may regard the homomorphism
$\psi_{5,6}\colon B_5\to{\mathbf S}(6)$
as the restriction of Artin's homomorphism
$\nu_6\colon B_6\to{\mathbf S}(6)$ to
the subgroup $B\cong B_5$ in $B_6$
generated by the first four canonical generators
$s_1,s_2,s_3,s_4$.
Since $\varepsilon_{\nu_6}=\rho\circ\nu_6$ and
$\varepsilon_{\psi_{5,6}}=\rho\circ \psi_{5,6}$, we have
$\varepsilon_{\nu_6}| {B_5}=\varepsilon_{\psi_{5,6}}$,
and thus the $B_5$-action $T_{\psi_{5,6}}$ coincides with
the restriction of the $B_6$-action $T_{\nu_6}$
to $B=B_5$.
It follows that {\sl the restriction $\mathbf z_{B_5}$
of any $1$-cocycle $\mathbf z\in \mathcal Z_{\nu_6}^1(B_6, A^{\oplus 6})$
to the subgroup $B_5=B\subset B_6$
belongs to $\mathcal Z_{\psi_{5,6}}^1(B_5, A^{\oplus 6})$.}
Moreover, {\sl if such a cocycle $\mathbf z$ is a coboundary, then
its restriction $\mathbf z_{B_5}$ is also a coboundary.}
So, in order to compute $H_{\nu_6}^1(B_6, A^{\oplus 6})$,
we may use some computations already made in the proof of 
Theorem \ref{Thm 6.4}.
\vskip0.2cm

\noindent Let $\mathbf z\in Z_{\nu_6}^1(B_6,A^{\oplus 6}$
be a cocycle with the values
$\mathbf z(s_i) = \mathbf h_i = (a_i^1,...,a_i^6)\in A^{\oplus 6}$, \
$1\le i\le 5$. Then the elements $a_i^j$ with $1\le i\le 4$ must satisfy
the system of equations $(\mathcal S)$. According to the proof of 
Theorem \ref{Thm 6.4},
they are of the form $({\mathcal H})$ with some $x_1,...,x_7\in A$.
The elements $a_i^j$ \ ($1\le i\le 4$) together with the elements
$a_5^1,...,a_5^6$ must satisfy the system of equations
$(s_i\rightleftarrows s_5)$, $(s_4\infty s_5)$
corresponding to the relations $s_is_5 = s_5s_i$ \ $(1\le i\le 3)$
and $s_4s_5s_4=s_5s_4s_5$. Using the formula
$$
T_{s_5}(a^1,a^2,a^3,a^4,a^5,a^6) = (a^4,a^3,a^2,a^1,a^6,a^5)
$$
for the transformation $T_{s_5}=(T_{\nu_6})_{s_5}$,
we can write down equations
$(s_i\rightleftarrows s_5)$, $(s_4\infty s_5)$ explicitly:
$$
\aligned
\ & a_1^1 - a_1^4 = a_5^1 - a_5^2\,;
\qquad a_1^2 - a_1^3 = a_5^2 - a_5^1\,;
\qquad a_1^3 - a_1^2 = a_5^3 - a_5^4\,;\\
& a_1^4 - a_1^1 = a_5^4 - a_5^3;
\qquad a_1^5 - a_1^6 = a_5^5 - a_5^6\,;
\qquad a_1^6 - a_1^5 = a_5^6 - a_5^5\,;
\endaligned
\eqnum{s_1\rightleftarrows s_5}
$$

$$
\aligned
\ & a_2^1 - a_2^4 = a_5^1 - a_5^5\,;
\qquad a_2^2 - a_2^3 = a_5^2 - a_5^3\,;
\qquad a_2^3 - a_2^2 = a_5^3 - a_5^2\,;\\
& a_2^4 - a_2^1 = a_5^4 - a_5^6\,;
\qquad a_2^5 - a_2^6 = a_5^5 - a_5^1\,;
\qquad a_2^6 - a_2^5 = a_5^6 - a_5^4\,;
\endaligned
\eqnum{s_2\rightleftarrows s_5}
$$

$$
\aligned
\ & a_3^1 - a_3^4 = a_5^1 - a_5^3\,;
\qquad a_3^2 - a_3^3 = a_5^2 - a_5^4\,;
\qquad a_3^3 - a_3^2 = a_5^3 - a_5^1\,;\\
& a_3^4 - a_3^1 = a_5^4 - a_5^2\,;
\qquad a_3^5 - a_3^6 = a_5^5 - a_5^6\,;
\qquad a_3^6 - a_3^5 = a_5^6 - a_5^5\,;
\endaligned
\eqnum{s_3\rightleftarrows s_5}
$$

$$
\aligned
\ &a_4^1 - a_4^4 + a_4^3 = a_5^1 - a_5^2 + a_5^6\,; \ \
a_4^2 - a_4^3 + a_4^4 = a_5^2 - a_5^1 + a_5^5\,; \\
&a_4^3 - a_4^2 + a_4^6 = a_5^3 - a_5^5 + a_5^1\,;\ \
a_4^4 - a_4^1 + a_4^5 = a_5^4 - a_5^6 + a_5^2\,; \\
&a_4^5 - a_4^6 + a_4^2 = a_5^5 - a_5^3 + a_5^4\,; \ \
a_4^6 - a_4^5 + a_4^1 = a_5^6 - a_5^4 + a_5^3\,.
\endaligned
\eqnum{s_4\infty s_5}
$$
Let us denote the system of equations
$(s_i\rightleftarrows s_5)$ \ ($1\le i\le 3$) and \ $(s_4\infty s_5)$
by $({\mathcal S}_{new})$.
\vskip0.2cm

\noindent By substitution of the expressions 
$({\mathcal H})$ for $a_i^j$ \ ($1\le i\le 4$) 
in terms of the parameters $x_i$ into the equations $({\mathcal S}_{new})$, 
we obtain the system of {\it nonhomogeneous} equations
$(\widetilde {{\mathcal S}_{new}})$ for
the remainder elements $a_5^j\in A$, \ $1\le j\le 6$.
This (very unpleasant!) procedure leads to the following result:
\vskip0.3cm

\noindent {\bfit Claim.} {\sl System of equations
$(\widetilde {{\mathcal S}_{new}})$
has a solution $a_5^j\in A$ \ ($1\le j\le 6$)
if and only if $x_1=0$. Clearly, the latter condition
means that the elements $a_i^j$ \ ($1\le i\le 4$) must be
chosen according to formulae $({\mathcal H}_0)$. If this is done,
the solution $(a_5^1,...,a_5^6)\in A^{\oplus 6}$
of $(\widetilde {{\mathcal S}_{new}})$
is unique and reads as follows:}
$$
\aligned
\ &a_5^1 = -2x_2+x_3+x_4+x_7,\\
&a_5^2 = \ \ x_7,\\
&a_5^3 = \ \ 2x_2-x_7,
\endaligned
\qquad
\aligned
\ &a_5^4 = \ \ 4x_2-x_3-x_4-x_7,\\
&a_5^5 = \ \ \ x_5,\\
&a_5^6 = \ \ 2x_2-x_5.
\endaligned
\eqnum{\mathbf h_5}
$$
\vskip0.3cm

\noindent Combined with formulae $({\mathcal H}_0)$,
this shows that
$$
\mathcal Z = \mathcal Z_{\nu_6}^1(B_6, A^{\oplus 6})\cong
A^{\oplus 6} = \{(x_1,...,x_7)\in A^{\oplus 7}|\, x_1=0\}.
$$

\noindent Now we must select the solutions corresponding to coboundaries.
To this end, we add the following new equations $(\mathbf h_5')$
$$
\aligned
\ &a_5^1 = -2x_2+x_3+x_4+x_7 = u_4 - u_1,\\
&a_5^2 = x_7 = u_3 - u_2,\\
&a_5^3 = 2x_2-x_7 = u_2 - u_3,\\
\endaligned
%\
\aligned
\ &a_5^4 = 4x_2-x_3-x_4-x_7 = u_1  - u_4,\\
&a_5^5 = x_5 = u_6 - u_5,\\
&a_5^6 = 2x_2-x_5 = u_5 - u_6,
\endaligned
\eqnum{\mathbf h_5'}
$$
to the equations $(\mathcal U_0)$; then we need to find out when 
the resulting system of equations
$(\mathcal U_0)$, $(\mathbf h_5')$ has a solution
$(u^1,...,u^6)\in A^{\oplus 6}$.
We must certainly assume that $x_2=0$ (this
is necessary for solvability of the equations
$(\mathcal U_0)$; see Remark \ref{Rmk 6.3}).
A straightforward computation shows that in fact $x_2=0$ is 
{\sl the only} condition for solvability of the system of equations
$(\mathcal U_0)$, $(\mathbf h_5')$. That is,
{\sl a cocycle $\mathbf z$ of the form $({\mathcal H}_0)$,$(\mathbf h_5)$ 
is a coboundary if and only if} $x_2=0$.
Hence 
$$
{\mathcal B} = {\mathcal B}_{\nu_6}^1(B_6, A^{\oplus 6})\cong A^{\oplus 5}
= \{\left.(x_2,...,x_7)\in A^{\oplus 6}\right|\, x_2=0\},
$$ 
and thus
$$
H_{\nu_6}^1(B_6, A^{\oplus 6}) = \mathcal Z/{\mathcal B}\cong A.
$$
Since the parameters $(x_3,...,x_7)\in A^{\oplus 5}$ 
are completely free both in cocycles and coboundaries,
any cohomology class may be represented by a cocycle of the form
$({\mathcal H}_0),(\mathbf h_5)$ with $x_3=\ldots =x_7=0$, which
gives formulae (\ref{6.23}) ($y=x_2\in A$).
\end{proof}

\begin{Remark}\label{Rmk 6.4} 
The homomorphisms
$$
\Omega = \mu\colon B_n\to{\mathbf S}(n) \ \
\text{and} \ \ 
\Omega = \nu_6\colon B_6\to{\mathbf S}(6)
$$ 
are surjective, and the image of the homomorphism
$\Omega = \psi_{5,6}\colon B_5\to{\mathbf S}(6)$
is isomorphic to ${\mathbf S}(5)$. In any of these cases, the centralizer of
the image $\Img \Omega$ in the corresponding symmetric group
is trivial, and thus the two equivalence relations $\approx$
and $\sim$ on the set of all $\Omega$-homomorphisms coincide
(see Proposition \ref{Prp: Omega-homomorphisms and cocycles}$(b)$).
Hence, to compute all the 
$\Omega$-homomorphisms up to conjugation, it is sufficient
to choose one cocycle $\mathbf z$ in each cohomology class
(which indeed was done in Theorems \ref{Thm 6.3}-\ref{Thm 6.5}),
and then to compute the corresponding $\Omega$-homomorphisms
$\varPhi_{\mathbf z}$ defined by (\ref{5.36}). For the latter step, we
must also compute the homomorphism $\varepsilon = \rho\circ\Omega$;
this is not a problem at all, as far as the splitting $\rho$ and
the homomorphism $\Omega$ are given explicitly.
\hfill $\bigcirc$
\end{Remark}

\noindent In the following three corollaries we consider only
the case when $A={\mathbb Z}/2{\mathbb Z}$, which is important for some
applications 
(see Sec. \ref{sec: 7. Homomorphisms B'(k) to S(n), B'(k) to B(n) (n<k) and
B(k) to S(n) (n<2k)}). We skip the proofs since they
follow immediately from the results stated in 
Theorems \ref{Thm 6.3}--\ref{Thm 6.5}
(see also Remarks \ref{Rmk 6.1} and \ref{Rmk 6.4}). 

\begin{Corollary}\label{Crl 6.6} 
The cohomology group
$$
H^1_{\mu}(B_n,({\mathbb Z}/2{\mathbb Z})^{\oplus n})
\cong ({\mathbb Z}/2{\mathbb Z})\oplus ({\mathbb Z}/2{\mathbb Z})
$$ 
consists of the cohomology classes of the following four cocycles:
%&
\begin{equation}\label{6.24}
\aligned
\ &\mathbf z_0(s_i) = \mathbf 0\qquad {\text {\rm {(zero cocycle)}}},\\
%\\
&\mathbf z_1(s_i) = \mathbf e_{i+1},\\
%\\
&\mathbf z_2(s_i) = \mathbf e_1 + \cdots +
\mathbf e_{i-1} + \mathbf e_{i+2} + \cdots + \mathbf e_n,\\
%\\
&\mathbf z_3(s_i) = \mathbf z_1(s_i) + \mathbf z_2(s_i) =
\mathbf e_1 + \cdots + \mathbf e_{i-1} + \mathbf e_{i+1}
+ \cdots + \mathbf e_n\\
\endaligned
%&\eqnum{6.24}
\end{equation}
%&
$($here $1\le i\le n-1)$.

Any $\mu_n$-homomorphism
$\varPhi\colon B_n\to G\subset {\mathbf S}(2n)$
is $H$-conjugate to one of
the following four $\mu$-homomorphisms
$\varphi_j, \ j = 0,1,2,3$ \ $($in each formula $1\le i\le n-1)$:
$$
\aligned
\ &\varphi_0 = \rho\circ\mu_n 
= \varepsilon_{\mu_n}\sim\mu_n\times\mu_n:\qquad
\varphi_0(s_i) = (2i-1,2i+1)(2i,2i+2); \\
%&\\
&\varphi_1(s_i)=\underbrace{(2i-1,2i+2,2i,2i+1)}_{4\text{\rm -cycle}}\,;\\
%&\\
&\varphi_2(s_i)
=(1,2)\cdots (2i-3,2i-2)\,
\underbrace{(2i-1,2i+1)(2i,2i+2)}_{\text {\rm two transpositions}}\times\\
&\hskip2.8in\times(2i+3,2i+4)\cdots (2n-1,2n);\\
%&\\
&\varphi_3(s_i)
= (1,2)\cdots (2i-3,2i-2)\,
\underbrace{(2i-1,2i+2,2i,2i+1)}_{4\text{\rm -cycle}}\times\\
&\hskip2.8in\times(2i+3,2i+4)\cdots (2n-1,2n)\,.\\
\endaligned
$$
\end{Corollary}

\begin{Corollary}\label{Crl 6.7} 
The cohomology group
$$
H^1_{\psi_{5,6}}(B_5,({\mathbb Z}/2{\mathbb Z})^{\oplus 6})
\cong ({\mathbb Z}/2{\mathbb Z})\oplus ({\mathbb Z}/2{\mathbb Z})
$$
consists of the cohomology classes of the four cocycles
$$
\mathbf z_{(x,y)}\,, \ 
(x,y)\in ({\mathbb Z}/2{\mathbb Z})\oplus ({\mathbb Z}/2{\mathbb Z})\,,
$$
that take on the generators $s_1,s_2,s_3,s_4$ the values
$$
\aligned
\ &\mathbf z_{(x,y)}(s_1) 
= (0,\ x, \ \ \ 0, \qquad \ x, \ \ \ \ 0, \ x),\qquad
\mathbf z_{(x,y)}(s_2) = (0,\ 0, \ x, \ x, \ x, \ 0),\\
&\mathbf z_{(x,y)}(s_3) = (y, \ y, \ x+y, \ x+y, \ x, \ 0),\qquad
\mathbf z_{(x,y)}(s_4) = (x, \ 0, \ 0, \ 0, \ x, \ x).
\endaligned
$$
Any $\psi_{5,6}$-homomorphism 
$\varPhi\colon B_5\to G\subset {\mathbf S}(12)$
is $H$-conjugate to one of
the following four $\psi_{5,6}$-homomorphisms
$\eta_i=\phi_{5,6;(x,y)}$, \
$(x,y)\in ({\mathbb Z}/2{\mathbb Z})\oplus ({\mathbb Z}/2{\mathbb Z})$, \ $i=0,1,2,3$:
$$
\aligned
\ &\eta_0=\phi_{5,6;(0,0)} = \varepsilon_{\psi_{5,6}}\colon\\
%&\\
&\qquad s_1\mapsto (1,3)(2,4)(5,7)(6,8)(9,11)(10,12)\,, \\
&\qquad s_2\mapsto (1,9)(2,10)(3,5)(4,6)(7,11)(8,12)\,, \\
&\qquad s_3\mapsto (1,5)(2,6)(3,7)(4,8)(9,11)(10,12)\,, \\
&\qquad s_4\mapsto (1,3)(2,4)(5,9)(6,10)(7,11)(8,12)\,; \\
\endaligned
$$
%&\\
$$
\aligned
\ &\eta_1=\phi_{5,6;(1,0)}\colon\\
%&
&\qquad s_1\mapsto (1,4,2,3)(5,8,6,7)(9,12,10,11)\,, \\
&\qquad s_2\mapsto (1,10,2,9)(3,6,4,5)(7,11,8,12)\,, \\
&\qquad s_3\mapsto (1,6,2,5)(3,8,4,7)(9,11,10,12)\,, \\
&\qquad s_4\mapsto (1,3,2,4)(5,10,6,9)(7,12,8,11)\,; \\
\endaligned
$$
%&\\
$$
\aligned
\ &\eta_2=\phi_{5,6;(0,1)}\colon\\
%&\\
&\qquad s_1\mapsto (1,3)(2,4)(5,7)(6,8)(9,11)(10,12)\,, \\
&\qquad s_2\mapsto (1,9)(2,10)(3,5)(4,6)(7,11)(8,12)\,, \\
&\qquad s_3\mapsto (1,6)(2,5)(3,8)(4,7)(9,11)(10,12)\,, \\
&\qquad s_4\mapsto (1,3)(2,4)(5,9)(6,10)(7,11)(8,12)\,; \\
\endaligned
$$
%&\\
$$
\aligned
\ &\eta_3=\phi_{5,6;(1,1)}\colon\\
%&\\
&\qquad s_1\mapsto (1,4,2,3)(5,8,6,7)(9,12,10,11)\,, \\
&\qquad s_2\mapsto (1,10,2,9)(3,6,4,5)(7,11,8,12)\,, \\
&\qquad s_3\mapsto (1,5,2,6)(3,7,4,8)(9,11,10,12)\,, \\
&\qquad s_4\mapsto (1,3,2,4)(5,10,6,9)(7,12,8,11)\,.
\endaligned
$$
\end{Corollary}

\begin{Corollary}\label{Crl 6.8} 
The group $H^1_{\nu_6}(B_6,({\mathbb Z}/2{\mathbb Z})^{\oplus 6})
\cong{\mathbb Z}/2{\mathbb Z}$
consists of the cohomology classes of the two cocycles
$\mathbf z_y$ \ ($y\in{\mathbb Z}/2{\mathbb Z}$)
that take on the generators $s_1,...,s_5$ the values
$$
\aligned
\ &\mathbf z_y(s_1)=(0,0,0,0,0,0), \ \
\mathbf z_y(s_2)=(0,0,0,0,0,0), \ \
\mathbf z_y(s_3)=(y,y,y,y,0,0),\\
&\qquad\qquad\qquad \mathbf z_y(s_4)=(0,0,0,0,0,0), \ \ \
\mathbf z_y(s_5)=(0,0,0,0,0,0).\qquad
\endaligned
$$
Any $\nu_6$-homomorphism
$\varPhi\colon B_6\to G\subset {\mathbf S}(12)$
is $H$-conjugate to one of
the following two $\nu_6$-homomorphisms
$\phi_y$, \ $y\in{\mathbb Z}/2{\mathbb Z}$:
$$
\aligned
\ &\phi_0 = \rho\circ\nu_6 = \varepsilon_{\nu_6}\colon\\
&\qquad s_1\mapsto (1,3)(2,4)(5,7)(6,8)(9,11)(10,12)\,, \\
&\qquad s_2\mapsto (1,9)(2,10)(3,5)(4,6)(7,11)(8,12)\,, \\
&\qquad s_3\mapsto (1,5)(2,6)(3,7)(4,8)(9,11)(10,12)\,, \\
&\qquad s_4\mapsto (1,3)(2,4)(5,9)(6,10)(7,11)(8,12)\,,\\
&\qquad s_5\mapsto (1,7)(2,8)(3,5)(4,6)(9,11)(10,12)\,;
\endaligned
$$

$$
\aligned
\ &\phi_1\colon \\
&\qquad s_1\mapsto (1,3)(2,4)(5,7)(6,8)(9,11)(10,12)\,, \\
&\qquad s_2\mapsto (1,9)(2,10)(3,5)(4,6)(7,11)(8,12)\,, \\
&\qquad s_3\mapsto (1,6)(2,5)(3,8)(4,7)(9,11)(10,12)\,, \\
&\qquad s_4\mapsto (1,3)(2,4)(5,9)(6,10)(7,11)(8,12)\,, \\
&\qquad s_5\mapsto (1,7)(2,8)(3,5)(4,6)(9,11)(10,12)\,.
\endaligned
$$
\end{Corollary}

\begin{Remark}\label{Rmk 6.5} 
For any natural $m$ and any group $G$, we denote by
$\mathbf 1_m$ the trivial homomorphism $G\to{\mathbf S}(m)$.
$\mathbf {\{2\}}$ denotes the unique non-trivial homomorphism
$B_n\to{\mathbf S}(2)$; so for every $i=1,...,n-1$,
$\mathbf {\{2\}}(s_i)$ is the unique transposition in ${\mathbf S}(2)$.
Given a homomorphism $\varphi\colon B_n\to{\mathbf S}(N)$,
we regard the disjoint products $\varphi\times\mathbf 1_m$ and
$\varphi\times\mathbf {\{2\}}$
as homomorphisms of $B_n$ to the groups 
${\mathbf S}(N+m)$ and ${\mathbf S}(N+2)$,
respectively (see Definition \ref{Def: symmetric groups}{\bf(d)}). 
For instance, the homomorphism
$\mu_n\times\mathbf 1_1\colon B_n\to{\mathbf S}(n+1)$
is defined by $(\mu_n\times\mathbf 1_1)(s_i)=(i,i+1)$, \ $1\le i\le n-1$.
\hfill $\bigcirc$
\end{Remark}
\vskip0.3cm

\noindent We skip the (trivial) proof of the next corollary;
$\varphi_j\colon B_n\to G\subset {\mathbf S}(2n)$
are the $\mu_n$-homomorphisms exhibited in Corollary \ref{Crl 6.6}.
 %%&&

\begin{Corollary}\label{Crl 6.9} 
Any $(\mu_n\times\mathbf 1_1)$-homomorphism
$\varPhi\colon B_n\to G\subset {\mathbf S}(2(n+1))$
is conjugate to one of the eight homomorphisms
$\varphi_j\times\mathbf 1_2$, $\varphi_j\times\mathbf {\{2\}}$, \ $j=0,1,2,3$.
\hfill $\square$
\end{Corollary}

\begin{Remark}\label{Rmk 6.6} 
Take any $n\ge 3$ and any $r\ge 2$.
Remark \ref{Rmk 6.1} and Theorem \ref{Thm 6.3} 
give rise to some non-cyclic homomorphisms
$B_n\to{\mathbf S}(rn)$. To simplify the form of
the final result, we identify the group ${\mathbf S}(n)$ with the group
${\mathbf S}({\mathbb Z}/n{\mathbb Z})$, and regard the group 
${\mathbf S}(rn)$ as the symmetric group of the direct product 
$\mathcal D(r,n)
=({\mathbb Z}/r{\mathbb Z})\times ({\mathbb Z}/n{\mathbb Z})$
via the identification
%&
\begin{equation}\label{6.25}
\aligned
\ &{\boldsymbol\Delta}_{rn}\ni a
\mapsto (R(a),N(a))\in ({\mathbb Z}/r{\mathbb Z})
\times ({\mathbb Z}/n{\mathbb Z}),\\
&\text{where} \ R(a)=|a-1|_r\in{\mathbb Z}/r{\mathbb Z} \ \ \text{and} \ \
N(a)=|(a-1-R(a))/r|_n\in{\mathbb Z}/n{\mathbb Z}
\endaligned
%\eqnum{6.25}
\end{equation}
%&
(here $|c|_r$ and $|N|_n$ denote the $r$-
and the $n$-residue of $N\in\mathbb N$ respectively).
Then the subgroup $H\cong ({\mathbb Z}/r{\mathbb Z})^{\oplus n}$ 
generated by the $r$-cycles
%&
\begin{equation}\label{6.26}
C_m = ((m-1)r+1, (m-1)r+2,\ldots , mr)\in{\mathbf S}(rn),
\ \ 1\le m\le n,
%\eqnum{6.26}
\end{equation}
%&
acts on $\mathcal D(r,n)$ by translations of the first argument: 
%&
\begin{equation}\label{6.27}
\aligned
\ &H\ni \mathbf h = (C_1^{a^0}\cdots C_n^{a^{n-1}})\colon
\,\mathcal D(r,n)\ni (R,N)\mapsto (R+a^N,N)\in \mathcal D(r,n), \\
&{\text{where}} \ \ a^0,...,a^{n-1}\in{\mathbb Z}/r{\mathbb Z}.
\endaligned
%\eqnum{6.27}
\end{equation}
%&
The subgroup $G\subset {\mathbf S}(rn)$ (the centralizer of the
element $\mathcal C = C_1\cdots C_n$) was already identified with
the semidirect product 
$({\mathbb Z}/r{\mathbb Z})^{\oplus n}\leftthreetimes {_\tau} {\mathbf S}(n)
=({\mathbb Z}/r{\mathbb Z})^{\oplus n}\leftthreetimes {_\tau} 
{\mathbf S}({\mathbb Z}/n{\mathbb Z})$; the
latter group acts on the set $\mathcal D(r,n)$ by permutations
as follows:
%&
\begin{equation}\label{6.28}
(\mathbf h,s)(R,N)=(R+b^{s(N)},s(N))\,,
%\eqnum{6.28}
\end{equation}
%&
where $\mathbf h=(b^0,...,b^{n-1})\in ({\mathbb Z}/r{\mathbb Z})^{\oplus n}$, \
$s\in{\mathbf S}({\mathbb Z}/n{\mathbb Z})$.
Clearly, {\sl this action is transitive and imprimitive} 
(any subset $({\mathbb Z}/r{\mathbb Z})\times \{N\}$ in $\mathcal D(r,n)$
is a set of imprimitivity). 
\vskip0.2cm

According to Theorem \ref{Thm 6.3}, there are $r^2$
cocycles $\mathbf z_{(x,y)}$, \ $(x,y)\in ({\mathbb Z}/r{\mathbb Z})^{\oplus 2}$,
representing all the cohomology classes. The $\mu$-homomorphism
$$
\varphi_{(x,y)}\colon B_n\to
({\mathbb Z}/r{\mathbb Z})^{\oplus n}\leftthreetimes {_\tau} {\mathbf S}(n)
\subset {\mathbf S}(\mathcal D(r,n))
$$
corresponding to the cocycle $\mathbf z_{(x,y)}$ is uniquelly 
defined by its values 
$\widehat s_{(x,y);i} = \varphi_{(x,y)}(s_i)\in{\mathbf S}(\mathcal D(r,n))$
on the canonical generators $s_i\in B_n$. The permutations
$\widehat s_i = \widehat s_{(x,y);i}$ act on the elements
$(R,N)\in \mathcal D(r,n)$ as follows:
%&
\begin{equation}\label{6.29}
\widehat s_i(R,N) = \left \{
\aligned 
\ &(R+y,N) \qquad\qquad {\text{if}} \ \ N\ne i-1,i;\\ 
&(R,N+1) \qquad\qquad {\text{if}} \ \ N=i-1;\\
&(R+x,N-1) \qquad \ {\text{if}} \ \ N=i;
\endaligned
\right .\qquad (1\le i\le n-1).
%\eqnum{6.29}
\end{equation}
%&
It is easily shown that the map $\sigma_i\to s_{x,y;i}$, $i=1,...,n-1$,
extends to a uniquely defined non-cyclic homomorphism
$\varphi_{x,y}\colon B_n\to{\mathbf S}({\mathcal D}(r,n))={\mathbf S}(rn)$.
This homomorphism $\varphi_{x,y}$ is transitive if and only
if the elements $x,y\in{\mathbb Z}/r{\mathbb Z}$ generate the whole
group ${\mathbb Z}/r{\mathbb Z}$, or, which is the same,
if and only if $x$,$y$ and $r$ are mutually co-prime. 
However $\varphi_{x,y}$ can never be
primitive.

It is easy to show that the homomorphism $\varphi_{(x,y)}$
defined by $(\ref{6.29})$ is transitive if and only if the elements
$x,y\in{\mathbb Z}/r{\mathbb Z}$ generate the whole group
${\mathbb Z}/r{\mathbb Z}$, or, which is the same, 
if and only if $x$,$y$ and $r$ are mutually co-prime.
However, this homomorphism can never be
primitive, since its image is contained in the imprimitive
permutation group $G=({\mathbb Z}/r{\mathbb Z})^{\oplus n}
\leftthreetimes{_\tau} {\mathbf S}(n)$.
\vskip0.3cm

\noindent Using the same approach and 
Theorems \ref{Thm 6.4}, \ref{Thm 6.5}, one can construct
``exceptional" non-cyclic homomorphisms
$B_5\to{\mathbf S}(6r)$ and 
$B_6\to{\mathbf S}(6r)$, \ $r\ge 2$.
\hfill $\bigcirc$
\end{Remark}

\begin{Remark}\label{Rmk 6.7} 
Assume that the homomorphism
$\Omega$ in diagram (\ref{6.1}) is a disjoint product of two
homomorphisms:
$$
\Omega = \Omega'\times\Omega''\colon B_n
\to{\mathbf S}(t')\times{\mathbf S}(t'')
\subset {\mathbf S}(t),\qquad t' + t''=t\,.
$$
Then we have the decomposition
$H=A^{\oplus t}=A^{\oplus t'}\oplus A^{\oplus t''}$,
the actions $\tau'$, $\tau''$ of the groups
${\mathbf S}(t')$ and ${\mathbf S}(t'')$ on
$A^{\oplus t'}$ and $A^{\oplus t''}$ respectively,
and the corresponding semidirect products
$G'=A^{\oplus t'}\leftthreetimes_{\tau'} {\mathbf S}(t')\subset G$
and $G''=A^{\oplus t''}\leftthreetimes_{\tau''}
{\mathbf S}(t'')\subset G$. Any two elements $g'\in G'$,
$g''\in G''$ commute in $G$ and $\pi(g' g'')\in {\mathbf S}(t')
\times{\mathbf S}(t'')$. It is readily seen that the image of any
$\Omega$-homomorphism $\varPhi\colon B_n\to G$ is contained in the subgroup
$G'\cdot G''\cong G'\times G''$. Hence $\varPhi$ is the direct product of
the two homomorphisms
$$
\varPhi'\colon B_n\to G'\qquad \text{and}\qquad \varPhi''\colon B_n\to G''.
$$
Each of the latter homomorphisms fits in its own
commutative diagram of the form (\ref{6.1}) and may be studied separately.
\hfill $\bigcirc$
\end{Remark}

\noindent Let us compute $\Omega$-cohomology for 
a {\sl cyclic} homomorphism
$\Omega\colon B_n\to{\mathbf S}(t)$.
In this case there is a permutation $S\in{\mathbf S}(t)$
such that $\Omega(s_i)=S$ for all $i=1,...,n-1$.
In view of Remark \ref{Rmk 6.7}, it suffices to handle the
following two cases: \ $(i)$ $t=1$; \
$(ii)$ $t\ge 2$ and $S$ is a $t$-cycle.

\begin{Theorem}
\label{Thm: if Omega is cyclic then every Omega-homomorphism is cyclic} 
Let $n>4$ and let
$\Omega\colon B_n\to{\mathbf S}(t)$ be a cyclic homomorphism.

$a)$ Every $\Omega$-homomorphism $\varPhi\colon B_n\to G$ is cyclic.

$b)$ If $(i)$ $t=1$ or $(ii)$ $t\ge 2$ and $S$ is a $t$-cycle,
then $H_\Omega^1(B_n,A^{\oplus t})\cong A$.
In case $(ii)$ each cohomology class contains a unique cocycle
$\mathbf z$ of the form
%&
\begin{equation}\label{6.30}
\mathbf z(s_i) = (a,\underbrace{0,...,0}_{t-1 \ \text{times}}),\qquad
a\in A, \ \ 1\le i\le n-1.
\end{equation}
%&
\end{Theorem}

\begin{proof} 
$a)$ Since $\pi\circ\varPhi=\Omega$ is cyclic,
we have $\pi[\varPhi(B'_n)]=\{1\}$; hence
$\varPhi(B'_n)$ is contained in the abelian
group $\Ker\pi = H$. Therefore, $\varPhi(B'_n)=\{1\}$
and $\varPhi$ is cyclic.

$b)$ In case $(i)$, the $B_n$-action on
the group $H=A$ is trivial and hence
$$
H_\Omega^1(B_n,A)\cong\Hom(B_n,A)\cong A\,.
$$
Consider case $(ii)$. By Lemma \ref{Lm 6.1}, the elements
$\mathbf h_i=(a_i^1,...,a_i^t)\in A^{\oplus t}$ are
the values $\mathbf z(s_i)$ of a cocycle $\mathbf z$ 
if and only if they satisfy relations (\ref{6.13}), 
(\ref{6.13'}), which may be written as
%&
\begin{equation}\label{6.31}
\mathbf h_i - T_S\mathbf h_i = 
\mathbf h_j - T_S\mathbf h_j ,\qquad 1\le i,j\le n-1, \ \
|i-j|\ge 2,
\end{equation}
%&
and
%&
\begin{equation}\label{6.31'}
\mathbf h_i - T_S\mathbf h_i + T_{S^2}\mathbf h_i = 
\mathbf h_{i+1} - T_S\mathbf h_{i+1} + T_{S^2}\mathbf h_{i+1},
\qquad 1\le i\le n-2,
%\eqnum{6.31'}
\end{equation}
%&
respectively. Since $n>4$, system (\ref{6.31}) includes the equations
$$
\mathbf h_1 - T_S\mathbf h_1 = \mathbf h_j - T_S\mathbf h_j, \ \ 3\le j\le n-1\,,
\qquad
\mathbf h_2 - T_S\mathbf h_2 = \mathbf h_j - T_S\mathbf h_j, \ \ 4\le j\le n-1\,;
$$
%$$
%\aligned
%&\mathbf h_1 - T_S\mathbf h_1 = \mathbf h_j - T_S\mathbf h_j, \ \ 3\le j\le n-1\,, \\
%&\mathbf h_2 - T_S\mathbf h_2 = \mathbf h_j - T_S\mathbf h_j, \ \ 4\le j\le n-1\,;
%\endaligned
%$$
thereby,
$$
\mathbf h_1 - T_S\mathbf h_1 = \mathbf h_2 - T_S\mathbf h_2
=\ldots = \mathbf h_{n-1} - T_S\mathbf h_{n-1}\,.
$$
Combined with (\ref{6.31'}), this shows that
$T_{S^2}\mathbf h_1 = T_{S^2}\mathbf h_2 = \ldots = T_{S^2}\mathbf h_{n-1}$.
However, $S^2$ is just a permutation of coordinates, and thus
%&
\begin{equation}\label{6.32}
\mathbf h_1 = \mathbf h_2 = \ldots = \mathbf h_{n-1}.
%\eqnum{6.32}
\end{equation}
%&
In turn, (\ref{6.32}) implies both (\ref{6.31}) and 
(\ref{6.31'}), which shows that
%&
\begin{equation}\label{6.33}
\mathcal Z=\mathcal Z_{\Omega}^1(B_n,A^{\oplus t})\cong A^{\oplus t}.
%\eqnum{6.33}
\end{equation}
%&
Furthermore, a cocycle $\mathbf z$ with the values
$\mathbf h_1 = \mathbf h_2 = \ldots = \mathbf h_{n-1} = \mathbf v
= (v^1,...,v^t)\in A^{\oplus t}$
is a coboundary if and only if there exists
$\mathbf h = (u^1,...,u^t)\in A^{\oplus t}$ such that
%&
\begin{equation}\label{6.34}
\mathbf v = T_S\mathbf h - \mathbf h.
%\eqnum{6.34}
\end{equation}
%&
Since $S$ is a $t$-cycle, it is readily seen that for a given
$\mathbf v$ equation (\ref{6.34}) has a solution $\mathbf h$ if and only if
%&
\begin{equation}\label{6.35}
\sum_{i=1}^t v^i = 0.
%\eqnum{6.35}
\end{equation}
%&
Thus, ${\mathcal B}\subset \mathcal Z=A^{\oplus t}$
consists of all the elements
$\mathbf v=(v^1,...,v^t)\in A^{\oplus t}$ that
satisfy (\ref{6.35}). Joined with (\ref{6.33}), this completes the proof
of statement $(b)$.
\end{proof}
\vfill

\newpage

%Sec. 7
\section{Homomorphisms $B'_k\to{\mathbf S}(n)$, \
$B'_k\to B_n$ and $B_k\to{\mathbf S}(n)$, \ $n\le 2k$}
\label{sec: 7. Homomorphisms B'(k) to S(n), B'(k) to B(n) (n<k) and
B(k) to S(n) (n<2k)}

\noindent Here we prove Theorems A$(c)$, E, F, and G.
Our first goal is Theorem E$(a)$.

\subsection{Homomorphisms $B_k\to{\mathbf S}(k+1)$}
\label{Ss: 7.1. Homomorphisms  B(k) to S(k+1)}
We start with the following obvious property of retractions $\Omega$.

\begin{Lemma}\label{Lm 7.1} 
Let $\psi\colon B_k
\to{\mathbf S}(n)$ be a homomorphism, and let 
$\mathfrak C = \{C_1,...,C_t\}$ be the $r$-component of the permutation 
$\widehat\sigma_1$. Assume that either $t < k-2\ne 4$ or $t\le 2$.
Then the homomorphism 
$\Omega  = \Omega_{\mathfrak C}$ is cyclic, 
i. e., there exists a permutation
$g\in{\mathbf S}(\mathfrak C)\cong {\mathbf S}(t)$ such that
%&
\begin{equation}\label{7.1}
\widehat\sigma_{i+2}C_m\widehat\sigma_{i+2}^{-1} = g(C_m) = C_{g(m)}
\end{equation}
%&
whenever $1\le i\le k-3$ and $1\le m\le t$.
\end{Lemma}

\begin{proof} 
The case $t\le 2$ is trivial.
If $t<k-2\ne 4$, $\Omega$ is cyclic by
Theorem \ref{Thm: homomorphisms B(k) to S(n), n<k}$(a)$.
\end{proof}

The next lemma might be proven by a straightforward (but rather long)
computation. Instead, we use the cohomology approach in order
to show how it works in the simplest case.

\begin{Lemma}\label{Lm 7.2} 
Assume that $3\le k\ne 4$. Let
$\psi\colon B_k\to{\mathbf S}(n)$ be a non-cyclic
homomorphism such that $[\widehat\sigma_1] = [2,2]$.
Then $n\ge k+2$. Moreover,
\vskip0.3cm

$a)$ if $n<2k$, then the homomorphism $\psi$ is conjugate to the
homomorphism
$$
\phi^{(1)}_{k,n}\colon \sigma_i\mapsto (1,2)(i+2,i+3), \ \
1\le i\le k-1, \ \ \ 
\phi^{(1)}_{k,n}\sim\mathbf {\{2\}}\times\mu_k\times\mathbf 1_{n-k-2};
$$

$b)$ if $n\ge 2k$, then $\psi$ is either conjugate to 
$\phi^{(1)}_{k,n}$ or conjugate to the homomorphism
$$
\phi^{(2)}_{k,n}\colon \sigma_i\mapsto
(2i-1,2i+1)(2i,2i+2), \ \ 1\le i\le k-1, \ \ \ 
\phi^{(2)}_{k,n}\sim\mu_k\times\mu_k\times\mathbf 1_{n-2k};
$$

$c)$ in any case the homomorphism $\psi$ is intransitive.
\end{Lemma}

\begin{proof} 
For $k = 3$, all the assertions follow
from Lemma \ref{Lm: p cycles in S(2p)}. Let $k>4$ and let $\widehat\sigma_1 = C_1 C_2$,
where $C_1=(1,3)$, $C_2=(2,4)$
\footnote{Here we choose this normalization of $\psi$
instead of the usual $C_1=(1,2)$, $C_2=(3,4)$.};
so $\mathfrak C=\{C_1, C_2\}$ 
is the only non-degenerate component of $\widehat\sigma_1$, with
$\Sigma=\supp\,\mathfrak C=\{1,2,3,4\}$. The corresponding retraction
$\Omega\colon B_{k-2}\to{\mathbf S}(2)$ is cyclic;
hence either $\Omega$ is trivial or $\Omega=\mathbf {\{2\}}$.
\vskip0.2cm

Suppose first that $\Omega=\mathbf {\{2\}}$. Then 
Theorem
\ref{Thm: if Omega is cyclic then every Omega-homomorphism is cyclic}
shows that the cocycles 
$z_0(s_i)=(0,0)\in ({\mathbb Z}/2{\mathbb Z})^{\oplus 2}$ 
and $z_1(s_i)=(1,0)\in ({\mathbb Z}/2{\mathbb Z})^{\oplus 2}$ ($1\le i\le k-3$) 
represent all the cohomology classes. 
It follows from Lemma \ref{Lm 5.5}
that the $\Omega$-homomorphism
$\varphi_{_\Sigma}\colon B_{k-2}\to G\subset {\mathbf S}(4)$
(up to $H$-conjugation) coincides with the
homomorphism $\varepsilon=\rho\circ\Omega$, \
$\varepsilon(s_i)=(1,2)(3,4)$ for all $i$ (the second possibility, namely,
$\varphi_{_\Sigma}(s_i)=C_1(1,2)(3,4)=(1,2,3,4)$
for all $i\le k-3$, cannot occur here, 
since $[\widehat\sigma_{i+2}]=[2,2]$). This means that
$\widehat\sigma_{i+2}| \Sigma=(1,3)(2,4)$, and the condition
$[\widehat\sigma_{i+2}]=[2,2]$ implies that
$\widehat\sigma_{i+2}=(1,3)(2,4)$ for all $i\le k-3$. Hence $\psi$
is cyclic, contradicting our assumption.
\vskip0.2cm

Thus, $\Omega$ is trivial and any $\Omega$-homomorphism
is just a homomorphism of $B_{k-2}$ to $H$.
There are precisely four such homomorphisms $\phi_j$, $0\le j\le 3$,
which are defined as follows:
$$
\phi_0(s_i)=1;\ \ \phi_1(s_i)=C_1;\ \ \phi_2(s_i)=C_2;
\ \ \phi_3(s_i)=C_1 C_2\ \ \ (1\le i\le k-3).
$$
Therefore we may assume that the homomorphism $\varphi_{_\Sigma}$
coincides with one of the homomorphisms $\phi_j$, \ $j=0,1,2,3$.
\vskip0.2cm

If $\varphi_{_\Sigma}=\phi_0$, then $\widehat\sigma_{i+2}|\Sigma
=\varphi_{_\Sigma}(s_i)=\phi_0(s_i)=1$. Hence all the permutations
$\widehat\sigma_3,...,\widehat\sigma_{k-1}$
are disjoint with $\widehat\sigma_1$,
and we may assume that $\widehat\sigma_3=(5,7)(6,8)$.
The relations $\widehat\sigma_2\infty \widehat\sigma_1$, \
$\widehat\sigma_2\infty \widehat\sigma_3$ and Lemma \ref{Lm: p cycles in S(2p)} imply
that (up to a $\widehat\sigma_1$- and
$\widehat\sigma_3$-admissible conjugation)
$\widehat\sigma_2 = (3,5)(4,6)$.
Since $\supp \widehat\sigma_4\cap\{1,2,3,4\}=\varnothing$ and
$\widehat\sigma_4\widehat\sigma_2=\widehat\sigma_2\widehat\sigma_4$,
Lemma \ref{Lm: braid-like triple of 3 cycles} shows that
$\supp\,\widehat\sigma_4\cap\supp\,\widehat\sigma_2=\varnothing$;
in particular, $5,6\not\in\supp \widehat\sigma_4$.
Since $\widehat\sigma_4\infty \widehat\sigma_3$,
it follows from Lemma \ref{Lm: p cycles in S(2p)} that $\widehat\sigma_4 = (7,9)(8,10)$
(up to conjugation that is $\widehat\sigma_i$-admissible for all $i\le 3$).
By induction, we obtain that $n\ge 2k$ and $\psi\sim \phi^{(2)}_{k,n}$.
\vskip0.2cm

The homomorphisms $\phi_1$, $\phi_2$ are not $H$-conjugate;
however they are $G$-conjugate; thus, it is sufficient to handle 
the case $\varphi_{_\Sigma}=\phi_1$.
In this case $C_1\preccurlyeq \widehat\sigma_i$ for all $i\ne 2$;
hence $\widehat\sigma_i = C_1 D_i$, where every $D_i$ is a transposition
disjoint with $C_1$ and $C_2$. Since
$\widehat\sigma_2\widehat\sigma_4=\widehat\sigma_4\widehat\sigma_2$, we have
$\widehat\sigma_2 C_1 D_4\widehat\sigma_2^{-1} = C_1 D_4$.
The relation $\widehat\sigma_1\infty\widehat\sigma_2$ implies that
$\widehat\sigma_2 C_1\widehat\sigma_2^{-1} = C_1$.
(For otherwise, $\widehat\sigma_2 C_1\widehat\sigma_2^{-1} = D_4$, and
the supports of $\widehat\sigma_1$ and
$\widehat\sigma_2$ have exactly two common
symbols belonging to the transposition $C_1$; however, this contradicts 
Lemma \ref{Lm: p cycles in S(2p)}.)
Hence the set $\Sigma_1=\supp\,\widehat\sigma_1$
is $\widehat\sigma_2$-invariant.
The relations $\widehat\sigma_2\infty\widehat\sigma_1$ and
$\widehat\sigma_3| \Sigma_1 = C_1$ imply that
$\widehat\sigma_2| \Sigma_1 = C_1$.
Thus, $C_1\preccurlyeq \widehat\sigma_2$.
Taking into account that $k>4$, it is easy to see that
$n\ge k+2$ and $\psi\sim \phi^{(1)}_{k,n}$.
\vskip0.2cm

Finally, if $\varphi_{_\Sigma}=\phi_3$, we have
$\widehat\sigma_{i+2}|_\Sigma=\varphi_{_\Sigma}(s_i)=\phi_3(s_i)=(1,3)(2,4)$
for all $i\le k-3$. But $[\widehat\sigma_{i+2}]=[2,2]$; hence
$\widehat\sigma_{i+2}=(1,3)(2,4)$ for all $i\le k-3$ and $\psi$ is cyclic,
contradicting our assumption.
Thereby, statements $(a)$ and $(b)$ of the lemma are
proven. The statement $(c)$ is a trivial corollary
of $(a)$ and $(b)$.
\end{proof}

\begin{Remark}\label{Rmk 7.1} 
{\sl Let $\psi\colon B_4\to
{\mathbf S}(n)$ be a non-cyclic homomorphism such that
$\widehat\sigma_{1} = [2,2]$. Then,
besides the possibilities described in statements
$(a)$ and $(b)$ of Lemma \ref{Lm 7.2}, only the following
four cases may occur:}
$$
\aligned
4c) \ &n\ge 5 \ \text{and}\\
&\psi \sim \phi ^{(3)}_{4,n}\colon \ \
\left \{
\aligned
\sigma _{1},\sigma _{3} &\mapsto (1,2)(3,4),\\
%\\
\sigma _{2} \ \ \ &\mapsto (1,2)(4,5);
\endaligned
\right .
\endaligned
\quad \ \ 
\aligned
4d) \ &n\ge 6 \ \text{and}\\
&\psi \sim \phi ^{(4)}_{4,n}\colon \ \
\left \{
\aligned
\sigma _{1},\sigma _{3} &\mapsto (1,2)(3,4),\\
%\\
\sigma _{2} \ \ \ &\mapsto (2,5)(4,6);
\endaligned
\right .
\endaligned
$$
\vskip0.3cm
$$
\aligned
4e) \ &n\ge 6 \ \text{and}\\
&\psi \sim \phi ^{(5)}_{4,n}\colon \ \
\left \{
\aligned
\sigma _{1} &\mapsto (1,2)(3,4),\\
%\\
\sigma _{2} &\mapsto (2,5)(4,6),\\
%\\
\sigma _{3} &\mapsto (1,4)(2,3);
\endaligned
\right .
\endaligned
\qquad \ \
\aligned
4f) \ &n\ge 7 \ \text{and}\\
&\psi \sim \phi ^{(6)}_{4,n}\colon \ \
\left \{
\aligned
\sigma _{1} &\mapsto (1,2)(3,4),\\
%\\
\sigma _{2} &\mapsto (2,5)(4,6),\\
%\\
\sigma _{3} &\mapsto (1,2)(6,7).
\endaligned
\right .
\endaligned \ \ \ \
$$
{\sl All these homomorphisms except of the homomorphism
$\phi^{(5)}_{4,6}$ are intransitive.} 

Notice that the homomorphism $\phi^{(5)}_{4,6}$ 
coincides with the homomorphism $\psi^{(3)}_{4,6}$
defined in (\ref{4.4}) (see also Proposition \ref{Prp 4.5}). 
\hfill $\bigcirc$
\end{Remark}
\vskip0.3cm

\begin{Theorem}\label{Thm: homomorphisms B(k) to S(k+1)} 
\index{Theorem!on homomorphisms $B_k\to{\mathbf S}(k+1)$\hfill}
\index{Homomorphisms!$B_k\to{\mathbf S}(k+1)$\hfill} 
$a)$ Any transitive homomorphism
$\psi\colon B_k\to{\mathbf S}(k+1)$ is cyclic
whenever $k>5$.
\vskip0.2cm

$b)$ Let $k>5$ and $\psi\colon B_k\to{\mathbf S}(k+1)$
be a non-cyclic homomorphism. Then either $\psi$ is
conjugate to the homomorphisms
$\mu_{k+1}^k = \mu_k\times\mathbf 1_1\colon B_k\to
{\mathbf S}(k+1)$ \ or $k=6$ and $\psi$ is conjugate to the homomorphism
$\nu_{7}^6 = \nu_6\times\mathbf 1_1\colon 
B_6\to{\mathbf S}(7)$, where $nu_6$ is the Artin homomorphism.
\end{Theorem}

\begin{proof} 
$a)$ Suppose first that for some $k\ge 7$ there exists
a non-cyclic transitive homomorphism
$\psi\colon B_k\to{\mathbf S}(k+1)$.
It is very well known that there is a prime $p$ such that
$(k+1)/2<p\le k-2$ (see Remark \ref{Rmk: existence of desired primes});
by Lemma \ref{Lm: Artin Fixed Point Lemma},
the permutation $\widehat\sigma_1$
has at least $k-2$ fixed points. Therefore
$\# \supp \widehat\sigma_1 \le 3$; thus,
either $[\widehat\sigma_1]=[2]$ or $[\widehat\sigma_1]=[3]$.
However this contradicts Lemma
\ref{Lm: if [psi(sigma(1))]=[3] or |Fix(psi(sigma(1))|>n-4 then psi is intransitive}.
\vskip0.2cm

Assume now that there is a non-cyclic transitive
homomorphism $\psi\colon B_6\to{\mathbf S}(7)$.
Let $T=\widehat\alpha_{3,6}
=\widehat\sigma_3\widehat\sigma_4\widehat\sigma_5\in{\mathbf S}(7)$.
We apply Corollary
\ref{Crl: if psi is non-cyclic then psi(sigma(i)) infty psi(sigma(i+1))}
with $i=3$, $j=6$ (so, $j-i+1=4$)
and obtain that $4$ divides $\ord T$. Hence the cyclic
decomposition of $T$ contains precisely one $4$-cycle $C$.
The permutation $\widehat\sigma_1$ commutes with $T$
and 
Lemma \ref{Lm: invariant sets and components of commuting permutations}$(b)$
implies that $\widehat\sigma_1|\supp C = C^q$
for some integer $q$, \ $0\le q\le 3$. Let us consider all
these possibilities for $q$.
\vskip0.2cm

If $q=0$, then $\supp C\subseteq \Fix \widehat\sigma_1$
and $\#\supp\,\widehat\sigma_1\le 3$, which contradicts Lemma
\ref{Lm: if [psi(sigma(1))]=[3] or |Fix(psi(sigma(1))|>n-4 then psi is intransitive}.
\vskip0.2cm

If $q=1$ or $q=3$, then $[C^q]=[4]$ and $C^q\preccurlyeq \widehat\sigma_1$,
contradicting
Lemma \ref{Lm: Artin Lemma on cyclic decomposition of hatsigma(1)}$(a)$ 
(with $k=6$, $n=7$, $r=4>7/2=n/2$).
\vskip0.2cm

Finally, let $q=2$. Then $[C^q]=[C^2]=[2,2]$
and $C^2\preccurlyeq \widehat\sigma_1$.
If $\widehat\sigma_1\ne C^2$, then either $[\widehat\sigma_1]=[2,2,2]$
and $\widehat\sigma_1$ has a unique fixed point, or
$[\widehat\sigma_1]=[2,2,3]$ and $\widehat\sigma_1$ has a unique
invariant set of length $3$ (the support of the $3$-cycle); 
however this contradicts Lemma
\ref{Lm: if |Inv(r,psi(sigma(1)))|=1 then psi is transitive}. 
Hence $\widehat\sigma_1 = C^2$
and $[\widehat\sigma_1]=[2,2]$, contradicting Lemma \ref{Lm 7.2}$(c)$. 
\vskip0.2cm

$b)$ Since $\psi$ is non-cyclic and (by the statement $(a)$) intransitive,
Theorem \ref{Thm: homomorphisms B(k) to S(n), n<k}$(a)$
shows that the group $G=\Img \psi\subset {\mathbf S}(k+1)$
has exactly one orbit $Q$ of length $k$ and one fixed point. 
Hence $\psi$ is the composition of its reduction 
$\psi_Q\colon B_k\to{\mathbf S}(Q)\cong {\mathbf S}(k)$
and the natural embedding ${\mathbf S}(Q)\hookrightarrow {\mathbf S}(k+1)$.
Clearly, $\psi_Q$ is a non-cyclic transitive homomorphism, and
Artin Theorem shows that (up to conjugation of $\psi$)
either $\psi_Q=\mu_k$ or $k=6$ and $\psi_Q=\nu_6$. 
This gives the desired result.
\end{proof}

\subsection{Certain homomorphisms of the commutator subgroup
$B'_k$}
\label{Ss: 7.2. Certain homomorphisms of the commutator subgroup B'(k)}
Our next goal is Theorem A$(c)$. 
For any $k\ge 4$, we have the embedding
$\lambda_k'\colon B_{k-2}\to B'_k$
defined by $s_i\mapsto c_i=\sigma_{i+2}\sigma_1^{-1}$, 
\ $1\le i\le k-3$ 
(Remark \ref{Rmk: embeddings lambda(k,m) of B(k-2) to B(k)}). 
Recall also that the 
{\sl multiple commutator subgroups} $H^{(n)}$ of a group $H$
are defined by $H^{(0)}=H$ and $H^{(n)}=(H^{(n-1)})'$ for
$n\ge 1$.

\begin{Lemma}
\label{Lm: composition B(k-2) to Bk' to H} 
Let $k>4$ and let $\psi\colon B'_k\to H$ be a group homomorphism;
consider the composition
$\phi
=\psi\circ\lambda_k'\colon B_{k-2}\stackrel{\lambda_k'}{\longrightarrow} B'_k
\stackrel{\psi}{\longrightarrow} H$.
\vskip0.2cm

$a)$ $\Img\phi\subseteq\Img\psi\subseteq H^{(n)}$ for any $n\ge 0$.
\vskip0.2cm

$b)$ Assume that either
\ $(i)$ $\phi(s_1)=\phi(s_q)$ for some $q$ that satisfies $2\le q\le k-3$
or \ \ $(ii)$ $\phi(s_1^{-1})=\phi(s_3)$.
Then $\psi$ is trivial. In particular, if $\phi$ is cyclic,
then $\psi$ is trivial.
\end{Lemma}

\begin{proof} 
$a)$ Since $k>4$, the group $B'_k$ 
is perfect; hence $B'_k=(B'_k)^{(n)}$ 
for any $n\ge 0$. Therefore,
$\Img\phi\subseteq\Img\psi = \psi(B'_k)
= \psi((B'_k)^{(n)})\subseteq H^{(n)}$
for any $n\ge 0$. 

$b)$ Clearly, $\phi(s_i)=\psi(\lambda_k'(s_i))=\psi(c_i)$;
thus, $(i)$ means $\psi(c_1)=\psi(c_q)$ and $(ii)$ means
$\psi(c_1^{-1})=\psi(c_3)$. Hence, to prove the lemma, it is sufficient
to show that {\sl the system of relations (\ref{1.14})-(\ref{1.21})
joined with one of the relations \ $(i_q)$ $c_1=c_q$,
\ $(ii_3)$ $c_1^{-1}=c_3$ defines a presentation of the trivial group.}
This is a simple exercise, and we only sketch the proof. 

Assume that $(i_q)$ is fulfilled. Then (\ref{1.19}) 
implies $v c_1 v^{-1} = c_1 u^{-1}$; by (\ref{1.16}), 
this shows that $c_1^2=wu$. Now (\ref{1.14}) 
implies $uc_1^2u^{-1} = w^2$, 
which leads to $w^2=uv$ and $u=w$. Using (\ref{1.14}) once again, we obtain
$uc_1u^{-1} = w = u$; hence $c_1=u=w$. In view of
(\ref{1.15}), this means that $c_1=u=w=1$. Braid relations
(\ref{1.20}), (\ref{1.21}) 
and the relation $c_1=1$ imply that $c_i=1$ for all $i$.
Finally, (\ref{1.18}) shows that $u=1$.

Assume  that $c_1^{-1}=c_3$. Then (\ref{1.19}) 
implies $vc_1v^{-1} = uc_1$; 
taking into account (\ref{1.16}), we obtain $uc_1 = c_1^{-1}w$. 
Then (\ref{1.14}) leads to $uwu^{-1} = wuw$. Using this and (\ref{1.15}), 
we obtain $c_1 = u^{-1}w$, and (\ref{1.14}) shows that $u=1$ and $w=c_1$. 
These relations and (\ref{1.15}) imply $c_1=w=1$,
and relations (\ref{1.20}), (\ref{1.21}) show that $c_i=1$ for all $i$.
Finally, from (\ref{1.18}) we obtain $v=1$.
\end{proof}

\begin{Lemma}\label{Lm 7.5} 
Consider a homomorphism
$\psi\colon B'_6\to{\mathbf S}(5)$
and assume that the composition
$\phi=\psi\circ\lambda_6'\colon B_4
\hookrightarrow B'_6\to{\mathbf S}(5)$ is intransitive. 
Then $\psi$ is trivial.
\end{Lemma}

\begin{proof} 
In view of Lemma \ref{Lm: composition B(k-2) to Bk' to H}$(b)$,
it is sufficient to show
that $\phi(s_1)=\phi(s_3)$. Set $G=\Img \phi\subseteq {\mathbf S}(5)$
and consider all $G$-orbits $\Sigma\subseteq {\boldsymbol\Delta}_5$ and all 
the reductions 
$$
\phi_\Sigma\colon B_4\to{\mathbf S}(\Sigma)
$$ 
of $\phi$ to these orbits. By our assumption, $\# \Sigma\le 4$
for any $G$-orbit $\Sigma$. If all the reductions are cyclic,
we are done. Assume that there is a $G$-orbit $\Sigma$ with
the non-cyclic reduction $\phi_\Sigma$; then $\Sigma$ is the
only orbit with this property and $\# \Sigma\ge 3$.
If $\# \Sigma=3$, then Theorem \ref{Thm 3.14}$(a)$ implies that 
$\phi_\Sigma(s_1)=\phi_\Sigma(s_3)$; in fact, we have
$\phi(s_1)=\phi(s_3)$ (since the reduction to any other
$G$-orbit is cyclic). Finally, assume that $\# \Sigma=4$.
By Lemma \ref{Lm: composition B(k-2) to Bk' to H}$(a)$, 
$G=\Img \phi\subseteq{\mathbf S}'(5)=\mathbf A(5)$;
hence $G$ contains only even permutations. 
The set ${\boldsymbol\Delta}_5-\Sigma$ consists of a single
point that is a fixed point of $G$. It follows
that {\sl the image of the non-cyclic transitive homomorphism
$\phi_\Sigma\colon B_4\to{\mathbf S}(\Sigma)
\cong {\mathbf S}(4)$ contains only even permutations.}
This property and the sentence $(c)$ of Artin Theorem imply that
$\phi_\Sigma$ is conjugate to the homomorphism
$\nu_{4,3}$; thus, $\phi_\Sigma(s_1)=\phi_\Sigma(s_3)$ 
and $\phi(s_1)=\phi(s_3)$.
\end{proof}

\begin{Theorem}
\label{Thm: Bk' has no homomorphisms to S(n) and Bn when k>4 and n<k} 
If $k>4$ and $n<k$, then the group
$B'_k$ does not possess non-trivial homomorphisms
to the groups ${\mathbf S}(n)$ and $B_n$.
\end{Theorem}

\begin{proof} 
Consider first a homomorphism 
$\psi\colon B'_k\to{\mathbf S}(n)$.
By Lemma \ref{Lm: composition B(k-2) to Bk' to H}$(a)$, we have 
$\Img \psi\subseteq {\mathbf S}'(n) = \mathbf A(n)$;
in particular, the homomorphism $\psi$ cannot be surjective.
The case $n<5$ is trivial, since for such $n$
the alternating group $\mathbf A(n)$ is solvable and the group 
$B'_k$ is perfect. So we may assume that
$4<n<k$ and $k>5$.  

Consider the composition
$\phi=\psi\circ\lambda_k'\colon B_{k-2}\stackrel{\lambda_k'}{\longrightarrow}
B'_k\stackrel{\psi}{\longrightarrow} {\mathbf S}(n)$.
By Lemma \ref{Lm: composition B(k-2) to Bk' to H}$(b)$,
it is sufficient to prove that $\phi$
is cyclic or at least satisfies the condition $\phi(s_1)=\phi(s_3)$.
If $k>6$ and $n<k-2$, then $\phi$ is cyclic
by Theorem \ref{Thm: homomorphisms B(k) to S(n), n<k}$(a)$.
Hence we must only consider the following
three cases: \ $i)$ $k=6$ and $n=5$; \ $ii)$ $k>6$ and $n=k-2$; \
$iii)$ $k>6$ and $n=k-1$.

$i)$ In this case we deal with the homomorphism
$\phi\colon B_4\to{\mathbf S}(5)$.
If $\phi$ is intransitive, then the conclusion follows by Lemma \ref{Lm 7.5}.
If $\phi$ is transitive, then Lemma \ref{Lm 4.2} 
shows that $\phi(s_1)=\phi(s_3)$.
\vskip0.3cm

$ii)$ In this case we deal with the homomorphism
$\phi\colon B_{k-2}\to{\mathbf S}(k-2)$.
As we noted above, the homomorphism $\psi$ cannot be surjective;
hence $\phi$ is non-surjective and
Lemma \ref{Lm: non-surjective homomorphism B(k) to S(k) is cyclic} 
implies that $\phi$ is cyclic.
\vskip0.3cm

$iii)$ In this case we deal with the homomorphism
$\phi\colon B_{k-2}\to{\mathbf S}(k-1)$.
We shall consider the following two cases:
\ $iii_1)$ the homomorphism $\phi$ is intransitive;
\ $iii_2)$ the homomorphism $\phi$ is transitive.
\vskip0.3cm
$iii_1)$ In this case we may also assume that the image
$G=\Img \phi\subset {\mathbf S}(k-1)$ has at least one $G$-orbit
$\Sigma\subset {\boldsymbol\Delta}_{k-1}$ such that the reduction
$\phi_\Sigma$ is non-cyclic.
Since $\phi$ is intransitive, it follows from
Theorem \ref{Thm: homomorphisms B(k) to S(n), n<k}$(a)$
that $\#\Sigma=k-2$; certainly, $\Sigma$ is the only orbit of such length.
By Lemma \ref{Lm: composition B(k-2) to Bk' to H}
$(a)$, $G=\Img \phi\subseteq{\mathbf S}'(k-1)=\mathbf A(k-1)$;
hence $G$ contains only even permutations.
The set ${\boldsymbol\Delta}_{k-1}-\Sigma$ consists of a single
point which is a fixed point of $G$. This implies
that the image of the non-cyclic homomorphism
$\phi_\Sigma\colon B_{k-2}\to{\mathbf S}(\Sigma)
\cong {\mathbf S}(k-2)$ contains only even permutations.
However this contradicts
Lemma \ref{Lm: non-surjective homomorphism B(k) to S(k) is cyclic}.
\vskip0.3cm

$iii_2)$ If $k>7$, then $k-2>5$ and $\phi$ is cyclic
by Theorem \ref{Thm: homomorphisms B(k) to S(k+1)}$(a)$.
Finally, if $k=7$, then $n=k-1=6$
and we deal with the transitive non-cyclic homomorphism
$\phi\colon B_5\to{\mathbf S}(6)$.
By Proposition \ref{Prp 4.9}, $\phi$ must be conjugate to the
homomorphism $\psi_{5,6}$. However, this is impossible, since
the $\psi_{5,6}$ is surjective and $\Img \phi\subseteq \mathbf A(6)$.
This concludes the proof for homomorphisms
$B'_k\to{\mathbf S}(n)$.
\vskip0.3cm

Consider now a homomorphism
$\varphi\colon B'_k\to B_n$.
As we have already proved, the composition
$\psi=\mu\circ\varphi\colon B'_k\stackrel{\varphi}{\longrightarrow} B_n
\stackrel{\mu}{\longrightarrow} {\mathbf S}(n)$
of the homomorphism $\varphi$ with the canonical epimorphism $\mu$
must be trivial. Therefore,
$\varphi(B'_k)\subseteq \Ker\, \mu = {PB_n}$.
By Corollary \ref{Crl: Perfect group has no homomorphisms to PB(k)}, 
the perfect group $B'_k$ does not
possess non-trivial homomorphisms to the pure braid group ${PB_n}$;
hence the homomorphism $\varphi$ is trivial.
\end{proof}

\begin{Remark}\label{Rmk 7.2} 
The groups $B'_3$
and $B'_4$ have many non-trivial homomorphisms
to any (non-trivial) group. Moreover, for any $k\ge 3$ and
any $n\ge k$ there exist non-trivial homomorphisms
$B'_k\to{\mathbf S}(n)$ and
$B'_k\to B_n$. This shows
that the conditions $k>4$ and $n<k$ in
Theorem \ref{Thm: Bk' has no homomorphisms to S(n) and Bn when k>4 and n<k}
are, in a sense, sharp. 
Theorem \ref{Thm: Bk' has no homomorphisms to S(n) and Bn when k>4 and n<k}
implies
Theorem \ref{Thm: homomorphisms B(k) to S(n), n<k} for $k>4$.
However, we could not skip the proof of
Theorem \ref{Thm: homomorphisms B(k) to S(n), n<k}
since the latter theorem was used essentially (and many times)
in the proof of
Theorem \ref{Thm: Bk' has no homomorphisms to S(n) and Bn when k>4 and n<k}. 
Note that Theorem \ref{Thm: homomorphisms B(k) to S(n), n<k} would hardly help
to prove very useful Proposition
\ref{Prp: transitive imprimitive homomorphisms Bk to S(n) is cyclic if n<2k}, 
while 
Theorem \ref{Thm: Bk' has no homomorphisms to S(n) and Bn when k>4 and n<k}
works perfectly.
\end{Remark}
\vskip0.3cm

Now we are almost ready to prove Theorem G (see Proposition
\ref{Prp: transitive imprimitive homomorphisms Bk to S(n) is cyclic if n<2k}
below).

\begin{Lemma}
\label{Lm: when transitive group homomorphism B to S(n) is primitive or not}  
Let $n<2k$ and let
$\psi\colon B\to{\mathbf S}(n)$ be a transitive group homomorphism.

$a)$ Assume that for every $m<k$ the commutator subgroup
$B'$ of $B$ does not admit non-trivial homomorphisms to ${\mathbf S}(m)$.
If $\psi$ is imprimitive, then it is abelian.

$b)$ Assume that for every $m<k$ the group $B$ itself does not
admit non-trivial homomorphisms to ${\mathbf S}(m)$.
Then the homomorphism $\psi$ is primitive.
\end{Lemma}

\begin{proof} 
The statement $(a)$ {\sl presumes} that $\psi$ is
imprimitive; to treat $(b)$ simultaneously with $(a)$,
we {\sl assume} that $\psi$ is imprimitive and show 
that this leads to a contradiction.

Let ${\boldsymbol\Delta}_n = Q_1\cup\cdots\cup Q_t$, \ $t\ge 2$,
be some decomposition of ${\boldsymbol\Delta}_n$ into 
imprimitivity sets of the group $G = \Img \psi\subset {\mathbf S}(n)$.
Since $\psi$ is transitive, $\# Q_1=\cdots = \# Q_t = r$,
where $r\ge 2$ and $rt=n$. Clearly,
%&
\begin{equation}\label{7.2}
2\le t<k\qquad \text{and}\qquad 2\le r<k.
\end{equation}
%&
Consider the normal subgroup $H\vartriangleleft G$ consisting
of all elements $h\in G$ such that every set $Q_i$ is
$h$-invariant. Thus, we have the exact sequence
%&
\begin{equation}\label{7.3}
1\to H\to G\stackrel{\pi}{\longrightarrow}\widetilde G\to 1,
\end{equation}
%&
where $\pi$ is the natural projection onto the quotient group
$\widetilde G=G/H$. This quotient group, in turn, possesses 
the natural embedding
$\widetilde G\hookrightarrow {\mathbf S}(\{Q_1,...,Q_t\})\cong{\mathbf S}(t)$. 
Consider the composition
$$
\widetilde\psi=\pi\circ\psi\colon B\stackrel{\psi}{\longrightarrow} G
\stackrel{\pi}{\longrightarrow} {\widetilde G},
\qquad {\widetilde G}\subseteq {\mathbf S}(t).
$$
Clearly, $\widetilde\psi$ is surjective.

Under the assumptions made in the statement $(b)$,
the homomorphism $\widetilde\psi$ must be trivial;
this means that the quotient group $\widetilde G$
is trivial and hence $H=G$. It follows that every
set $Q_i$ is $G$-invariant. However this contradicts
the transitivity of the homomorphism $\psi$.

Under the assumptions made in the statement $(a)$, 
the restriction of $\widetilde\psi$ to the
commutator subgroup $B'$ must be trivial; hence,
$\widetilde\psi\colon B\to \widetilde G$ is
a surjective abelian homomorphism, and the group
$\widetilde G$ is abelian. Thereby, the exact sequence (\ref{7.3}) shows that
the commutator subgroup $G'$ of $G$ is contained in $H$.
In particular, we have
%&
\begin{equation}\label{7.4}
\psi(B')\subseteq G'\subseteq H.
\end{equation}
%&
Every $Q_i$ is $H$-invariant, and (\ref{7.4}) 
shows that the restriction 
$\psi'=\psi| B'\to G'\subseteq H$
may be regarded as the disjoint product of the reductions
$$
\psi'_{Q_i}\colon B'\to{\mathbf S}(Q_i)\cong {\mathbf S}(r),
\ \ 1\le i\le t, \qquad
\psi'_{Q_i}(b)=(\psi'(b))| Q_i
\ \ \text{for all} \ \ b\in B'.
$$
In view of (\ref{7.2}), each homomorphism $\psi'_{Q_i}$ is trivial.
Hence the homomorphism $\psi'=\psi| B'$ is
trivial, and our original homomorphism $\psi$ is abelian.
\end{proof}

\begin{Proposition}
\label{Prp: transitive imprimitive homomorphisms Bk to S(n) is cyclic if n<2k}
Let $k>4$ and $n<2k$. Then

$a)$ Any transitive imprimitive homomorphism
$\psi\colon B_k\to{\mathbf S}(n)$ is cyclic.

$b)$ Any transitive homomorphism
$\psi'\colon B'_k\to{\mathbf S}(n)$ is
primitive.
\end{Proposition}

\begin{proof} 
Both $(a)$ and $(b)$ are immediate consequences of  
Theorem \ref{Thm: Bk' has no homomorphisms to S(n) and Bn when k>4 and n<k}
and the corresponding statements of
Lemma \ref{Lm: when transitive group homomorphism B to S(n) is primitive or not}
(with $B=B_k$ when proving $(a)$ and $B=B'_k$ when proving $(b)$).
\end{proof}

\subsection{Homomorphisms $B_k\to{\mathbf S}(k+2)$}
\label{Ss: 7.3. Homomorphisms  B(k) to S(k+2)}
Here we prove Theorem E$(b)$
(see Theorem \ref{Thm: homomorphisms B(k) to S(k+2)} below). 
To this end, we need some preparation.
In the following lemma, $\varphi_1$ is the homomorphism
exhibited in Corollary \ref{Crl 6.6} (with $n=k$).

\begin{Lemma}\label{Lm: homomorphisms Bk to S(n) for n<=2k and
supp(hatsigma(1))<6}   
Let $6<k<n\le 2k$ and let
$\psi\colon B_k\to{\mathbf S}(n)$ be a transitive
non-cyclic homomorphism. If \ $\# \supp \widehat\sigma_1<6$, then
$n=2k$, \ $\widehat\sigma_1$ is a $4$-cycle, and the homomorphism
$\psi$ is conjugate to the homomorphism $\varphi_1$.
\end{Lemma}
 %%&&

\begin{proof} 
A priori, we have the following $6$ possibilities
for the cyclic type of $\widehat\sigma_1$:
$$
[\widehat\sigma_1] = [2]; \ [2,2]; \ [2,3]; \
[3]; \ [4]; \ [5].
$$ 
The types $[2]$, $[3]$, and $[2,2]$ are forbidden by Lemma
\ref{Lm: if [psi(sigma(1))]=[3] or |Fix(psi(sigma(1))|>n-4 then psi is intransitive}
and Lemma \ref{Lm 7.2}$(c)$ respectively.
Let us note that 
$E(n/E(k/2))\le E(2k/E(k/2))\le 4$ for any $k\ge 6$;
hence the inequality of
Lemma \ref{Lm: condition no cycles of the same length restricts support}
eliminates the types $[2,3]$ and $[5]$. 

So, we are left with the type $[\widehat\sigma_1]=[4]$; in this case
$[\widehat\sigma_i]=[4]$ for all $i$. Put $\Sigma_i=\supp \widehat\sigma_i$.
The first statement of
Lemma \ref{Lm: condition no cycles of the same length restricts support}
says that $\Sigma_i\cap\Sigma_j=\varnothing$ for $|i-j|\ge 2$;
in particular, $\Sigma_i\cap\Sigma_{i+2}=\varnothing$ for $1\le i\le k-3$.
Since $\psi$ is non-cyclic, we have
$\widehat\sigma_{i+1}\infty\widehat\sigma_i$, \
$\widehat\sigma_{i+1}\infty\widehat\sigma_{i+2}$,
and all three permutations are $4$-cycles. The set $\Sigma_{i+1}$
cannot coincide with one of the sets $\Sigma_i$, $\Sigma_{i+2}$ 
(for otherwise, $\widehat\sigma_{i+1}$ would be disjoint with
one of the permutations $\widehat\sigma_i$, $\widehat\sigma_{i+2}$).
Hence Lemma \ref{Lm: braid-like couples of 4-cycles}
implies that $\# (\Sigma_i\cap\Sigma_{i+1})=2$
and $\# (\Sigma_i\cap\Sigma_{i+2})=2$.
It follows immediately that for any $m\le k-1$ the union
$\Sigma_1\cup\Sigma_2\cup\ldots\cup\Sigma_m\subseteq {\boldsymbol\Delta}_n$ 
consists of $4+2(m-1)$ points. In particular, for $m=k-1$ 
this union consists of $4+2(k-2)=2k\ge n$ points; so, $n=2k$.
Without loss of generality, we may assume that
$\widehat\sigma_1=(1,4,2,3)$.
It follows from Lemma \ref{Lm: braid-like couples of 4-cycles}
and from what has been proved above that %RRRRRRRRRRRRRRRR
we may assume that $\widehat\sigma_2=(3,6,4,5)$
(any of other possibilities can be reduced to this case by a
$\widehat\sigma_1$-admissible conjugation of $\psi$).
Taking into account the above arguments and 
the property $\Sigma_1\cap\Sigma_3=\varnothing$,
we obtain that (up to an admissible conjugation)
$\widehat\sigma_3=(5,8,6,7)$, and so on. Hence $\psi\sim\varphi_1$.
\end{proof}

\begin{Lemma}
\label{Lm1: homomorphism with all components of hatsigma(1) of length<=k-3} 
Let $k>6$ and let $\psi\colon B_k\to{\mathbf S}(n)$ be
a homomorphism such that all components
of $\widehat\sigma_1$ $($including the degenerate 
component $\Fix \widehat\sigma_1)$
are of lengths at most $k-3$. Then $\psi$ is cyclic.
\end{Lemma}

\begin{proof} 
If $\widehat\sigma_1=\id$, then $\psi$ is trivial.
So, we may assume that for some $r\ge 2$ the 
permutation $\widehat\sigma_1$ has the $r$-component $\mathfrak C$
of some length $t\ge 1$. 
Put $\Sigma_{\mathfrak C} = \supp \mathfrak C$ and consider the retraction
$\Omega_{\mathfrak C}\colon B_{k-2}\to{\mathbf S}(\mathfrak C)
\cong {\mathbf S}(t)$ of $\psi$ to $\mathfrak C$
(see Sec. \ref{sec: Retraction of homomorphisms; homomorphisms and cohomology}). 
According to our assumptions, we have $k-2>4$ and $t<k-2$;
by Lemma \ref{Lm 7.1}, $\Omega_{\mathfrak C}$ is cyclic.
It follows from
Theorem
\ref{Thm: if Omega is cyclic then every Omega-homomorphism is cyclic}$(a)$
that the $\Omega_{\mathfrak C}$-homomorphism 
$\varphi_{{\Sigma_{\mathfrak C}}}\colon B_{k-2}
\to G_{\mathfrak C}\subset {\mathbf S}(rt)$
is also cyclic. This means that
$\varphi_{{\Sigma_{\mathfrak C}}}(s_1)=\ldots
=\varphi_{{\Sigma_{\mathfrak C}}}(s_{k-3})$,
and thus
%&
\begin{equation}\label{7.5}
\widehat\sigma_3| {\Sigma_{\mathfrak C}}
=\ldots =\widehat\sigma_{k-1}| {\Sigma_{\mathfrak C}}.
\end{equation}
%&
Put $\Sigma=\bigcup_{\mathfrak C} \Sigma_{\mathfrak C}$,
where $\mathfrak C$ runs over all the non-degenerate components of
$\widehat\sigma_1$; clearly, $\Sigma=\supp\,\widehat\sigma_1$. The sets
$\Sigma$ and $\Sigma' = \Fix \widehat\sigma_1={\boldsymbol\Delta}_n-\Sigma$
are invariant under all the permutations 
$\widehat\sigma_3,...,\widehat\sigma_{k-1}$. Since (\ref{7.5}) holds for
{\sl every} non-degenerate component $\mathfrak C$ of $\widehat\sigma_1$,
it follows that there is a permutation $S\in{\mathbf S}(\Sigma)$ such that
\ \ $\widehat\sigma_3|\Sigma=\ldots
=\widehat\sigma_{k-1}|\Sigma = S$. \
By our assumption, the degenerate component 
$\Sigma' = \Fix \widehat\sigma_1$ contains at most $k-3$ points:
%&
\begin{equation}\label{7.6}
\#\Sigma' = \# \Fix \widehat\sigma_1 \le k-3.
\end{equation}
%&
Set $S_i' = \widehat\sigma_i|\Sigma'$, \ $i=3,...,k-1$.
Clearly, 
%&
\begin{equation}\label{7.7}
\widehat\sigma_i=S\cdot S_i' \qquad {\text{for all}} \ \ i=3,...,k-1.
\end{equation}
%&
For any $i=3,...,k-1$, we have $\supp S\cap\supp S_i'=\varnothing$;
hence it follows from (\ref{7.7}) that the permutation 
$S_3',...,S_{k-1}'$ satisfy the standard braid relations
$S_i' S_j'=S_j' S_i'$ for
$|i-j|>1$ and $S_i' S_{i+1}' S_i'
= S_{i+1}' S_i' S_{i+1}'$ for
$3\le i<k-1$. This means that we can define a group homomorphism
$\phi\colon B_{k-2}\to{\mathbf S}(\Sigma')$
by \ $\phi(s_i)=S_{i+2}'$, \ \ $i=1,...,k-3$.
By Theorem \ref{Thm: homomorphisms B(k) to S(n), n<k}$(a)$,
condition (\ref{7.6}) implies that
this homomorphism $\phi$ is cyclic.
Hence $S_3'=\ldots =S_{k-1}'$.
In view of (\ref{7.7}), this shows that
$\widehat\sigma_3=\ldots =\widehat\sigma_{k-1}$
and the homomorphism $\psi$ is cyclic.
\end{proof}

\noindent The following lemma supplies the upper bound $t\le k-2$
for the length $t$ of any non-degenerate component of every permutation
$\widehat\sigma_i$ (provided $6<k<n<2k$ and $\psi$ is non-cyclic).
In fact Theorem \ref{Thm: homomorphisms B(k) to S(n), 6<k<n<2k}
shows that $t\le (k+1)/2$;
however, we are not ready to prove the latter statement now.

\begin{Lemma}\label{Lm: no components of length t>k-3}  
Let $6<k<n<2k$ and let
$\psi\colon B_k\to{\mathbf S}(n)$ be a homomorphism
such that the permutation $\widehat\sigma_1$ has
a non-degenerate component $\mathfrak C$ of length $t>k-3$.
Then $\psi$ is cyclic.
\end{Lemma}

\begin{proof} 
Suppose to the contrary that $\psi$ is non-cyclic.
The assumptions $n<2k$ and $t>k-3$ imply that $\mathfrak C$ is the
$2$-component of $\widehat\sigma_1$. Put $\Sigma=\supp \mathfrak C$ and
$\Sigma' = {\boldsymbol\Delta}_n - \Sigma$; so $\# \Sigma = 2t$
and $\# \Sigma' = n-2t < 2k - 2(k-2) = 4$.
Since $k>6$, any homomorphism $B_{k-2}\to S(\Sigma')$
is cyclic (Theorem \ref{Thm: homomorphisms B(k) to S(n), n<k}$(a)$).
In particular, the homomorphism
$\varphi_{_{\Sigma'}}$ is cyclic, and Lemma
\ref{Lm: if psi: B_k to S(n) is non-cyclic and varphi(Sigma') is abelian then phi, varphi(Sigma) and Omega are non-abelian} implies
that {\sl the homomorphisms
$\Omega\colon B_{k-2}\to{\mathbf S}(t)$
and $\varphi_{_{\Sigma}}\colon B_{k-2}
\to G\subset {\mathbf S}(2t)$
must by non-cyclic}.

We may assume that the homomorphism $\psi$ is normalized;
this means that
$$
\Sigma = \{1,2,\ldots,2t\},\ \ \widehat\sigma_1|\Sigma = C_1\cdots C_t, \
\text{where} \ C_m = (2m-1,2m) \ \text{for} \ m = 1,\ldots ,t.
$$
Clearly, $t\le k-1$; so, either $t=k-2$ or $t=k-1$.
Therefore, we must consider the following three cases:
\vskip0.3cm

\begin{itemize}

\item[$i)$] $t=k-2$ and $\widehat\sigma_1$
is a disjoint product of $k-2$ transpositions;

\item[$ii)$] $t=k-2$ and $\widehat\sigma_1$
is a disjoint product of $k-2$ transpositions and a $3$-cycle;

\item[$iii)$] $t=k-1$ and $\widehat\sigma_1$
is a disjoint product of $k-1$ transpositions.
\end{itemize}
\vskip0.3cm

$i)$ In this case we deal with the non-cyclic homomorphisms
$\Omega\colon B_{k-2}\to{\mathbf S}(k-2)$ and
$\varphi_{_{\Sigma}}\colon B_{k-2}\to G\subset {\mathbf S}(2k-4)$.
By Artin Theorem and Theorem \ref{Thm: Improvement of Artin Theorem}, either
\vskip0.3cm

\ \ \ $i_a)$ $\Omega$  {\sl is conjugate to the canonical epimorphism}
$\mu\colon B_{k-2}\to{\mathbf S}(k-2)$

\noindent or

\ \ \ $i_b)$ $k = 8$ {\sl and $\Omega$ is conjugate to Artin's homomorphism}
$\nu_6\colon B_6\to{\mathbf S}(6)$.
\vskip0.3cm
\noindent We shall show that these cases are impossible.
\vskip0.2cm %%&&

$i_a)$ In this case, without loss of generality
we may assume that $\Omega =\mu$ (see Remark \ref{Rmk 5.1}).
Then the homomorphism $\varphi_{_{\Sigma}}$ must be $H$-conjugate
to one of the four $\mu$-homomorphisms
$\varphi_j\colon B_{k-2}\to G \subseteq {\mathbf S}(2k-4)$,
\ $j = 0,1,2,3$, listed in Corollary \ref{Crl 6.6} (with $n=k-2$).
\vskip0.2cm

The homomorphisms $\varphi_0,\varphi_1,\varphi_3$ may be
eliminated by trivial reasons. Indeed, if $\varphi_{_{\Sigma}}$
is conjugate to
one of $\varphi_0,\varphi_1$, then the support of
$\widehat\sigma_{i+2}|_\Sigma=\varphi_{_{\Sigma}}(s_i)$ consists of $4$
points; however,
$\# \supp \widehat\sigma_{i+2} = \# \supp \widehat\sigma_1 = 2k-4$, and
the rest $2k-8$ points of $\supp \widehat\sigma_{i+2}$
must be situated in the set $\Sigma'$ containing at most $3$ points,
which is impossible. If $\varphi_{_{\Sigma}}\sim \varphi_3$,
then the cyclic decomposition of
$\widehat\sigma_{i+2}|_\Sigma=\varphi_{_{\Sigma}}(s_i)$
contains a $4$-cycle; but this is not the case,
since $\widehat\sigma_{i+2}\sim\widehat\sigma_1$.
\vskip0.2cm

By condition (!!!) (see Declaration that follows Remark \ref{Rmk 5.1}),
we may assume that $\varphi_{_{\Sigma}}$ {\sl coincides}
with $\varphi_2$. Then any permutation
$\widehat\sigma_{i+2}|_\Sigma=\varphi_{{\Sigma}}(s_i)=\varphi_2(s_i)$
is a product of $k-2$ disjoint transpositions.
However $\widehat\sigma_{i+2}$ itself is such a product; therefore,
$\widehat\sigma_{i+2}=\varphi_2(s_i)$.
In particular, for $i=1$ and $i=k-3$, we obtain
\vskip0.3cm
$$
\aligned
&\widehat\sigma_3 \ \ \ =
\underbrace{(1,3)(2,4)}(5,6)\cdots (2k-9,2k-8)(2k-7,2k-6)(2k-5,2k-4)\quad
\text{and}\\
&\widehat\sigma_{k-1} = (1,2)(3,4)(5,6)\cdots (2k-9,2k-8)
\underbrace{(2k-7,2k-5)(2k-6,2k-4)}
\endaligned
$$
respectively. Using these formulae and computing the permutations
$\widehat\sigma_3\cdot C\cdot \widehat\sigma_3^{-1}$
for some particular transpositions $C\preccurlyeq \widehat\sigma_{k-1}$,
we obtain:
$$
\aligned
&\widehat\sigma_3 (1,2)\widehat\sigma_3^{-1}=(3,4),\qquad \
\widehat\sigma_3 (2k-7,2k-5) \widehat\sigma_3^{-1}=(2k-6,2k-4),\\
&\widehat\sigma_3 (3,4) \widehat\sigma_3^{-1}=(1,2),\qquad \
\widehat\sigma_3 (2k-6,2k-4) \widehat\sigma_3^{-1}=(2k-7,2k-5).
\endaligned
$$
In view of our definition of the homomorphism $\Omega^*$,
these formulae show that the cyclic decomposition of the permutation
$\Omega^*(s_1)\in{\mathbf S}(\mathfrak C^*)\cong {\mathbf S}(k-2)$
must contain at least two disjoint transpositions. On the other hand,
since $\Omega\sim\mu$, the permutation $\Omega(s_1)\sim\mu(s_1)$
is a transposition; hence $\Omega^*(s_1)\ne\Omega(s_1)$,
contradicting Lemma \ref{Lm: Omega=Omega*}$(a)$.
\vskip0.3cm

$i_b)$ We may assume that $\Omega =\nu_6$ (see Remark \ref{Rmk 5.1}).
Then the homomorphism $\varphi_{_{\Sigma}}$ must be $H$-conjugate to one of
the two $\nu_6$-homomorphisms
$\phi_j\colon B_6\to G\subseteq {\mathbf S}(12)$, \ $j = 0,1$,
listed in Corollary \ref{Crl 6.8}. We may assume that
$\varphi_{_{\Sigma}}$ {\sl coincides} with one of the homomorphisms $\phi_j$.
Then any permutation 
$\widehat\sigma_{i+2}|_\Sigma=\varphi_{_{\Sigma}}(s_i)$, \ $i=1,...,5$,
is a product of $6$ disjoint transpositions.
Since $\widehat\sigma_{i+2}$ itself is such a product,
this means that either $\widehat\sigma_{i+2}=\phi_0(s_i)$
or $\widehat\sigma_{i+2}=\phi_1(s_i)$ for all $i=1,...,5$.

In each of these two cases the permutations $\phi_j(s_1)$ and
$\phi_j(s_5)$ contain two common transpositions, namely,
$D_5=(9,11)$ and $D_6=(10,12)$. Hence the conjugation
of the six transpositions $D_m\preccurlyeq \widehat\sigma_7=\phi_j(s_5)$
by the permutation $\widehat\sigma_3=\phi_j(s_1)$ does not
change these transpositions $D_5$ and $D_6$. Consequently,
the permutation $\Omega^*(s_1)\in{\mathbf S}(\mathfrak C^*)
\cong {\mathbf S}(6)$ 
(defined by this conjugation) has at least two fixed points. 
However this contradicts Lemma \ref{Lm: Omega=Omega*}$(a)$,
since the permutation $\Omega(s_1)=\nu_6(s_1)
=(C_1,C_2)(C_3,C_4)(C_5,C_6)\in{\mathbf S}(\mathfrak C)
\cong {\mathbf S}(6)$ has no fixed points.
This concludes the proof in case $(i)$.
\vskip0.2cm %%&&

$ii)$ In this case $2k>n\ge\#\supp\,\widehat\sigma_1 = 2(k-2)+3 = 2k-1$;
so, $n=2k-1$ and $\widehat\sigma_1$ has no fixed points. Hence
$\Sigma'$ (the support of the $3$-cycle) is
the only $\widehat\sigma_1$-invariant set of length $3$
and, by Lemma \ref{Lm: if |Inv(r,psi(sigma(1)))|=1 then psi is transitive}, 
$\psi$ is intransitive. Clearly, $\psi$ is the disjoint
product of its reductions $\psi_\Sigma$ and $\psi_{\Sigma'}$ to
the $(\Img \psi)$-invariant sets $\Sigma$ and $\Sigma'$ respectively.

The homomorphism $\psi_{\Sigma'}$ is cyclic ($k>\#\Sigma'$).
As to the homomorphism 
$$
\widetilde\psi = \psi_\Sigma\colon B_k
\to{\mathbf S}(\Sigma)\cong {\mathbf S}(2k-4),
$$
it must be non-cyclic (since $\psi$ is non-cyclic). Clearly,
the permutation
$$
\widetilde\psi(\sigma_1)=\psi_\Sigma(\sigma_1)=\psi(\sigma_1)|\Sigma
$$ 
is a product of $k-2$ disjoint transposition.
However, it has been already proven that this is impossible (see case $(i)$).
\vskip0.3cm

$iii)$ In this case either $n=2k-2$ or $n=2k-1$, and the set
$\Sigma'$ contains at most $1$ point.
We deal with the non-cyclic homomorphism
$\Omega\colon B_{k-2}\to{\mathbf S}(k-1)$.
It follows from Theorem \ref{Thm: homomorphisms B(k) to S(k+1)}
and Proposition \ref{Prp 4.9} 
that either
\vskip0.3cm

\ \ \ $iii_a)$ $\Omega$ {\sl is conjugate to the homomorphism}
$$
\mu_{k-1}^{k-2} = \mu_{k-2}\times\mathbf 1_1
\colon B_{k-2}\to{\mathbf S}(k-1),
$$
\noindent or

\ \ \ $iii_b)$ $k = 7$ {\sl and $\Omega$ 
is conjugate to the homomorphism}
$\psi_{5,6}\colon B_5\to{\mathbf S}(6)$.
\vskip0.3cm

\noindent Let us show that these cases are impossible.
\vskip0.3cm

$iii_a)$ Again, we may assume that $\Omega=\mu_{k-1}^{k-2}$.
Then the homomorphism $\varphi_{_{\Sigma}}$ must be $H$-conjugate
to one of the eight homomorphisms $\psi_{0;j}$, \ $\psi_{1;j}$,
$$
\psi_{0;j}=\varphi_j\times\mathbf 1_2, \ \
\psi_{1;j}=\varphi_j\times\mathbf {\{2\}},\ \ \ j=0,1,2,3,
$$
listed in Corollary \ref{Crl 6.9} (with $n=k-2$).

The homomorphisms $\psi_{x;0}$, $\psi_{x;1}$, $\psi_{x;3}$ ($x=0,1$) may be
eliminated as in case $(i_a)$.
Indeed, if $\varphi_{_{\Sigma}}$ is conjugate to one of $\psi_{x;1}$,
$\psi_{x;3}$, then the cyclic decomposition of
$\widehat\sigma_{i+2}|_\Sigma=\varphi_{_{\Sigma}}(s_i)$
contains a $4$-cycle, which is impossible 
(for $\widehat\sigma_{i+2}\sim\widehat\sigma_1$).
If $\varphi_{_{\Sigma}}\sim\psi_{x;0}$, then either
$[\widehat\sigma_{i+2}|_\Sigma]=[2,2]$ or
$[\widehat\sigma_{i+2}|_\Sigma]=[2,2,2]$.
Since $\widehat\sigma_{i+2}$ is a disjoint product of $k-1$ transpositions,
at least $(k-1)-3=k-4$ of them must be situated on the set
$\Sigma'$ containing at most $1$ point, which is impossible.

So, we may assume that the homomorphism $\varphi_{_{\Sigma}}$ coincides
with one of the homomorphisms $\psi_{x;2}$, \ $x=0,1$. 
In any of these two cases
$$               
\widehat\sigma_{i+2}|\Sigma
=\varphi_{_{\Sigma}}(s_i)=\psi_{x;2}(s_i)=\varphi_2(s_i)\cdot T^x,
$$ 
where $T=(2k-3,2k-2)$. Hence $\widehat\sigma_{i+2}=\varphi_2(s_i)\cdot T^x$
for all $i=1,...,k-2$.
For $i=1$ and $i=k-3$, we have respectively,
$$
\aligned
\ &\widehat\sigma_3 \ \ \ =
\underbrace{(1,3)(2,4)}(5,6)\cdots (2k-9,2k-8)(2k-7,2k-6)(2k-5,2k-4)
\cdot T^x,\\
&\widehat\sigma_{k-1} = (1,2)(3,4)(5,6)\cdots (2k-9,2k-8)
\underbrace{(2k-7,2k-5)(2k-6,2k-4)}\cdot T^x.
\endaligned
$$
As in case $(i_a)$,
$$
\aligned
&\widehat\sigma_3 (1,2) \widehat\sigma_3^{-1}=(3,4),\qquad \
\widehat\sigma_3 (2k-7,2k-5) \widehat\sigma_3^{-1}=(2k-6,2k-4),\\
&\widehat\sigma_3 (3,4) \widehat\sigma_3^{-1}=(1,2),\qquad \
\widehat\sigma_3 (2k-6,2k-4) \widehat\sigma_3^{-1}=(2k-7,2k-5).
\endaligned
$$
This implies that the cyclic decomposition
of $\Omega^*(s_1)$ contains at least two disjoint
transpositions, which is impossible, since 
$\Omega(s_1)=\mu_{k-1}^{k-2}(s_1)=(1,2)$.
\vskip0.3cm

$iii_b)$ In this case, we can assume that
$\Omega=\psi_{5,6}\colon B_5\to{\mathbf S}(6)$
and $\varphi_{_{\Sigma}}$ coincides with one of the four
homomorphisms $\eta_j\colon B_5\to S(12)$,
\ $j=0,1,2,3$, \ exhibited in Corollary \ref{Crl 6.7}. 
We may certainly exclude the homomorphisms $\eta_1,\eta_3$
because of the $4$-cycles presence.
If $\varphi_{_{\Sigma}}=\eta_j$, \ $j=0$ or $j=2$,
then, in fact, $\widehat\sigma_{i+2}=\eta_j(s_i)$ for all
$i=1,2,3,4$. In both cases the permutations $\eta_j(s_1)$ and
$\eta_j(s_4)$ contain the two common transpositions
$D_1=(1,3)$ and $D_2=(2,4)$; as in case $(i_b)$, it follows that
the permutation $\Omega^*(s_1)\in{\mathbf S}(\mathfrak C^*)
\cong {\mathbf S}(6)$ has at least two fixed points.
However this contradicts
Lemma \ref{Lm: Omega=Omega*}$(a)$, since the permutation 
$\Omega(s_1)=\psi_{5,6}(s_1)=(C_1,C_2)(C_3,C_4)(C_5,C_6)\in 
{\mathbf S}(\mathfrak C)
\cong {\mathbf S}(6)$ has no fixed points.
\end{proof}

\begin{Corollary}\label{Crl: k-2 fixed points of psi(sigma(1)) for 6<k<n<2k} 
\index{Fixed points!of $\widehat\sigma_1$ for $6<k<n<2k$\hfill}
\index{Artin Fixed Point Lemma\hfill|phantom\hfill}
\index{Artin Lemma!on fixed points of $\widehat\sigma_1$\hfill}
\index{Artin Fixed Point Lemma!analog of\hfill}
Let $6<k<n<2k$ and let $\psi\colon B_k\to{\mathbf S}(n)$ be
a non-cyclic homomorphism. Then the permutation
$\widehat\sigma_1$ has at least $k-2$ fixed points.
\end{Corollary}

\begin{proof} 
Lemma \ref{Lm: no components of length t>k-3}
implies that $\widehat\sigma_1$ does not
possess non-degenerate components of length $t>k-3$; 
combined with 
Lemma \ref{Lm1: homomorphism with all components of hatsigma(1) of length<=k-3}, 
this shows that $\#\Fix \widehat\sigma_1\ge k-2$.
\end{proof}

\begin{Remark}\label{Rmk: about Corollary on fixed points without primes} 
Corollary \ref{Crl: k-2 fixed points of psi(sigma(1)) for 6<k<n<2k}
is a useful and powerful
partner of Artin Fixed Point Lemma \ref{Lm: Artin Fixed Point Lemma}. 
\index{Fixed points!of $\widehat\sigma_1$ for $6<k<n<2k$\hfill}
\index{Artin Fixed Point Lemma\hfill|phantom\hfill}
\index{Artin Lemma!on fixed points of $\widehat\sigma_1$\hfill}
\index{Artin Fixed Point Lemma!analog of\hfill}
If $6<k<n<2k$ and $n$ is ``far" from $k$, 
we cannot be sure that there is a prime number $p$ between $n/2$ and
$k-2$ and  Lemma \ref{Lm: Artin Fixed Point Lemma} does not work,
while Corollary \ref{Crl: k-2 fixed points of psi(sigma(1)) for 6<k<n<2k}
applies. One should however take into account that 
Artin Fixed Point Lemma \ref{Lm: Artin Fixed Point Lemma}
played important part in the proof of
Corollary \ref{Crl: k-2 fixed points of psi(sigma(1)) for 6<k<n<2k}
(by means of Artin Theorem,
Theorem \ref{Thm: homomorphisms B(k) to S(n), n<k}$(a)$,
Theorem \ref{Thm: homomorphisms B(k) to S(k+1)}, 
Lemma \ref{Lm1: homomorphism with all components of hatsigma(1) of length<=k-3}
and Lemma \ref{Lm: no components of length t>k-3}).
\end{Remark}

\noindent In order to study homomorphisms 
$B_k\to{\mathbf S}(k+2)$,
we start with the exceptional cases $k=5$ and $k=6$.

\begin{Remark}\label{Rmk 7.4} 
Note that for $k>4$ there are the following
evident non-cyclic homomorphisms
$\psi\colon B_k\to{\mathbf S}(k+2)$:
$$
\aligned
\ &\psi^k_{k+2}=\mu_k\times\mathbf 1_2\colon B_k
\to{\mathbf S}(k+2);\qquad
\widetilde\psi^k_{k+2}=\mu_k\times\mathbf {\{2\}}\colon B_k
\to{\mathbf S}(k+2);\\
&\phi^6_8 \ \ \ = \ \nu_6\times\mathbf 1_2
\colon B_6\to{\mathbf S}(8);
\qquad\qquad 
\widetilde\phi^6_8 \ \ \ =\nu_6\times\mathbf {\{2\}}\colon 
\,B_6\to{\mathbf S}(8);\\
&\qquad\qquad\qquad\qquad\qquad \ \ \
\psi^5_7=\psi_{5,6}\times\mathbf 1_1\colon B_5
\to{\mathbf S}(7),\\
\endaligned
$$
\noindent where $\psi_{5,6}$ is defined by (\ref{4.6}) or (\ref{4.6'}),
and $\nu_6$ is Artin's homomorphism. All 
these homomorphisms have the
following property: {\sl the image of any generator $\sigma_i$
is a product of disjoint transpositions.}
\end{Remark}

\begin{Proposition}\label{Prp: B(5) to S(7)} 
Let
$\psi\colon B_5\to{\mathbf S}(7)$
be a non-cyclic homomorphism. Then $\psi$ is intransitive
and conjugate to one of the homomorphism $\psi^5_7$,
\ $\widetilde\psi^5_7$, $\phi^5_7$. In particular,
every permutation $\widehat\sigma_i$,
\ $1\le i\le 4$, is a product of disjoint transpositions.
\end{Proposition}

\begin{proof}
Suppose to the contrary that $\psi$ is transitive.
Then Lemma \ref{Lm: condition no cycles of the same length restricts support}
implies that if all the cycles 
$C_\nu\preccurlyeq \widehat\sigma_1$ are of the distinct
lengths $r_\nu$, then $\sum r_\nu \le 3$. This eliminates
all cyclic types for $\widehat\sigma_1$ but $[2]$, $[3]$, $[2,2]$,
$[3,3]$, $[2,2,2]$, $[2,2,3]$. However, Lemma
\ref{Lm: if [psi(sigma(1))]=[3] or |Fix(psi(sigma(1))|>n-4 then psi is intransitive}
eliminates $[2]$ and $[3]$,
Lemma \ref{Lm 7.2}$(c)$ eliminates $[2,2]$,
Lemma \ref{Lm: if |Inv(r,psi(sigma(1)))|=1 then psi is transitive} 
eliminates the rest types, and a contradiction ensues. 

Further, since $\psi$ is non-cyclic and (as we already have proved)
intransitive, the group $G=\Img \psi\subset {\mathbf S}(7)$
has exactly one orbit $Q$ of length $L$, \ $5\le L\le 6$;
put $Q'={\boldsymbol\Delta}_7 - Q$. 
Clearly, either $Q'$ is a $G$-orbit
of length $7-L$ or $Q'$ consists of $7-L$ fixed points.
The homomorphism $\psi$ is the disjoint product of its reductions
$\psi_Q$ and $\psi_{Q'}$.
The reduction $\psi_Q$ is a non-cyclic transitive homomorphism
$B_5\to{\mathbf S}(Q)\cong {\mathbf S}(L)$. 
By Artin Theorem and Proposition \ref{Prp 4.9}, we obtain that 
(up to conjugation of $\psi$) either $L=5$ and $\psi_Q=\mu_5$ 
or $L=6$ and $\psi_Q=\psi_{5,6}$. The reduction 
$\psi_{Q'}$ is either trivial or takes all $\sigma_i$ to
the same transposition. This concludes the proof.
\end{proof}

In the sequel we use the following theorem of C. Jordan
(see \cite{Wiel}, Theorem 13.3, Theorem 13.9} or
\cite{Hall}, Theorem 5.6.2, Theorem 5.7.2):
\vskip0.2cm

{\scaps Jordan Theorem.} 
{Let $G\subseteq {\mathbf S}(n)$ be a primitive group of permutations.
If $G$ contains a transposition, then $G={\mathbf S}(n)$;
if $G$ contains a $3$-cycle, then either $G={\mathbf S}(n)$
or $G=\mathbf A(n)$. Moreover, if $n=p+r$, where $p$ is prime, $r\ge 3$,
and $G$ contains a $p$-cycle, then either $G={\mathbf S}(n)$ 
or $G=\mathbf A(n)$.}
%\endproclaim
\vskip0.3cm

\begin{Remark}\label{Rmk 7.5} 
The following fact (certainly, well known
to experts in permutation groups)
follows trivially from Jordan Theorem:
\vskip0.3cm

{\sl Let $G\subseteq {\mathbf S}(n)$ be a primitive permutation group.
Assume that either \ $i)$ $k<n$ and $G$ contains
a $p$-cycle of length $p\le 3$, or
\ $ii)$ $k=6$ and $8\le n\le 9$, or \ $iii)$ $k=7$ and $10\le n\le 13$. 
Then the group $G$ cannot be isomorphic to ${\mathbf S}(k)$.}
\vskip0.3cm

Indeed, suppose to the contrary that $G\cong S(k)$; then
$\# G = k! < n!/2$. Let us show 
that either $G={\mathbf S}(n)$ or $G=\mathbf A(n)$
(clearly, this will contradict the above inequality).
In case $(i)$ our statement follows immediately from Jordan Theorem.
In all other cases the assumption $G\cong S(k)$ implies
that there is an element $g\in G$ of order $p$, 
where $p=5$ in case $(ii)$ and $p=7$ in case $(iii)$.
It follows from the constraints on $n$ that $g$ is a $p$-cycle
and $r=n-p\ge 3$. Hence the last sentence of Jordan Theorem
shows that either $G={\mathbf S}(n)$ or $G=\mathbf A(n)$.
\end{Remark}

\begin{Proposition}\label{Prp: B(6) to S(8)} 
Any transitive homomorphism
$\psi\colon B_6\to{\mathbf S}(8)$ is cyclic.
\end{Proposition}

\begin{proof} 
Suppose that $\psi$ is non-cyclic. We claim that then
{\sl any cycle $C\preccurlyeq \widehat\sigma_1$
must be a transposition, that is, $\widehat\sigma_1^2=1$.} 
If this is the case,
then $\widehat\sigma_i^2=1$ for all $i$ and hence
$PB_6\subseteq\Ker\psi$. Therefore, there is a homomorphism
$\phi\colon {\mathbf S}(6)\to{\mathbf S}(8)$
such that $\psi=\phi\circ\mu_6$. Since $\mu_6$ is surjective, we
have $G=\Img \psi=\Img \phi\subset {\mathbf S}(8)$. The group
$G$ is transitive and non-abelian,
and the group ${\mathbf S}(6)$ has no proper non-abelian quotient groups;
hence $G\cong {\mathbf S}(6)$. Moreover, by Proposition
\ref{Prp:
transitive imprimitive homomorphisms Bk to S(n) is cyclic if n<2k}$(a)$,
the group $G$ is primitive. These properties contradict Jordan Theorem 
(see Remark \ref{Rmk 7.5}). 

To justify the above claim, assume, on the contrary, that
$\widehat\sigma_1^2\ne 1$, and let us find out the possible cyclic types of
$\widehat\sigma_1$. Certainly, $\widehat\sigma_1$ cannot contain a cycle of
length $>4$ (say by
Lemma \ref{Lm: Artin Lemma on cyclic decomposition of hatsigma(1)}$(a)$).
Since $E(8/E(6/2))=2$,
Lemma \ref{Lm: condition no cycles of the same length restricts support} shows 
that if all the cycles $\preccurlyeq\widehat\sigma_1$
are of distinct lengths, then $\widehat\sigma_1$ is a transposition
(in fact, this is forbidden by Lemma
\ref{Lm: if [psi(sigma(1))]=[3] or |Fix(psi(sigma(1))|>n-4 then psi is intransitive}); 
so, we may assume that \ 
(*) {\sl $\widehat\sigma_1$ contains at least two cycles of the same length}.
Under this assumption,
Lemma \ref{Lm: Artin Lemma on cyclic decomposition of hatsigma(1)}$(b)$
shows that if 
$\widehat\sigma_1$ contains a $4$-cycle, then either $[\widehat\sigma_1]=[4,4]$
or $[\widehat\sigma_1]=[4,2,2]$. Moreover, (*) and
Lemma \ref{Lm: if |Inv(r,psi(sigma(1)))|=1 then psi is transitive}
eliminate all cyclic types with a $3$-cycle. This
shows that \ (**) {\sl either $[\widehat\sigma_1]=[4,4]$ or
$[\widehat\sigma_1]=[4,2,2]$}. Consider now the permutation
$A=\widehat\alpha_{3,5}=\widehat\sigma_3\widehat\sigma_4$.
By Corollary
\ref{Crl: if psi is non-cyclic then psi(sigma(i)) infty psi(sigma(i+1))} 
(with $i=3$ and $j=5$), $\ord A$ is divisible by $3$;
hence only the following cyclic types of $A$ may occur:
%$[3]$, $[3,2]$, $[3,2,2]$, $[3,3]$, $[3,3,2]$, $[3,4]$, $[3,5]$,
%$[6]$, $[6,2]$.
$$
[3],\, [3,2],\, [3,2,2],\, [3,3],\, [3,3,2],\, [3,4],\, [3,5],\,
[6],\, [6,2]\,. %%&&&
$$
Since $\widehat\sigma_1$ satisfies (**) and commutes with $A$,
Lemma \ref{Lm: invariant sets and components of commuting permutations}$(b)$
eliminates all cyclic types but $[3,3]$ and $[3,3,2]$ 
(were $A$ of any other of the listed above types,
$\widehat\sigma_1$ would contain either a power
of a $3$-cycle or a power of a $6$-cycle; however this contradicts (**)). 
The same argument eliminates the type $[\widehat\sigma_1]=[4,2,2]$
(a power of a $4$-cycle cannot occure in the cyclic decomposition of
a permutation of cyclic type $[3,3]$ or $[3,3,2]$). Finally, the case
$[\widehat\sigma_1]=[4,4]$ is also impossible; indeed, $A$ has
only one invariant set of length $2$ (two fixed points
for $[A]=[3,3]$ and the support of a transposition 
for $[A]=[3,3,2]$); the latter set must be 
$\widehat\sigma_1$-invariant, which cannot happen if
$[\widehat\sigma_1]=[4,4]$). This completes the proof.
\end{proof}

\begin{Theorem}\label{Thm: homomorphisms B(k) to S(k+2)} 
$a)$ Any transitive homomorphism
$\psi\colon B_k\to{\mathbf S}(k+2)$ is cyclic
whenever $k>4$.
\vskip0.2cm

\noindent $b)$ Let $k>4$ and let $\psi\colon B_k\to{\mathbf S}(k+2)$
be a non-cyclic homomorphism. Then either $\psi$ is conjugate
to one of the homomorphism $\psi^k_{k+2}$, \ $\widetilde\psi^k_{k+2}$
$($this may happen for any $k)$ or $k=5$ and $\psi$ is conjugate
to the homomorphism $\phi^5_7$, or, finally, $k=6$ and $\psi$
is conjugate to one of the homomorphisms $\phi^6_8$, \ $\widetilde\phi^6_8$.
\end{Theorem}

\begin{proof} 
$a)$ The cases $k=5$ and $k=6$ are already
considered in Propositions \ref{Prp: B(5) to S(7)} and
\ref{Prp: B(6) to S(8)}; assume now that $k>6$
and that the homomorphism $\psi$ is non-cyclic.
By Corollary \ref{Crl: k-2 fixed points of psi(sigma(1)) for 6<k<n<2k} 
$(a)$, the permutation $\widehat\sigma_1$ has at least
$k-2$ fixed points, and thus $\# \supp \widehat\sigma_1\le (k+2)-(k-2)=4$.
Hence only the following cyclic types of $\widehat\sigma_1$ may
occur: $[2]$, $[2,2]$, $[3]$, $[4]$. However, $[2]$, $[3]$ are forbidden
by Lemma
\ref{Lm: if [psi(sigma(1))]=[3] or |Fix(psi(sigma(1))|>n-4 then psi is intransitive},
\ $[2,2]$ is forbidden by Lemma \ref{Lm 7.2}
$(c)$, and $[4]$ is forbidden by
Lemma \ref{Lm: condition no cycles of the same length restricts support},
since the numbers $R_k=(k+2)/E(k/2)$ satisfy $E(R_7)=3$ for $k>6$
and $E(R_k)=2$ for $k>7$.
\vskip0.3cm

$b)$ Since $\psi$ is non-cyclic and (by the statement $(a)$) intransitive,
Theorem \ref{Thm: homomorphisms B(k) to S(n), n<k}
$(a)$ shows that the group $G=\Img \psi\subset {\mathbf S}(k+2)$
has exactly one orbit $Q$ of length $L$, \ $k\le L\le k+1$;
set $Q'={\boldsymbol\Delta}_{k+2} - Q$.
Clearly, either $Q'$ is a $G$-orbit
of length $k+2-L$ or $Q'$ consists of $k+2-L$ fixed points.
The homomorphism $\psi$ is the disjoint product of its reductions
$\psi_Q$ and $\psi_{Q'}$.
The reduction $\psi_Q$ is a non-cyclic transitive homomorphism
$B_k\to{\mathbf S}(Q)\cong {\mathbf S}(L)$.
By Artin Theorem, Theorem \ref{Thm: homomorphisms B(k) to S(k+1)}
and Proposition \ref{Prp 4.9}, 
we obtain that
(up to conjugation of $\psi$) either $\psi_Q=\mu_k$, or $k=5$
and $\psi_Q=\psi_{5,6}$, or, finally, $k=6$ and $\psi_Q=\nu_6$.
The reduction $\psi_{Q'}$ is either trivial or takes all $\sigma_i$
to the same transposition. This concludes the proof.
\end{proof}

\subsection{Homomorphisms $B_k\to {\mathbf S}(n)$, \
$6<k<n<2k$}
\label{Ss: Homomorphisms B(k) to S(n), 6<k<n<2k}
Our main goal now is to prove Theorem F$(a)$ 
(see Theorem \ref{Thm: homomorphisms B(k) to S(n), 6<k<n<2k} below).
We prove this theorem by induction on $k$.
Lemma \ref{Lm: no components of length t>k-3}
and Lemma
\ref{Lm2: homomorphism with all components of hatsigma(1) of length<=k-3}
enable us to pass from $k$ to $k+2$ (the inductive step);
Lemma \ref{Lm: B(7) to S(n) for 7<n<14} and
Lemma \ref{Lm: B(8) to S(n), 8<n<16}
provide a base of induction.
\vskip0.3cm

\noindent{\bfit Convention.} Given a homomorphism
$\psi\colon B_k\to{\mathbf S}(n)$,
we use the following notation. We set
%&
\begin{equation}\label{7.8}
\aligned
\ &\widehat\sigma_i = \psi(\sigma_i), \ 1\le i\le k-1;\qquad
\Sigma_i = \supp \widehat\sigma_i; \qquad \Sigma_i' = \Fix \widehat\sigma_i;\\
&N = \# \Sigma_i; \qquad
N' = \# \Sigma_i';\qquad G=\Img \psi \subseteq {\mathbf S}(n).
\endaligned
\end{equation}
%&
Obviously, $N$, $N'$ {\sl do not depend on $i$, and $N+N'=n$.}
Moreover, {\sl if $k>4$, \ $n<2k$, and the homomorphism $\psi$ is non-cyclic,
then $G$ is a non-abelian primitive permutation group of degree $n$}
(see Proposition \ref{Prp: transitive imprimitive homomorphisms Bk to S(n) is cyclic if n<2k}). We denote by
%&
\begin{equation}\label{7.9}
\phi\colon B_{k-2}\stackrel{\cong}{\longrightarrow} B
\stackrel{\psi|{B}}{\longrightarrow} {\mathbf S}(n)
\end{equation}
%&
the restriction of $\psi$ to the subgroup
$B\cong B_{k-2}\subset B_k$
generated by $\sigma_3,...,\sigma_{k-1}$, and set
$$
H=\Img \phi\subseteq G\subseteq {\mathbf S}(n).
$$
Since $\widehat\sigma_1$ commutes with
$\widehat\sigma_3,...,\widehat\sigma_{k-1}$,
{\sl the sets $\Sigma_1$ and $\Sigma_1'$ are $H$-invariant, and
the homomorphism $\phi$
is the disjoint product $\varphi\times \varphi'$ of its reductions
%&
\begin{equation}\label{7.10}
\varphi=\phi_{\Sigma_1}\colon B_{k-2}\to
{\mathbf S}(\Sigma_1)\cong {\mathbf S}(N)
\end{equation} 
%&
and
%&
\begin{equation}\label{7.11}
\varphi'=\phi_{\Sigma_1'}
\colon B_{k-2}\to
{\mathbf S}(\Sigma_1')\cong {\mathbf S}(N')
\end{equation}
%&
to the sets $\Sigma_1$ and $\Sigma_1'$ respectively.}
We consider also the restriction 
%&
\begin{equation}\label{7.9t}
\widetilde\phi\colon 
\,B_{k-2}\stackrel{\cong}{\longrightarrow} {\widetilde B}
\stackrel{\psi|{\widetilde B}}{\longrightarrow} {\mathbf S}(n)
\end{equation}
%&
of the homomorphism $\psi$ to the subgroup 
$\widetilde{B}\cong B_{k-2}\subset B_k$
generated by $\sigma_1,...,\sigma_{k-3}$.
{\sl The homomorphism $\widetilde\phi$ is the disjoint product
$\widetilde\varphi\times \widetilde\varphi'$ of its reductions
%&
\begin{equation}\label{7.10t}
\widetilde\varphi=\widetilde\phi_{\Sigma_{k-1}}
\colon B_{k-2}\to
{\mathbf S}(\Sigma_{k-1})\cong {\mathbf S}(N)
\end{equation}
%&
and
%&
\begin{equation}\label{7.11t}
\widetilde\varphi'=\widetilde\phi_{\Sigma_{k-1}'}
\colon B_{k-2}\to
{\mathbf S}(\Sigma_{k-1}')\cong {\mathbf S}(N')
\end{equation}
%&
to the sets $\Sigma_{k-1}$ and $\Sigma_{k-1}'$ respectively.}
\hfill $\square$

\begin{Lemma} 
\label{Lm2: homomorphism with all components of hatsigma(1) of length<=k-3}
Let $k>6$ and let $\psi\colon B_k\to{\mathbf S}(n)$ be
a non-cyclic homomorphism. Assume that every non-degenerate component
of $\widehat\sigma_1$ is of length at most $k-3$
{\rm (this is certainly the case if $6<k<n<2k$;
see Lemma \ref{Lm: no components of length t>k-3})}.
Then the following statements hold true:
\vskip0.2cm

\noindent $a)$ The homomorphism $\varphi=\phi_{\Sigma_1}$
is cyclic and the homomorphism
$\varphi'=\phi_{\Sigma_1'}$ is non-cyclic.
In particular,
%&
\begin{equation}\label{7.12}
\widehat\sigma_i = S S_i' \qquad \text{for all} \ \ i=3,...,k-1\,,
\end{equation}
%&
where the permutation
%&
\begin{equation}\label{7.13}
S=\varphi(\sigma_i)=\widehat\sigma_i|{\Sigma_1}
\in{\mathbf S}(\Sigma_1)\cong {\mathbf S}(N)
\end{equation}
%&
does not depend on $i$, and all the permutations
%&
\begin{equation}\label{7.14}
S_i'=\varphi'(\sigma_i)
=\widehat\sigma_i|{\Sigma_1'}\in{\mathbf S}(\Sigma_1')
\cong {\mathbf S}(N')
\end{equation}
%&
are disjoint with $S$ and $\widehat\sigma_1$.
\vskip0.2cm

\noindent $b)$ The homomorphism
$\widetilde\varphi=\widetilde\phi_{\Sigma_{k-1}}$
is cyclic and the homomorphism
$\widetilde\varphi'=\widetilde\phi_{\Sigma_{k-1}'}$
is non-cyclic. In particular,
%&
\begin{equation}\label{7.12t}
\widehat\sigma_i = \widetilde S \widetilde {S_i'}
\qquad \text{for all} \ \ i=1,...,k-3,
\end{equation}
%&
where the permutation
%&
\begin{equation}\label{7.13t}
\widetilde S=\widetilde\varphi(\sigma_i)=\widehat\sigma_i|{\Sigma_{k-1}}
\in{\mathbf S}(\Sigma_{k-1})\cong {\mathbf S}(N)
\end{equation}
%&
does not depend on $i$, and all the permutations
%&
\begin{equation}\label{7.14t}
\widetilde {S_i'}=\widetilde\varphi'(\sigma_i)
=\widehat\sigma_i|{\Sigma_{k-1}'}\in{\mathbf S}(\Sigma_{k-1}')
\cong {\mathbf S}(N')
\end{equation}
%&
are disjoint with $\widetilde S$ and $\widehat\sigma_{k-1}$.
\vskip0.2cm

\noindent $c)$ Actually, $S=\widetilde S$ and
$S_i'=\widetilde {S_i'}$ for all $i=3,...,k-3$.
Moreover, if $\psi$ is transitive,
then $S=\widetilde S=1$, \ $\widehat\sigma_i=S_i'$
for $i\ge 3$ $($so every such $\widehat\sigma_i$ is disjoint with
$\widehat\sigma_1)$, the homomorphism $\varphi=\phi_{\Sigma_1}$ is trivial,
and the homomorphism $\varphi'=\phi_{\Sigma_1'}$ coincides with
the restriction $\phi$ of $\psi$ to the subgroup
$B_{k-2}\subset B_k$ generated by $\sigma_3,...,\sigma_{k-1}$.
\end{Lemma}

\begin{proof} 
$a)$ By Lemma \ref{Lm 7.1}, for any non-degenerate component
$\mathfrak C=\{C_1,...,C_t\}$ 
(of some length $t$, \ $1\le t\le k-3$) the retraction
$$
\Omega=\Omega_{\mathfrak C}\colon B_{k-2}\to
{\mathbf S}(\mathfrak C)\cong {\mathbf S}(t)
$$
is cyclic. By Theorem
\ref{Thm: if Omega is cyclic then every Omega-homomorphism is cyclic}$(a)$,
any $\Omega$-homomorphism
$B_{k-2}\to G_{\mathfrak C}\subseteq {\mathbf S}(\supp \mathfrak C)$
is cyclic; in particular, the homomorphism
$\varphi_{\supp \mathfrak C}\colon B_{k-2}\to G_{\mathfrak C}
\subseteq {\mathbf S}(\supp \mathfrak C)$ is cyclic
(recall that $G_{\mathfrak C}$ is the centralizer of the element
$\mathcal C=C_1\cdots C_t$ in ${\mathbf S}(\supp \mathfrak C)$).
Thereby,
$$
\varphi_{\supp \mathfrak C}(s_1)=\ldots = \varphi_{\supp \mathfrak C}(s_{k-3}).
$$
By the definition of $\varphi_{\supp \mathfrak C}$, we have
$$
\widehat\sigma_{i+2}|{\supp \mathfrak C}
=\psi(\sigma_{i+2})|{\supp \mathfrak C}
=\varphi_{\supp \mathfrak C}(s_i) = \varphi_{\supp \mathfrak C}(s_1),\qquad
i=1,...,k-3;
$$
hence the reduction 
$\phi_{\supp \mathfrak C}\colon B_{k-2}
\to{\mathbf S}(\supp \mathfrak C)$
of $\phi$ to the $H$-invariant set $\supp \mathfrak C$ is a cyclic homomorphism.
The homomorphism $\varphi=\phi_{\Sigma_1}$ 
is the disjoint product of all the
reductions $\phi_{\supp \mathfrak C}$, where $\mathfrak C$ runs over
all the non-degenerate components of $\widehat\sigma_1$.
Hence $\varphi$ is cyclic and the permutation
$S$ defined by (\ref{7.13}) does not depend on $i$.
Clearly, the permutations $S_i'$ \ ($i=3,...,k-1$) defined by (\ref{7.14})
are disjoint with $S$ and $\widehat\sigma_1$.
Finally, the homomorphism $\varphi'$ must be non-cyclic
(for otherwise, $S_3'=...=S_{k-1}'$ and
(\ref{7.12}) shows that $\psi$ is cyclic).

The statement $(b)$ follows by the same argument
(one has just to work with the permutation $\widehat\sigma_{k-1}$
and the sets $\Sigma_{k-1}$, $\Sigma_{k-1}'$ 
instead of $\widehat\sigma_1$, $\Sigma_1$ and
$\Sigma_1'$ respectively).

$c)$ The proof of this statement is an exercise in elementary set theory. 
By $(a)$ and $(b)$, we have
$$
\aligned
\ &\supp S\subseteq \Sigma_1 \ \ \ 
=\supp \widehat\sigma_1,\qquad \ \,
\supp S_i'\subseteq \Sigma_1' \ \ \ =\Fix \widehat\sigma_1,
\qquad \ 3\le i\le k-1,\\
&\supp \widetilde S\subseteq \Sigma_{k-1}
=\supp \widehat\sigma_{k-1},\quad \,\,
\supp \widetilde{S_i'}
\subseteq \Sigma_{k-1}'=\Fix \widehat\sigma_{k-1},
\quad \, 1\le i\le k-3.
\endaligned
$$
Taking into account the representation for $\widehat\sigma_1$
given by 
(\ref{7.12t}), 
we see that
$$
\supp S\subseteq \Sigma_1=\supp \widehat\sigma_1
=(\supp \widetilde S)\cup (\supp \widetilde {S_1'})
$$
and  
$$
\supp S_i'\subseteq \Sigma_1'=
{\boldsymbol\Delta}_n-\supp \widehat\sigma_1
={\boldsymbol\Delta}_n-((\supp \widetilde S)
\cup (\supp \widetilde {S_1'})),
\ \ \ 3\le i\le k-1.
$$
Hence 
%&
\begin{equation}\label{*}
(\supp S_i')\cap (\supp \widetilde S)=\varnothing\qquad
\text{for} \ \ 3\le i\le k-1.
\end{equation}
%&
Completely analogously, using the representation for $\widehat\sigma_{k-1}$
given by (\ref{7.12}), we get
%&
\begin{equation}\label{*t}
(\supp \widetilde {S_i'})\cap (\supp S)=\varnothing\qquad
\text{for} \ \ 1\le i\le k-3.
\end{equation}
%& 
There are at least two $i$'s such that 
$3\le i\le k-3$ (since $k>6$); for any such $i$,
formulae (\ref{7.12}), (\ref{7.12t})
provide us with the two representations of 
$\widehat\sigma_i$ in the form of disjoint products:
%&
\begin{equation}\label{**}
S S_i'=\widehat\sigma_i = \widetilde S \widetilde{S_i'}.
\end{equation}
%&
Obviously, (\ref{*}), (\ref{*t}), and (\ref{**}) imply $S=\widetilde S$ and 
$S_i'=\widetilde{S_i'}$ for $3\le i\le k-3$. 
In view of (\ref{7.12}) and (\ref{7.12t}), the set 
$Q=\supp S=\supp \widetilde S$ is invariant under all the
permutations $\widehat\sigma_1,...,\widehat\sigma_{k-1}$. Clearly,
$Q\ne {\boldsymbol\Delta}_n$ (for $\widehat\sigma_1$
has at least $k-2$ fixed points). 
If $\psi$ is transitive, the set $Q$ must be empty, and
we get $S=\widetilde S=1$. Hence $\widehat\sigma_i=S_i$ for $3\le i\le k-1$,
\ $\phi=\varphi'$, and $\widehat\sigma_1$ is disjoint
with $\widehat\sigma_3,...,\widehat\sigma_{k-1}$.
\end{proof}
%end of 2nd pt of sec 7

To get a base for induction, we study some homomorphisms
of $B_7$ and $B_8$.

\begin{Lemma}\label{Lm: B(7) to S(n) for 7<n<14} 
Let $7<n<14$. Then any transitive homomorphism \
$\psi\colon B_7\to{\mathbf S}(n)$ is cyclic.
\end{Lemma}

\begin{proof} 
Suppose to the contrary that $\psi$ is non-cyclic.
By Lemma \ref{Lm: homomorphisms Bk to S(n) for n<=2k and supp(hatsigma(1))<6},
we have $N=\#\Sigma_1=\#\supp\,\widehat\sigma_1\ge 6$.
Corollary \ref{Crl: k-2 fixed points of psi(sigma(1)) for 6<k<n<2k} shows that
$N'=\#\Sigma_1'=\#\Fix\,\widehat\sigma_1\ge 5$.
Hence $11\le N+N'=n\le 13$ and $5\le N'\le 7$.
\vskip0.2cm

By Lemma
\ref{Lm2: homomorphism with all components of hatsigma(1) of length<=k-3}$(a,c)$,
all permutations $\widehat\sigma_3,...,\widehat\sigma_6$
are supported in the set $\Sigma_1'$,
and the non-cyclic homomorphism $\phi$ (that is, the
restriction of $\psi$ to the subgroup
in $B_7$ generated by $\sigma_3,...,\sigma_6$)
coincides with its reduction
$$
\varphi'=\phi_{\Sigma_1'}\colon B_5\to
{\mathbf S}(\Sigma_1')\cong {\mathbf S}(N').
$$
%\vskip0.2cm

{\bfit Claim.} {\sl Every permutation $\widehat\sigma_i$, \
$1\le i\le 6$, is a product of disjoint transpositions.}
\vskip0.2cm

It suffices to prove this claim for $i\ge 3$;
let us deal with such $i$'s.
We consider the following three cases:
$N'=5$, \ $N'=6$ and $N'=7$. If $N'=5$, then
$\varphi'$ must be transitive (otherwise, all orbits are
of length $<5$ and $\phi=\varphi'$ is cyclic);
by Artin Theorem, any $\widehat\sigma_i$ is a transposition.
If $N'=7$, Claim follows from Proposition \ref{Prp: B(5) to S(7)}.
Let $N'=6$. If $\varphi'$ is intransitive,
then its image in ${\mathbf S}(\Sigma_1')\cong {\mathbf S}(6)$ has
exactly one orbit $Q$ of length $5$ and one fixed point;
Artin Theorem applies to the reduction of $\varphi'$
to $Q$, and we see that any $\widehat\sigma_i$ is a transposition.
Finally, if $\varphi'$ is transitive,
then, by Proposition \ref{Prp 4.9}, $\varphi'\sim\psi_{5,6}$
and every $\widehat\sigma_i$ is a product of $3$ disjoint transpositions.
\vskip0.2cm

To conclude the proof of the lemma, we use the approach that
was already used in the proof of Proposition \ref{Prp: B(6) to S(8)}.
Claim shows that $\widehat\sigma_i^2=1$ for all $i$ and hence
the non-abelian primitive group $G=\Img \psi\subset {\mathbf S}(n)$ \ 
($7<n<14$) is isomorphic to ${\mathbf S}(7)$; in view of Remark \ref{Rmk 7.5},
this contradicts Jordan Theorem.
\end{proof}

To treat homomorphisms $B_8\to{\mathbf S}(n)$, \
$8<n<16$, we need the following fact.

\begin{Proposition}\label{Prp 7.18} 
\noindent $a)$ Any transitive homomorphism
$\psi\colon B_6\to{\mathbf S}(9)$ is cyclic.

$b)$ Any non-cyclic homomorphism
$\psi\colon B_6\to{\mathbf S}(9)$ is conjugate
to a disjoint product $\psi_1\times\psi_2$,
where $\psi_1\colon B_6\to{\mathbf S}(6)$
is either $\mu_6$ or $\nu_6$, and
$\psi_2\colon B_6\to{\mathbf S}(3)$
is a cyclic homomorphism. In particular,
either every $\widehat\sigma_i$ is a disjoint product of
transpositions or every $\widehat\sigma_i$ is a disjoint product of
transpositions and a $3$-cycle that does not depend on $i$.
\end{Proposition}

\begin{proof} 
$a)$ Suppose to the contrary that $\psi$ is non-cyclic. By Proposition
\ref{Prp: transitive imprimitive homomorphisms Bk to S(n) is cyclic if n<2k}
$(a)$ and Remark \ref{Rmk 7.5}, \ $\widehat\sigma_1^2\ne 1$.
Using this and Lemmas
\ref{Lm: if |Inv(r,psi(sigma(1)))|=1 then psi is transitive},
\ref{Lm: Artin Lemma on cyclic decomposition of hatsigma(1)},
\ref{Lm: condition no cycles of the same length restricts support},
we can eliminate all cyclic types of $\widehat\sigma_1$ but the following
three:
\ $i)$ $[2,2,3]$; \ $ii)$ $[3,3]$; \ $iii)$  $[3,3,3]$. Let us consider
these cases.
\vskip0.2cm

$i)$ For $i\ge 3$, let $\widehat\sigma_i'$ be the restriction
of $\widehat\sigma_i$ to the support of the $3$-cycle $C$ in $\widehat\sigma_1$;
then, by
Lemma \ref{Lm: invariant sets and components of commuting permutations},
$\widehat\sigma_i=C^{q_i}$, \ $0\le q_i\le 2$.
Clearly, $q_i\ne 0$ (for $\widehat\sigma_i$ has only $2$ fixed points);
hence all $C^{q_i}$ are $3$-cycles with the same support
$\supp C$. Now, $C^{q_5}$ is the {\sl only} $3$-cycle in the
cyclic decomposition of $\widehat\sigma_5$, and $\widehat\sigma_2$ commutes
with $\widehat\sigma_5$. Hence the set $\supp C$ is $(\Img \psi)$-invariant,
which contradicts the transitivity of $\psi$.
\vskip0.3cm

$ii)$ In this case the $3$-component $\mathfrak C=\{C_0, C_1\}$
is the only non-degenerate component of $\widehat\sigma_1$ 
and we have the corresponding retraction
$\Omega=\Omega_{\mathfrak C}\colon B_4\to{\mathbf S}(2)$.
Clearly, either $\Omega$ is trivial or all
$\Omega(s_i)$ coincide with the transposition $(C_0,C_1)$. In any
case $\Omega(s_i^2)=1$, which means that
$\widehat\sigma_{i+2}^2 C_j \widehat\sigma_{i+2}^{-2}=C_j$ whenever
$j=0,1$ and $i=1,2,3$; thus,
$\widehat\sigma_{i+2}^2|{\supp \mathfrak C}=C_0^{q_{0,i}} C_1^{q_{0,i}}$
with some $q_{j,i}$, \ $0\le q_{j,i}\le 2$.
Because of $[\widehat\sigma_{i+2}]=[3,3]$, this implies
that for some $p_{j,i}$, \ $0\le p_{j,i}\le 2$,
the permutations $\widehat\sigma_{i+2}$ themselves satisfy
%&
\begin{equation}\label{7.15}
\widehat\sigma_{i+2}|{\supp \mathfrak C}=C_0^{p_{0,i}} C_1^{p_{1,i}},
\qquad i=1,2,3.
\end{equation}
%&
Since $\# ({\boldsymbol\Delta}_9-\supp \mathfrak C)=3$, the conditions
$[\widehat\sigma_{i+2}]=[3,3]$ and (\ref{7.15}) show that the permutations
$\widehat\sigma_{i+2}$, $i=1,2,3$, commute with each other, which is impossible.
\vskip0.3cm

$iii)$ In this case the only non-degenerate component of $\widehat\sigma_1$
is the $3$-component $\mathfrak C=\{C_0, C_1, C_2\}$,
and we have the retraction
$\Omega=\Omega_{\mathfrak C}\colon B_4\to{\mathbf S}(3)$.
We consider the following two cases: \ $iii_1)$ $\Omega$ is non-cyclic;
\ $iii_2)$ $\Omega$ is cyclic.
\vskip0.3cm

$iii_1)$ In this case, by Theorem \ref{Thm 3.14}, $\Omega\sim\mu_3\circ\pi$,
where $\pi\colon B_4\to B_3$
is the canonical epimorphism. This means that
all $\Omega(s_i)$ are transpositions; hence also
$\Omega(s_i^3)$ are transpositions. Therefore, there is a value $j$,
\ $j=0,1,2$, such that $\widehat\sigma_3^3 C_j \widehat\sigma_3^{-3}\ne C_j$.
However this contradicts the relation $\widehat\sigma_3^3=1$.
\vskip0.2cm

$iii_2)$ In this case all $\Omega(s_i)=A$, where $A\in{\mathbf S}(3)$
does not depend on $i$. If $A^2=1$, then
$\widehat\sigma_{i+2}^2 C_j \widehat\sigma_{i+2}^{-2}=C_j$ for
$i=1,2,3$, \ $j=0,1,2$; combined with the condition
$[\widehat\sigma_{i+2}]=[3,3,3]$, this shows that the permutations
$\widehat\sigma_{i+2}$, $i=1,2,3$, commute with each other, which
is impossible. Finally, if $A^2\ne 1$, then $A$ is a $3$-cycle,
and we may assume that
%&
\begin{equation}\label{7.16}
\widehat\sigma_i C_j \widehat\sigma_i^{-1}=C_{|j+1|_3},
\qquad i=1,2,3, \ \ j=0,1,2,
\end{equation}
%&
where $|N|_3\in{\mathbb Z}/3{\mathbb Z}$ denotes the $3$-residue of $N\in\mathbb N$. Let
$C_j=(a_j^0,a_j^1,a_j^2)$, \ $j=0,1,2$.
It follows from (\ref{7.16}) that there exist $t(j,i)\in{\mathbb Z}/3{\mathbb Z}$
such that
%&
\begin{equation}\label{7.17}
\widehat\sigma_i (c_j^k)=c_{j+1}^{k+t(j,i)},
\ \ \ j,k\in{\mathbb Z}/3{\mathbb Z}, \ \ i=3,4,5.
\end{equation}
%&

\noindent The condition $\widehat\sigma_i^3=1$ implies that
%&
\begin{equation}\label{7.18}
t(0,i)+t(1,i)+t(2,i)=0,\qquad i=3,4,5
\end{equation}
%&
(here and below all the equalities are in ${\mathbb Z}/3{\mathbb Z}$).
Using $\widehat\sigma_3\infty\widehat\sigma_4$, we obtain that
%&
\begin{equation}\label{7.19}
\aligned
\ &t(0,3)+t(1,4)+t(2,3)=t(0,4)+t(1,3)+t(2,4),\\
&t(0,3)+t(2,4)+t(1,3)=t(0,4)+t(2,3)+t(1,4).
\endaligned
\end{equation}
%&
Relations (\ref{7.18}), (\ref{7.19}) show that $t(j,3)=t(j,4)$
for all $j=0,1,2$; it follows that
$$
\widehat\sigma_3 (c_j^k)=c_{j+1}^{k+t(j,3)}
=c_{j+1}^{k+t(j,4)}=\widehat\sigma_4 (c_j^k)
$$
and $\widehat\sigma_3=\widehat\sigma_4$.
This contradiction concludes the proof of the statement $(a)$.
The proof of $(b)$ follows immediately from $(a)$, 
Theorem \ref{Thm: homomorphisms B(k) to S(k+1)}
and Theorem \ref{Thm: homomorphisms B(k) to S(k+2)}.
\end{proof}

\begin{Remark}\label{Rmk 7.6} 
By Theorem \ref{Thm: homomorphisms B(k) to S(k+1)},
Proposition \ref{Prp: B(6) to S(8)} and Proposition \ref{Prp 7.18},
any transitive homomorphism
$B_6\to{\mathbf S}(n)$ is cyclic whenever $7\le n\le 9$.
However, there is a non-cyclic transitive
homomorphism $B_6\to{\mathbf S}(10)$. To see this,
consider all $10$ partitions of ${\boldsymbol\Delta}_6$ 
into two (disjoint) subsets
consisting of $3$ points. The group ${\mathbf S}(6)$ acts transitively
on the family $\mathfrak P\cong {\boldsymbol\Delta}_{10}$ 
of all these partitions;
this action defines the transitive homomorphism
${\mathbf S}(6)\to{\mathbf S}(10)$; the composition of the canonical
projection $\mu_6$ with this homomorphism is a non-cyclic transitive
homomorphism $B_6\to{\mathbf S}(10)$. Under suitable
notations, this homomorphism looks as follows:
$$
\aligned
\widehat\sigma_1=(1,2)&(3,4)(5,6); \ \ \widehat\sigma_2=(1,7)(3,8)(5,9);
\ \ \widehat\sigma_3=(3,6)(4,5)(7,10);\\
&\widehat\sigma_4=(1,3)(2,4)(7,8);\qquad
\widehat\sigma_5=(3,5)(4,6)(8,9).
\endaligned
$$
Instead of the canonical epimorphism $\mu_6$, one could use Artin's
homomorphism $\nu_6$.
\end{Remark}
\vskip0.3cm

\begin{Lemma}\label{Lm: B(8) to S(n), 8<n<16}  
Let $8<n<16$. Then any transitive homomorphism
$\psi\colon B_8\to{\mathbf S}(n)$ is cyclic.
\end{Lemma}

\begin{proof} 
Suppose to the contrary that $\psi$ is non-cyclic.
By Lemma \ref{Lm: homomorphisms Bk to S(n) for n<=2k and supp(hatsigma(1))<6},
$N=\#\Sigma_1=\#\supp\,\widehat\sigma_1\ge 6$.
Corollary \ref{Crl: k-2 fixed points of psi(sigma(1)) for 6<k<n<2k}
shows that $N'=\Sigma_1'=\#\Fix\,\widehat\sigma_1\ge 6$.
Hence $12\le N+N'=n\le 15$ and $6\le N'\le 9$.

By Lemma
\ref{Lm2: homomorphism with all components of hatsigma(1) of length<=k-3}$(a,c)$,
all the permutations $\widehat\sigma_3,...,\widehat\sigma_7$
are supported in $\Sigma_1'$, and the restriction 
$$
\phi=\varphi'\colon B_6\to
{\mathbf S}(\Sigma_1')\cong {\mathbf S}(N')
$$
of $\psi$ to the subgroup in $B_8$ 
generated by $\sigma_3,...,\sigma_7$ is a non-cyclic homomorphism.
As usual, we set
$H=\Img \phi\subseteq {\mathbf S}(\Sigma_1')\cong {\mathbf S}(N')$.
\vskip0.3cm

\noindent{\bfit Claim.} {\sl There is exactly one $H$-orbit
$Q\subseteq \Sigma_1'$ of length $6$. The reduction
$$
\phi_Q\colon B_6\to{\mathbf S}(Q)\cong {\mathbf S}(6)
$$
is conjugate to one of the homomorphisms $\mu_6$, $\nu_6$.
The complement $Q'=\Sigma'\setminus Q$ contains at most $3$
points, and there is a permutation $A\in{\mathbf S}(Q')$
such that every $\widehat\sigma_i$ \ $(i=3,...,7)$ is a disjoint
product of some transpositions and this permutation $A$.}
%\vskip0.3cm

Indeed, since $\#\Sigma_1'=N'\le 9$ and the homomorphism
$\phi$ is non-cyclic, theorems \ref{Thm: homomorphisms B(k) to S(n), n<k}$(a)$,
\ref{Thm: homomorphisms B(k) to S(k+1)},
\ref{Thm: homomorphisms B(k) to S(k+2)} and Proposition \ref{Prp 7.18}
show that there is only one $H$-orbit $Q$ of length $6$.
In view of Artin Theorem and
Theorem \ref{Thm: homomorphisms B(k) to S(n), n<k}$(a)$, other
statements of Claim follow immediately from this fact.
\vskip0.3cm

We have the following cases:
\ $i)$ $N'= 6$; \ $ii)$ $N'= 7$; \
$iii)$ $N'= 8$; \ $iv)$ $N'= 9$.
\vskip0.3cm

Let us show that in all these cases the primitive group
$G=\Img \psi\subseteq {\mathbf S}(n)$ is isomorphic to
${\mathbf S}(8)$ and, besides, contains a $3$-cycle; 
this will contradict Jordan Theorem. 

$i)$ In this case $\Sigma_1'=Q$; hence either
$\phi=\phi_Q\sim\mu_6$ or to $\phi=\phi_Q\sim\nu_6$. Therefore,
$\widehat\sigma_i^2=1$ for all $i$ and $G\cong {\mathbf S}(8)$.

If $\phi\sim \mu_6$, 
then any $\widehat\sigma_i$ is a transposition,
and the product $(\widehat\sigma_3 \widehat\sigma_4)^2$ is a $3$-cycle in $G$
(in fact, the element $\widehat\sigma_3 \widehat\sigma_4$ 
itself is a $3$-cycle;
we take its square only to unify the proofs for all cases
$(i)-(iv)$).

If $\phi\sim \nu_6$,
then the permutation $(\widehat\sigma_3 \widehat\sigma_4\cdots \widehat\sigma_7)^2$
is a $3$-cycle in $G$ (here the square is essential, since 
$\widehat\sigma_3 \widehat\sigma_4\cdots \widehat\sigma_7$
is of cyclic type $[3,2]$).

$ii)$ In this case $Q'$ consists of one point
that is a fixed point of $H$. Applying the same arguments
as in case $(i)$, we obtain the desired result.

$iii)$ The only difference with the previous cases is that
all the permutations $\widehat\sigma_i$, \ $i\ge 3$, may contain one
additional disjoint transposition $A$.
However, this does not change anything (the square kills this transposition).

$iv)$ Here $\# Q'=3$. Hence either
$A=1$ or $[A]=[2]$ or, finally, $[A]=[3]$.
In the first two cases we follow the same arguments as above.
Let us show that the third case cannot occur.
Indeed, $N'=9$ and $N\ge 6$; thus, $N=6$ 
(for $N+N'=n\le 15$). That is, the support of
any permutation $\widehat\sigma_i$ consists of $6$ points. If $A$
is a $3$-cycle, it takes $3$ points from the $6$, and the rest
three places cannot be filled by transpositions.
This concludes the proof.
\end{proof}

Now we are ready to prove Theorem F$(a)$. Actually, the proof
is simple, since the main work was already done.

\begin{Theorem}\label{Thm: homomorphisms B(k) to S(n), 6<k<n<2k}
Let $6<k<n<2k$. Then:
\vskip0.2cm

\noindent $a)$ Any transitive homomorphism
$B_k\to{\mathbf S}(n)$ is cyclic.
\vskip0.2cm

\noindent $b)$ Any non-cyclic homomorphism
$\psi\colon B_k\to{\mathbf S}(n)$ is
conjugate to a homomorphism of the form $\mu_k\times\widetilde\psi$,
where $\widetilde\psi\colon B_k\to{\mathbf S}(n-k)$
is a cyclic homomorphism.
\end{Theorem}

\begin{proof} 
$a)$ Let us call $\Gamma(m)$ the following conjecture: 
\vskip0.3cm

\noindent{\it Conjecture $\Gamma(m)$.} 
{\sl Every transitive homomorphism
$\psi\colon B_k\to{\mathbf S}(n)$ is cyclic whenever
$6<k\le m$ and $6<k<n<2k$.}
\vskip0.3cm

We have already proved $\Gamma(m)$ for
$m=7$ and $m=8$ (Lemma \ref{Lm: B(7) to S(n) for 7<n<14},
Lemma \ref{Lm: B(8) to S(n), 8<n<16}).
Suppose that $\Gamma(m)$ is fulfilled for some $m\ge 7$.
We shall show that then $\Gamma(m+2)$ is fulfilled.
By Induction Principle, this will prove the statement $(a)$.

Suppose to the contrary that $\Gamma(m+2)$
is wrong. That is, for some $k$ and $n$ that satisfy
$k\le m+2$ and $6<k<n<2k$ there exists
a transitive non-cyclic homomorphism
$\psi\colon B_k\to{\mathbf S}(n)$.
It follows from Lemma \ref{Lm: B(7) to S(n) for 7<n<14}
and Lemma \ref{Lm: B(8) to S(n), 8<n<16}
that $k>8$; hence $6<k-2\le m$ and $\Gamma(k-2)$ is fulfilled.

By Lemmas \ref{Lm: homomorphisms Bk to S(n) for n<=2k and supp(hatsigma(1))<6}
and Corollary \ref{Crl: k-2 fixed points of psi(sigma(1)) for 6<k<n<2k}, we have
$N\ge 6$ and $N'\ge k-2$; thus
%&
\begin{equation}\label{7.20}
6<k-2\le N' \le n-6<2k-6<2(k-2).
\end{equation}
%&
By Lemma
\ref{Lm2: homomorphism with all components of hatsigma(1) of length<=k-3}$(a,c)$,
the restriction $\phi\colon B_{k-2}\to{\mathbf S}(n)$
of $\psi$ to the subgroup $B_{k-2}\subset B_k$
generated by $\sigma_3,...,\sigma_{k-1}$
coincides with its non-cyclic reduction
$$
\varphi'=\phi_{\Sigma_1'}
\colon B_{k-2}\to
{\mathbf S}(\Sigma_1')\cong {\mathbf S}(N')
$$
to the $H$-invariant set $\Sigma_1'=\Fix \widehat\sigma_1$.
To conclude the proof of the statement $(a)$, it is sufficient
to prove the following
\vskip0.3cm

\noindent{\bfit Claim.} $\widehat\sigma_i^2=1$ for all $i$
and hence
$G=\Img \psi\cong {\mathbf S}(k)$. Moreover, the primitive permutation group
$G\subseteq {\mathbf S}(n)$ contains a $3$-cycle.
\vskip0.3cm

Indeed, as we know, these properties are incompatible; hence
the assumption that $\Gamma(m+2)$ is wrong leads to a contradiction.

To justify Claim, assume first that $N' = k-2$.
Since $\varphi'$ is non-cyclic,
Theorem \ref{Thm: homomorphisms B(k) to S(n), n<k}$(a)$ shows
that $\varphi'$ is transitive; by Artin Theorem,
$\phi=\varphi'\sim \mu_{k-2}$ (for $k-2>6$).
Hence any $\widehat\sigma_i$ is a transposition
and $\widehat\sigma_1\widehat\sigma_2$ is a $3$-cycle containing in $G$.

Assume now that $N' > k-2$.
Since $\Gamma(k-2)$ is fulfilled, any transitive homomorphism
$B_{k-2}\to{\mathbf S}(N')$
is cyclic; therefore, {\sl the homomorphism $\phi=\varphi'$ must be
intransitive.}
The reduction $\phi_Q$ of the non-cyclic intransitive
homomorphism $\phi=\varphi'$ to any $H$-orbit
$Q\subset \Sigma_1'$ is a transitive homomorphism
$B_{k-2}\to{\mathbf S}(Q)$.
Clearly, $\# Q<N'<2(k-2)$. If $\# Q\ne k-2$, then
$\phi_Q$ is cyclic (this follows from
Theorem \ref{Thm: homomorphisms B(k) to S(n), n<k}$(a)$
whenever $\# Q < k-2$; if $\# Q > k-2$, then $\phi_Q$
must be cyclic by our assumption that $\Gamma(k-2)$ is fulfilled).
Hence there exists a unique $H$-orbit $Q$ of length $k-2$, and
the reduction $\phi_Q$ of $\phi$ to this orbit
is non-cyclic and transitive. Since $k-2>6$, Artin Theorem shows
that $\phi_Q\sim \mu_{k-2}$.
Let $Q' = \Sigma' - Q$; clearly,
$\phi$ is the disjoint product of the reductions
$\phi_Q$ and $\phi_{Q'}$, and $\phi_{Q'}$ is cyclic.
This means that there is
a permutation $A\in{\mathbf S}(Q')$ such that
for every $i$, \ $3\le i\le k-1$, the permutation
$\widehat\sigma_i$ is the disjoint product of $A$ and the transposition
$A_i=\phi_Q(\sigma_i)$.

Let us show that $A^2=1$. Indeed, if this is not the case,
then for some $r>2$ the cyclic decomposition of $A$ contains an $r$-cycle.
Let $\mathfrak C_r(A)$ be the $r$-component of $A$.
Since any $\widehat\sigma_i$, \ $3\le i\le k-1$, is the disjoint
product of $A$ and the transposition $A_i$,
we obtain that $\mathfrak C_r(A)$ is, actually,
the $r$-component of every $\widehat\sigma_i$, \ $3\le i\le k-1$.
Thereby, the support $\Sigma_{\mathfrak C_r(A)}$ of this component
$\mathfrak C_r(A)$ is invariant under all the permutations
$\widehat\sigma_3,...,\widehat\sigma_{k-1}$. Moreover, the permutations
$\widehat\sigma_1$, $\widehat\sigma_2$
commute with $\widehat\sigma_{k-1}$
and hence the set $\Sigma_{\mathfrak C_r(A)}$ is
invariant under $\widehat\sigma_1$ and $\widehat\sigma_2$. However
this contradicts the transitivity of $\psi$.
%\vfill\pagebreak
\vskip0.3cm

Since $\widehat\sigma_i=A_i A$ for $i\ge 3$ and all $A_i$ are transpositions,
the property $A^2=1$ implies that $\widehat\sigma_i^2=1$ for
$i\ge 3$ and hence for all $i$. Moreover,
$\widehat\sigma_3 \widehat\sigma_4 = A_3 A A_4 A = A_3 A_4$
is a $3$-cycle containing in $G$. This concludes the proof of Claim
and proves the statement $(a)$.
In view of Theorem \ref{Thm: homomorphisms B(k) to S(n), n<k}$(a)$
and Artin Theorem, $(a)$ implies $(b)$.
\end{proof}
 %%&&

\begin{Remark}\label{Rmk 7.7} 
There is a non-cyclic transitive homomorphism
$B_6\to{\mathbf S}(10)$ (see Remark \ref{Rmk 7.6}).
On the other hand, for any $k\ge 3$, Corollary \ref{Crl 6.6} 
provides us with the four non-cyclic homomorphisms
$\varphi_j\colon B_k\to{\mathbf S}(2k)$, \ $j=0,1,2,3$.
The homomorphism $\varphi_0$ is intransitive, and
the homomorphisms $\varphi_j$, \ $j=1,2,3$, are transitive 
(see also Remark \ref{Rmk 6.6}).
These remarks show that the conditions $6<k<n<2k$ 
of Theorem \ref{Thm: homomorphisms B(k) to S(n), 6<k<n<2k}$(a)$ are,
in a sense, sharp.
\end{Remark}
\vskip0.3cm

\subsection{Homomorphisms $B_k\to{\mathbf S}(2k)$}
\label{Ss: 7.4. Homomorphisms B(k) to S(2k)} %@
Here we prove Theorem F$(b)$. In the following lemma 
(which is similar to Lemma \ref{Lm: no components of length t>k-3}) 
we make use of the model homomorphisms
$\varphi_j\colon B_k\to{\mathbf S}(2k)$ exhibited 
in Definition \ref{Def: model and standard homomorphisms B(k) to S(2k)}
(see also Corollary \ref{Crl 6.6}, with $k$ instead of $n$).

\begin{Lemma}
\label{Lm: psi B(k) to S(2k) with psi(sigma(1)) having component of length>k-3} 
Let $k>6$ and let 
$\psi\colon B_k\to{\mathbf S}(2k)$ be a non-cyclic homomorphism
such that the permutation $\widehat\sigma_1$ has
a non-degenerate component $\mathfrak C$ of length $t>k-3$. 
Then either $t=k-2$ and $\psi\sim\varphi_3$ or
$t=k$ and $\psi\sim\varphi_2$.
\end{Lemma}

\begin{proof} 
We follow the proof of Lemma \ref{Lm: no components of length t>k-3}.
Clearly, $\mathfrak C$
must be the $2$-component of $\widehat\sigma_1$, and $t\le k$;
hence either $t=k-2$ or $t=k-1$ or $t=k$. 

Set 
%&
\begin{equation}\label{7.21}
\aligned
\Sigma&=\supp \mathfrak C, \ \ \ \Sigma'
= {\boldsymbol\Delta}_{2k} - \Sigma,\qquad
Q=\supp \widehat\sigma_1, \ \ \ Q'={\boldsymbol\Delta}_{2k}-Q,\\
\# \Sigma &= 2t, \ \qquad \# \Sigma' = 2k-2t\le 2k - 2(k-2) = 4.
\endaligned
\end{equation}
%&
Since $k>6$, any homomorphism 
$B_{k-2}\to{\mathbf S}(\Sigma')$
is cyclic (Theorem \ref{Thm: homomorphisms B(k) to S(n), n<k}$(a)$).
In particular, the homomorphism $\varphi_{_{\Sigma'}}$ is cyclic and
Lemma
\ref{Lm: if psi: B_k to S(n) is non-cyclic and varphi(Sigma') is abelian then phi, varphi(Sigma) and Omega are non-abelian} implies
that {\sl the homomorphisms 
$\Omega\colon B_{k-2}\to{\mathbf S}(\mathfrak C)\cong{\mathbf S}(t)$
and $\varphi_{_{\Sigma}}\colon B_{k-2}\to G\subset {\mathbf S}(2t)$
must by non-cyclic}.

We may assume that the homomorphism $\psi$ is normalized;
this means that
$$
\Sigma = \{1,2,\ldots,2t\}, \ \
\widehat\sigma_1|\Sigma = C_1\cdots C_t, \ \text{where} \
C_m = (2m-1,2m) \ \text{for} \ m = 1,\ldots ,t.
$$
We must consider the following five cases:
\vskip0.3cm

\begin{itemize} 

\item[$i)$] $t=k-2$ and $\widehat\sigma_1$
is a disjoint product of $k-2$ transpositions;

\item[$ii)$] $t=k-2$ and $\widehat\sigma_1$
is a disjoint product of $k-2$ transpositions and a $3$-cycle;

\item[$ii')$] $t=k-2$ and $\widehat\sigma_1$
is a disjoint product of $k-2$ transpositions and a $4$-cycle $F_1$;

\item[$iii)$] $t=k-1$ and $\widehat\sigma_1$
is a disjoint product of $k-1$ transpositions;

\item[$iv)$] $t=k$ and $\widehat\sigma_1$
is a disjoint product of $k$ transpositions.
\end{itemize}
\vskip0.3cm

First we prove that cases $(i)$, $(ii)$, $(iii)$ are
impossible.
\vskip0.2cm

Case $(i)$ may be eliminated by the same argument
that were used in the proof of Lemma \ref{Lm: no components of length t>k-3}
(the only difference is that now $\Sigma'$ consists of four points;
actually, this does not change anything).
\vskip0.2cm

$ii)$ In this case $\widehat\sigma_1$ has exactly one fixed point.
Clearly, this point is also the only fixed point of any $\widehat\sigma_i$
(see, for instance,
Lemma \ref{Lm: if |Inv(r,psi(sigma(1)))|=1 then psi is transitive}).
Hence $\psi$ is the disjoint product $\psi_Q\times{\mathbf 1}_1$, where
$\psi_Q\colon B_k\to{\mathbf S}(Q)\cong{\mathbf S}(2k-1)$
is the reduction of $\psi$ to the $(\Img \psi)$-invariant
set $Q=\supp \widehat\sigma_1$. Since $\psi$ is non-cyclic,
$\psi_Q$ is non-cyclic too, and the permutation
$\psi_Q(\sigma_1)=\widehat\sigma_1$ has a $2$-component of length $k-2$.
However this contradicts Lemma \ref{Lm: no components of length t>k-3}.
\vskip0.2cm

$iii)$ In this case $\Sigma=Q$, \ $\Sigma'=Q'$, and we deal
with the non-cyclic homomorphisms
$\varphi_{_\Sigma}\colon B_{k-2}\to G\subset{\mathbf S}(2k-2)$
and $\Omega\colon B_{k-2}\to{\mathbf S}(k-1)$.
By Remark \ref{Rmk 5.1}, Proposition \ref{Prp 4.9}, 
and Theorem \ref{Thm: homomorphisms B(k) to S(k+1)},
we must consider the following two cases:
\ $iii_a)$ $\Omega=\mu_{k-2}\times\mathbf 1_1$;
$iii_b)$ $k=7$ and
$\Omega=\psi_{5,6}\colon B_5\to{\mathbf S}(6)$.
 %%&&

In case $(iii_a)$, as in Lemma \ref{Lm: no components of length t>k-3}, 
$\varphi_{_\Sigma}$ must be conjugate
to one of the eight homomorphisms $\psi_{x;j}$ listed in
Corollary \ref{Crl 6.9} (with $n=k-2$).
All these cases may be eliminated in the same way as in
Lemma \ref{Lm: no components of length t>k-3}
(the only difference is that now the set $\Sigma'$ 
consists of two points, which does not change anything).
Case $(iii_b)$ is impossible by the same reasons as in
Lemma \ref{Lm: no components of length t>k-3}.
\vskip0.3cm

Now we must handle cases $(ii')$ and $(iv)$.
\vskip0.2cm

$ii')$ We prove that in this case $\psi\sim\varphi_3$.
We may assume that
$$
\Sigma=\{1,...,2k-4\},\qquad
\Sigma'=\{2k-3,2k-2,2k-1,2k\}.
$$
We deal with the non-cyclic homomorphisms
$$
\varphi_{_\Sigma}\colon B_{k-2}\to G\subset{\mathbf S}(2k-4),
\qquad
\Omega\colon B_{k-2}\to{\mathbf S}(k-2).
$$
Clearly, either \ $ii_a')$ $\Omega\sim\mu_{k-2}$ or
$ii_b')$ $k=8$ and $\Omega\sim\nu_6$.
\vskip0.2cm

Case $(ii_b')$ is actually impossible; this may be
proven by the argument used in
Lemma \ref{Lm: no components of length t>k-3}
in case $(i_b)$
(the only difference is that now the cyclic decomposition
of $\widehat\sigma_1$ contains the additional $4$-cycle
$F_1\in{\mathbf S}(\{13,14,15,16\})$,
and $\widehat\sigma_3,...,\widehat\sigma_7$ contain
the additional $4$-cycle $F=F_1^{\pm 1}$;
however, this does not change anything).
\vskip0.2cm

$ii_a')$ We may assume that $\Omega=\mu_{k-2}$.
Then, by Corollary \ref{Crl 6.6}, $\varphi_{_\Sigma}$ is conjugate
to one of the homomorphisms $\varphi_j$, \ $j=0,1,2,3$ (with $n=k-2$).
For $j=0,1$ we have
$\widehat\sigma_{i+2}|_\Sigma=\varphi_{_\Sigma}(s_i)=\varphi_j(s_i)$;
hence $\# \supp (\widehat\sigma_{i+2}|\Sigma)=4$
and there is no room in $\Sigma'$ for the
rest $2k-4$ points of $\supp \widehat\sigma_{i+2}$. For $j=2$ we have
$\widehat\sigma_{i+2}|\Sigma=\varphi_2(s_i)$, and
hence $\widehat\sigma_{i+2}=\varphi_2(s_i)F$ for all $i\ge 1$, where
$F\in{\mathbf S}(\{2k-3,2k-2,2k-1,2k\})$ is a $4$-cycle;
the argument used in Lemma \ref{Lm: no components of length t>k-3}
for case $(i_a)$
show that this is impossible.
So, we are left with the case $j=3$, i. e.
$\varphi_{_\Sigma}=\varphi_3\colon B_{k-2}\to{\mathbf S}(2k-4)$.
Without loss of generality, we may assume that
$F_1=(2k-3,2k,2k-2,2k-1)$;
to simplify notation, put $a=2k-3$, $b=2k-2$, $c=2k-1$, $d=2k$.
Since $\psi$ is normalized and
$\widehat\sigma_{i+2}|\Sigma=\varphi_3(s_i)$, we have
%&
\begin{equation}\label{7.22}
\widehat\sigma_1=(1,2)(3,4)
\cdots 
(2k-7,2k-6)(2k-5,2k-4)\underbrace{(a,d,b,c)}_{4\text{-cycle}}
\end{equation}
%&
and
%&
\begin{equation}\label{7.23}
\aligned
\widehat\sigma_i=(1,2)(3,4)&\cdots (2i-7,2i-6)\,
\underbrace{(2i-5,2i-2,2i-4,2i-3)}_{4\text{-cycle}}\,(2i-1,2i)\\
&\qquad\qquad\qquad\qquad\qquad
\times\cdots\times (2k-5,2k-4)(a,b)(c,d)
\endaligned
\end{equation}
%&
for $3\le i\le k-1$. To recover the homomorphism $\psi$, we
need to compute $\widehat\sigma_2$. This permutation commutes
with $\widehat\sigma_i$, \ $4\le i\le k-1$; the cyclic decomposition
of $\widehat\sigma_i$ contains only one $4$-cycle, namely,
$F_i=(2i-5,2i-2,2i-4,2i-3)$; therefore, each of the sets
$\{2i-5,2i-2,2i-4,2i-3\}$,\ $4\le i\le k-1$, must be
$\widehat\sigma_2$-invariant. It follows that each of the sets
$\{3,4\},\{5,6\},...,\{2k-5,2k-4\}$ is $\widehat\sigma_2$-invariant.
Since the permutation $\widehat\sigma_2\sim\widehat\sigma_1$
has no fixed points, its cyclic decomposition contains the
product
$$
A=(3,4)(5,6)\cdots (2k-5,2k-4)
$$
of $k-3$ disjoint transpositions;
it must also contain one more transposition
$T$ and some $4$-cycle $F_2$. By Lemma \ref{Lm: Omega=Omega*}$(a)$,
$\Omega^*=\Omega=\mu_{k-2}$; it follows that exactly $k-4$
from the $k-2$ transpositions
$(1,2)$, $(3,4)$,..., $(2k-9,2k-8)$, $(a,b)$, $(c,d)$ that occur
in the cyclic decomposition of $\widehat\sigma_{k-1}$
must be $\Omega^*(s_2)$-invariant (that is, invariant
under conjugation by $\widehat\sigma_2$), and the rest
two transpositions must mutually interchange. Evidently,
the $k-3$ transpositions $(3,4),(5,6),\ldots ,(2k-5,2k-4)$
are the fixed points of $\Omega^*(s_2)$; hence,
exactly one of the transpositions $(1,2)$, $(a,b)$, $(c,d)$
must be a fixed point of $\Omega^*(s_2)$; denote
this transposition by $T$.
\vskip0.2cm

Let us show that $T\ne (1,2)$.
Indeed, if $T=(1,2)$, then the cyclic decomposition of $\widehat\sigma_2$
contains the product $P=(1,2)A$ and a $4$-cycle
$F_2$ supported on $\{a,b,c,d\}$. Since $\widehat\sigma_{k-1}$ commutes
with $\widehat\sigma_2$, and the product $(a,b)(c,d)$ is contained in
$\widehat\sigma_{k-1}$ (see (\ref{7.22})),
we have $(a,b)(c,d)=F_2^2$; thereby $F_2=(a,c,b,d)^{\pm 1}$.
It is easy to check that in this case $\widehat\sigma_2$,
$\widehat\sigma_1$ cannot be a braid-like couple.

So, either $T=(a,b)$ or $T=(c,d)$. If $T=(a,b)$, then
either $F_2=(1,c,2,d)$ or $F_2=(1,d,2,c)$; however,
conjugation by $(a,b)(c,d)$ does not change
$\widehat\sigma_1,\widehat\sigma_3,...,\widehat\sigma_{k-1}$,
$T$ and transforms
$(1,d,2,c)$ into $(1,c,2,d)$. Similarly, if $T=(c,d)$, then
either $F_2=(1,a,2,b)$ or $F_2=(1,b,2,a)$, and the same
conjugation transforms $(1,b,2,a)$ into $(1,a,2,b)$.
Moreover, conjugation by $(a,d,b,c)$
does not change $\widehat\sigma_1,\widehat\sigma_3,...,\widehat\sigma_{k-1}$
and transforms $(c,d)$ into $(a,b)$ and $(1,a,2,b)$ into
$(1,d,2,c)$. Hence, without loss of generality, we may assume
that $T=(a,b)$ and $F_2=(1,c,2,d)$, and thus
$$
\widehat\sigma_2=(1,c,2,d)(3,4)(5,6)\cdots (2k-5,2k-4)(a,b).
$$
Finally, we conjugate the original homomorphism $\psi$ by the permutation
\begin{equation}\label{7.24}
B=\begin{pmatrix}
1 & 2 & 3 &... & 2k-4 & a & b & c & d \\
5 & 6 & 7 &... & 2k   & 1 & 2 & 3 & 4
\end{pmatrix}
\end{equation}
%&
and obtain the homomorphism
$\widetilde\psi\colon B_k\to{\mathbf S}(2k)$,
$$
\widetilde\psi(\sigma_i)=B\psi(\sigma_i)B^{-1}=B\widehat\sigma_iB^{-1},
\qquad 1\le i\le k-1,
$$
that coincides with $\varphi_3$. This concludes the proof
in case $(ii')$.
\vskip0.3cm

$iv)$ We prove that in this case $\psi\sim\varphi_2$.
We deal with non-cyclic homomorphisms
$\Omega\colon B_{k-2}\to{\mathbf S}(k)$ and
$\phi\colon B_{k-2}\to{\mathbf S}(2k)$ (we write $\phi$ instead %%@@
of $\varphi_{_\Sigma}$, for $\Sigma={\boldsymbol\Delta}_{2k}$). According
to Theorem \ref{Thm: homomorphisms B(k) to S(k+2)},
we must consider the following cases:
\vskip0.3cm

$iv_a)$ $k=7$ and $\Omega=\psi_{5,6}\times \mathbf 1_1$;

$iv_b)$ $k=8$ and $\Omega=\nu_6\times \omega$;

$iv_c)$ $\Omega=\mu_{k-2}\times \omega$;
\vskip0.2cm

\noindent in cases $(iv_b)$, $(iv_c)$, \ 
$\omega\colon B_6\to{\mathbf S}(2)$ is some (cyclic)
homomorphism.
\vskip0.3cm

\noindent With help of Corollary \ref{Crl 6.7}, Corollary \ref{Crl 6.8}
and Theorem
\ref{Thm: if Omega is cyclic then every Omega-homomorphism is cyclic}
we show below that $(iv_a)$ and $(iv_b)$ cannot actually occur.

Since $[\widehat\sigma_i]=[\underbrace{2,...,2}_{k\ \text{times}}]$,
in case $(iv_a)$ we may assume
that $\phi$ coincides with one of the four homomorphisms 
$\eta_j\times\bold{[2]}\colon B_5\to
{\mathbf S}(12)\times{\mathbf S}(2)\subset{\mathbf S}(14)$, \ $j=0,1,2,3$
(see Corollary \ref{Crl 6.7}), where
$\bold{[2]}\colon B_5\to{\mathbf S}(2)$
is the cyclic homomorphism sending each generator $s_i$
to the transposition $(13,14)$.

For $j=1,3$ the cyclic decompositions of $\eta_j(s_i)$
contain $4$-cycles, which cannot occur here.

For $j=0,2$ the permutation
$A=\widehat\sigma_4\widehat\sigma_5$ is a product of four disjoint $3$-cycles
supported on ${\boldsymbol\Delta}_{12}=\{1,...,12\}$.
Since $\widehat\sigma_2$
commutes with $A$, ${\boldsymbol\Delta}_{12}$ is $\widehat\sigma_2$-invariant.
Hence ${\boldsymbol\Delta}_{12}$ is $(\Img \psi)$-invariant, and the reduction
$\psi_{{\boldsymbol\Delta}_{12}}\colon B_7\to{\mathbf S}(12)$
is a non-cyclic homomorphism; it follows from
Theorem \ref{Thm: homomorphisms B(k) to S(n), 6<k<n<2k}
that $\Img \psi_{{\boldsymbol\Delta}_{12}}$ must have an invariant set
$E\subset{\boldsymbol\Delta}_{12}$ of cardinality $12-7=5$. In particular,
$E$ must be invariant under all the permutations
$\widehat\sigma_{i+2}|{{\boldsymbol\Delta}_{12}}=\eta_j(s_i)$, \ $1\le i\le 4$,
which is not the case.
\vskip0.2cm

\noindent In case $(iv_b)$ we may assume that $\phi$ is of the form
$$
\phi=\phi_j\times\varphi_\omega\colon
\,B_6\to{\mathbf S}(12)\times{\mathbf S}(4)\subset{\mathbf S}(16)
$$
(see Corollary \ref{Crl 6.8}). Here
$\varphi_\omega\colon B_6\to{\mathbf S}(4)$
is a cyclic homomorphism defined by the following conditions:
{\em if the homomorphism
$\omega\colon B_6\to{\mathbf S}(2)$ is trivial,
then $\varphi_\omega(s_i)=(13,14)(15,16)$ \ $(i=1,...,5)$;
otherwise, i. e., for $\omega$ sending all $s_i$ to $(1,2)$, \
$\varphi_\omega(s_i)=(13,15)(14,16)$ \ $(i=1,...,5)$}.
As in case $(iv_a)$, the permutation $A=\widehat\sigma_4\widehat\sigma_5$
is a product of four disjoint $3$-cycles supported on
${\boldsymbol\Delta}_{12}$,
the set ${\boldsymbol\Delta}_{12}$ is $(\Img \psi)$-invariant, and the reduction
$\psi_{{\boldsymbol\Delta}_{12}}\colon B_8\to{\mathbf S}(12)$
is a non-cyclic homomorphism. It follows from
Theorem \ref{Thm: homomorphisms B(k) to S(n), 6<k<n<2k} that
there is an $(\Img \psi)$-invariant set $E\subset{\boldsymbol\Delta}_{12}$
of cardinality $4$. However, the formulae for $\phi_0,\phi_1$
show that even the permutations
$\widehat\sigma_{i+2}|{{\boldsymbol\Delta}_{12}}=\phi_j(s_i)$,
\ $1\le i\le 5$, do not have a common invariant set of such cardinality
(in fact, $\phi_1$ is transitive and $\Img \phi_0$ has in
${\boldsymbol\Delta}_{12}$ exactly two orbits, each of length $6$).
Hence case $(iv_b)$ is impossible.
\vskip0.3cm

\noindent We are left with case $(iv_c)$. By Corollary \ref{Crl 6.6}
and Theorem
\ref{Thm: if Omega is cyclic then every Omega-homomorphism is cyclic},
we may assume that $\phi$ is of the form
%&
\begin{equation}\label{7.25}
\phi=\varphi_2\times\varphi_\omega\colon
\,B_{k-2}\to{\mathbf S}(2k-4)\times{\mathbf S}(4)\subset{\mathbf S}(2k)
\end{equation}
%&
(see Corollary \ref{Crl 6.8}); here
$\varphi_\omega\colon B_{k-2}\to{\mathbf S}(4)$
is a cyclic homomorphism defined as follows: if the homomorphism
$\omega\colon B_{k-2}\to{\mathbf S}(2)$ is trivial,
then $\varphi_\omega(s_i)=(2k-3,2k-2)(2k-1,2k)$, \ $i=1,...,k-3$;
and $\varphi_\omega(s_i)=(2k-3,2k-1)(2k-2,2k)$, \ $i=1,...,k-3$,
for the only non-trivial $\omega$. (In fact, there is one more possibility,
namely, $\varphi_\omega(s_i)=(2k-3,2k)(2k-2,2k-1)$, \ $i=1,...,k-3$;
if so, we conjugate $\psi$ by the transposition
$(2k-1,2k)$ and reduce this case to the previous one).

First, we note that the set $R=\{2k-3,2k-2,2k-1,2k\}$ cannot be
$\widehat\sigma_2$-invariant. Otherwise, $R$ would be
$(\Img \psi)$-invariant and the reduction $\psi_S$ of $\psi$
to the complement $S={\boldsymbol\Delta}_{2k}\setminus R$
would be a {\sl non-cyclic transitive} homomorphism
$B_k\to{\mathbf S}(2k-4)$, which contradicts
Theorem \ref{Thm: homomorphisms B(k) to S(n), 6<k<n<2k}
($\psi_S$ must be transitive, since $\varphi_2$ is so).

Now we show that $\omega$ must be non-trivial.
Indeed, if $\omega$ is trivial,
then $\Omega=\mu_{k-2}\times\mathbf 1_2$ and the action of $\Omega(s_2)$
on the $2$-component $\mathfrak C$ of $\widehat\sigma_1$
(that is, the conjugation by $\widehat\sigma_4$)
interchanges some two transpositions and does not move
the rest $k-2$. By Lemma \ref{Lm: Omega=Omega*},
the action of $\Omega^*(s_2)$ on the
$2$-component $\mathfrak C^*$ of $\widehat\sigma_{k-1}$
(that is, the conjugation by $\widehat\sigma_2$) is of the same type.
This means that $\widehat\sigma_2$ and $\widehat\sigma_{k-1}$ have exactly
$k-2$ common transpositions. Any of these common transpositions is
neither $(2k-3,2k-2)$ nor $(2k-1,2k)$, since otherwise the reduction
of $\psi$ to at least one of the complements
${\boldsymbol\Delta}_{2k}\setminus\{2k-3,2k-2\}$, \
${\boldsymbol\Delta}_{2k}\setminus\{2k-1,2k\}$,
${\boldsymbol\Delta}_{2k}\setminus\{2k-3,2k-2,2k-1,2k\}$
would be a non-cyclic transitive homomorphism
$$
B_k\to{\mathbf S}(n),\qquad 6<k<n, \ \ n=2k-2 \ \
\text{or} \ \ n=2k-4,
$$
which contradicts Theorem \ref{Thm: homomorphisms B(k) to S(n), 6<k<n<2k}.
Hence the conjugation by $\widehat\sigma_2$ interchanges the transpositions
$(2k-3,2k-2)$, $(2k-1,2k)$. It follows that the set
$\{2k-3,2k-2,2k-1,2k\}$ is $\widehat\sigma_2$-invariant;
however, we have already proved that this is impossible.

Taking into account that $\psi$ is normalized and using
(\ref{7.25}) and what has been proven above, we see that
$$
\widehat\sigma_1=(1,2)(3,4)(5,6)\cdots (2k-5,2k-4)(a,b)(c,d)
$$
and for $3\le i\le k-1$
$$
\aligned
\widehat\sigma_i=(1,2)(3,4)\cdots (2i-7,2i-6)&
\underbrace{(2i-5,2i-3)(2i-4,2i-2)}(2i-1,2i)\\
&\qquad\qquad\qquad\times\cdots\times (2k-5,2k-4)(a,c)(b,d),
\endaligned
$$
where $a=2k-3$, $b=2k-2$, $c=2k-1$, $d=2k$.
Now it is convenient to conjugate the original homomorphism $\psi$
by the permutation
$$
C=\begin{pmatrix}
1 & 2 & 3 &... & 2k-4 & a & b & c & d \\
5 & 6 & 7 &... & 2k   & 1 & 3 & 2 & 4
\end{pmatrix};
$$
we denote the new homomorphism by $\widetilde\psi$, but
preserve the notations $\widehat\sigma_i$ for all the permutations
$\widetilde\psi(\sigma_i)$, \ $1\le i\le k-1$. Clearly,
\begin{equation}\label{7.26}
\ \ \
\widehat\sigma_1=\underbrace{(1,3)(2,4)}(5,6)\cdots (2k-3,2k-2)(2k-1,2k)
\end{equation}
%&
and for $3\le i\le k-1$
%&
\begin{equation}\label{7.27}
\aligned
\widehat\sigma_i=(1,2)(3,4)\cdots (2i-3,2i-2)&
\underbrace{(2i-1,2i+1)(2i,2i+2)}(2i+3,2i+4)\\
&\qquad\times\cdots\times (2k-3,2k-2)(2k-1,2k).
\endaligned
\end{equation}
%&
\noindent{\bfit Claim.} {\sl The set ${\boldsymbol\Delta}_6=\{1,2,3,4,5,6\}$
is $\widehat\sigma_2$-invariant and the restriction of $\widehat\sigma_2$
to the complement ${\boldsymbol\Delta}_{12}\setminus{\boldsymbol\Delta}_6$
coincides with $(7,8)(9,10)\cdots (2k-1,2k)$.}
\vskip0.2cm

\noindent Indeed, it follows from (\ref{7.27}) that
$$
\widehat\sigma_i\widehat\sigma_{i+1}=(2i-1,2i+2,2i+3)(2i,2i+1,2i+4)\,;
$$
if $4\le i\le k-2$, this product commutes with $\widehat\sigma_2$;
hence for such $i$ every set $\{2i-1,2i,2i+1,2i+2,2i+3,2i+4\}$ is
$\widehat\sigma_2$-invariant. In particular,
the set ${\boldsymbol\Delta}_6=\{1,2,3,4,5,6\}$ 
is $\widehat\sigma_2$-invariant.
If $k>7$ then each of the sets
$\{7,8\}, \{9,10\},..., \{2k-1,2k\}$ is $\widehat\sigma_2$-invariant;
since is $\widehat\sigma_2$ has no fixed points, this shows that
the cyclic decomposition of $\widehat\sigma_2$
contains the disjoint product 
$$
(7,8)(9,10)\cdots (2k-1,2k)\,.
$$
\noindent Consider the case $k=7$. We still have the two
$\widehat\sigma_2$-invariant sets
$$
\{7,8,9,10,11,12\}\quad \text{and}\quad \{9,10,11,12,13,14\}\,; 
$$
hence each of
the sets $\{7,8\}$, $\{9,10,11,12\}$, $\{13,14\}$ is
$\widehat\sigma_2$-invariant. Since $\widehat\sigma_2$ has no fixed points
and is a product of disjoint transpositions, it must contain
the transpositions $(7,8)$, $(13,14)$ and some two transpositions
that are supported in $\{9,10,11,12\}$. The product
$\widehat\sigma_5\widehat\sigma_6=(9,12,13)(10,11,14)$ commutes with
$\widehat\sigma_2$, and we already know that $\widehat\sigma_2$ contains
the transposition $(13,14)$; hence the restriction
of $\widehat\sigma_2$ to $\{9,10,11,12\}$ coincides with
$(9,10)(11,12)$. This completes the proof of Claim.
\vskip0.3cm

To complete the whole proof, we consider the restrictions
$$
A=\widehat\sigma_1|{{\boldsymbol\Delta}_6}=(1,3)(2,4)(5,6),\ \
B=\widehat\sigma_2|{{\boldsymbol\Delta}_6},\ \
C=\widehat\sigma_4|{{\boldsymbol\Delta}_6}=(1,2)(3,4)(5,6).
$$
Claim shows that the restrictions of $\widehat\sigma_1$ and
$\widehat\sigma_2$ to the complement
${\boldsymbol\Delta}_{2k}\setminus{\boldsymbol\Delta}_6$ coincide;
hence $A\infty B$. Clearly, $AC=CA$; we know also that $B$ must
be a product of $3$ disjoint transpositions supported
in ${\boldsymbol\Delta}_6$.
There exist exactly $4$ permutations $B$ that satisfy all these
conditions:
$$
\aligned
B_1&=(1,2)(3,5)(4,6);\qquad B_2=(1,2)(3,6)(4,5);\\
B_3&=(1,5)(2,6)(3,4);\qquad B_4=(1,6)(2,5)(3,4).
\endaligned
$$
If $B=B_1$, then $\widetilde\psi$ coincides with the homomorphism
$\varphi_2$ from Corollary \ref{Crl 6.6} (with $n=k$). 
Any of the other three
possibilities leads to a conjugate homomorphism. Indeed,
the conjugation by the permutation $(5,6)(7,8)\cdots (2k-1,2k)$
interchanges $B_1$ with $B_2$, and $B_3$ with $B_4$; the conjugation
by the permutation $(1,3)(2,4)$ interchanges
$B_1$ with $B_3$. Further, these two conjugations preserve
formulae (\ref{7.26}), (\ref{7.27}). They preserve also the
form of the
restriction of $\widehat\sigma_2$ to
${\boldsymbol\Delta}_{12}\setminus{\boldsymbol\Delta}_6$
exhibited in Claim. This concludes the proof.
\end{proof}

The following statement is similar to Corollary
\ref{Crl: k-2 fixed points of psi(sigma(1)) for 6<k<n<2k}.

\begin{Corollary}\label{Crl 7.22} 
Let $k>6$ and let
$\psi\colon B_k\to{\mathbf S}(2k)$ be a non-cyclic homomorphism
such that $\widehat\sigma_1$ has at most $k-3$ fixed points.
Then either $\psi\sim\varphi_2$ or $\psi\sim\varphi_3$.
\end{Corollary}

\begin{proof} 
It follows from Lemma
\ref{Lm1: homomorphism with all components of hatsigma(1) of length<=k-3}
that $\widehat\sigma_1$ must have
a non-degenerate component of length at least $k-2$;
Lemma 
\ref{Lm: psi B(k) to S(2k) with psi(sigma(1)) having component of length>k-3}
completes the proof.
\end{proof}

\begin{Theorem}
\label{Thm: non-cyclic transitive B(k) to S(2k) is standard for k>6} 
Any non-cyclic transitive homomorphism
$\psi\colon B_k\to{\mathbf S}(2k)$ is standard
whenever $k>6$. In other words, $\psi$ is
conjugate to one of the model homomorphisms
$\varphi_1$, $\varphi_2$, $\varphi_3$.
\end{Theorem}

\begin{proof} 
Here we prove the theorem for $k>8$;
the cases $k=7,8$ will be considered separately.
We use the notation introduced in
Sec. \ref{sec: 7. Homomorphisms B'(k) to S(n), B'(k) 
to B(n) (n<k) and B(k) to S(n) (n<2k)} (see Convention therein).
Suppose to the contrary that there is
a non-cyclic transitive homomorphism $\psi$ that
is conjugate neither to $\varphi_1$ nor to
$\varphi_2$ nor to $\varphi_3$.
By Lemma
\ref{Lm: psi B(k) to S(2k) with psi(sigma(1)) having component of length>k-3},
every non-degenerate component
of $\widehat\sigma_1$ is of length at most $k-3$. Hence,
Lemma
\ref{Lm2: homomorphism with all components of hatsigma(1) of length<=k-3}$(a,c)$
applies to the homomorphism $\psi$;
this leads to the following conclusions:

$a)$ the reduction $\varphi=\phi_{\Sigma_1}\colon B_{k-2}
\to{\mathbf S}(\Sigma_1)\cong {\mathbf S}(N)$
of the homomorphism
$\phi=\psi|{B_{k-2}}\colon B_{k-2}\to{\mathbf S}(2k)$
to the set $\Sigma_1=\supp \widehat\sigma_1$ is trivial
(here and below $B_{k-2}\subset B_k$
is the subgroup generated by $\sigma_3,...,\sigma_{k-1}$);

$b)$ the reduction
$\varphi'=\phi_{\Sigma_1'}\colon B_{k-2}
\to{\mathbf S}(\Sigma_1')\cong {\mathbf S}(N')$
of $\phi$ to the set $\Sigma_1'=\Fix \widehat\sigma_1$ is non-cyclic;

$c)$ since
$\phi=\phi_{\Sigma_1}\times\phi_{\Sigma_1'}\colon
\,B_{k-2}\to{\mathbf S}(\Sigma_1)\times{\mathbf S}(\Sigma_1')
\subset {\mathbf S}(2k)$
and $\phi_{\Sigma_1}$ is trivial, we see
that eventually $\phi$ coincides with $\phi_{\Sigma_1'}$.
 %%&

By Lemma \ref{Lm: homomorphisms Bk to S(n) for n<=2k and supp(hatsigma(1))<6},
we have $N\ge 6$, and Corollary \ref{Crl 7.22} shows that
$N'\ge k-2$. It follows that
$6<k-2\le N' = 2k-N\le 2k-6<2(k-2)$.
Since the homomorphism $\phi_{\Sigma_1'}$ is non-cyclic,
Theorem \ref{Thm: homomorphisms B(k) to S(n), n<k}$(a)$,
Artin Theorem and Theorem \ref{Thm: homomorphisms B(k) to S(n), 6<k<n<2k}
imply that $\phi_{\Sigma_1'}$ is conjugate to a homomorphism
of the form $\mu_{k-2}\times\nu$,
where $\nu\colon B_{k-2}\to{\mathbf S}(N'-k+2)$
is some cyclic homomorphism.
Hence, without loss of generality we may assume
that $\Sigma_1'=\{1,...,N'\}$, \
$\Sigma_1=\{N'+1,...,2k\}$ and $\widehat\sigma_i=(i,i+1)\cdot S$
(disjoint product) for $3\le i\le k-1$, where $S$ is some permutation
not depending on $i$ and supported on a set
$Q\subseteq\Sigma_1'-\{3,...,k\}$.
We have $\# Q\le N' -(k-2)\le 2k-6-k+2=k-4$ and
$\# Q + 2 = \# \supp \widehat\sigma_i = N\ge 6$;
hence {\sl the set $Q=\supp S$ is nonempty
and does not coincide with the whole set ${\boldsymbol\Delta}_{2k}$}.
Clearly, $Q$ is $\widehat\sigma_i$-invariant
for any $i\ne 2$ (since $Q\subseteq\Fix\,\widehat\sigma_1$ and
$S\preccurlyeq\widehat\sigma_i$ for every $i\ge 3$).
On the other hand, $\widehat\sigma_2$ commutes with any
$\widehat\sigma_i=(i,i+1)\cdot S$, \ $4\le i\le k-1$, and thus
each of the sets $\{4,5\}\cup Q$,...,$\{k-1,k\}\cup Q$
is $\widehat\sigma_2$-invariant.
Hence their intersection $Q$ is $\widehat\sigma_2$-invariant;
this contradicts the transitivity of $\psi$ and
concludes the proof.
\end{proof}

\newpage

%Sec. 7
\section{Homomorphisms $B'_k\to{\mathbf S}(k)$ and endomorphisms of $B'_k$}
\label{sec: Homomorphisms B'(k) to S(k) and endomorphisms of B'(k)}

In this section we apply the results of 
Sec. \ref{sec: 7. Homomorphisms B'(k) to S(n), B'(k) to B(n) (n<k) and
B(k) to S(n) (n<2k)} to prove Theorem C and Theorem D. 
We start with some preparations to the proof of Theorem C.
\vskip0.2cm

Assume that $k>4$ and consider a {\sl non-trivial} homomorphism
$\psi\colon B'_k\to{\mathbf S}(k)$. Taking into
account the presentation of the commutator subgroup
$B'_k$ given by (\ref{1.14})-(\ref{1.21}), we denote the 
$\psi$-images of the generators $u,v,w,c_i$
by $\widehat u,\widehat v,\widehat w,\widehat c_i$ respectively. The latter
permutations satisfy the system of equations
%&
\begin{equation}\label{8.1}
\ \widehat u\widehat c_1\widehat u^{-1} = \widehat w,
\qquad\qquad\qquad\qquad\qquad \ \ \ \ \ \ \ \
\end{equation}
%&
%&
\begin{equation}\label{8.2}
\ \widehat u\widehat w\widehat u^{-1} = \widehat w^2\widehat c_1^{-1}\widehat w,
\qquad\qquad\qquad\qquad\qquad
\end{equation}
%&
\begin{equation}\label{8.3}
\ \widehat v\widehat c_1\widehat v^{-1} = \widehat c_1^{-1}\widehat w,
\qquad\qquad\qquad\qquad\qquad \ \ \
\end{equation}
%&

\begin{equation}\label{8.4}
\widehat v\widehat w\widehat v^{-1} = (\widehat c_1^{-1}\widehat w)^3\widehat c_1^{-2}\widehat w,
\qquad\qquad\qquad \ \ \ \
\end{equation}
%&
\begin{equation}\label{8.5}
\qquad \ \ \widehat u\widehat c_i = \widehat c_i\widehat v
\qquad\qquad\qquad \ \ \ (2\le i\le k-3),\qquad \ \ 
\end{equation}
%&
\begin{equation}\label{8.6}
\qquad \ \ \widehat v\widehat c_i = \widehat c_i\widehat u^{-1}\widehat v
\qquad\qquad \ \ \ \ (2\le i\le k-3),\qquad \ \
\end{equation}
%&
\begin{equation}\label{8.7}
\qquad \ \ \ \ \ \widehat c_i\widehat c_j = \widehat c_j\widehat c_i
\qquad\qquad\qquad \ (1\le i<j-1\le k-4),
\end{equation}
%&
\begin{equation}\label{8.8}
\widehat c_i\widehat c_{i+1}\widehat c_i = \widehat c_{i+1}\widehat c_i\widehat c_{i+1}
\qquad (1 \le i\le k-4).
\end{equation}
%&
Consider the embedding 
$$
\lambda_k'\colon B_{k-2}\hookrightarrow B'_k, \ \ \
\lambda_k'(s_i)=c_i,\qquad 1\le i\le k-3,
$$
and the composition
$$
\phi=\psi\circ\lambda_k'\colon B_{k-2}\stackrel{\lambda_k'}{\longrightarrow}
B'_k\stackrel{\psi}{\longrightarrow} {\mathbf S}(k),\ \ \ 
\phi(s_i)={\widehat c}_i, \qquad 1\le i\le k-3.
$$

\begin{Definition}[\caps tame homomorphisms $B'_k\to{\mathbf S}(k)$]
\label{Def: tame homomorphism psi: B'(k) to S(k)} 
For a non-trivial
homomorphism $\psi\colon B'_k\to{\mathbf S}(k)$,
set $G=\Img \psi\subseteq {\mathbf S}(k)$ and
$H=\Img \phi\subseteq {\mathbf S}(k)$.
For any $H$-orbit $Q\subseteq {\boldsymbol\Delta}_k$ we put
$Q' ={\boldsymbol\Delta}_k-Q$ and denote by
$\phi_Q\colon B_{k-2}\to{\mathbf S}(Q)$
and $\phi_{Q'}\colon B_{k-2}\to{\mathbf S}(Q')$
the reductions of $\phi$ to the $H$-invariant sets $Q$ and $Q'$, 
respectively; $\phi$ is the disjoint product
$\phi_Q\times\phi_{Q'}$. A non-trivial homomorphism $\psi$
is called {\em tame} if there is an $H$-orbit 
$Q\subset {\boldsymbol\Delta}_k$ of length $k-2$.
This orbit $Q$ (if it exists) is the only $H$-orbit
of length $\ge k-2$; we call it the {\em tame orbit} of $\psi$;
clearly $\# Q' =2$ and $\psi$ coincides with 
the disjoint product $\phi_Q\times\phi_{Q'}$ of
the non-cyclic transitive homomorphism 
$\phi_Q\colon B_{k-2}\to{\mathbf S}(Q)\cong {\mathbf S}(k-2)$
and the cyclic homomorphism
$\phi_{Q'}\colon B_{k-2}\to{\mathbf S}(Q')\cong {\mathbf S}(2)$
(a priori $\phi_{Q'}$ might be trivial; we shall see that
this cannot happen).

A group homomorphism $K\to{\mathbf S}(k)$ is said to be
{\it even} if its image is contained in the alternating subgroup 
$\mathbf A(k)={\mathbf S}'(k)$.
\hfill $\bigcirc$
\end{Definition}

\noindent By Lemma \ref{Lm: composition B(k-2) to Bk' to H},
for any non-trivial homomorphism
$\psi\colon B'_k\to{\mathbf S}(k)$ \ ($k>4$) we have:
\vskip0.3cm

\begin{itemize}

\item[$(*)$] {\sl the homomorphisms $\psi$ and $\phi$ are even, that is,
$H\subseteq G\subseteq \mathbf A(k)$; moreover,
$\phi$ is non-cyclic, $\phi(s_1)\ne \phi(s_3)$,
and $\phi(s_1^{-1})\ne \phi(s_3)$}.

\end{itemize}
\vskip0.3cm

%RRRRRRRRRRRRRRRRRRRRRR

\noindent In the next lemma we establish some properties of the
permutations $\widehat u$, $\widehat v$, $\widehat w$, $\widehat c_i$
corresponding to a non-trivial homomorphism $\psi$. 

\begin{Lemma}\label{Lm 8.2} 
$a)$ $[\widehat w]=[\widehat c_1^{-1}\widehat w]
=[\widehat c_1]=\ldots =[\widehat c_{k-3}]$ and
$[\widehat u] = [\widehat v]=[\widehat u^{-1}\widehat v]$.

$b)$ All the permutations $\widehat u,\widehat v,\widehat w,\widehat c_i$
are non-trivial $($and even$)$.

$c)$ $\widehat u$ commutes with all the permutations
$\widehat c_{i,j}=\widehat c_i\widehat c_j^{-1}$, \ $2\le i,j\le k-3$.

$d)$ If $\widehat c_1^2=1$, then $\widehat u^3=\widehat v^3=1$
and $\widehat v=\widehat u^{-1}$.

$e)$ If $k=5$, then $[\widehat c_1]\ne [3]$ and $[\widehat c_1]\ne [5]$.
\end{Lemma}

\begin{proof} 
$a)$ Follows immediately from relations (\ref{8.1}), (\ref{8.3}), 
(\ref{8.5}) (with $i=2$), (\ref{8.6}) (with $i=2$),  
and (\ref{8.8}), which shows that all the elements $\widehat c_i$ 
are conjugate to each other.
\vskip0.2cm

$b)$ Since $\psi$ is non-trivial, it is sufficient to show
that if one of the permutations
$\widehat u,\widehat v,\widehat w,\widehat c_i$
is trivial, then all of them are trivial.
If some $\widehat c_i=1$ or $\widehat w=1$, then it follows from $(a)$ that
$\widehat w=\widehat c_1=\ldots =\widehat c_{k-3}=1$, and (\ref{8.5}), 
(\ref{8.6}) imply that $\widehat u =\widehat v =1$. 
If $\widehat u = 1$ or $\widehat v = 1$,
then, by $(a)$, we have $\widehat u = \widehat v = 1$, 
and (\ref{8.1}), (\ref{8.2})
imply $\widehat w = \widehat c_1 = \widehat w^2$; hence $\widehat w = 1$.

$c)$ Relations (\ref{8.5}) may be written in the form
$$
\widehat v = \widehat c_2^{-1}\widehat u\widehat c_2
= \widehat c_3^{-1}\widehat u\widehat c_3
=\ldots = \widehat c_{k-3}^{-1}\widehat u\widehat c_{k-3};
$$
this shows that $\widehat u = (\widehat c_i\widehat c_j^{-1})\cdot\widehat u
\cdot (\widehat c_i\widehat c_j^{-1})^{-1}$ for all $2\le i,j\le k-3$.

$d)$ By $(a)$, the condition $\widehat c_1^2=1$ implies
that $\widehat c_i^2=1$ for all $i$; hence $c_i^{-1}=c_i$ and
relation (\ref{8.5}) for $i=2$ may be written in the form
$$
\widehat u = \widehat c_2\widehat v\widehat c_2^{-1}
= \widehat c_2^{-1}\widehat v\widehat c_2\,.
$$
In view of (\ref{8.6}) with $i=2$, 
the right hand side of the latter relation
is equal to $\widehat u^{-1}\widehat v$, and we get $\widehat u
= \widehat u^{-1}\widehat v$;
hence $\widehat v = \widehat u^2$.
Using the same relations (\ref{8.5}), (\ref{8.6}) with $i=2$
and $\widehat c_2^2=1$, we have
$$
\widehat v = \widehat c_2 \widehat u^{-1}\widehat v \widehat c_2^{-1} =
(\widehat c_2 \widehat u \widehat c_2^{-1})^{-1}\cdot \widehat c_2 \widehat v \widehat c_2^{-1}
=(\widehat c_2^{-1} \widehat u \widehat c_2)^{-1}\cdot \widehat u
= \widehat v^{-1}\widehat u;
$$
thus $\widehat u = \widehat v^2$. Thereby $\widehat u^3=\widehat v^3=1$
and $\widehat v=\widehat u^{-1}$.
\vskip0.2cm

$e)$ Assume that $[\widehat c_1]=[p]$, where $p=3$ or $p=5$.
Then also $[\widehat c_1^{-1}]=[p]$; by $(a)$,
we have $[\widehat w]=[\widehat c_1^{-1}\widehat w]=[p]$.
If $p=3$, then $(\widehat c_1^{-1}\widehat w)^3=1$, $\widehat c_1^{-2}=\widehat c_1$,
and (8.4) shows that
$[\widehat c_1\widehat w]=[\widehat c_1^{-2}\widehat w]
=[\widehat w]=[3]$; however this contradicts 
Lemma \ref{Lm: 3- and 5-cycles}$(a)$
(with $A=\widehat c_1^{-1}$ and $B=\widehat w$).
\vskip0.2cm

Consider the case $p=5$. Then
$[\widehat c_1]=[\widehat c_1^{-1}]
=[\widehat w]=[\widehat c_1^{-1}\widehat w]=[5]$.
It follows from (\ref{8.1}) and (\ref{8.2}) that 
$$
\widehat w^2\widehat c_1^{-1}\widehat w
=\widehat u\widehat w\widehat u^{-1}\qquad \text{and}\qquad
\widehat u\widehat c_1^2\widehat w\widehat u^{-1} 
= \widehat w^4\widehat c_1^{-1}\widehat w
=\widehat w^{-1}\widehat c_1^{-1}\widehat w\,;
$$
hence $[\widehat w^2\widehat c_1^{-1}\widehat w]=[5]$ and
$[\widehat c_1^2\widehat w]=[5]$.
Moreover, from (\ref{8.3}), (\ref{8.4}) we have
$$
\widehat v\widehat c_1\widehat w\widehat v^{-1}
= (\widehat c_1^{-1}\widehat w)^4\widehat c_1^{-2}\widehat w
=(\widehat c_1^{-1}\widehat w)^{-1}\cdot\widehat c_1^{-1}
\cdot (\widehat c_1^{-1}\widehat w)\,;
$$
therefore $[\widehat c_1\widehat w]=[\widehat c_1^{-1}]=[5]$. 
Taking $A=\widehat c_1^{-1}$ and
$B=\widehat w$, we see that all the permutations $A$, $B$, $AB$,
$A^{-1}B$, $A^{-2}B$, and $B^2AB$ are $5$-cycles in ${\mathbf S}(5)$, 
which contradicts Lemma \ref{Lm: 3- and 5-cycles}$(b)$.
\end{proof}

\noindent The following lemma brings us essentially closer to
the desired result.

\begin{Lemma}
\label{Lm: tame and non-tame homomorphisms B'(k) to S(k)} 
$a)$ The homomorphism $\psi$ is tame whenever
$k\ne 6$.

$b)$ If $k=6$ and $\psi$ is non-tame, then the homomorphism 
$\phi$ is transitive and conjugate to the homomorphism \
$\widetilde\nu_6'\colon B_4\to{\mathbf S}(6)$ defined
in Remark \ref{Rmk 4.2}.

$c)$ If $\psi$ is tame, then the reduction  
$\phi_Q\colon B_{k-2}\to{\mathbf S}(Q)\cong {\mathbf S}(k-2)$
to the tame orbit $Q$ is conjugate to the canonical epimorphism 
$\mu_{k-2}\colon B_{k-2}\to{\mathbf S}(k-2)$,
and the reduction $\phi_{Q'}\colon B_{k-2}
\to{\mathbf S}(Q')\cong {\mathbf S}(2)$ is a non-trivial
homomorphism. In particular, 
$\widehat c_i=\psi(c_i)=\phi(s_i)=S_i T$, \ $1\le i\le k-3$, where every
$S_i=\phi_Q(s_i)$ is a transposition supported in $Q$, 
and $T$ is the $($only$)$ transposition supported on $Q'$.
\end{Lemma}

\begin{proof} 
We start with the following claim, which is true for any
$k>4$ and any non-trivial homomorphism $\psi$:
\vskip0.2cm

\noindent{\bfit Claim 1.} {\sl There exists {\rm {(exactly one)}} $H$-orbit
of length $q\ge k-2$.}
\vskip0.2cm

\noindent For $k\ne 6$ this follows immediately from the property $(*)$ 
and Theorem \ref{Thm: homomorphisms B(k) to S(n), n<k}$(a)$.
For $k=6$, we deal with the non-cyclic even 
homomorphism $\phi\colon B_4\to{\mathbf S}(6)$ that
satisfies $(*)$. In this case there exists (exactly one)
$H$-orbit of length $q\ge 4$. Indeed,
let $Q$ be an $H$-orbit of some length $q$. If $q\le 3$ and $\phi_Q$ 
is non-cyclic, then, by Theorem \ref{Thm 3.14}, $\phi_Q(s_1)=\phi_Q(s_3)$.
Hence, if $\# Q\le 3$ for all $H$-orbits, then
the homomorphism $\phi$ cannot satisfy $(*)$. 
\vskip0.3cm

\noindent{\bfit Claim 2.} {\sl If $\psi$ is non-tame, then $k=6$.}
\vskip0.2cm

\noindent Taking into account Claim 1, we may assume that there is an
$H$-orbit $Q$ with $\# Q=q>k-2$. Clearly, either $q=k-1$ or $q=k$;
in any case, $\# Q'\le 1$ and $\phi=j_Q\circ\phi_Q$, where 
$j_Q\colon {\mathbf S}(Q)\hookrightarrow {\mathbf S}(k)$
is the natural embedding.
Since $\phi_Q$ is non-cyclic and transitive,
Theorem \ref{Thm: homomorphisms B(k) to S(k+1)} and
Theorem \ref{Thm: homomorphisms B(k) to S(k+2)} 
show that this could only happen in one of the following five cases:
\vskip0.3cm

\begin{itemize}

\item[$i)$] $k=5$, $k-2=3$, $q=4$, \ 
$\phi=j_Q\circ\phi_Q\colon B_3\stackrel{\phi_Q}{\longrightarrow}
{\mathbf S}(4)\stackrel{j_Q}{\hookrightarrow} {\mathbf S}(5)$,
\ $\phi_Q$ is transitive and non-cyclic, $\phi$ is even;

\item[$ii)$] $k=6$, $k-2=4$, $q=5$, \ 
$\phi=j_Q\circ\phi_Q\colon B_4\stackrel{\phi_Q}{\longrightarrow }
{\mathbf S}(5)\stackrel{j_Q}{\hookrightarrow} {\mathbf S}(6)$,
\ $\phi_Q$ is transitive and non-cyclic, $\phi$ is even;

\item[$iii)$] $k=7$, $k-2=5$, $q=6$, \
$\phi=j_Q\circ\phi_Q\colon B_5\stackrel{\phi_Q}{\longrightarrow}
{\mathbf S}(6)\stackrel{j_Q}{\hookrightarrow} {\mathbf S}(7)$,
\ $\phi_Q$ is transitive and non-cyclic, $\phi$ is even;

\item[$iv)$] $k=5$, $k-2=3$, $q=5$, \ 
$\phi=\phi_Q\colon B_3\to{\mathbf S}(5)$, \ $\phi$ is transitive, 
non-cyclic and even;

\item[$v)$] $k=6$, $k-2=4$, $q=6$, \
$\phi=\phi_Q\colon B_4\to{\mathbf S}(6)$, \ $\phi$ is transitive,
non-cyclic and even.

\end{itemize}

\vskip0.3cm

\noindent In fact all these cases but $(v)$ are impossible.
Indeed, in case $(i)$, applying Proposition \ref{Prp 4.1}$(a)$, we
see that the (even!) homomorphism $\phi_Q$ must be conjugate
to the homomorphism $\psi_{3,4}^{(2)}$; clearly
$\widehat c_1\sim\psi_{3,4}^{(2)}(s_1)=(2,3,4)$; 
hence $[\widehat c_1]=[3]$,
which contradicts Lemma \ref{Lm 8.2}$(e)$.
In case $(ii)$, by Lemma \ref{Lm 4.2}, the homomorphism $\phi$ would
satisfy $\phi(s_1)=\phi_Q(s_1)=\phi_Q(s_3)=\phi(s_3)$, which
contradicts property $(*)$. In case $(iii)$, by Proposition \ref{Prp 4.9},
the homomorphism $\phi_Q$ must be conjugate to the homomorphism
$\psi_{5,6}$ that sends any $s_i$ into an {\sl odd} permutation; 
clearly, $\phi$ makes the same (for $\# Q' = 1$), 
which is impossible (since $\phi$ must be even).
To eliminate $(iv)$, we use Proposition \ref{Prp 4.1} $(b)$, which shows 
that the homomorphism $\phi$ must be conjugate to the homomorphism
$\psi_{3,5}$; so, $\widehat c_1\sim \psi_{3,5}(s_1)=(1,4,3,2,5)$; however
this contradicts Lemma \ref{Lm 8.2}$(e)$. This proves Claim 2 and the
statement $(a)$ of the lemma.
\vskip0.2cm

$b)$ If $k=6$ and $\psi$ is non-tame, the proof of Claim 2 shows that
we are in the situation of case $(v)$.
By Proposition \ref{Prp 4.5} and condition $(*)$, the homomorphism $\phi$
must be conjugate to one of the homomorphisms $\psi_{4,6}^{(i)}$
defined by (\ref{4.4}). However, $\psi_{4,6}^{(1)}$ and $\psi_{4,6}^{(2)}$
are not even, and for $\psi_{4,6}^{(4)}$ we have
$\psi_{4,6}^{(4)}(s_1^{-1})=\psi_{4,6}^{(4)}(s_3)$, which is uncompatible
with $(*)$; hence $\psi\sim \psi_{4,6}^{(3)}$. By Remark \ref{Rmk 4.2},
$\psi_{4,6}^{(3)}$ is conjugate to the homomorphism
$\widetilde\nu_6'$. 

$c)$ Since $\psi$ is tame, the reduction 
$\phi_Q\colon B_{k-2}\to{\mathbf S}(Q)\cong{\mathbf S}(k-2)$
is a non-cyclic transitive homomorphism. If this homomorphism
is conjugate to $\mu_{k-2}$, the other assertions
of the statement $(c)$ are evident (note that if $\phi_Q\sim\mu_{k-2}$,
then the "complementary" reduction $\phi_{Q'}\colon B_{k-2}
\to{\mathbf S}(Q')\cong {\mathbf S}(2)$
{\sl must} be non-trivial, since the homomorphism $\psi$ is even). 

Let us assume that $\phi_Q$ is not conjugate to $\mu_{k-2}$; by Artin
Theorem, this may only happen if $k=6$ or $k=8$. The complementary
reduction $\phi_{Q'}$ is either trivial or takes each
$s_i$ to the only transposition $T$ supported on $Q'$;
in any case, we have $\phi_{Q'}(s_1)=\phi_{Q'}(s_3)$
and $\phi_{Q'}(s_1^{-1})=\phi_{Q'}(s_3)$.
If $k=6$, the reduction $\phi_Q$ must be conjugate 
to one of Artin's homomorphisms $\nu_{4,j}$, \ $1\le j\le 3$; 
however, in each of these cases we have either
$\psi(s_1)=\psi(s_3)$ or $\psi(s_1^{-1})=\psi(s_3)$,
which contradicts $(*)$. 

Finally, we must show that
the case when $k=8$ and $\phi_Q\sim \nu_6$ is impossible. Since $\nu_6(s_1)$
is the product of three disjoint transpositions and $\phi$ must be even,
the complementary reduction $\phi_{Q'}$ sends each $s_i$ 
to the only transposition $T$ supported on $Q'$.
Without loss of generality, we may assume that $T=(1,2)$ and
%&
\begin{equation}\label{8.9}
\phi\colon \left\{
\aligned
\ &s_1\mapsto \widehat c_1=(1,2)(3,4)(5,6)(7,8), \ \
s_2\mapsto \widehat c_2=(1,2)(3,7)(4,5)(6,8),\\
&s_3\mapsto \widehat c_3=(1,2)(3,5)(4,6)(7,8), \ \
s_4\mapsto \widehat c_4=(1,2)(3,4)(5,7)(6,8),\\
&\qquad\qquad\qquad\qquad s_5\mapsto \widehat c_5=(1,2)(3,6)(4,5)(7,8).
\endaligned
\right .
\end{equation}
%&
By Lemma \ref{Lm 8.2}$(d)$, \ $\widehat u^3=1$;
since $\widehat u$ is even and non-trivial,
we see that either $[\widehat u]=[3]$ or $[\widehat u]=[3,3]$.
By Lemma \ref{Lm 8.2} $(c)$, $\widehat u$ commutes with all the permutations
$\widehat c_{i,j}= \widehat c_i\widehat c_j^{-1}$, \ $i,j\ge 2$;
in particular, this is the case for $\widehat c_{2,3}=(3,4,8)(5,7,6)$.
Since $\Fix \widehat c_{2,3}=\{1,2\}$, this set
is $\widehat u$-invariant. It follows that $\{1,2\}\subseteq \Fix \widehat u$
(the cyclic decomposition of $\widehat u$ cannot contain
a transposition). Hence $\supp \widehat u\subseteq \{3,4,5,6,7,8\}$.
Further, $\widehat c_{3,5}=(3,4)(5,6)$. The set $\{7,8\}$ is the fixed points
set of the permutation $(3,4)(5,6)$ acting on $\{3,4,5,6,7,8\}$; 
therefore, it must be $\widehat u$-invariant; as above, this shows that
$\{7,8\}\subseteq \Fix \widehat u$ and $\supp \widehat u\subseteq \{3,4,5,6\}$.
Therefore, $\widehat u$ must be a $3$-cycle supported in $\{3,4,5,6\}$;
however such a permutation cannot commute with $(3,4)(5,6)$ and a
contradiction ensues.
\end{proof}

\noindent Recall that we denote by $\mu_k'$ the restriction of the canonical
projection 
$$
\mu_k\colon B_k\to{\mathbf S}(k)
$$ 
to the commutator subgroup $B'_k$; similarly, $\nu_6'$
denotes the restriction to $B'_6$ of Artin's
homomorphism $\nu_6$. If $\psi=\mu_k'$, then 
%&
\begin{equation}\label{8.10}
\aligned
\ &\widehat u=(1,3,2), \qquad \widehat v=(1,2,3), \qquad \widehat w=(1,3)(2,4), \\
&\widehat c_i=(1,2)(i+2,i+3),\qquad 1\le i\le k-3.
\endaligned
\end{equation}
%&
Moreover, if $k=6$ and $\psi=\nu_6'$, then
%&
\begin{equation}\label{8.11}
\aligned
\ &\widehat u=(1,3,6)(2,5,4), \ \ \widehat v=(1,6,3)(2,4,5), \ \
\widehat w=(2,3)(5,6), \\
&\widehat c_1=(1,4)(2,3), \ \ \widehat c_2=(3,6)(4,5), \ \ \widehat c_3=(1,3)(2,4).
\endaligned
\end{equation}
%&
\begin{Remark}\label{Rmk 8.1} 
Let $k>4$. In view of
Lemma \ref{Lm: tame and non-tame homomorphisms B'(k) to S(k)}, in order
to classify non-trivial homomorphisms
$\psi\colon B'_k\to{\mathbf S}(k)$
up to conjugation, it suffices to study the following
two cases:
\vskip0.3cm

\begin{itemize}

\item[$i)$] The homomorphism $\psi$ is tame, with the tame
$H$-orbit $Q=\{3,4,...,k\}$. The reduction 
$\phi_Q\colon B_{k-2}\to{\mathbf S}(k-2)$
coincides with the ``shifted" canonical epimorphism 
$$
\widetilde\mu_{k-2}\colon B_{k-2}\to{\mathbf S}(Q),
\ \ \widetilde\mu_{k-2}(s_i)=(i+2,i+3), \ \
1\le i\le k-3.
$$
$Q'=\{1,2\}$ and the complementary reduction
$\phi_{Q'}\colon B'_{k-2}
\to{\mathbf S}(Q')\cong {\mathbf S}(2)$
is of the form $\phi_{Q'}(s_i)=(1,2)$, \ $1\le i\le k-3$.
The homomorphism $\phi$ is the disjoint product
$\phi_Q\times\phi_{Q'}$ and
%&
\begin{equation}\label{8.12}
\widehat c_i=\phi(s_i)=(1,2)(i+2,i+3)\qquad \text{for all} \ \
i=1,...,k-3.
\end{equation}
%&
%\vskip0.3cm

\item[$ii)$] $k=6$ and the homomorphism $\psi$ is non-tame, 
with the only $H$-orbit $Q={\boldsymbol\Delta}_6$. The homomorphism
$\phi\colon B_4\to{\mathbf S}(6)$
coincides with the homomorphism $\widetilde\nu_6'$ and
%&
\begin{equation}\label{8.13}
\aligned
\ &\phi(s_1)= \widehat c_1=(1,4)(2,3), \ \ \
\phi(s_2)= \widehat c_2=(3,6)(4,5), \\
&\qquad\qquad\qquad\phi(s_3)= \widehat c_3=(1,3)(2,4).
\endaligned
\end{equation}
%&
\end{itemize}

\noindent Let us say that $\psi$ is {\it reduced} 
if it is either of type $(i)$ or of type $(ii)$.
\end{Remark}

\begin{Lemma}\label{Lm 8.4} 
Let $\psi$ be a reduced homomorphism of type $(i)$. 

$a)$ If $k\ge 6$, then $\widehat u(\{1,2,3\})=\{1,2,3\}$
and $\widehat u(\{4,5,6\})=\{4,5,6\}$.

$b)$ If $k\ge 7$, then $4,5,...,k\in \Fix \widehat u$ and $\widehat u$
is a $3$-cycle supported on $\{1,2,3\}$.
\end{Lemma}

\begin{proof} 
By Lemma \ref{Lm 8.2}$(c)$ and relation (\ref{8.12}), 
$\widehat u$ commutes with each of the permutations 
$$
\widehat c_{i,i+1}=\widehat c_i\widehat c_{i+1}^{-1}=(i+2,i+3,i+4), \ \
2\le i\le k-4.
$$
Hence each of the sets $\{4,5,6\}$, 
$\{5,6,7\}$,...,$\{k-2,k-1,k\}$ is $\widehat u$-invariant.
The union, the intersection, and the difference 
of two $\widehat u$-invariant sets are $\widehat u$-invariant.
This implies $(a)$. Moreover, if $k\ge 7$, we have
$$
j+2,j+5\in \Fix \widehat u \ \ \text{and} \ \
\widehat u(\{j+3,j+4\})=\{j+3,j+4\} \ \text{whenever} \ 2\le j\le k-5;
$$
by Lemma \ref{Lm 8.2} $(d)$, all the cycles in the cyclic
decomposition of $\widehat u$ are of length $3$; hence,
$j+3,j+4\in \Fix \widehat u$ \ and \ $\supp \widehat u=\{1,2,3\}$.
\end{proof}

\begin{Theorem}
\label{Thm: homomorphisms B'(k) to S(k)} 
Let $k>4$ and let $\psi\colon B'_k\to{\mathbf S}(k)$
be a non-trivial homomorphism. Then
either $\psi\sim \mu_k'$ or $k=6$ and
$\psi\sim \nu_6'$. In any case $\Img\psi=\mathbf A(k)$
and $\Ker\psi=J_k={PB_k}\cap B'_k$.
\end{Theorem}

\begin{proof} 
By Remark \ref{Rmk 8.1}, we may assume that $\psi$ is 
reduced. Let us start with case $(i)$. By Lemma \ref{Lm 8.2}$(a)$,
$[\widehat c_1^{-1}\widehat w]=[\widehat w]=[\widehat c_1]=[2,2]$;
hence $\widehat w$ and $\widehat c_1$ cannot be disjoint.
\vskip0.3cm

\noindent{\bfit Claim 1.} $\supp \widehat w = \{1,2,3,4\}$ and
either $\widehat w=(1,3)(2,4)$ or $\widehat w=(1,4)(2,3)$.
\vskip0.3cm

\noindent Let $m=\# (\{1,2,3,4\}\cap\supp \widehat w)$. 
We already know that
$m\ge 1$. The values $m=1$ and $m=3$ cannot occur by trivial reasons 
($m=1$ implies $[\widehat c_1^{-1}\widehat w]=[3,2,2]$; and if $m=3$, then
either $[\widehat c_1^{-1}\widehat w]=[5]$ or $[\widehat c_1^{-1}\widehat w]=[3]$).
Assume that $m=2$, that is, $\supp \widehat w=\{a,b,p,q\}$,
where $a,b\in \{1,2,3,4\}$ and $p,q\ge 5$. Then
$k\ge 6$. By (\ref{8.1}) and (\ref{8.12}), 
$\widehat w=\widehat u(1,2)(3,4)\widehat u^{-1}$;
hence $\widehat u(\{1,2,3,4\})=\supp \widehat w = \{a,b,p,q\}$.
In view of Lemma \ref{Lm 8.4} $(a)$, this shows that
$\{1,2,3\}=\widehat u(\{1,2,3\})\subset\widehat u(\{1,2,3,4\})=\{a,b,p,q\}$,
which contradicts the condition $p,q\ge 5$. Thus, $\widehat w$ is a
product of two disjoint transpositions supported on $\{1,2,3,4\}$,
and the condition
$[(1,2)(3,4)\cdot\widehat w]=[\widehat c_1^{-1}\widehat w]=[2,2]$ implies the desired
result.
\vskip0.3cm 

If $\widehat w=(1,4)(2,3)$, we conjugate the homomorphism $\psi$
by the transposition $(1,2)$ and obtain a homomorphism that
sends any $c_i$ into $\widehat c_i$ and sends $w$ into $(1,3)(2,4)$;
therefore, without loss of generality we may assume that the
original homomorphism $\psi$ itself satisfies the condition
%&
\begin{equation}\label{8.14}
\widehat w=\psi(w)=(1,3)(2,4).
\end{equation}
%&
Then relation (\ref{8.1}) takes the form
%&
\begin{equation}\label{8.15}
\widehat u(1,2)(3,4)\widehat u^{-1}=(1,3)(2,4);
\end{equation}
%&
in particular,
%&
\begin{equation}\label{8.16}
\widehat u(\{1,2,3,4\})=\{1,2,3,4\}.
\end{equation}
%&
Taking into account (\ref{8.10}), (\ref{8.12}), and (\ref{8.14}),
we conclude the proof of the theorem
in case $(i)$ by proving the following claim:
\vskip0.3cm

\noindent{\bfit Claim 2.} $a)$ Any $i\ge 4$ is a fixed point of $\widehat u$,
and thus $\widehat u$ is a $3$-cycle supported on $\{1,2,3\}$.
\ $b)$ $\widehat u=(1,3,2)$ and $\widehat v=(1,2,3)$.
\vskip0.3cm

In view of Lemma \ref{Lm 8.4}$(b)$, we need to prove $(a)$ only for $k=5,6$.
For $k=6$, Lemma \ref{Lm 8.4}$(a)$ 
shows that $\widehat u(\{4,5,6\})=\{4,5,6\}$;
by (\ref{8.16}), 
we have $\widehat u(4)=4$ and $\widehat u(\{5,6\})=\{5,6\}$.
In fact, $\{5,6\}\subset \Fix \widehat u$
(since $\widehat u$ cannot contain a transposition);
this proves $(a)$ for $k=6$. If $k=5$, (\ref{8.16}) shows that
$\widehat u(5)=5$. Relation (\ref{8.5}) (with $i=2$)
and Lemma \ref{Lm 8.2}$(d)$ imply 
$(\widehat u\widehat c_2)(5)=(\widehat c_2\widehat u^{-1})(5)$.
Since $\widehat u(5)=5$ and (by (\ref{8.12}) 
$\widehat c_2(5)=4$, this means that
$\widehat u(4)=4$, which concludes the proof of $(a)$.
To prove $(b)$, we note that $\widehat u = (1,3,2)$ is the only $3$-cycle
supported on $\{1,2,3\}$ that satisfies (\ref{8.15}).
\vskip0.2cm

Case $(ii)$ may be treated by straightforward computations;
however they are too long, and we prefer to use a simple trick.
Namely, instead of the original homomorphism $\psi$ of type $(ii)$,
we consider its composition $\widetilde\psi=\varkappa\circ\psi$ with
the outer automorphism $\varkappa$ of the group ${\mathbf S}(6)$.
(see (\ref{4.3})).
It is completely clear that $\widetilde\psi$ is
a tame homomorphism of type $(i)$; it follows
from what has been proven above that $\widetilde\psi\sim \mu_6'$.
The automorphism $\varkappa$ is involutive and
$\nu_6'=\varkappa\circ\mu_6'$; therefore,
$\psi\sim \nu_6'$.
\end{proof}

\begin{Corollary}
\label{Crl: every homomorphism B'(k) to S(k) extends to B(k)} 
For $k>4$ each homomorphism $\psi\colon B'_k\to{\mathbf S}(k)$
extends to $B_k$. For a non-trivial $\psi$ an extension
$\Psi\colon B_k\to{\mathbf S}(k)$ is unique.
\end{Corollary}

\begin{proof} 
The existence of an extension $\Psi$ follows immediately from
Theorem \ref{Thm: homomorphisms B'(k) to S(k)}.
The uniqueness follows from the facts that
$\mu_k(B'_k)={\mathbf A}(k)$,
\ $\nu_6(B'_6)={\mathbf A}(6)$
and for any $k\ge 3$ the centralizer
of ${\mathbf A}(k)$ in ${\mathbf S}(k)$ is trivial.
\end{proof}

\begin{Remark}\label{Rmk 8.2} 
In view of Artin Theorem,
Corollary \ref{Crl: every homomorphism B'(k) to S(k) extends to B(k)}
implies Theorem \ref{Thm: homomorphisms B'(k) to S(k)}.
However, I have no idea how to extend non-trivial homomorphisms
$\psi\colon B'_k\to{\mathbf S}(k)$ to
homomorphisms $\Psi\colon B_k\to{\mathbf S}(k)$
{\sl without} Theorem \ref{Thm: homomorphisms B'(k) to S(k)}.
\end{Remark}

\begin{Theorem}
\label{Thm: J(k) is a completely characteristic subgroup of B'(k)} 
Let $k>4$. Then the pure commutator subgroup $J_k = {PB_k}\cap B'_k$
is a completely characteristic subgroup of the group $B'_k$, that is,
$\phi(J_k)\subseteq J_k$ for any endomorphism $\phi\colon B'_k\to B'_k$.
Moreover, $\phi^{-1}(J_k) = J_k$ and $\Ker\, \phi \subset J_k$
for every non-trivial endomorphism $\phi$.
\end{Theorem}

\begin{proof} 
The case of trivial $\phi$ is trivial.
Given a non-trivial $\phi$, consider the composition
$$
\psi=\mu_k'\circ\phi\colon B'_k\stackrel{\phi}{\longrightarrow} B'_k
\stackrel{\mu_k'}{\longrightarrow} {\mathbf S}(k).
$$
This homomorphism $\psi$ must be non-trivial, since otherwise
$\Img\phi\subseteq \Ker\mu_k' \subset {PB_k}$
and Markov Theorem implies that $\phi$ is trivial.
By Theorem \ref{Thm: homomorphisms B'(k) to S(k)},
either $\psi\sim\mu_k'$ or $k=6$ and
$\psi\sim\nu_6'$; in any of these cases,
$\Ker\psi={PB_k}\cap B'_k = J_k$ and
we have
$$
J_k = \Ker\psi = \Ker\,(\mu_k'\circ\phi)
= \phi^{-1}(\Ker\mu_k') = \phi^{-1}(J_k).
$$
Certainly, this shows also that $\phi(J_k)\subseteq J_k$
and $\Ker\phi \subset J_k$.
\end{proof}

\begin{Remark}\label{Rmk 8.3} 
For a non-trivial endomorphism $\phi$ the inclusion
$\Ker\, \phi \subset J_k$ must be strict,
since $B'_k/J_k\cong \mathbf A(k)$
and $B'_k$ is torsion free. It seems
that for $k>4$ no examples of non-trivial endomorphisms
$B'_k\to B'_k$ 
with non-trivial kernels are known. I conjectured that
for $k>4$ a proper quotient group of the commutator
subgroup $B'_k$ cannot be torsion free (this would imply
that any non-trivial endomorphism of $B'_k$ must
be injective). I was told that D. Goldsmith's braid group
(which is a proper non-abelian quotient group of $B_k$)
is torsion free. For sure, this is true if $k=3$, 
but I newer saw any proof for $k>4$. 
If so, this would disprove my conjecture.

E. Artin \cite{Art47b} proved that the pure braid group ${PB_k}$
is a characteristic subgroup of the braid group $B_k$, that is,
$\phi({PB_k})={PB_k}$ for any {\sl automorphism} $\phi$ of 
the {\sl whole} braid group $B_k$ (see also Theorem
\ref{Thm: for k ne 4 psi(PB(k))<PB(k) for any non-integral psi}).
Formally, for $k>4$
Theorem \ref{Thm: J(k) is a completely characteristic subgroup of B'(k)} 
is essentially stronger than this Artin theorem
(and also essentially stronger than Theorem
\ref{Thm: for k ne 4 psi(PB(k))<PB(k) for any non-integral psi}, 
which, in turn, is an improvement of Artin's result). However I do not know
any non-trivial endomorphism of $B'_k$ ($k>4$) that
is not an automorphism. Seemingly, nobody knows whether there is 
an automorphism of $B'_k$ that cannot be extended to
an automorphism of the whole braid group $B_k$.
In view of these remarks, it may actually happen that
Theorem \ref{Thm: J(k) is a completely characteristic subgroup of B'(k)}
does not say more than Artin's result says. Nevertheless,
Theorem \ref{Thm: J(k) is a completely characteristic subgroup of B'(k)}
works in some situations when Artin Theorem and Theorem
\ref{Thm: for k ne 4 psi(PB(k))<PB(k) for any non-integral psi}
(in their present forms) are useless.
\end{Remark}
\vfill

\newpage

%Sec. 9
\section{Special homomorphisms $B_k\to B_n$}
\label{sec: Special homomorphisms B(k) to B(n)}
Let $\sigma_1,...,\sigma_k$ be the canonical generators
in $B_k$, $\alpha=\sigma_1\cdots\sigma_k$, \
$\beta=\alpha\sigma_1=\sigma_1\cdots\sigma_k\sigma_1$.
\vskip0.2cm

Recall the notions and notation introduced
in Definitions \ref{Def: special system of generators} 
and \ref{Def: special homomorphism}.
Given a special system of generators $\{a,b\}$ in 
the braid group $B_m$, we denote by  
${\mathcal H}_m(a,b)\subset B_m$ consisting of all the elements 
$g^{-1}a^q g$ and $g^{-1}b^q g$, where $g$ runs over 
$B_m$ and $q$ runs over ${\mathbb Z}$. 
Furthermore, a homomorphism
$\varphi\colon B_k\to B_n$ is said to be 
{\it special} if there exists a special system of generators
$\{a,b\}$ in $B_n$ such that 
$\varphi({\mathcal H}_k(\alpha,\beta))\subseteq {\mathcal H}_n(a,b))$.
\vskip0.2cm

According Murasugi Theorem mentioned in 
Sec. \ref{Ss: Main results},
the latter property is equivalent
to the following condition: 
\vskip0.2cm

\begin{itemize}

\item {\sl for every element $h\in B_k$ 
of finite order modulo the center $CB_k$ of %the group
\, $B_k$, its image $\varphi(h)$ has the same property
in the group $B_n$, that is, $\varphi(h)$
is an element of finite order modulo the center 
$CB_n$ of the group $B_n$}. 

\end{itemize}

\vskip0.2cm

\noindent It was stated in \cite{Lin71} (see also \cite{Lin79}
and \cite{Lin03,Lin04a}, Sec. 9, for the proof) that
\vskip0.2cm

\begin{itemize}

\item {\sl for every holomorphic mapping 
$f\colon{\mathbf G}_k\to{\mathbf G}_n$,
every point $z^\circ\in{\mathbf G}_k$, and any choice 
of isomorphisms $B_k\cong\pi_1({\mathbf G}_k,z^\circ)$
and $\pi_1({\mathbf G}_n,f(z^\circ))\cong B_n$
the induced homomorphism
$$
f_*\colon B_k\cong\pi_1({\mathbf G}_k,z^\circ)\to
\pi_1({\mathbf G}_n,f(z^\circ))\cong B_n
$$
is special}.

\end{itemize}
\vskip0.2cm

\noindent This property motivates the study of special 
braid homomorphisms.
\vskip0.3cm

\noindent According to Definition-Notation \ref{Def: four progressions}, 
$P(k)$ denotes the union
of the four increasing infinite arithmetic progressions
$$
P^i(k)=\left\{\left. p^i_j(k)=p^i_1(k)+(j-1)d(k)\right|\
j\in\mathbb N\right\}\,,
\quad 1\le i\le 4\,,
$$
with the same difference $d(k) = k(k-1)$ and starting with
the initial terms
$$
\quad\qquad 
p^1_1(k)= k, \ \ p^2_1(k)= k(k-1), \ \ p^3_1(k)= k(k-1)+1, \ \
{\text {and}} \ \ p^4_1(k)= (k-1)^{2}
$$
respectively
\vskip0.2cm

\noindent Let $\varphi\colon B_k\to B_n$ be
a special homomorphism. Then there exist
a special system of generators $\{a,b\}$ in $B_n$,
an element $g\in B_n$, and integers $p,q$ such that
at least one of the following four conditions is fulfilled:
$$
\aligned
\ &1) \ \varphi(\alpha)=a^p \ \ \text{and} \ \ \varphi(\beta)=gb^q g^{-1};
\qquad
2) \ \varphi(\alpha)=a^p \ \ \text{and} \ \ \varphi(\beta)=ga^q g^{-1};\\
&3) \ \varphi(\alpha)=b^p \ \ \text{and} \ \ \varphi(\beta)=gb^q g^{-1};
\qquad
4) \ \varphi(\alpha)=b^p \ \ \text{and} \ \ \varphi(\beta)=ga^q g^{-1}.
\endaligned
$$
We denote by ${\mathcal S}_i(k,n)$ the set
of all special homomorphisms
$\varphi\colon B_k\to B_n$
that satisfy condition $(i)$ ($1\le i\le 4$).

\begin{Theorem}
\label{Thm: arithmetical properties of special homomorphisms B(k) to B(n)} 
Assume that for some $k\ne 4$ and some $n$ there
exists a non-integral special homomorphism
$\varphi\colon B_k\to B_n$.
Then $n\in P(k)$. More precisely, there exist integers $l,t$ and
an element $g\in B_n$ such that
\vskip0.2cm

\noindent $a)$ if $\varphi\in {\mathcal S}_1(k,n)$, 
then $l\ge 0$, \ $n=k+lk(k-1)\in P^1(k)$,
\ $(t,k(k-1))=1$, and
$$
\varphi(\alpha)=a^p, \ \ \varphi(\beta)=gb^q g^{-1}, \ \
\text{where} \ \  p=t(l(k-1)+1), \ q=t(lk+1).
$$
\vskip0.2cm

\noindent $b)$ If $\varphi\in {\mathcal S}_2(k,n)$, 
then $l\ge 1$, \ $n=lk(k-1)\in P^2(k)$,
\ $g$ commutes with $a^{tk(k-1)}$, and
$$
\varphi(\alpha)=a^p, \ \ \varphi(\beta)=ga^q g^{-1}, \ \
\text{where} \ \ where \ \ p=t(k-1), \ q=tk.
$$
\vskip0.2cm

\noindent 
$c)$ if $\varphi\in {\mathcal S}_3(k,n)$, then $l\ge 1$, \ $n=lk(k-1)+1\in P^3(k)$,
\ $g$ commutes with $b^{tk(k-1)}$, and
$$
\varphi(\alpha)=b^p, \ \ \varphi(\beta)=ga^q g^{-1}, \ \
\text{where} \ \ p=t(k-1), \ q=tk;
$$
\vskip0.2cm

\noindent $d)$ if $\varphi\in {\mathcal S}_4(k,n)$, then $l\ge 1$,
\ $n=(k-1)(lk-1)\in P^4(k)$,
\ $(t,k(k-1))=1$, and
$$
\varphi(\alpha)=b^p, \ \ \varphi(\beta)=ga^q g^{-1}, \ \
\text{where} \ \ p=t(l(k-1)-1), \ q=t(lk-1).
$$
In particular, if $k\ne 4$ and $n\notin P(k)$, then
any special homomorphism $B_k\to B_n$
is integral.
\end{Theorem}

\begin{proof} There exist a homomorphism
$\delta\colon B_n\to{\mathbb Z}$ such that
$\delta(a)=n-1$ and $\delta(b)=n$
(since $B_n/B'_n\cong{\mathbb Z}$
and $a^n=b^{n-1}$). Moreover, since $\alpha^k=\beta^{k-1}$,
we have
%&
\begin{equation}\label{9.46}
k\delta(\varphi(\alpha))=(k-1)\delta(\varphi(\beta)).
%\eqnum{9.46}
\end{equation}
%&
The element $a^n=b^{n-1}$ is central in $B_n$.
Hence if $\varphi\in {\mathcal S}_1(k,n)\cup {\mathcal S}_2(k,n)$,
then the element $\varphi(\alpha^n)=a^{pn}$ commutes with
$\varphi(\beta)$; since $k\ne 4$, Lemma
\ref{Lm: k|m if psi(alpha m)<->psi(beta),(k-1)|m if psi(beta m)<->psi(alpha)}$(a)$
implies that
$k$ divides $n$. If $\varphi\in {\mathcal S}_1(k,n)\cup {\mathcal S}_3(k,n)$, then
$\varphi(\beta^{n-1})=gb^{q(n-1)}g^{-1}=b^{q(n-1)}$ commutes with
$\varphi(\alpha)$, and Lemma
\ref{Lm: k|m if psi(alpha m)<->psi(beta),(k-1)|m if psi(beta m)<->psi(alpha)}$(b)$ 
shows that $k-1$ divides $n-1$.
Completely similar, if $\varphi\in {\mathcal S}_2(k,n)\cup {\mathcal S}_4(k,n)$,
then $\varphi(\beta^n)=ga^{qn}g^{-1}=a^{qn}$ commutes with $\varphi(\alpha)$
and $k-1$ divides $n$, and $\varphi\in {\mathcal S}_3(k,n)\cup {\mathcal S}_4(k,n)$
implies that $\varphi(\alpha^{n-1})=b^{p(n-1)}$ commutes
with $\varphi(\beta)$ and $k$ divides $n-1$.

Assume that $\varphi\in{\mathcal S}_1(k,n)$, that is,
$\varphi(\alpha)=a^p$, \ $\varphi(\beta)=gb^q g^{-1}$.
It follows from the above consideration that $k$ divides $n$
and $k-1$ divides $n-1$. Hence there exists an integer
$l\ge 0$ such that $n=k+lk(k-1)$. Relation (\ref{9.46}) shows that
$k(n-1)p=(k-1)nq$; therefore, $(lk+1)p=(l(k-1)+1)q$. Since the
numbers $lk+1$, $l(k-1)+1$ are co-prime, there exists an integer $t$
such that $p=t(l(k-1)+1)$ and $q=t(lk+1)$. Let us show that
$(t,k(k-1))=1$. Indeed, if $m=(t,k)>1$, then the ratios
$t'=t/m$ and $k'=k/m$ are integral and
$1\le k'<k$, \ $pk'=tk'(l(k-1)+1)
=mt' k'(l(k-1)+1)=t' k(l(k-1)+1)=t' n$.
Hence the element $\varphi(\alpha^{k'})=a^{pk'}=a^{t' n}$
commutes with $\varphi(\beta)$;
by Lemma
\ref{Lm: k|m if psi(alpha m)<->psi(beta),(k-1)|m if psi(beta m)<->psi(alpha)}$(a)$,
$k$ must be a divisor of $k'$, which is impossible. Similarly, one can
check that the inequality $(t,k-1)>1$ leads to a contradiction;
this completes the proof in the case when $\varphi\in {\mathcal S}_1(k,n)$.

Assume now that $\varphi\in{\mathcal S}_2(k,n)$, that is,
$\varphi(\alpha)=a^p$, \ $\varphi(\beta)=ga^q g^{-1}$.
Then $k$ and $k-1$ divide $n$; hence $n=lk(k-1)$ for some
integer $l\ge 1$. Relation (\ref{9.46}) shows that $kp=(k-1)q$;
consequently, there is an integer $t$ such that $p=t(k-1)$
and $q=tk$. Taking into account the relations $kp=(k-1)q=tk(k-1)$ and
$\alpha^k=\beta^{k-1}$, we obtain $a^{tk(k-1)}=ga^{tk(k-1)}q^{-1}$;
thus, $g$ commutes with $a^{tk(k-1)}$. This concludes the proof
in the case when $\varphi\in{\mathcal S}_2(k,n)$.

We skip the proofs for the cases $\varphi\in{\mathcal S}_3(k,n)$
and $\varphi\in{\mathcal S}_4(k,n)$, which are very similar to the
cases considered above.
\end{proof}

\begin{Remark}\label{Rmk 9.2} 
$B_3$ possesses non-integral special homomorphisms
$B_3\to B_n$ for every $n$ that is not
forbidden by Theorem
\ref{Thm: arithmetical properties of special homomorphisms B(k) to B(n)}.
Moreover, the conditions
$\varphi(\alpha)=b$, \ $\varphi(\beta)=a^2$
define the special non-integral (actually, surjective) homomorphism
$\varphi\colon B_4\to B_3$;
therefore, if there exists a non-integral special homomorphism
$B_3\to B_n$, then there
is a non-integral special homomorphism $B_4\to B_n$.
For $k>4$ and $n\in P^3(k)\cup P^4(k)$ I do not know any example of
a non-integral special homomorphism $B_k\to B_n$;
however, for any $k$ and any $n\in P^1(k)\cup P^2(k)$
such homomorphisms do exist (see examples below).
\hfill $\bigcirc$
\end{Remark}

\begin{Examples}\label{Ex: Examples}
Let $k$ and $m$ be natural numbers.
There is a natural way to construct some interesting
embeddings $B_k\hookrightarrow B_{mk}$.
Let $\sigma_1,...,\sigma_{k-1}$
be the canonical generators in $B_k$
and let $\chi\colon B_k\to{\mathbb Z}$
be the canonical integral projection 
(see Sec. \ref{Ss: Commutator subgroup B'(k)}).
Any geometric $m$-braid (that is, a braid on $m$ strings)
may be considered as a ``thin rope"; take some $m$-braid $v$
and consider the $mk$-braids
$\sigma_i\circ v$ obtained by replacing of every string
in $\sigma_i$ by the thin rope $v$.
The correspondence $\sigma_i\mapsto\sigma_i\circ v\in B_{mk}$,
\ $i=1,...,k-1$, defines a homomorphism
$B_k\to B_{mk}$.
For any $k$-braid $g$ its image in $B_{mk}$ is the
$mk$-braid $g\circ v^{\chi(g)}$ obtained by replacing of every string
in $g$ by the thin rope $v^{\chi(g)}$.
\vskip0.2cm

\noindent This construction may be modified as follows. Let $\sigma_i*v$
be the $mk$-braid obtained from $\sigma_i$ by replacing of 
the $i$'th string
by the thin rope $v$ and all the rest strings by the thin ropes
corresponding to the trivial $m$-braid.
The correspondence $\sigma_i\mapsto\sigma_i*v\in B_{mk}$,
\ $i=1,...,k-1$, defines a homomorphism
$\phi_v\colon B_k\to B_{mk}$.
\vskip0.2cm

\noindent More formally, the homomorphism $\phi_v$ may be described as
follows. Take the canonical generators $s_1,...,s_{mk-1}$ of
$B_{mk}$ and define the elements $a_{ij}$
according to (\ref{alpha beta}), that is,
%&
\begin{equation}\label{9.47}
a_{ii}=1 \ \ \text{for all} \ i, \qquad \text{and} \ \
a_{ij} = s_i s_{i+1}\cdots s_{j-1} \ \ \text{for} \ \ 1\le i<j\le mk.
%\eqnum{9.47}
\end{equation}
%&
Put
%&
\begin{equation}\label{9.48}
u_i=a_{im,(i+1)m}a_{im-1,(i+1)m-1}\cdots a_{(i-1)m+1,im+1},
\ \ 1\le i\le k-1.
%\eqnum{9.48}
\end{equation}
%&
Let $v=v(x_1,...,x_{m-1})$ be any word in variables
$x_1,x_1^{-1},...,x_{m-1},x_{m-1}^{-1}$. Define the elements
$v_i\in B_{mk}$ by
%&
\begin{equation}\label{9.49}
v_i=v\left(s_{(i-1)m+1},s_{(i-1)m+2},...,s_{(i-1)m+m-1}\right),
\ \ 1\le i\le k.
%\eqnum{9.49}
\end{equation}
%&
The elements $u_i$, $v_j$ satisfy the following relations:
%&
\begin{equation}\label{9.50}
\aligned
&u_i u_j=u_j u_i \ \ \text{if} \ i,j=1,...,k-1, \ \ |i-j|>1;
\qquad v_i v_j=v_j v_i \ \ \text{for} \ i,j=1,...,k;\\
&u_i v_j=v_j u_i \ \ \text{if} \ j<i \ \ \text{or} \ \ j>i+1;
\qquad\qquad u_i v_{i+1}=v_i u_i \ \ \text{for} \ 1\le i<k-1;\\
&u_i v_i =v_{i+1} u_i \ \ \text{and} \ \
u_i u_{i+1} u_i=u_{i+1} u_i u_{i+1} \ \ \text{for} \ 1\le i<k-1;\\
&u_i u_{i+1}\cdots u_j = a_{im,(j+1)m}a_{im-1,(j+1)m-1}\cdots
a_{(i-1)m+1,jm+1}.
\endaligned
\end{equation}
%&
These relations imply that the correspondence
$\sigma_i\mapsto \sigma_i*v=v_i u_i$, \ $1\le i\le k-1$,
defines a homomorphism
$\phi_v\colon B_k\to B_{mk}$.
\vskip0.2cm

\noindent If $k>2$ and $v\ne 1$, then the homomorphism $\phi_v$ 
is non-abelian.
It turns out that for suitable $m$ and $v\in B_m$
the homomorphism $\phi_v$ is special.
\end{Examples}

\newpage

\bibliographystyle{amsplane}
%\ifx\undefined\bysame
%\newcommand{\bysame}{\leavevmode\hbox to3em{\hrulefill}\,}
%\fi

%\printindex
\begin{theindex}

  \item $E(x)$\hfill 12, 23
  \item $J_k$\hfill 8, 9
  \item $PB_k$\hfill 5
  \item $PB_k^{(r)}$\hfill 5
  \item $P^i(k)$\hfill 10
  \item $\GCD(m,n)$\hfill 17
  \item $\Inv H$; \ $\Inv_r H$\hfill 2
  \item $\lambda_k'$\hfill 8
  \item $\lambda_{k,m}$\hfill 8
  \item $\modo{N}_r$\hfill 17
  \item $\nu_6$\hfill 6
  \item $\nu_6'$\hfill 9
  \item $\nu_{4,1}$, $\nu_{4,2}$, $\nu_{4,3}$\hfill 6
  \item $\preccurlyeq$\hfill 2
  \item ${\mathbb F}_m$;\, $\mathbb F$\hfill 1
  \item ${\mathfrak C}_r(A)$\hfill 2

  \indexspace

  \item Arithmetic progressions $P^i(k)$\hfill 10
  \item Artin Fixed Point Lemma
    \subitem analog of\hfill 79
  \item Artin Fixed Point Lemma\hfill \phantom\hfill{11}, 
		\phantom\hfill{79}
  \item Artin homomorphism $\nu_6$\hfill 6
  \item Artin homomorphisms $\nu_{4,1}$, $\nu_{4,2}$, $\nu_{4,3}$\hfill 
		6
  \item Artin Lemma
    \subitem on cyclic decomposition of $\widehat\sigma_1$\hfill 12, 22
    \subitem on fixed points of $\widehat\sigma_1$\hfill 11, 79
  \item Artin Theorem
    \subitem improvement of\hfill 11, 30
  \item Artin Theorem\hfill \phantom\hfill{6}, 15

  \indexspace

  \item Braid-like couples\hfill 2

  \indexspace

  \item Canonical epimorphism $\mu$\hfill 5
  \item Canonical homomorphism $B'_k\to{\mathbf S}(k)$\hfill 8
  \item Canonical presentation of $B_k$\hfill 3
  \item Center\hfill 6
  \item Chebyshev Theorem\hfill 11
  \item Cohomology\hfill \phantom\hfill{12}
  \item Conjugate homomorphisms\hfill 2
  \item Cyclic type\hfill 2

  \indexspace

  \item Disjoint product\hfill 3

  \indexspace

  \item Embedding $\lambda_k'$\hfill 8, 11
  \item Embeddings $\lambda_{k,m}$ and $\lambda_k'$\hfill 8
  \item Endomorphisms of $B'_k$\hfill 15

  \indexspace

  \item Fadell-Neuwirth Theorem\hfill 3
  \item Fixed points
    \subitem of $\widehat\sigma_1$ for $6<k<n<2k$\hfill 79
  \item Fixed points and primes\hfill 11

  \indexspace

  \item Gorin's relation\hfill 7
  \item Gorin-Lin Theorem\hfill 6, 15
  \item Group
    \subitem Hopfian\hfill 2
    \subitem perfect\hfill 1
    \subitem residually finite\hfill 1

  \indexspace

  \item Homomorphism
    \subitem $\Omega$\hfill 12
    \subitem abelian\hfill 2
    \subitem cyclic\hfill 2
    \subitem integral\hfill 2
    \subitem primitive\hfill 3
    \subitem special\hfill 10
    \subitem transitive\hfill 3
  \item Homomorphism $\nu_6$\hfill 6
  \item Homomorphism $\nu_6'$\hfill 9
  \item Homomorphisms
    \subitem $B'_k\to{\mathbf S}(k)$\hfill 15
    \subitem $B'_k\to{\mathbf S}(n)$\hfill 11
    \subitem $B_k\to B_n$, $k>n$\hfill 15
    \subitem $B_k\to{\mathbf S}(2k)$
      \subsubitem model\hfill 9
      \subsubitem standard\hfill 9
    \subitem $B_k\to{\mathbf S}(k)$\hfill 6
    \subitem $B_k\to{\mathbf S}(k+1)$\hfill 11, 71
    \subitem $B_k\to{\mathbf S}(n)$ for $n<k$\hfill 11
    \subitem $B_k\to{\mathbf S}(n)$, $k<n\le 2k$\hfill 14
    \subitem $\nu_{4,1}$, \ $\nu_{4,2}$, \ $\nu_{4,3}$\hfill 6
  \item Homomorphisms $B_k\to{\mathbf S}(n)$, $n<k$\hfill 27

  \indexspace

  \item Improvement of Artin Theorem\hfill 11, 30

  \indexspace

  \item Lemma
    \subitem of Artin on fixed points \hfill 11

  \indexspace

  \item Malcev Theorem\hfill 1
  \item Markov Theorem\hfill 5, 15
  \item Model homomorphisms $B_k\to{\mathbf S}(2k)$\hfill 9
  \item Murasugi Theorem\hfill 10

  \indexspace

  \item Non-abelian special homomorphisms\hfill 10
  \item Normal series in ${PB_k}$\hfill 5

  \indexspace

  \item Presentation of $B_k'$\hfill 7
  \item Pure braid group\hfill 5
  \item Pure commutator subgroup $J_k$\hfill 8, 9

  \indexspace

  \item r-component of a permutation\hfill 2
  \item Reduction of homomorphism\hfill 3

  \indexspace

  \item Special homomorphisms
    \subitem non-abelian\hfill 10
  \item Special homomorphisms\hfill 10
  \item Special presentation of $B_k$\hfill 4
  \item Special system of generators in $B_k$\hfill 4
  \item Standard homomorphisms $B_k\to{\mathbf S}(2k)$\hfill 9
  \item Standard system of generators in $B_k$\hfill 5

  \indexspace

  \item Theorem
    \subitem of Artin
      \subsubitem improvement of\hfill 11, 30
    \subitem of Artin\hfill 6, 15
    \subitem of Chebyshev\hfill 11
    \subitem of Fadell-Neuwirth\hfill 3
    \subitem of Gorin-Lin\hfill 6, 15
    \subitem of Malcev\hfill 1
    \subitem of Markov\hfill 5, 15
    \subitem of Murasugi\hfill 10
    \subitem on endomorphisms of $B_k$ for $k\ne 4$\hfill 8
    \subitem on homomorphisms
      \subsubitem $B_k'\to B_n$ for $n<k$\hfill 8
      \subsubitem $B_k'\to{\mathbf S}(k)$\hfill 9
      \subsubitem $B_k'\to{\mathbf S}(n)$ for $n<k$\hfill 8
      \subsubitem $B_k\to B_n$ for $n<k$\hfill 8
      \subsubitem $B_k\to{\mathbf S}(2k)$ for $k>6$\hfill 9
      \subsubitem $B_k\to{\mathbf S}(k)$\hfill 6
      \subsubitem $B_k\to{\mathbf S}(n)$ for $6<k<n<2k$\hfill 9
      \subsubitem $B_k\to{\mathbf S}(n)$ for $n<k$\hfill 8
    \subitem on homomorphisms $B_k\to{\mathbf S}(k+1)$\hfill 11, 71
    \subitem on homomorphisms $B_k\to{\mathbf S}(n)$, $n<k$\hfill 27
    \subitem on pure commutator subgroup $J_k$\hfill 9
  \item Torsion\hfill 3

\end{theindex}

\end{document}